%% file: document.tex
\documentclass[11pt]{elsarticle}
\usepackage{amsmath,amssymb,amsthm}
\usepackage{mathrsfs}
\usepackage[bookmarks=true,
  bookmarksnumbered=true,
  bookmarksopen=true,
  colorlinks,
  pdfborder=001,
linkcolor=black]{hyperref}
\usepackage{bm}
\usepackage{array}
\usepackage{float}
\usepackage{color}
\usepackage{lineno}
\usepackage{subcaption}
\usepackage{stmaryrd}
\usepackage{extarrows}
\usepackage{subcaption}
\usepackage{multirow}

\newtheorem{rmk}{Remark}[section]
\newtheorem{example}{Example}[section]
\newtheorem{proposition}{Proposition}[section]
\newtheorem{definition}{Definition}[section]
\newproof{pf}{Proof}
\numberwithin{equation}{section}
\numberwithin{figure}{section}
\numberwithin{table}{section}
\allowdisplaybreaks

\newcommand\diag{\mathrm{diag}}

\newcommand\dd{\mathrm{d}}

\newcommand\bF{\bm{F}}

\newcommand\bv{\bm{v}}

\newcommand\bx{\bm{x}}

\newcommand\bV{\bm{V}}
\newcommand\bU{\bm{U}}

\newcommand\bT{\bm{T}}

\newcommand\tbF{\widetilde{\bm{F}}}
\newcommand\tbU{\widetilde{\bm{U}}}
\newcommand\thF{\widehat{\bm{F}}}
\newcommand\thU{\widehat{\bm{U}}}

\newcommand\pd[2]{\dfrac{\partial {#1}}{\partial {#2}}}

\newcommand\abs[1]{\lvert #1 \rvert}
\newcommand\jump[1]{\llbracket #1 \rrbracket}
\newcommand\mean[1]{\{\!\!\{ #1 \}\!\!\}}
\newcommand\meanln[1]{\{\!\!\{ #1 \}\!\!\}^{\text{ln}}}

\begin{document}

	\begin{frontmatter}
	
	\title{High-order accurate entropy stable adaptive moving mesh finite difference schemes for (multi-component) compressible Euler equations with the stiffened equation of state}
	
	            \author{Shangting Li}
            	\ead{shangtl@pku.edu.cn}
              \author{Junming Duan\corref{cor2}}
	          \ead{duanjm@pku.edu.cn}
            	\address{Center for Applied Physics and Technology, HEDPS and LMAM,
				School of Mathematical Sciences, Peking University, Beijing 100871, P.R. China}
			\cortext[cor2]{Current address: 
\'Ecole polytechnique F\'ed\'erale de Lausanne, 1015 Lausanne, Switzerland.}
		   \cortext[cor1]{Corresponding author. Fax:~+86-10-62751801.}
	         \author{Huazhong Tang\corref{cor1}}
	         \ead{hztang@math.pku.edu.cn}
			\address{Nanchang Hangkong University, Jiangxi Province, Nanchang 330000, P.R. China;  Center for Applied Physics and Technology, HEDPS and LMAM,
				School of Mathematical Sciences, Peking University, Beijing 100871, P.R. China}

	
	\begin{abstract}
		This paper extends the  high-order entropy stable (ES) adaptive moving mesh finite difference schemes developed in 
		\cite{duan2021highorder}
		to the two- and three-dimensional (multi-component) compressible Euler equations
    	with   the stiffened equation of state.
	The two-point entropy conservative (EC) flux is  first  constructed in the curvilinear coordinates.
	The high-order semi-discrete EC schemes are given with the aid of the two-point EC flux and the high-order  discretization of  the geometric conservation laws, 
and then the high-order semi-discrete ES schemes  satisfying the entropy inequality
 are derived by adding the high-order  dissipation term based on the multi-resolution weighted essentially non-oscillatory (WENO) reconstruction for the scaled entropy variables to
 the EC schemes.
 The explicit strong-stability-preserving Runge-Kutta methods are used for the time discretization and the  mesh points are adaptively redistributed  by iteratively solving the mesh redistribution equations with an
  appropriately chosen monitor function.  Several  2D and 3D numerical tests  are conducted on the parallel computer system with the MPI programming to validate the accuracy and the ability to capture  effectively  the localized structures of the proposed schemes.
	\end{abstract}
	
	
\end{frontmatter}
\input{Intro}

\input{NumSches}

\input{MovingMesh}

\input{NumTests}

\input{Conc}

\input{WENO}



\input{References}
\end{document}

%% file: Intro.tex
\section{Introduction}\label{section:Intro}
This paper is concerned with the high-order accurate {entropy stable (ES)} adaptive moving mesh finite difference schemes for the $d$-dimensional (multi-component) compressible Euler equations \cite{Larrouturou1989}
\begin{equation}\label{eq:ConserLaw}
\begin{aligned}
&\frac{\partial \bU}{\partial t}+\sum_{k=1}^{d} \frac{\partial \bF_{k}(\bU)}{\partial x_{k}}=0, \\
&\boldsymbol{U}= \left(
\rho_1, \cdots,\rho_N, \rho \bm{v}^\mathrm{T},  E
\right)^\mathrm{T}, \\
&\boldsymbol{F}_{k}=\left(
\rho_1 v_{k}, \cdots, \rho_N v_{k},\rho v_{k}\bm{v}^\mathrm{T} + p\bm{e}_{k}^\mathrm{T}, (E+p) v_{k}
\right)^\mathrm{T},
\end{aligned}
\end{equation}
where $d=2$ or 3, $\rho_{\ell}$
denotes the $\ell$th species density,  $ \ell = 1,\cdots,N$, $\rho=\sum_{\ell=1}^{N}\rho_{\ell}$ is the total density, $\bm{v} = (v_{1}, \cdots, v_{d})^{\mathrm{T}}$  denotes the velocity  vector, $\bm{e}_{k}$ is  the $k$th column of the $d \times d$ unit matrix, and $E=\rho e  +\rho|\bm{v}|^{2} /2$ is the total energy with   the specific internal energy $e$.
Assume that  the $N$ species fluids are in thermal equilibrium and
the pressure mixture $p$ is governed by the stiffened equation of state (EOS) \cite{SHYUE1998208}
$$
\dfrac{p+ p_{\infty}}{\Gamma -1} + p_{\infty}= \rho e=\sum_{\ell=1}^{N} \rho_{\ell} e_{\ell},\quad \rho_{\ell} e_{\ell} =c_{v, \ell}\rho_{\ell} T +{p_{\infty, \ell}},
$$
where  $T$ is the  temperature,  $p_{\infty,\ell}$ is the pressure constant   related to the material,
and $\Gamma :=\left(\sum_{\ell=1}^{N}  \Gamma_{\ell}c_{v,\ell}\rho_{\ell}\right)/\left(\sum_{\ell=1}^{N}  c_{v,\ell}\rho_{\ell}\right)$
with the specific heat  at constant volume $c_{v, \ell}$ and the   $\ell$th species adiabatic index $\Gamma_{\ell}$, see e.g.  \cite{Larrouturou1989}.  If taking $p_\infty = \sum_{\ell=1}^{N}  p_{\infty,\ell}$,
then the above EOS can be rewritten as follows
\begin{align}\label{eq:StiffenedEOS}
p= \sum_{\ell  = 1}^{N} \left(\rho_{\ell} R_\ell T-
{p_{\infty, \ell}}\right),
\end{align}
with $R_\ell =  c_{v, \ell}\left(\Gamma_{\ell}-1\right)$.  
The stiffened EOS is often used when considering water under very high pressures (typical applications are underwater nuclear explosions and sonic shock lithotripsy  etc.).
If $p_{\infty, \ell} = 0$, $\ell = 1,\cdots, N$, then   \eqref{eq:StiffenedEOS}  reduces to the ideal gas EOS.
The multi-component ($N$ species) compressible Euler equations \eqref{eq:ConserLaw} can also be  viewed as the ``one-component'' compressible Euler equations complemented with $N-1$ species mass-conservation equations  \cite{GOUASMI2020112912}. 
In view of this, one may discretize  the multi-component ($N$ species) compressible Euler equations \eqref{eq:ConserLaw} wholly, see e.g. \cite{LARROUTUROU199159,MR991358,SHYUE1998208}, or discretize the ``one-component'' compressible Euler equations and
the complemented $N-1$   species  equations  separately, see e.g. \cite{1987Explicit}.
 Besides,  the level-set methods \cite{davis1992interface, MULDER1992209},   the volume-of-fluid methods \cite{miller1996high,puckett19923d}, and the BGK-based method \cite{xu1997bgk} were also studied for the multi-component Euler equations.
Even if the initial data are sufficiently smooth, the quasi-linear hyperbolic conservation laws such as \eqref{eq:ConserLaw}
may have discontinuous  solutions
 so that
one should consider the weak solutions which are not  unique in general and single out the physically relevant solution among all the weak solutions by the entropy condition etc.

\begin{definition}[Entropy function]
	A scalar function $\eta(\bU)$ is called an entropy function  for the system \eqref{eq:ConserLaw} if there exist associated entropy fluxes $q_{k}(\bU)$ satisfying
	$$
	q_{k}^{\prime}(\bU)=\bV^{\mathrm{T}} \bF_{k}^{\prime}(\bU), \quad k=1, \cdots, d,
	$$
	where $\bV=\eta^{\prime}(\bU)^{\mathrm{T}}$ is called the entropy variables and $\left(\eta, q_{k}\right)$ forms an entropy pair. Further, one can define the entropy potential $\phi$ and corresponding flux $\psi_k$   by using the conjugate variables as follows
	\begin{align*}
	\phi:=\bV^\mathrm{T}\bU-\eta,\quad \psi_k:=\bV^\mathrm{T}\bF_k - q_k.
	\end{align*}
\end{definition}

If there exists a strictly convex entropy pair  for the hyperbolic conservation laws
\eqref{eq:ConserLaw}, $\eta''(\bU)>0$,
then the entropy solution should satisfy the following entropy condition
\begin{equation}\label{eq:entrosolution}
\frac{\partial \eta(\bU)}{\partial t}+\sum_{k  = 1}^{d} \frac{\partial q_{k}(\bU)}{\partial x_{k}} \leq 0,
\end{equation}
where the equality  holds for the smooth solutions, while
 the inequality is for the nonsmooth solutions in the sense of distributions.

Integrating \eqref{eq:entrosolution} in space with
periodic or zero entropy flux
boundary conditions, the total amount of {the} entropy $\int 	\eta(\bU(\cdot,t)){\rm{d}}\bm{x}$ decreases in time. This is a generalization of the $L^2$-energy bound encountered in the linear case \cite{Tadmor2003Entropy}.
To select the physically relevant solution, it is important to construct the high-order accurate  entropy conservative (EC)  or ES schemes which satisfy a discrete or semi-discrete version of  the entropy condition \eqref{eq:entrosolution}. 
The framework of the second-order EC scheme which satisfies the semi-discrete entropy identity was  established in  \cite{Tadmor1987The, Tadmor2003Entropy}, and the higher-order extension was introduced in \cite{Lefloch2002Fully}.
It should be noticed that the EC schemes may produce oscillations near the discontinuities so that some numerical dissipative terms need to be added to obtain the ES schemes which  suppress possible oscillations \cite{Fjordholm2012Arbitrarily}.
With the help of the summation-by-parts (SBP) operators \cite{Carpenter2014Entropy, Gassner2013A}, the ES discontinuous Galerkin (DG) schemes {were} developed, including the space-time DG method \cite{Hiltebrand2014Entropy},
the DG spectral element methods
 \cite{Gassner2013A,Carpenter2014Entropy} and the DG methods on the unstructured simplex meshes \cite{Chen2020Review}.
Recently, the EC or ES schemes {were}  extended to 
the relativistic hydrodynamic
equations {\cite{Bhoriya2020Entropy,Duan2020RHD,DUAN2021109949}},  the relativistic magnetohydrodynamic equations \cite{duan2021highorder, Duan2020RMHD,  Wu2020Entropy},
the multi-component Euler equations 
 \cite{GOUASMI2020112912,Renac2021},
 and so on.

Adaptive moving mesh methods  have been playing an  important role in solving partial differential equations due to improving the efficiency and quality of the numerical computations, including the grid redistribution approaches \cite{Brackbill1993An,Brackbill1982Adaptive,Ren2000An,Wang2004A,Winslow1967Numerical}, the moving finite element methods \cite{Davis1982,Miller1981} and the moving mesh PDEs methods \cite{CAO1999221,CENICEROS2001609,Stockie2001}.
 The readers  are  referred to the
 review articles \cite{Budd2009Adaptivity,Tang2005Moving} and references therein.
This paper focuses on the high-order ES adaptive moving mesh finite difference schemes for the  (multi-component) compressible Euler equations $(N = 1, 2)$ with the stiffened EOS.
The   two-point EC fluxes for the (multi-component) compressible Euler equations with the stiffened EOS are 
explicitly derived in curvilinear coordinates,
and then  are combined with the high-order discrete geometric conservation laws to give the high-order EC fluxes.
The high-order ES fluxes in curvilinear coordinates are obtained by adding appropriate dissipation terms into the high-order EC fluxes, which are built on the jump of the high-order accurate  multi-resolution WENO reconstruction values of the scaled entropy variables.
The mesh adaptation is implemented by iteratively solving the Euler-Lagrange equations of the mesh adaptation functional in the computational domain with appropriate  monitor function.
Compared to the single-component case, the monitor function needs to contain more information on the solutions of the multi-component compressible Euler equations to produce high-quality mesh.
The semi-discrete schemes are equipped with
the third-order accurate explicit strong-stability preserving (SSP) Runge-Kutta (RK) schemes to obtain the
fully-discrete schemes.

%
This paper is organized as follows. Section \ref{section: EntCon} introduces the form of the  system \eqref{eq:ConserLaw} in the curvilinear coordinates and corresponding entropy conditions.
 Section \ref{section:Num} gives a sufficient condition for the EC fluxes, two-point EC fluxes,
the high-order  discrete geometric conservation laws, and the high-order EC and ES schemes for the (multi-component) compressible Euler equations with the stiffened EOS in curvilinear coordinates.
Adaptive moving mesh strategy is presented in Section \ref{section:MM}. Several  2D and 3D numerical results are presented in Section \ref{section:Result} to validate the effectiveness and performance of our schemes  on the parallel computer system with the MPI communication. Section \ref{section:Conc} gives some conclusions.

%% file: NumSches.tex
\section{Entropy conditions in curvilinear coordinates}\label{section: EntCon}
 This section introduces the entropy conditions in curvilinear coordinates similar to that in \cite{duan2021highorder}.
%
The adaptive moving meshes in the physical domain $\Omega_p$
with coordinates $\bx=(x_1,\cdots,x_d)$  can be generated as the images of a reference mesh in the computational domain $\Omega_c$ with coordinates $\bm{\xi}=(\xi_1,\cdots,\xi_d)$  by a time dependent,
differentiable, one-to-one coordinate mapping {$\bx = \bx(\bm{\xi},t)$},
which can be expressed as
\begin{align}\label{eq:transf}
t=\tau,\ \ \bx=\bx(\bm{\xi},\tau),\ \
\bm{\xi}=(\xi_1,\cdots,\xi_d)\in\Omega_c,
\end{align}
under which the system \eqref{eq:ConserLaw} can be transformed as the following conservative form 
\begin{align}\label{eq:ConserLaw_curv}
\pd{\left(J\bU\right)}{\tau}+\sum_{k=1}^d\dfrac{\partial}{\partial\xi_k}\left[{\left(J\pd{\xi_k}{t}\bU\right)}
+\sum_{j=1}^d{\left(J\pd{\xi_k}{x_j}\bF_j\right)}\right]=0,
\end{align}
where  $J=\det\left(\pd{(t,\bx)}{(\tau,\bm{\xi})}\right)$.
For  \eqref{eq:transf}, one has
the following geometric conservation laws (GCLs)
\begin{equation}\label{eq:GCL}
\begin{aligned}
&\text{VCL:}\quad \pd{J}{\tau}+\sum_{k=1}^d\dfrac{\partial}{\partial\xi_k}{\left(J\pd{\xi_k}{t}\right)}=0,\\
&\text{SCLs:}\quad \sum_{k=1}^d\dfrac{\partial}{\partial\xi_k}{\left(J\pd{\xi_k}{x_j}\right)}=0,~ j=1,\cdots,d,
\end{aligned}
\end{equation}
where the
 volume conservation law (VCL)  implies that the volumetric increment of a moving cell is equal to the sum of the changes along the surfaces that enclose the cell,
while the  surface conservation laws (SCLs)  indicate that the cell volume should be closed by its surfaces \cite{Zhang1993Discrete}.

Utilizing  \eqref{eq:ConserLaw_curv} and the GCLs \eqref{eq:GCL} can derive the entropy condition in curvilinear coordinates   
\begin{align*}
\pd{\left(J\eta\right)}{\tau}+\sum_{k=1}^d\dfrac{\partial}
{\partial\xi_k}\left[{\left(J\pd{\xi_k}{t}\eta\right)}
+\sum_{j=1}^d{\left(J\pd{\xi_k}{x_j}q_j\right)}\right]\leq0,
\end{align*}
where the equality  holds
for the smooth solutions of  \eqref{eq:ConserLaw}, and
the  inequality  is in the sense of {distributions}
for the nonsmooth solutions.

\section{Numerical schemes}\label{section:Num}
This section only presents  the 3D moving mesh EC and ES schemes for the system \eqref{eq:ConserLaw_curv} on the structured hexahedral mesh
 following \cite{duan2021highorder},
 	because the 2D schemes can be considered as the degenerative case, see the appendices in \cite{duan2021highorder} for more details.

Let us choose $\Omega_c$ as a cuboid
$[a_1,b_1]\times[a_2,b_2]\times[a_3,b_3]$ and divide it into a   fixed orthogonal uniform mesh
$\{(\xi_{1,i_1},\xi_{2,i_2},
\xi_{3,i_3})$:
$a_k=\xi_{k,1}<\cdots<\xi_{k,i_k}
<\cdots<\xi_{k,N_k}=b_k$, $k=1,2,3\}$
with the constant mesh size $\Delta \xi_k=\xi_{k,i_k+1}-\xi_{k,i_k}$.
For the sake of brevity, the index $\bm{i}=(i_1,i_2,i_3)$   denotes the point
$(\xi_{1,i_1}, \xi_{2,i_2}, \xi_{3,i_3})$ and the notation $\{\bm{i},k,m\}$ means that the index $\bm{i}$ increases $m$ along $i_k$-direction, e.g. $\{\bm{i},1,\frac12\}$ is $(i_1+\frac12, i_2, i_3)$.

Consider the following semi-discrete conservative {$2w$}-order ($w\geq1$) finite difference schemes
for
\eqref{eq:ConserLaw_curv}
and the first equation in \eqref{eq:GCL}
\begin{align}
\label{eq:semi_U}
&\dfrac{\dd}{\dd t}(J\bU)_{\bm{i}}=
-\sum_{k=1}^3\dfrac{1}{\Delta \xi_k}\left(\left(\widehat{\bF}_k\right)_{\bm{i},k,+\frac12}^{2w\rm{th}}-\left(\widehat{\bF}_k\right)_{\bm{i},k,-\frac12}^{2w\rm{th}}\right),
\\
\label{eq:semi_J}
&\dfrac{\dd}{\dd t}J_{\bm{i}}=
-\sum_{k=1}^3\dfrac{1}{\Delta \xi_k}\left(\left(\widehat{J\pd{\xi_k}{t}}\right)_{\bm{i},k,+\frac12}^{2w\rm{th}}-\left(\widehat{J\pd{\xi_k}{t}}\right)_{\bm{i},k,-\frac12}^{2w\rm{th}}\right),
\end{align}
where
$J_{\bm{i}}(t)$ and $(J\bU)_{\bm{i}}(t)$ approximate the point values of
$J\left(t,\bm{\xi}\right)$ and $(J\bU)(t,\bm{\xi})$ at $\bm{i}$, respectively, { $\left(\widehat{\bF}_k\right)_{\bm{i},k,\pm\frac12}^{2w\rm{th}}(t)$}
is the numerical flux approximating the continuous flux
{$\left(J\pd{\xi_k}{t}\bU+\sum\limits_{j=1}^3 J\pd{\xi_k}{x_j}\bF_j\right)(t,\bm{\xi})$} at $\{\bm{i}, k, \pm\frac{1}{2}\}$, $k=1,2,3$, and
 {
 $\left(\widehat{J\pd{\xi_k}{t}}\right)_{\bm{i},k,\pm\frac12}^{2w\rm{th}}$}
is the
flux   approximating the  metric $J\pd{\xi_k}{t}$ at $\{\bm{i}, k, \pm\frac{1}{2}\}$, which is used to update the metric Jacobian $J_{\bm{i}}$,
see \eqref{eq:GCL_flux} in Section \ref{sec:sufficon}.
 Corresponding discrete version of the SCLs in  \eqref{eq:GCL}
  will be of the form
\begin{equation}\label{eq:SCL_dis}
\sum_{k=1}^3\dfrac{1}{\Delta \xi_k}\left(\left(\widehat{J\pd{\xi_k}{x_j}}\right)_{\bm{i},k,+\frac12}^{2w\rm{th}}-\left(\widehat{J\pd{\xi_k}{x_j}}\right)_{\bm{i},k,-\frac12}^{2w\rm{th}}\right)=0,~j=1,2,3.
\end{equation}

\begin{definition}
The semi-discrete  scheme \eqref{eq:semi_U}-\eqref{eq:semi_J} or its flux is EC,
if  its solution   satisfies the semi-discrete entropy identity
	\begin{equation}\label{eq:NumEntropyID}
	\dfrac{\dd}{\dd t}J_{\bm{i}}\eta(\bU_{\bm{i}}(t))
	+\sum_{k=1}^3\dfrac{1}{\Delta \xi_k}\left(\left(\widehat{q}_k\right)_{\bm{i},k,+\frac12}^{2w\rm{th}}(t)-\left(\widehat{q}_k\right)_{\bm{i},k,-\frac12}^{2w\rm{th}}(t)\right)=0,
	\end{equation}
	where the numerical entropy flux
	$\left(\widehat{q}_k\right)_{\bm{i},k,+\frac12}^{2w\rm{th}}$ is consistent with the entropy flux
	$J\pd{\xi_k}{t}\eta+\sum\limits_{j=1}^3J\pd{\xi_k}{x_j}q_j$.
 \end{definition}

\subsection{A sufficient condition for {the} EC  fluxes}\label{sec:sufficon}

%
%

 Similar to the special relativistic (magneto)hydrodynamics in \cite{duan2021highorder}, one can  deduce the following sufficient condition for the two-point EC  fluxes of \eqref{eq:semi_U}.

 \begin{proposition}\label{prop:ECFlux_curv}\rm
	If   a  two-point flux  $\thF_k\left(\bU_l, \bU_r, \left(J\pd{\xi_k}{\zeta}\right)_l, \left(J\pd{\xi_k}{\zeta}\right)_r \right)$, $\zeta=t,x_1,x_2,x_3$,
	 being consistent with
 	$J\pd{\xi_k}{t}\bU+\sum\limits_{j=1}^3 J\pd{\xi_k}{x_j}\bF_j$,
	 satisfies
	\begin{align}\label{eq:ECCondition_curv}
	\left(\bV(\bU_r)-\bV(\bU_l)\right)^\mathrm{T}\thF_k=&	~\dfrac12\left(\left(J\pd{\xi_k}{t}\right)_l+\left(J\pd{\xi_k}{t}\right)_r\right)\left(\phi(\bU_r)-\phi(\bU_l)\right) \nonumber
\\	&+\sum_{j=1}^3\dfrac12\left(\left(J\pd{\xi_k}{x_j}\right)_l+\left(J\pd{\xi_k}{x_j}\right)_r\right)\left(\psi_{j}(\bU_r)-\psi_{j}(\bU_l)\right),
	\end{align}
	then the scheme \eqref{eq:semi_U} is EC, 	where the subscripts $l$ and $r$ represent two
	states, denoted  respectively by the left and right states.
 \end{proposition}

  If the two-point EC flux  $\thF_k\left(\bU_l, \bU_r, \left(J\pd{\xi_k}{\zeta}\right)_l, \left(J\pd{\xi_k}{\zeta}\right)_r \right)$ satisfying \eqref{eq:ECCondition_curv}
  is symmetric, $\zeta=t,x_1,x_2,x_3$,
then one can further derive the semi-discrete  $2w$th-order EC schemes \eqref{eq:semi_U}-\eqref{eq:semi_J}
	 with the following $2w$th-order EC  fluxes  
		\begin{align}\label{eq:ECFlux_curv_2p}
	&{\left(\thF_k\right)}_{\bm{i},k,+\frac12}^{{2w\rm{th}}}
	=~\sum_{m=1}^w\alpha_{w,m}\sum_{s=0}^{m-1}\thF_k\left(\bU_{\bm{i},k,-s}, \bU_{\bm{i},k,-s+m},
	\left(J\pd{\xi_k}{\zeta}\right)_{\bm{i},k,-s}, \left(J\pd{\xi_k}{\zeta}\right)_{\bm{i},k,-s+m}\right),\\
	\label{eq:GCL_flux}
	&\left(\widehat{J\pd{\xi_k}{\zeta}}\right)_{\bm{i},k,+\frac12}^{2w\rm{th}}
	=\sum_{m=1}^w\alpha_{w,m}\sum_{s=0}^{m-1}\dfrac12\left(\left(J\pd{\xi_k}{\zeta}\right)_{\bm{i},k,-s}
	+\left(J\pd{\xi_k}{\zeta}\right)_{\bm{i},k,-s+m}\right) ,
	\end{align}
	where
	the constants $\{\alpha_{w,m}\}$  satisfy the conditions  \cite{Lefloch2002Fully}
	$$
	\sum\limits_{m=1}^{w} m\alpha_{w, m}=1, \quad \sum\limits_{m=1}^{w} m^{2 s-1} \alpha_{w, m}=0, \
 \ s=2, \ldots, w.
	$$
 In this case,  corresponding  numerical entropy fluxes can be chosen
   as follows
	\begin{subequations}
		\begin{align*}
		\left(\widehat{q_k}\right)_{\bm{i},k,+\frac12}^{2w{\rm{th}}}=\sum_{m=1}^w\alpha_{w,m}\sum_{s=0}^{m-1}\widehat{q}_k\bigg(
		\bU_{\bm{i},k,-s}, &\bU_{\bm{i},k,-s+m}, \left(J\pd{\xi_k}{\zeta}\right)_{\bm{i},k,-s}, \left(J\pd{\xi_k}{\zeta}\right)_{\bm{i},k,-s+m}
		\bigg),\\
		\label{eq:NumEntropyFlux_curv}
		\widehat{q}_k\left(\bU_l, \bU_r, \left(J\pd{\xi_k}{\zeta}\right)_l, \left(J\pd{\xi_k}{\zeta}\right)_r \right)=&~\dfrac12\left(\bV(\bU_l)+\bV(\bU_r)\right)^\mathrm{T}\thF_k\left(\bU_l, \bU_r, \left(J\pd{\xi_k}{\zeta}\right)_l, \left(J\pd{\xi_k}{\zeta}\right)_r \right)\nonumber\\
		&-\dfrac14\left(\left(J\pd{\xi_k}{t}\right)_l+\left(J\pd{\xi_k}{t}\right)_r\right)\left(\phi(\bU_l)+\phi(\bU_r)\right)\nonumber\\
		&-\sum_{j=1}^3\dfrac14\left(\left(J\pd{\xi_k}{x_j}\right)_l+\left(J\pd{\xi_k}{x_j}\right)_r\right)\left(\psi_j(\bU_l)+\psi_j(\bU_r)\right).
		\end{align*}
	\end{subequations}

\subsection{Discrete GCLs}\label{subsection:GCLs}
 This section gives  {the} discrete GCLs, which are essential in the proof of the EC or ES property of the schemes. Failing to satisfy discrete GCLs may lead to a misrepresentation of the convective velocities and extra sources or sinks in the physically conservative media  \cite{Zhang1993Discrete}.

	To achieve the discrete SCLs, following \cite{duan2021highorder},
the $2w$th-order accurate discretizations for $J\pd{\xi_k}{x_j},~ k,j = 1,2,3,$ can be constructed.
For example, when $j = 1$, one has
\begin{equation}\label{eq:SCLCoeff}
\begin{aligned}
&\left(J\pd{\xi_1}{x_1}\right)_{\bm{i}}
=\left(\pd{x_2}{\xi_2}\pd{x_3}{\xi_3}-\pd{x_2}{\xi_3}\pd{x_3}{\xi_2}\right)_{\bm{i}}
=\dfrac{1}{\Delta\xi_2\Delta\xi_3}
\left(\delta_3\left[\delta_2\left[x_2\right]x_3\right]-\delta_2\left[\delta_3\left[x_2\right]x_3\right]\right),\\
&\left(J\pd{\xi_2}{x_1}\right)_{\bm{i}}
=\left(\pd{x_2}{\xi_3}\pd{x_3}{\xi_1}-\pd{x_2}{\xi_1}\pd{x_3}{\xi_3}\right)_{\bm{i}}
=\dfrac{1}{\Delta\xi_3\Delta\xi_1}
\left(\delta_1\left[\delta_3\left[x_2\right]{x_3}\right]-\delta_3\left[\delta_1\left[x_2\right]{x_3}\right]\right),\\
&\left(J\pd{\xi_3}{x_1}\right)_{\bm{i}}
=\left(\pd{x_2}{\xi_1}\pd{x_3}{\xi_2}-\pd{x_2}{\xi_2}\pd{x_3}{\xi_1}\right)_{\bm{i}}
=\dfrac{1}{\Delta\xi_1\Delta\xi_2}
\left(\delta_2\left[\delta_1\left[x_2\right]{x_3}\right]-\delta_1\left[\delta_2\left[x_2\right]{x_3}\right]\right),\\
\end{aligned}
\end{equation}
with the $2w$th-order central difference operator  in the $\xi_k$-direction
\begin{align*}
\delta_k[a_{\bm{i}}]=\dfrac12\sum_{m=1}^w\alpha_{w,m}\left(a_{\bm{i},k,+m} - a_{\bm{i},k,-m}\right).
\end{align*}
Combining \eqref{eq:SCLCoeff} with the $2w$th-order discretizations of the fluxes in \eqref{eq:GCL_flux}  easily  gets the discrete SCLs \eqref{eq:SCL_dis}.

	Regarding the discrete VCL, for the transformation \eqref{eq:transf},  one has
\begin{equation*}
\begin{aligned}
J\pd{\xi_k}{t}=-\sum_{j=1}^3\pd{x_j}{t}\left(J\pd{\xi_k}{x_j}\right),~k=1,2,3,
\end{aligned}
\end{equation*}
which can be approximated efficiently and easily by
\begin{equation}\label{eq:VCLCoeff}
\left(J\pd{\xi_k}{t}\right)_{\bm{i}}=-\sum_{j=1}^3(\dot{x}_j)_{\bm{i}}\left(J\pd{\xi_k}{x_j}\right)_{\bm{i}},
\end{equation}
where $\left(J\pd{\xi_k}{x_j}\right)_{\bm{i}}$ is given by \eqref{eq:SCLCoeff},
 and $(\dot{x}_j)_{\bm{i}}$, $j=1,2,3$, are the mesh velocities at $\bm{i}$ { and} will be determined in Section \ref{section:MM}.   Combining  \eqref{eq:VCLCoeff} with the  fluxes
 \eqref{eq:GCL_flux} can yield  the semi-discrete VCL \eqref{eq:semi_J}.

\begin{rmk}\rm
 It is known that violating  the free-stream condition may cause large errors and even lead to numerical instabilities for the high-order schemes \cite{VISBAL2002155}. It is proved  \cite{duan2021highorder} that
 the free-stream condition is satisfied by our high-order accurate fully-discrete adaptive moving mesh finite difference schemes derived by integrating \eqref{eq:semi_U}-\eqref{eq:semi_J} with the third-order accurate explicit SSP
RK  schemes \cite{Gottlieb2001Strong}.
\end{rmk}

\subsection{Two-point EC flux}

This subsection focuses on the construction of  a two-point
EC flux satisfying \eqref{eq:ECCondition_curv}.
One can verify that the following flux, similar to that in \cite{DUAN2021109949}, meets  the requirement
\begin{align}\label{eq:ECfluxMM}
\thF_k\left(\bU_l,\bU_r,\left(J\pd{\xi_k}{\zeta}\right)_l,\left(J\pd{\xi_k}{\zeta}\right)_r\right)
=&~\dfrac12\left(\left(J\pd{\xi_k}{t}\right)_l + \left(J\pd{\xi_k}{t}\right)_r\right)\widetilde{\bU} \nonumber\\
&+\sum_{j=1}^3\dfrac12\left(\left(J\pd{\xi_k}{x_j}\right)_l + \left(J\pd{\xi_k}{x_j}\right)_r\right)\widetilde{\bF}_j,
\end{align}
where $\widetilde{\bU}$ and $\widetilde{\bF}_j$
satisfy the following conditions, respectively,
\begin{align}\label{eq:ECCondition_comp}
\left(\bV_r-\bV_l\right)^\mathrm{T}\widetilde{\bU}=\phi_r - \phi_l,\quad
\left(\bV_r-\bV_l\right)^\mathrm{T}\widetilde{\bF}_j=(\psi_j)_r - (\psi_j)_l.
\end{align}

 In the following, we will give  the explicit expressions  of  the symmetric two-point EC fluxes
 for  the single- and two-component compressible Euler equations ($N = 1, 2$) with the stiffened EOS, separately.

\subsubsection{Single-component compressible Euler equations $(N=1)$}
This subsection   begins to construct the symmetric  two-point EC flux for the single-component compressible Euler equations $(N=1)$.

Assume that the numerical solutions satisfy $\rho_1, T >0$, and define the thermodynamic entropy as  $S_1 = c_{v,1}\ln{T} - R_1\ln{\rho_1}$, see  \cite{etde_20457592}. It is easy to
prove that the smooth solutions of \eqref{eq:ConserLaw} with $N=1$ satisfy
	\begin{align*}
		\frac{\partial(\rho_1  S_1)}{\partial t}+\sum_{k=1}^{3} \frac{\partial\left(\rho_1 v_{k} S_1\right)}{\partial x_{k}}=0.
\end{align*}
	If define
\begin{equation}\label{EQ:entropypair-1D}
		\eta(\bU) = -{\rho_1 S_1},  \quad q_k(\bU) = \eta v_k,
\end{equation} and
	$$\bV= \eta'(\bU)^{\mathrm{T}} = \left(- S_1 -
	\dfrac{|\bv|^2 }{2 T} + c_{v,1}\Gamma_{1}, \dfrac{ \bv^{\mathrm{T}}}{T},
	-\dfrac{1}{T}
	\right)^{\mathrm{T}},$$
then
	one can verify that  $\partial{\bU}/\partial{\bV}$ is symmetric positive definite, and the matrix $\pd{\bF_{k}}{\bU}\pd{\bU}{\bV}$ is symmetric,
	so that \eqref{eq:ConserLaw}
	can be symmetrized with the change of variables $\bU\to \bV$,
	and $(\eta, q_k)$ forms a convex entropy pair of \eqref{eq:ConserLaw} with $N=1$.
In this case,
  the entropy potential $\phi$ and entropy potential flux $\psi_k$ are explicitly given by
	\begin{align}\label{eq:potentialN1}
		\phi = R_1 \rho_1 -  \dfrac{p_{\infty,1}}{T}, ~~~\psi_k =\phi v_k,~~~k = 1,2,3.
	\end{align}

 If choosing the parameter vector as
$$
\bm{z}  = \left({z}_1, {z}_2, {z}_3,  {z}_4, {z}_5\right)^{\mathrm{T}} = \left(
\rho_1,  \bv^{\mathrm{T}},
1/T
\right)^{\mathrm{T}},
$$
 and using
	  the identity
	$\jump{ab} = \mean{a}\jump{b} + \mean{b}\jump{a}$,
where $\jump{a}$ and $\mean{a}$ are the jump and mean of $a$,  respectively,
then the jumps of $\bV, \phi,$ and $\psi_1$ can be rewritten
as the following linear combinations of the jumps of   $\bm{z}$
 	\begin{align}\label{eq:jump}
 \left\{\begin{array}{l}
 \jump{\bV_1} = c_{v,1} \dfrac{\jump{z_5}}{\meanln{z_5}} + R_1\dfrac{\jump{z_1}}{\meanln{z_1}} - \mean{z_5}\sum\limits_{m=2}^{4}\mean{z_m}\jump{z_m} - \dfrac{1}{2}\sum\limits_{m=2}^{4}\mean{z_m^2}\jump{z_5}, \\
 \jump{\bV_m} = \mean{z_m}\jump{z_5} +\mean{z_5} \jump{z_m},
 ~~~m = 2,3,4, \\
 \jump{\bV_5} = -\jump{z_5},\\
 \jump{\phi} =R_1\jump{z_1} -  p_{\infty, 1}\jump{z_5},\\
 \jump{\psi_1} = \left(R_1\jump{z_1} -  p_{\infty, 1}\jump{z_5}\right)\mean{z_2} + \left(R_1\mean{z_1} -  p_{\infty, 1}\mean{z_5}\right)\jump{z_2},
 \end{array}\right.
 \end{align}
where  $\meanln{a}:=\jump{a}/\jump{\ln{a}}$, $a >0$ is the logarithmic mean, see \cite{Ismail2009Affordable}.
 If substituting \eqref{eq:jump} into \eqref{eq:ECCondition_comp}
 and equating the coefficients of the same jump terms on each side of the identity \eqref{eq:ECCondition_comp}, then
	\begin{align*}
\left\{\begin{array}{l}
\dfrac{R_1}{\meanln{z_1}}\tbU_1 = R_1,\\
-\mean{z_m}\mean{z_5}\tbU_1    + \mean{z_5}\tbU_{m}    = 0, ~~~m= 2,3,4,\\
\dfrac{c_{v,1}} {\meanln{z_5}}\tbU_1 - \dfrac{1}{2}\sum\limits_{m=2}^{4}\mean{z_m^2} \tbU_1    + \sum\limits_{m=2}^{4}\mean{z_m} \tbU_{m}    - \tbU_5    =
-p_{\infty,1},
\end{array}\right.
\end{align*}
and
	\begin{align*}
\left\{\begin{array}{l}
\dfrac{R_1}{\meanln{z_1}}\tbF_{1,1}    = R_1\mean{z_2},\\
-\mean{z_2}\mean{z_5}\tbF_{1,1}     + \mean{z_5}\tbF_{1,2}    = R_1\mean{z_1} -  p_{\infty, 1}\mean{z_5}, \\
-\mean{z_m}\mean{z_5} \tbF_{1,1}     + \mean{z_5}\tbF_{1,m}    = 0, ~~~m = 3,4,\\
\dfrac{c_{v,1}} {\meanln{z_5}}\tbF_{1,1}  - \dfrac{1}{2}\sum\limits_{m=2}^{4}\mean{z_m^2}
 \tbF_{1,1}    + \sum\limits_{m=2}^{4}\mean{z_m} \tbF_{1,m}    - \tbF_{1,5}   =
-\mean{z_2}p_{\infty,1},
\end{array}\right.
\end{align*}
where $\tbU_{m}$ and $\tbF_{1,m}$ denote the $m$th component of $\tbU$ and $\tbF_{1}$ with $m = 1,\cdots, 5$, respectively.
Solving the above two systems of the linear equations yields the expressions  of $\widetilde{\bU}$ and $\widetilde{\bF}_{1}$ as follows
\begin{align*}
&\widetilde{\bU}   =
\begin{pmatrix}
\meanln{z_1}\\
\mean{z_2}\meanln{z_1}\\
\mean{z_3}\meanln{z_1}\\
\mean{z_4}\meanln{z_1}\\
\left(\dfrac{c_{v,1}}{\meanln{z_5}} - \dfrac{1}{2}\sum\limits_{m=2}^{4}\mean{z_m^2}
\right)\meanln{z_1} +
\meanln{z_1}\sum\limits_{m=2}^{4}\mean{z_m}^2  + {p_{\infty,1}}
\end{pmatrix},\\
&\tbF_1   =
\begin{pmatrix}
\meanln{z_1}\mean{z_2},\\
\mean{z_2}\tbF_{1,1}   +\dfrac{1}{\mean{z_5}} R_1\mean{z_1}-
{ p_{\infty,1}}\\
\mean{z_3}\tbF_{1,1}   \\
\mean{z_4}\tbF_{1,1}   \\
\left(\dfrac{c_{v,1}}{\meanln{z_5}} - \dfrac{1}{2}\sum\limits_{m=2}^{4}\mean{z_m^2}
\right)\tbF_{1,1}
+ \sum\limits_{m=2}^{4}\left( \mean{z_m}
\tbF_{1,m}    \right) +{p_{\infty,1} \mean{z_2} }
\end{pmatrix}.
\end{align*}
 For $k=2,3$, $\tbF_k$ may be similarly gotten.

\subsubsection{Two-component compressible Euler equations ($N = 2$)}\label{section3.3.2}
Similarly,  the smooth solutions of the two-component compressible Euler equations \eqref{eq:ConserLaw}
with
the  stiffened EOS satisfy
\begin{align*}
\frac{\partial(\rho  S)}{\partial t}+\sum_{k=1}^{3} \frac{\partial\left(\rho v_{k} S\right)}{\partial x_{k}}=0,
\end{align*}
with the entropy of the mixture
\begin{align*}
\rho S:=\sum_{\ell=1}^{2} \rho_{\ell} S_{\ell},\quad S_{\ell}:=c_{v,\ell} \ln (T)-R_{\ell} \ln \left(\rho_{\ell}\right),\end{align*}
where $S_{\ell}$ is the thermodynamic entropy of species $\ell$.
With the help of the thermodynamic entropy,
 the mathematical entropy pair of \eqref{eq:ConserLaw}  may be defined by
\begin{align}\label{EQ:entropypair}
\eta(\bU) = -{\rho S}, \quad q_k(\bU) = \eta v_k, \quad k = 1,2,3.
\end{align}
Because
	for $\rho_1>0, \rho_2>0, T>0,$    $\pd{\bU}{\bV}$ is  symmetric positive definite, and  $\pd{\bF_{k}}{\bU}\pd{\bU}{\bV}$ is symmetric,
	so that the equations \eqref{eq:ConserLaw}
	can be symmetrized with   $\eta(\bU), q_k(\bU)$. 
In this case, the entropy potential $\phi$ and the entropy potential flux $\psi_k$ can be explicitly given by
\begin{align}\label{eq:potentialN2}
\phi =\sum\limits_{\ell  = 1}^{ 2} \left(R_\ell \rho_\ell -  \dfrac{p_{\infty,\ell}}{T} \right), \quad \psi_k =\phi v_k,
\end{align}
with the entropy {variables}
$\bV = \left(- S_1 -
\dfrac{|\bv|^2 }{2 T} + c_{v,1}\Gamma_{1}, - S_2 -
\dfrac{|\bv|^2 }{2 T} + c_{v,2}\Gamma_{2}, \dfrac{ \bv^{\mathrm{T}}}{T},
-\dfrac{1}{T}
\right)^{\mathrm{T}}$.

If choosing the parameter vector
$
\bm{z}  = \left({z}_1, {z}_2, {z}_3,  {z}_4, {z}_5, z_6\right)^{\mathrm{T}} = \left(
\rho_1, \rho_2,  \bv^{\mathrm{T}},
1/T
\right)^{\mathrm{T}}$, then
  the jumps of the entropy variables $\bV$,
  the  entropy potential $\phi$ and the entropy potential flux $\psi_1$ can be rewritten as
\begin{align}\label{eq:jump_MC}
\left\{\begin{array}{l}
\jump{\bV_1} = c_{v,1} \dfrac{\jump{z_6}}{\meanln{z_6}} + R_1\dfrac{\jump{z_1}}{\meanln{z_1}} - \mean{z_6}\sum\limits_{m=3}^{5}\mean{z_m}\jump{z_m} - \dfrac{1}{2}\sum\limits_{m=3}^{5}\mean{z_m^2}\jump{z_6}, \\
\jump{\bV_2} = c_{v,2} \dfrac{\jump{z_6}}{\meanln{z_6}} + R_2\dfrac{\jump{z_2}}{\meanln{z_2}} - \mean{z_6}\sum\limits_{m=3}^{5}\mean{z_m}\jump{z_m} - \dfrac{1}{2}\sum\limits_{m=3}^{5}\mean{z_m^2}\jump{z_6}, \\
\jump{\bV_m} = \mean{z_m}\jump{z_6} +\mean{z_6} \jump{z_m},
~~~m = 3,4,5, \\
\jump{\bV_6} = -\jump{z_6},\\
\jump{\phi} = \sum\limits_{\ell=1}^{2} \left(R_\ell\jump{z_\ell} -  p_{\infty, \ell}\jump{z_6}\right),\\
\jump{\psi_1} = \sum\limits_{\ell=1}^{2} \left(R_\ell\jump{z_\ell} -  p_{\infty, \ell}\jump{z_6}\right)\mean{z_3} + \sum\limits_{\ell=1}^{2} \left(R_\ell\mean{z_\ell} -  p_{\infty, \ell}\mean{z_6}\right)\jump{z_3}.
\end{array}\right.
\end{align}
 Substituting it into \eqref{eq:ECCondition_comp} gives
\begin{align*}
\left\{\begin{array}{l}
\dfrac{R_1}{\meanln{z_1}}\tbU_1 = R_1,\\
\dfrac{R_2}{\meanln{z_2}}\tbU_2    = R_2,\\
-\mean{z_m}\mean{z_6}\sum\limits_{\ell  = 1}^{2}\tbU_\ell    + \mean{z_6}\tbU_{m}    = 0, ~~~m= 3,4,5,\\
\sum\limits_{\ell=1}^{2}\left(\dfrac{c_{v,\ell}} {\meanln{z_6}}\tbU_\ell   \right)  - \dfrac{1}{2}\sum\limits_{m=3}^{5}\mean{z_m^2} \sum\limits_{\ell  = 1}^{2}\tbU_\ell    + \sum\limits_{m=3}^{5}\mean{z_m} \tbU_{m}    - \tbU_6    =
-\sum\limits_{\ell=1}^{2}p_{\infty,\ell},
\end{array}\right.
\end{align*}
and
\begin{align*}
\left\{\begin{array}{l}
\dfrac{R_1}{\meanln{z_1}}\tbF_{1,1}    = R_1\mean{z_3},\\
\dfrac{R_2}{\meanln{z_2}}\tbF_{1,2}    = R_2\mean{z_3},\\
-\mean{z_3}\mean{z_6}\sum\limits_{\ell=1}^{2} \tbF_{1,\ell}     + \mean{z_6}\tbF_{1,3}    =\sum\limits_{\ell=1}^{2} \left( R_\ell\mean{z_\ell} -  p_{\infty, \ell}\mean{z_6}\right), \\
-\mean{z_m}\mean{z_6}\sum\limits_{\ell=1}^{2} \tbF_{1,\ell}     + \mean{z_6}\tbF_{1,m}    = 0, ~~~m = 4,5,\\
\sum\limits_{\ell=1}^{2}\left(\dfrac{c_{v,\ell}} {\meanln{z_6}}\tbF_{1,\ell}   \right)  - \dfrac{1}{2}\sum\limits_{m=3}^{5}\mean{z_m^2}
\sum\limits_{\ell=1}^{2} \tbF_{1,\ell}    + \sum\limits_{m=3}^{5}\mean{z_m} \tbF_{1,m}    - \tbF_{1,6}   =
-\mean{z_3}\sum\limits_{\ell=1}^{2}p_{\infty,\ell},
\end{array}\right.
\end{align*}
where  $\tbF_{1,m}   $ and $\tbU_{m}   $ denote the $m$th component of $\tbF_{1}   $ and  $\tbU   $ with $ m=  1,\cdots, 6,$ respectively.
Solving those linear {systems} gives
$\widetilde{\bU}   $ and $\widetilde{\bF}_1   $ as follows
\begin{align*}
&\widetilde{\bU}   =
\begin{pmatrix}
\meanln{z_1}\\\meanln{z_2}\\
\mean{z_3}\sum\limits_{\ell=1}^{2}\meanln{z_\ell}\\
\mean{z_4}\sum\limits_{\ell=1}^{2}\meanln{z_\ell}\\
\mean{z_5}\sum\limits_{\ell=1}^{2}\meanln{z_\ell}\\
\sum\limits_{\ell=1}^{2}\left[\meanln{z_\ell}
\left(\dfrac{c_{v,\ell}}{\meanln{z_6}} - \dfrac{1}{2}\sum\limits_{m=3}^{5}\mean{z_m^2}
\right)\right]
+
 \sum\limits_{\ell=1}^{2}\meanln{z_\ell}\sum\limits_{m=3}^{5}\mean{z_m}^2  + {\sum\limits_{\ell  = 1}^{2}p_{\infty,\ell}}
\end{pmatrix},\\
&\tbF_1   =
\begin{pmatrix}
\meanln{z_1}\mean{z_3},\\
\meanln{z_2}\mean{z_3},\\
\mean{z_3}\sum\limits_{\ell=1}^{2}\tbF_{1,\ell}   +\dfrac{1}{\mean{z_6}}\left(\sum\limits_{\ell=1}^{2} R_\ell\mean{z_\ell}\right) -
{ \sum\limits_{\ell  = 1}^{2}p_{\infty,\ell}}\\
\mean{z_4}\sum\limits_{\ell=1}^{2}\tbF_{1,\ell}   \\
\mean{z_5}\sum\limits_{\ell=1}^{2}\tbF_{1,\ell}   \\
\sum\limits_{\ell=1}^{2}
\left[\tbF_{1,\ell}
\left(\dfrac{c_{v,\ell}}{\meanln{z_6}} - \dfrac{1}{2}\sum\limits_{m=3}^{5}\mean{z_m^2}
\right)\right]
+ \sum\limits_{m=3}^{5}\left( \mean{z_m}
\tbF_{1,m}    \right) +{\sum\limits_{\ell = 1}^{2}p_{\infty,\ell} \mean{z_3} }
\end{pmatrix}.
\end{align*}
 For $k=2,3$, $\tbF_k$ may be similarly derived.
For the entropy pair $(\eta,q_k)$ in  \eqref{EQ:entropypair},
$\tbF_k$  can also be obtained by choosing respectively  the angles
$\varphi=0,\theta=\pi/2$ and $\varphi=\pi/2,\theta=0$
in
%
%
%
	\begin{align}\label{eq:rotation}
\tbF_{\varphi,\theta}({\bU_l}, {\bU_r}):= \bT^{-1} \tbF_{1}
	\left(\thU_l,\thU_r\right),
	\end{align}
which is the EC flux approximating the flux $\cos\varphi\cos\theta \bF_{1}({\bU})
+ \cos\varphi\sin\theta \bF_{2}({\bU}) + \sin\varphi \bF_{3}({\bU})$,
	where  $\varphi\in[0,2\pi)$, $\theta\in[0,\pi]$, $\thU_l := \bT\bU_l$,  $\thU_r := \bT\bU_r$,
 $\bT$ is  the expanded rotational matrix defined by
	\begin{align}\label{eq:rotamatrix}
	&\bT =
	\begin{bmatrix}
	\bm{I}_N & 0                      & 0                      & 0           & 0 \\
	0 & \cos\varphi\cos\theta  & \cos\varphi\sin\theta  & \sin\varphi & 0 \\
	0 & -\sin\theta            & \cos\theta             & 0           & 0 \\
	0 & -\sin\varphi\cos\theta & -\sin\varphi\sin\theta & \cos\varphi & 0 \\
	0 & 0                      & 0                      & 0           & 1
	\end{bmatrix},
	\end{align}
with	 the unit $N\times N$ matrix $\bm{I}_N$, $N = 1, 2$.
In fact,   $\tbF_{\varphi,\theta}({\bU_l}, {\bU_r})$ satisfies
the sufficient condition of the EC flux
$$
\jump{\bV(\bU)} ^\mathrm{T} \tbF_{\varphi,\theta}({\bU_l}, {\bU_r})
=
\cos\varphi\cos\theta \jump{\psi_{1}(\bU)} + \cos\varphi\sin\theta\jump{\psi_{2}(\bU)} + \sin\varphi
\jump{\psi_{3}(\bU)},
$$
because
	\begin{align*}
	&\begin{aligned}
	&\thU := \bT\bU= \left(\rho_{1}, \cdots, \rho_{N}, \hat{v}_{1}, \hat{v}_{2}, \hat{v}_{3}, E\right)^{\mathrm{T}},\quad\hat{v}_1=\cos\varphi\cos\theta v_1+ \cos\varphi\sin\theta v_2 + \sin\varphi v_3,\\
	&\hat{v}_2= -\sin\theta v_1+\cos \theta v_2, \quad \hat{v}_3 =  -\sin\varphi\cos\theta v_1 -\sin\varphi\sin\theta v_2 + \cos\varphi v_3,
	\end{aligned}
\\
&	{\bV(\thU)} = {\left(- S_1 -
	\dfrac{|\bv|^2 }{2 T} + c_{v,1}\Gamma_{1},\cdots,  - S_N -
	\dfrac{|\bv|^2 }{2 T} + c_{v,N}\Gamma_{N}, ~\dfrac{ \widehat{\bm{v}}}{T},
	-\dfrac{1}{T}
	\right)^{\mathrm{T}}},
	\end{align*}
so that one has
${\bV(\thU)} ={\bV(\bT\bU)} = {\bT \bV(\bU)}$
and
%
%
	\begin{align*}
	\jump{\bV(\bU)} ^\mathrm{T} \tbF_{\varphi,\theta}({\bU_l}, {\bU_r})&= \jump{\bV({\bU})}^\mathrm{T}\bT^{-1}\widetilde{\bF}_1(\widehat{\bU}_l, \widehat{\bU}_r)
=\jump{\bV(\thU)} ^\mathrm{T}{\tbF}_1({\thU_l}, {\thU_r})
\\
&    \xlongequal[]{\eqref{eq:ECCondition_comp}}
\jump{\psi_1(\widehat{\bU})}=\jump{\sum\limits_{\ell  = 1}^{ N} \left(R_\ell \rho_\ell -  \dfrac{p_{\infty,\ell}}{T} \right) \hat{v}_1}\\
	& =  \cos\varphi\cos\theta \jump{\psi_1({\bU})}+ \cos\varphi\sin\theta \jump{\psi_2({\bU})} + \sin\varphi \jump{\psi_3({\bU})}. 
%
	\end{align*}

\subsection{ES schemes}\label{sec:es_scheme}
It is known that 
the EC schemes   work well for the smooth solutions, but
they may produce severe nonphysical oscillations if  the solutions contain discontinuities. In order to suppress those numerical oscillations,  a suitable dissipation term should be added to the EC flux \eqref{eq:ECFlux_curv_2p} to make the schemes satisfy the semi-discrete entropy inequality for the given entropy pair. Similar to \cite{duan2021highorder},  the high-order accurate ES flux may be  given by
\begin{equation}\label{eq:HOstable}
\left({{\thF}_k}\right)_{\bm{i}, k,\pm\frac12}
=\left({\thF_k}\right)_{\bm{i}, k,\pm\frac12}^{2w{\rm th}}
-\dfrac12\left( \bT^{-1}\bm{R}(\bT\bU)\left|\widetilde{\bm{\Lambda}}(\bT\bU)\right|
\right)_{\bm{i}, k,\pm\frac12}\bm{Y}_{\bm{i},k,\pm\frac{1}{2}}
\jump{\bm{W}}_{\bm{i}, k,\pm\frac12}^{\tt {WENOMR}},
\end{equation}
where
$\left|\widetilde{\bm{\Lambda}}(\bT\bU)\right|:=\max\limits_{m}\left\{ \left| J\pd{\xi_k}{t}+L_k\lambda_m(\bT\bU) \right|\right\}\bm{I},
$
  $L_{k}=\sqrt{\sum\limits_{j=1}^{3}\left(J \frac{\partial \xi_{k}}{\partial x_{j}}\right)^{2}}$,
  the rotational matrix $\bT$ is given in \eqref{eq:rotamatrix}
with  
\begin{align*}
&\theta = \arctan\left(\left(J\pd{\xi_k}{x_2}\right)\Big/\left(J\pd{\xi_k}{x_1}\right)\right),\\
&\varphi = \arctan\left(\left(J\pd{\xi_k}{x_3}\right)\Bigg/\sqrt{\left(J\pd{\xi_k}{x_1}\right)^2
+\left(J\pd{\xi_k}{x_2}\right)^2}\right),
\end{align*}
and $\bm{R}$ is the scaled right eigenvector matrix  satisfying
\begin{equation}\label{eq:eigen}
\pd{\bU}{\bV}=\bm{R}\bm{R}^\mathrm{T},\quad
\pd{\bF_1}{\bU}=\bm{R}\bm{\Lambda}\bm{R}^{-1},
\end{equation}
here  $\bm{\Lambda}$ is the diagonal matrix, whose diagonal elements are  the eigenvalues of the matrix  $\pd{\bF_1}{\bU}$.
	The high-order accurate jump terms in \eqref{eq:HOstable} are defined as
$\jump{\bm{W}}_{\bm{i}, k, \pm\frac12}^{\tt {WENOMR}} := \bm{W}_{\bm{i},k,\pm\frac12}^{\tt {WENOMR},+} -
\bm{W}_{\bm{i},k,\pm\frac12}^{\tt {WENOMR},-}
$ with the left and right limit values $\bm{W}_{\bm{i},k,\pm\frac12}^{\tt {WENOMR},-}$ and  $\bm{W}_{\bm{i},k,\pm\frac12}^{\tt {WENOMR},+}$ obtained by  the high-order multi-resolution WENO reconstruction   \cite{WANG2021105138}.
The diagonal matrix $\bm{Y}_{\bm{i},k,\pm\frac12}$ is  chosen as
\begin{align*}
&\left(\bm{Y}_{m,m}\right)_{\bm{i},k,\pm\frac12} =
\left\{\begin{array}{ll}1, & \text{if} ~ {\rm{sign}}\left(\jump{\bm{W}_m}_{\bm{i}, k, \pm\frac12}^{\tt {WENOMR}} \right)
	{\rm{sign}}\left(
	\jump{\bm{W}_m}_{\bm{i}, k, \pm\frac12}
	\right)
	>0, \\ 0, & \text {otherwise},\end{array}\right.
\end{align*}
in order to ensure  the ``sign''  property \cite{Biswas2018Low}.

\begin{rmk}\rm
Our computations will take  the fifth-order multi-resolution WENO reconstruction, which
	  uses unequal-sized stencils  and arbitrary positive linear weights whose sum is one \cite{WANG2021105138}, see \ref{section:WENO}. The ES adaptive moving mesh schemes based on the multi-resolution WENO reconstruction  can better capture the localized structures for the (multi-component) flow problems, and outperform their  counterparts based on the classical WENO reconstruction \cite{Jiang1996Efficient}  with a slight increase in the computational cost, see Section \ref{section:Result}.
\end{rmk}

Before ending this section, we give the scaled eigenvector matrix $\bm{R}$
for the multi-component compressible Euler equations ($N = 1,2$) with the stiffened EOS.

For the single-component compressible Euler equations ($N = 1$),
the diagonal matrix $\bm{\Lambda}  $ is given by
$$\bm{\Lambda} =\mbox{diag}\{\lambda_1,\dots,\lambda_5\}  =\mbox{diag}\{v_1 -c_s,v_1, v_1, v_1, v_1+c_s\},$$
where $c_s$ is the speed of sound given by $c_s^2 = R_{1}\Gamma_{1}T $. After some algebraic manipulations,
the scaled eigenvector matrix $\bm{R}$ can be expressed as
\begin{align*}
\begin{bmatrix}
1 & 1                     & 0                      & 0           & 1 \\
v_1-c_s& v_1 & 0& 0  & v_1+c_s\\
v_2 & v_2 & 1& 0 & v_2 \\
v_3      & v_3            & 0           & 1 & v_3\\
H-c_sv_1 & \frac12|\bv|^2     & v_2     & v_3   & H+c_sv_1
\end{bmatrix}
\begin{bmatrix}
\dfrac{\rho_1}{2\Gamma_1R_1}& 0                   & 0                      & 0           & 0 \\
0& \dfrac{\rho_1}{c_{v,1}\Gamma_{1}}  & 0& 0  & 0\\
0 & 0 & \rho_1 T& 0 & 0\\
0& 0          & 0           &  \rho_1 T& 0\\
0& 0     & 0 & 0  &  \dfrac{\rho_1}{2\Gamma_1R_1}
\end{bmatrix}^{\frac12}
,
\end{align*}
where $ H  = \left(E+p\right)/{\rho_1}$ is the total enthalpy.
In practice, the values  of  $\bm{R}_{\bm{i}, k, +\frac{1}{2}}$  and $\left|\widetilde{\bm{\Lambda}}\right|_{\bm{i}, k,+\frac12}$
are calculated by using some ``averaged" values of the  primitive variables as follows
\begin{align*}
\overline{\rho_1} = \meanln{\rho_1}_{\bm{i},k,+\frac{1}{2}}, ~ \overline{\bv} = \mean{\bv}_{\bm{i},k,+\frac{1}{2}}, ~ \overline{p} +p_{\infty,1} = \left(\dfrac{ \meanln{\rho_1}}{ \meanln{\rho_1/(p+p_{\infty,1})}}\right)_{\bm{i},k,+\frac{1}{2}}.
\end{align*}

For the two-component compressible Euler equations ($N = 2$) with
the stiffened EOS, utilizing the  similar procedure for the ideal EOS in \cite{GOUASMI2020112912}, the diagonal matrix  is given by
$$\bm{\Lambda} =\mbox{diag}\{\lambda_1,\dots,\lambda_6\}   =\mbox{diag}\{v_1 -c_s,v_1, v_1, v_1, v_1, v_1+c_s\},$$
where $c_s ^2= R \Gamma T$ with  $R: = \left(\sum \left(\Gamma_{\ell}-1\right)c_{v,\ell}\rho_{\ell}\right)/\sum \rho_{\ell}$.
	The scaled eigenvector matrix $\bm{R}$ can be obtained by scaling the right eigenvectors $\widetilde{\bm{R}}$ using
	a symmetric block diagonal matrix $\bm{D}$, i.e. $\bm{R} = \widetilde{\bm{R}}\bm{D}$,
	where $\widetilde{\bm{R}}$ is
	\begin{align*}
		\begin{bmatrix}
			Y_1& 1                  &0              & 0       & 0           & Y_1 \\
			Y_2&  0                  &1              &0        &0                                    & Y_2\\
			v_1-c_s                              & v_1                &v_1         & 0       & 0  & v_1+c_s\\
			v_2                                      & v_2                &v_2           & c_s   & 0 & v_2 \\
			v_3                                     & v_3                 &v_3            & 0        & c_s & v_3\\
			H-c_sv_1                          & \frac12|\bv|^2  - \dfrac{d_1}{\Gamma-1}  &  \frac12|\bv|^2  -  \dfrac{d_2}{\Gamma-1} & c_sv_2     & c_sv_3   & H+c_sv_1
		\end{bmatrix},
	\end{align*}
where  $d_\ell = h_\ell - \Gamma c_{v,\ell}T$, and $ H  = {\sum \rho_\ell h_\ell}/{\rho} + \frac12 |\bv|^2 $ with $h_\ell := c_{v,\ell}T + R_\ell T$. 	
According to \eqref{eq:eigen},  the explicit expression for {the} matrix $\bm{D}\bm{D}^{\mathrm{T}}$ is provided by
\begin{align*}
&\bm{D}\bm{D}^{\mathrm{T}}=  \widetilde{\bm{R}}^{-1} \dfrac{\partial\bU}{\partial\bV}\widetilde{\bm{R}}^{-\mathrm{T}} = \dfrac{\rho}{\Gamma R}\diag\left(1/2,\bm{D}^{2Y}, 1,1, 1/2\right),\\
&\bm{D}^{2 Y}:=Y_1Y_2\left[\begin{array}{cc}
(\Gamma-1) Y_{1}/Y_2+\left(\Gamma R_{2} / R_{1}\right) & -1 \\
-1 & (\Gamma-1) Y_{2}/Y_{1}+\left(\Gamma R_{1} / R_{2}\right)
\end{array}\right],
\end{align*}
with $Y_\ell = {\rho_\ell}/{\rho}$.
If decomposing $\bm{D}^{2Y}$ as
$$
\bm{D}^{2 Y}=\bm{D}^{Y}\left(\bm{D}^{Y}\right)^{T},\quad \bm{D}^{Y}:=
\sqrt{Y_{1} Y_{2}}
\left[\begin{array}{cc}
- \sqrt{\Gamma R_{2} / R_{1}} & \sqrt{(\Gamma-1)Y_1/Y_2} \\
 \sqrt{\Gamma R_{1} / R_{2}} &  \sqrt{(\Gamma-1)Y_2/Y_1}
\end{array}\right],
$$
then $\bm{D} = \dfrac{\sqrt{\rho}}{\sqrt{\Gamma R}}\diag\left(1/\sqrt{2},\bm{D}^{Y}, 1,1, 1/\sqrt{2}\right)$.
Thus the scaled eigenvector matrix $\bm{R}$ is obtained by
\begin{align*}
\widetilde{\bm{R}} \times
\sqrt{\dfrac{\rho}{\Gamma{R}}}
\begin{bmatrix}
\dfrac{1}{\sqrt{2}}& 0        &0           & 0                      & 0           & 0 \\
0& -\sqrt{Y_1Y_2}\sqrt{\Gamma R_2/R_1} & Y_1\sqrt{\Gamma-1}& 0  &0 & 0\\
0 & \sqrt{Y_1Y_2}\sqrt{\Gamma R_1/R_2}         & Y_2\sqrt{\Gamma-1}& 0     & 0 &0\\
0& 0          & 0 &1     & 0 & 0\\
0& 0          & 0  &0     & 1 & 0\\
0& 0           & 0 & 0      &0 &  \dfrac{1}{\sqrt{2}}
\end{bmatrix}.
\end{align*}
Similarly, the values of $\bm{R}_{\bm{i}, k,+\frac12}$  and $\left|\widetilde{\bm{\Lambda}}\right|_{\bm{i}, k,+\frac12}$ are calculated
by using some averaged values at $\{\bm{i},k,+\frac12\}$
\begin{align*}
\overline{\rho_1} &= \meanln{\rho_1},~ \overline{\rho_2} = \meanln{\rho_2},
\overline{\rho} = \meanln{\rho}, ~ \overline{\bv} = \mean{\bv}, ~ \overline{T} = 1/\meanln{1/T},\\
\overline{R} &= \mean{R},~ \overline{\Gamma} = \mean{\Gamma},~\sum_{\ell=1}^{2}\overline{ \rho_\ell h_\ell} = \sum_{\ell=1}^{2} \left({c_{v,\ell} \overline{\rho_\ell}}{\overline{T}} + p_{\infty,
	\ell} \right) + \mean{p}.
\end{align*}

%% file: MovingMesh.tex
\section{Adaptive moving mesh strategy}\label{section:MM}
This section  introduces the adaptive moving mesh strategy \cite{duan2021highorder},
 but will omit {the} dependence of the variables on $t$ for convenience, unless otherwise stated.

The mesh is moved adaptively, which is equivalent to finding
   the coordinate transformation $\bx=\bx(\bm{\xi})$ from $\Omega_c$ to   $\Omega_p$ mentioned in Section \ref{section: EntCon}.
Such transformation
can be obtained by
solving the mesh {redistribution equations}
\begin{equation}\label{eq:mesh_EL}
\nabla_{\bm{\xi}}\cdot\left(\bm{G}_k\nabla_{\bm{\xi}}x_k\right)=0,
~\bm{\xi}\in\Omega_c,~k=1,2,3,
\end{equation}
which may be the Euler-Lagrange equations or the stationary variational conditions
for  minimizing the  mesh adaptation functional, 
where $\bm{G}_k$ is the symmetric positive definite matrix  depending on the solutions
	of the underlying governing equations \eqref{eq:ConserLaw} or their derivatives.
The simplest choice of $\bm{G}_k$ is
\begin{align*}
  \bm{G}_k={\Theta}\bm{I}_3,
\end{align*}
where the monitor function  {$\Theta$}  is   positive  and taken in this paper as
\begin{align}\label{eq:monitor}
\Theta=\Big({1+\sum_{k=1}^\kappa\alpha_k\Big(
\dfrac{\abs{\nabla_{\bm{\xi}}\sigma_k}}
{\max\abs{\nabla_{\bm{\xi}}\sigma_k}}\Big)^2}\Big)^{1/2},
\end{align}
here $\sigma_k$ is a physical variable, $\alpha_k$ is {a} non-negative parameter, and $\kappa$ is the number of the chosen physical {variables}.
 Using the {second-order accurate} central difference scheme and the Jacobi iteration,
the mesh equations \eqref{eq:mesh_EL} are approximated by
\begin{equation*}
\begin{aligned}
&\sum\limits_{k =1}^{d} \left[\Theta_{\bm{i}, k, +\frac12 }\left(\bx_{\bm{i}, k,1}^{[\nu]}-\bx_{\bm{i}}^{[\nu+1]}\right)
-\Theta_{\bm{i}, k, -\frac12}\left(\bx_{\bm{i}}^{[\nu+1]}-\bx_{\bm{i}, k, -1}^{[\nu]}\right)\right] {/\Delta \xi_k^2}=0, \ \nu=0,1,\cdots,\mu,
\end{aligned}
\end{equation*}
where ${\bx}^{[0]}_{\bm{i}}
  :=\bx^n_{\bm{i}}$, and
  $\Theta_{\bm{i}, k, \pm\frac12 } :=\frac12\left(
 \Theta_{\bm{i},k}+\Theta_{\bm{i},k,\pm1}
 \right).$
The total iteration number $\mu$ is taken as $10$ in our numerical tests.
The final adaptive mesh  is given by
$
\bx^{n+1}_{\bm{i}}
:=\bx^n_{\bm{i}}
+  
{\Delta_\tau}   (\delta_\tau{\bx})^{n}_{\bm{i}},
%
%
%
$
where
$
 (\delta_\tau {\bx})^{n}_{\bm{i}}:=
   {\bx}^{[\mu]}_{\bm{i}}
  -\bx^n_{\bm{i}},
$
and the parameter ${\Delta_\tau}$
 is  the limiter of the movement of mesh points satisfying
\begin{equation*}
  {\Delta_\tau}\leqslant
  \begin{cases}
  	-\frac{1}{2(\delta_\tau{x_k})_{\bm{i}}^n}\left[(x_k)^n_{\bm{i}}-(x_k)^n_{\bm{i}, k, -1}\right],
  	~ (\delta_\tau{x_k})_{\bm{i}}^n<0, \\
   \quad\frac{1}{2(\delta_\tau{x_k})_{\bm{i}}^n}\left[(x_k)^n_{\bm{i}, k, +1}-(x_k)^n_{\bm{i}}\right],
  		~ (\delta_\tau{x_k})_{\bm{i}}^n>0,
  \end{cases}
  k = 1,\cdots, d.
\end{equation*}
Finally, the mesh velocity at $t=t_n$ in \eqref{eq:VCLCoeff} is determined by
$
\dot{\bx}^{n}_{\bm{i}}:=
{\Delta_\tau}   (\delta_\tau{\bx})^{n}_{\bm{i}}/\Delta t_n
$
with the time stepsize $\Delta t_n$, obtained by \eqref{eq:cfl} in Section \ref{section:Result}.

\begin{rmk}\rm
	In order to weaken the singularity of {the} monitor function near the strong discontinuity, it is useful to apply the following low pass filter
  \begin{align*}
  \Theta_{i_1,i_2,i_3}\leftarrow&\sum_{j_1,j_2,j_3=0,\pm 1}\left(\dfrac{1}{2}\right)^{\abs{j_1}+\abs{j_2}+\abs{j_3}+3}
  \Theta_{i_1+j_1,i_2+j_2,i_3+j_3},
  \end{align*}
  to smooth {the} monitor function $3\sim 10$ times.
\end{rmk}

%% file: NumTests.tex
\section{Numerical results}\label{section:Result}
    This section conducts several 2D and 3D numerical tests to validate the accuracy and the ability
    in capturing the localized structures of   {the} previous fifth-order adaptive moving mesh methods. The fully-discrete schemes
	are derived by using the third-order accurate explicit SSP
	RK time discretization \cite{Gottlieb2001Strong}, and   implemented in parallel with  the MPI parts of the PLUTO code \cite{2007PLUTO}. All computations are performed on the CPU nodes of the High-performance Computing Platform of Peking University (Linux Redhat environment, two Intel Xeon E5-2697A V4  per node, and core frequency of 2.6GHz).
	The time stepsize $\Delta t_{n}$ is determined by the CFL condition
	\begin{align} \label{eq:cfl}
	\Delta t_{n} \leqslant
\frac{C_{\text{\tiny \tt CFL}}}{\max\limits_{\bm{i}}\left\{\sum_{k =1}^{3}\varrho_{\bm{i}, k}^{n} / \Delta \xi_{k} \right\}},
	\end{align}
but it will be taken as $C_{\text{\tiny \tt CFL}} (\min\Delta\xi_{k})^{5/3}$ in all accuracy tests in order to make the spatial error dominant,
	where $\varrho_{\bm{i}, k}$  is the spectral radius of the eigen-matrix in the $\xi_k$-direction, and
$C_{\text{\tiny \tt CFL}}$ is taken as 0.4 and 0.3 in 2D and 3D examples, respectively,  unless otherwise stated.
For the sake of convenience,
the {fully-discrete} fifth-order finite difference  schemes with the multi-resolution WENO reconstruction and the ES fluxes on the uniform   and   moving meshes are denoted by
 ``{\tt UM-WENOMR}" and ``{\tt MM-WENOMR}", respectively.
For a comparison,  the counterparts of ``{\tt MM-WENOMR}" with the classical WENO reconstruction \cite{Jiang1996Efficient} denoted by ``{\tt MM-WENOJS}" are also implemented.

	\subsection{Single-component compressible Euler equations $(N = 1)$}
	This section {considers} some numerical experiments  on the
2D and 3D single-component compressible Euler equations $(N = 1)$ with the ideal and stiffened EOS.
The adiabatic index $\Gamma_1$ and the specific heat at constant volume $c_{v,1}$ are respectively
 taken as $1.4$ and 1, unless otherwise stated.
	
	\begin{example}[2D isentropic vortex]\label{ex:2DVortex}\rm
This example is to test the accuracy of {\tt MM-WENOMR} by solving	the problem describing an isentropic vortex propagating periodically at a constant speed in the 2D domain $\Omega_p=[-10, 10]^{2}$.
Initially, the domain
is  divided into uniformly $N_1\times N_1$  rectangular cells,
and
the vortex perturbation
	$$
	(\delta v_1, \delta v_2)=\frac{\epsilon}{2 \pi} e^{0.5\left(1-r^{2}\right)}((x_2+2), -(x_1+2)), \quad \delta
 T=-\frac{ \epsilon^{2}}{8 \Gamma_1 \pi^{2}} e^{\left(1-r^{2}\right)}, \quad \delta S_1=0,
	$$
is added to the mean flow
$\rho_1=1$, $T=1$ and $(v_1, v_2)=(1,1)$,
	where $p_{\infty,1} = 0$, $r^{2}=(x_1+2)^{2}+(x_2+2)^{2}$ with the vortex strength  $\epsilon=5$.
The monitor function is chosen as
	$$
	\Theta=\left({1+\frac{20 |\nabla_{\bm{\xi}} \rho_1|}{\max |\nabla_{\bm{\xi}} \rho_1|}+\frac{10 |\Delta_{\bm{\xi}} \rho_1|}{\max |\Delta_{\bm{\xi}} \rho_1|}}\right)^{1/2}.
	$$

Figure \ref{fig:2DVortex} presents the adaptive mesh  with $N_1 = 80$ and the density contours ($10$ equally spaced contour lines) at $t=0,2,4$.
Figure \ref{fig:2DVortex2} 
 shows the $\ell^1$- and $\ell^\infty$-errors in $\rho_1$ at $t = 4$ versus $N_1$, the orders of convergence  and  the {discrete} total entropy $\sum_{i_{1}, i_{2}} J_{i_{1}, i_{2}} \eta\left(\boldsymbol{U}_{i_{1}, i_{2}}\right) \Delta \xi_{1} \Delta \xi_{2}$ with respect to time
 by	 using the EC   and ES schemes with $N_1=160$, respectively.
It is easy to see that
the  mesh points well and adaptively concentrate near the large gradient area of the density  as expected,
  {\tt MM-WENOMR} gets the fifth-order accuracy, and the
	EC scheme almost keeps the total entropy  conservative, while the total entropy of the ES  scheme decays in time.


	\begin{figure}[!ht]
	\centering
	\begin{subfigure}[b]{0.31\textwidth}
		\centering
		\includegraphics[width=1.0\textwidth]{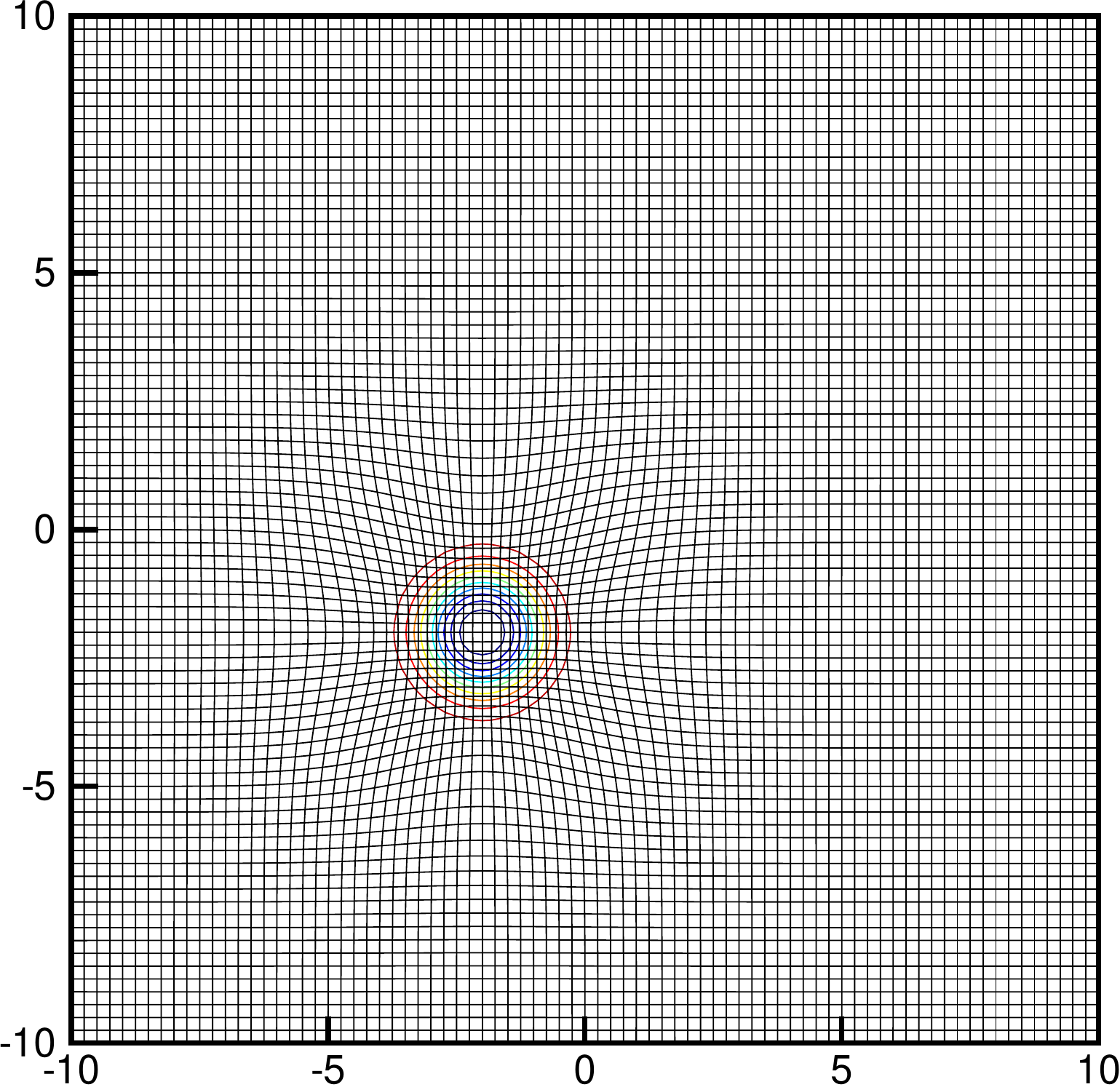}
		\caption{$t=0$}
	\end{subfigure}
	\begin{subfigure}[b]{0.31\textwidth}
		\centering
		\includegraphics[width=1.0\textwidth]{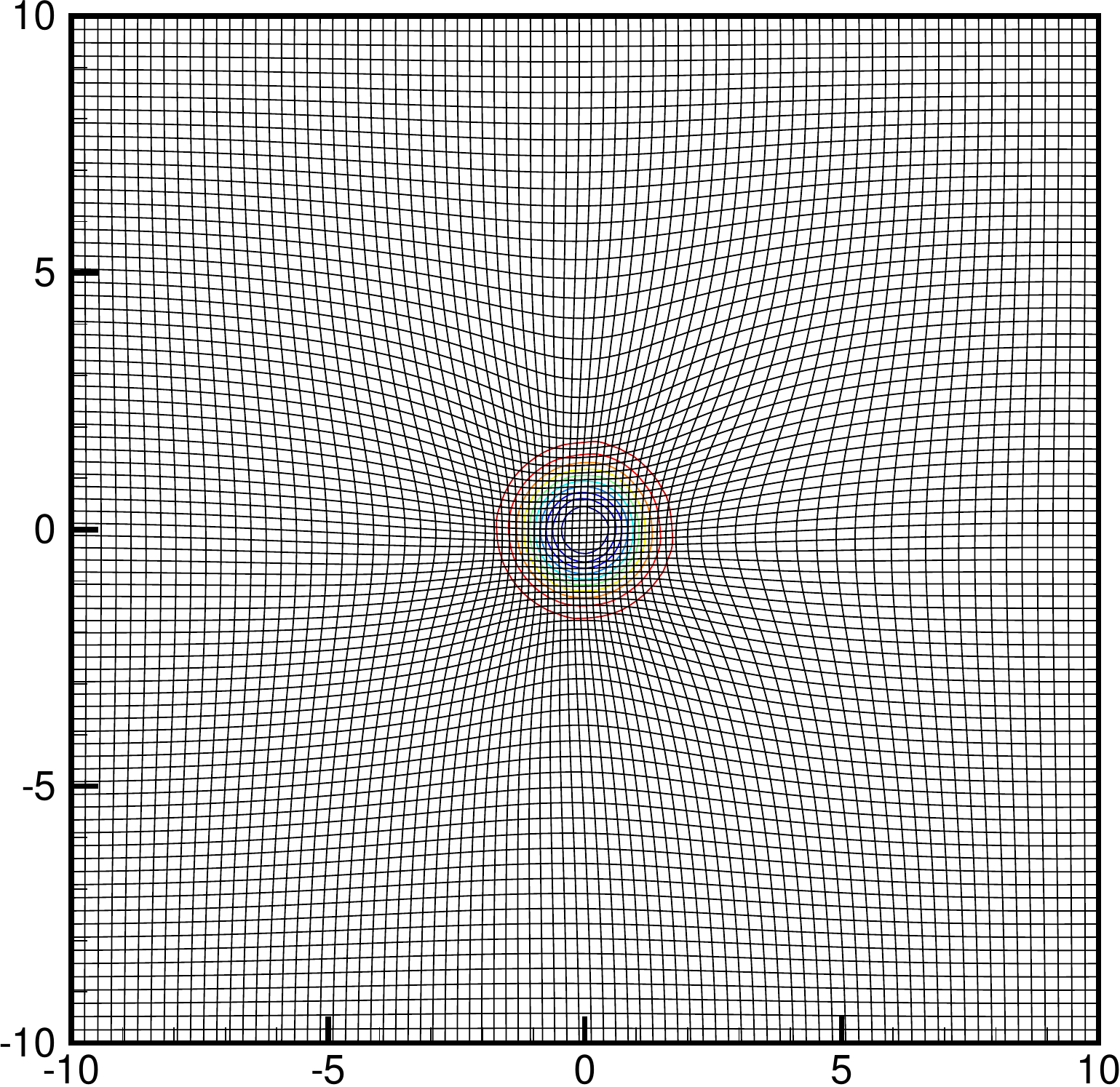}
		\caption{$t=2$}
	\end{subfigure}
	\begin{subfigure}[b]{0.31\textwidth}
		\centering
		\includegraphics[width=1.0\textwidth]{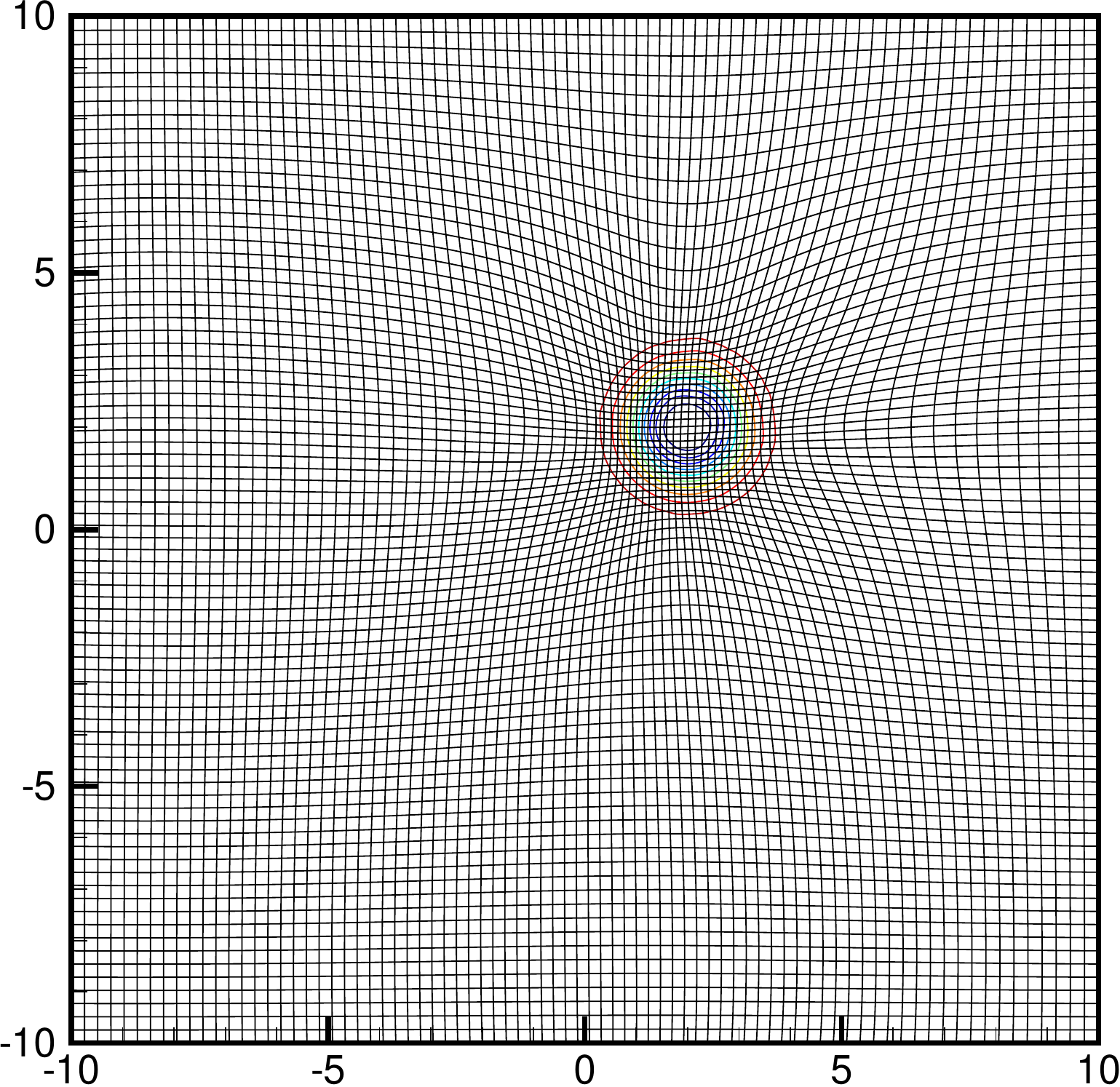}
		\caption{$t=4$}
	\end{subfigure}
	\caption{Example \ref{ex:2DVortex}.  Adaptive meshes and density contours at  $t=0,2,4$ with  $80\times 80$ cells.}
	\label{fig:2DVortex}
\end{figure}

\begin{figure}[!ht]
	\centering
	\begin{subfigure}[b]{0.36\textwidth}
		\centering
\includegraphics[width=1.0\textwidth]{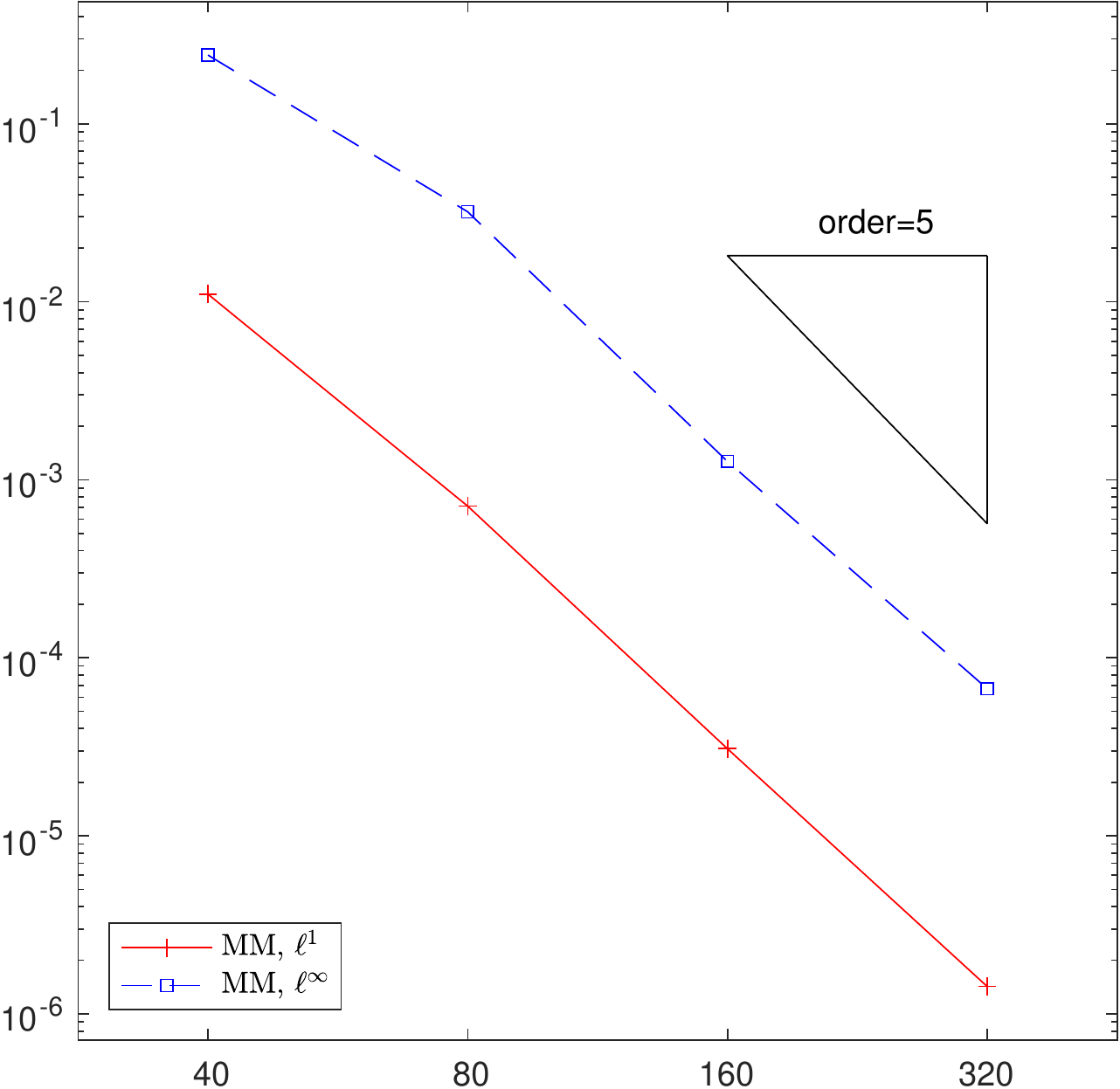}
	\end{subfigure}
	\begin{subfigure}[b]{0.36\textwidth}
		\centering
		\includegraphics[width=1.0\textwidth]{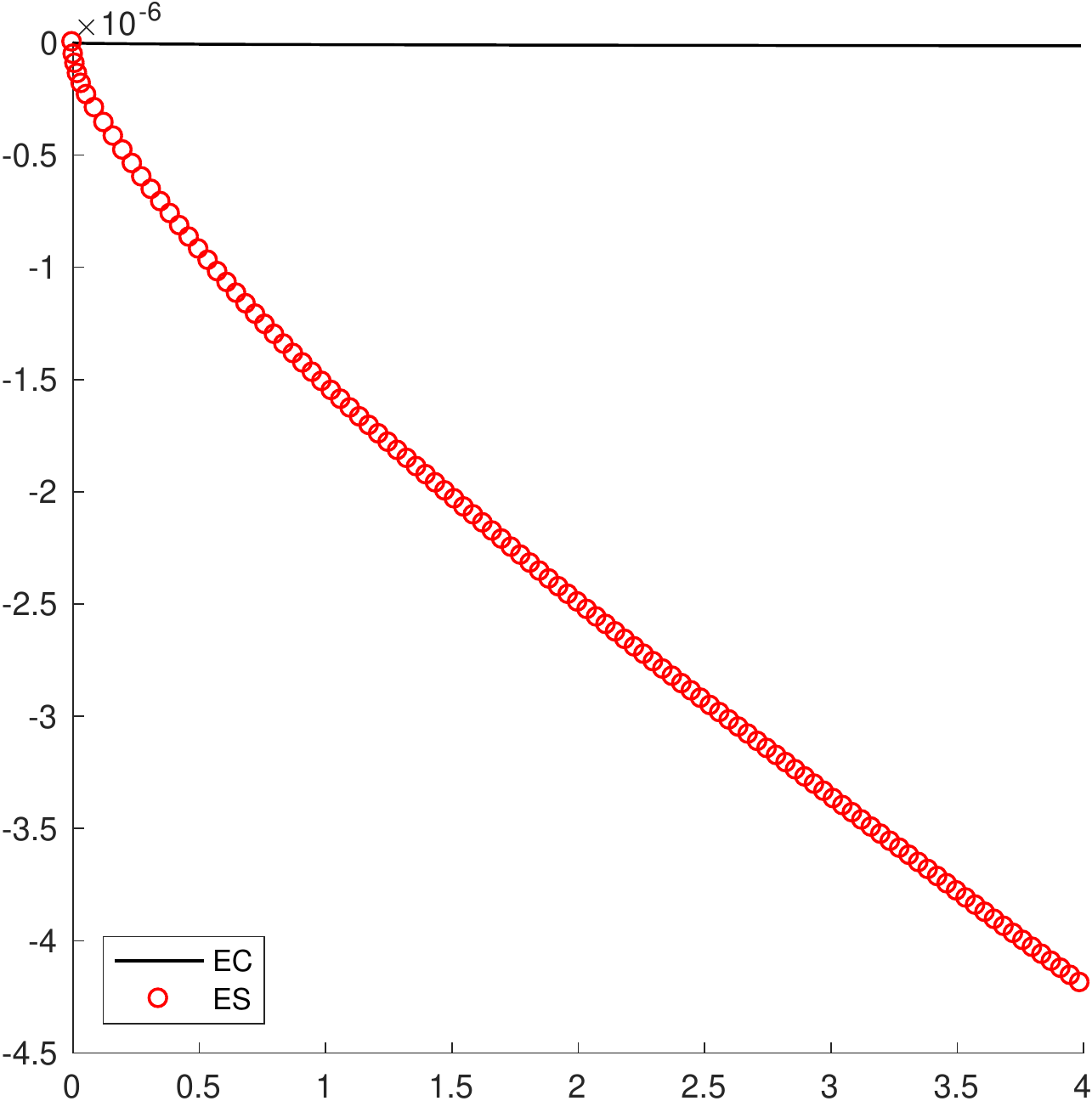}
	\end{subfigure}
	\caption{Example \ref{ex:2DVortex}.
Left:  $\ell^1$- and $\ell^\infty$-errors in $\rho_1$ at $t=4$ versus $N_1$ and orders of convergence;
right:   total entropy with respect to $t$  with $N_1 = 160$.
}
	\label{fig:2DVortex2}
\end{figure}

\end{example}

	\begin{example}[Quasi 2D  shock tube]\label{ex:2DQP}\rm
The initial data
	are
	$$
	\label{QRP}
	\left(\rho_1, v_{1}, v_{2}, p\right)=\left\{\begin{array}{ll}
	(1,0,0,1), & x_1<0.5, \\
	(0.75,0,0,0.05), & x_1>0.5,
	\end{array}\right.
	$$
	with $p_{\infty,1} = 1$ and $\Gamma_1 = 3$, see   \cite{wu2008general}. The exact solution consists of  a narrow rarefaction wave, a contact discontinuity and a right  moving  shock wave.
	The monitor function is chosen as \eqref{eq:monitor}
		with $\kappa = 1, \sigma_1 = \rho_1$ and $\alpha_1 = 1200.$
	
	\begin{figure}[!ht]
		\centering
		\begin{subfigure}[b]{0.402\textwidth}
			\centering
			\includegraphics[width=1.0\linewidth]{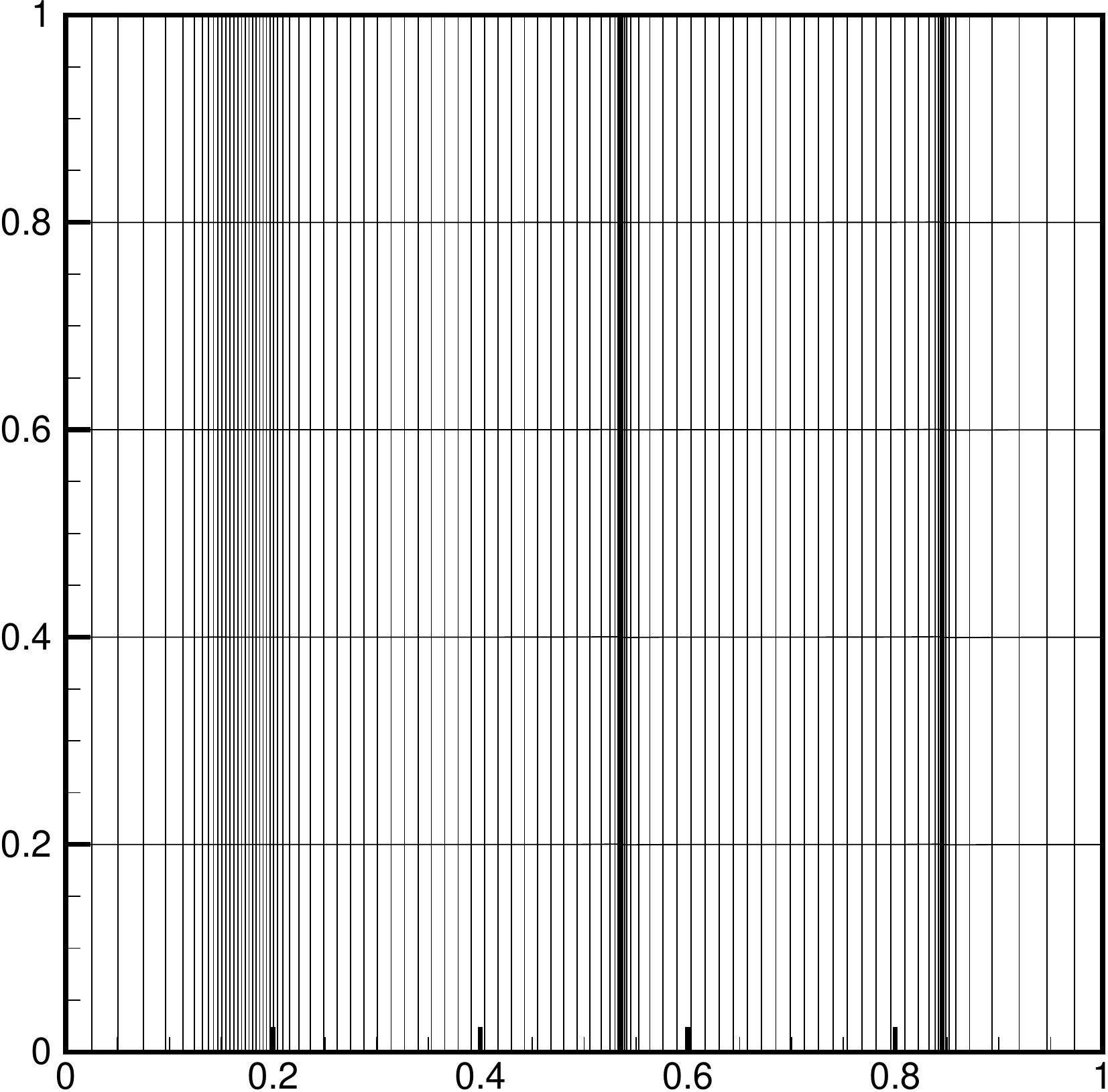}
		\end{subfigure}
		\begin{subfigure}[b]{0.4\textwidth}
			\centering
			\includegraphics[width=1.0\linewidth,  clip]{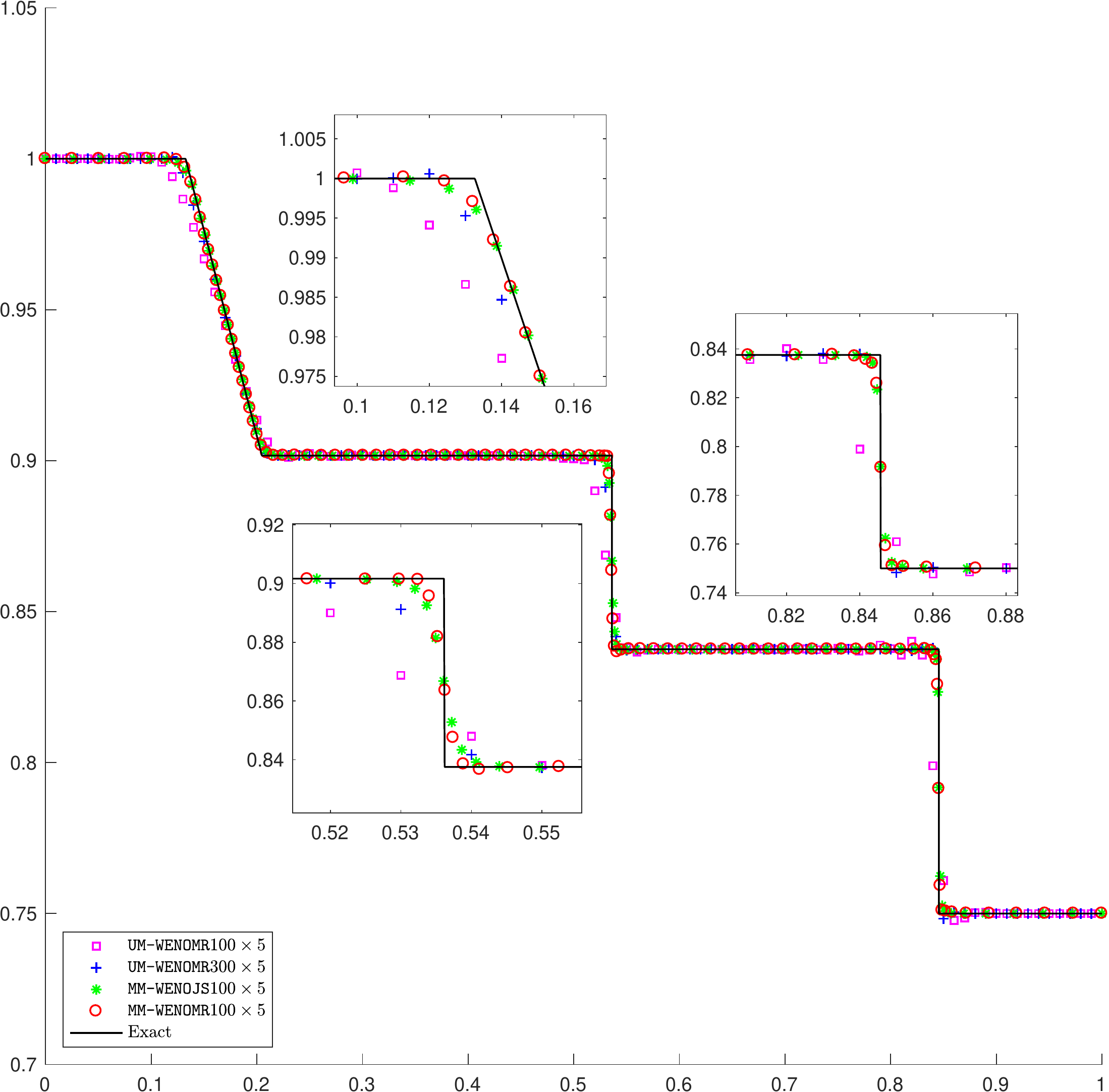}
		\end{subfigure}
		\caption{Example {\ref{ex:2DQP}}.
Adaptive mesh of {$100 \times 5$} cells (Left) and densities (Right) at
$t = 0.15$.}
		\label{fig:QRP}
	\end{figure}
	
	Figure \ref{fig:QRP} shows the   adaptive mesh  and   the densities at $t = 0.15$ obtained respectively by the schemes on the adaptive moving mesh  and the uniform mesh.
We see that	the mesh points adaptively concentrate near the large gradient area of the density, 
		 and   {\tt MM-WENOMR} with $N_{1}=100$ is superior  to   {\tt UM-WENOMR}
		 with $N_{1}=100$, better than {\tt MM-WENOJS} with $N_1 = 100$ near  the contact discontinuity and {\tt UM-WENOMR} with $N_{1}=300$ near
		 the head of the rarefaction wave, the contact discontinuity. 
		
\end{example}

	 \begin{example}[2D Riemann problem \uppercase\expandafter{\romannumeral1}]\label{ex:RP1}\rm
	 The initial data are \cite{BRIO2001177}
	 $$
	\left(\rho_1, v_{1}, v_{2}, p\right)=\left\{\begin{array}{ll}
	(0.5313,0,0,0.4)), & x_1>0.5, x_2>0.5,\\
	(1,0.7276,0,1), & x_1<0.5, x_2>0.5, \\
	(0.8,0,0,1), & x_1<0.5, x_2<0.5,\\
	(1,0,0.7276,1), &\text {otherwise, }
	\end{array}\right.
	$$
	with $p_{\infty,1} = 0$.
	The initial
	 discontinuities are two shock waves and two contact discontinuities.

	  \begin{table}
	 	\centering
	 	\resizebox{.95\columnwidth}{!}{
	 	\begin{tabular}{l|cccc}
	 		\hline & {\tt MM-WENOMR} & {\tt MM-WENOJS}& {\tt UM-WENOMR}  &{\tt UM-WENOMR}  \\
	 		\hline Example \ref{ex:RP1} & $2 \mathrm{m} 56\mathrm{s}$  ($200 \times 200$ cells)
	 		& $2 \mathrm{m} 37\mathrm{s}$  ($200 \times 200$ cells)
	 		& $28 \mathrm{s}$  ($200 \times 200$ cells)& $7 \mathrm{m} 9 \mathrm{s}$ ($500 \times 500$ cells)\\
	 		Example \ref{ex:RP2}  &$2\mathrm{m} 50 \mathrm{s}$ ($200 \times 200$ cells) & $2\mathrm{m}38 \mathrm{s}$ ($200 \times 200$ cells)  & $36 \mathrm{s}$ ($200 \times 200$ cells)& $13\mathrm{m} 10 \mathrm{s}$ ($600 \times 600$ cells)\\
	 		\hline
	 	\end{tabular}
 	}
	 	\caption{CPU times of Examples \ref{ex:RP1}-\ref{ex:RP2} ($4$ cores).}
	 		\label{CPU}
	 \end{table}
The monitor function is  the same as that used in Example \ref{ex:2DQP}, and the linear weights of the multi-resolution WENO reconstruction are taken as $\chi_1 = 0.95, \chi_2 = 0.045$ and $\chi_3 = 0.005$. %
 Figure \ref{fig:RP1} gives the adaptive mesh of {\tt MM-WENOMR} with $200 \times 200$ cells,  the density contours ($40$
 equally spaced contour lines) and the densities along $x_1 = x_2$ at $t = 0.25$.
 The schemes can capture important flow structures such as the {Mach} reflection resulting from the initial two shock waves  and the following emerged jet moving towards the lower left direction near the center of the domain.
One can see that {\tt MM-WENOMR} with $200 \times 200$ cells gives sharper transitions near the shock waves  than {\tt UM-WENOMR}  with $500 \times 500$ cells and {\tt MM-WENOJS} with $200 \times 200$ cells,
and	  the resolution of {\tt MM-WENOMR} is better than {\tt MM-WENOJS}, see Figure \ref{fig:RP1_Comp},  although its CPU time    is slightly larger than {\tt MM-WENOJS}, see Table \ref{CPU}.
From Table \ref{CPU}, we can also see that {\tt MM-WENOMR} only takes $41.0\%$  CPU time of {\tt UM-WENOMR}
	 with $500 \times 500$ cells, showing the high efficiency of the adaptive moving mesh scheme.

	 \begin{figure}[!ht]
	 	 \centering
	 		\begin{subfigure}[b]{0.32\textwidth}
	 			\centering
 	 			\includegraphics[width=1.0\linewidth]{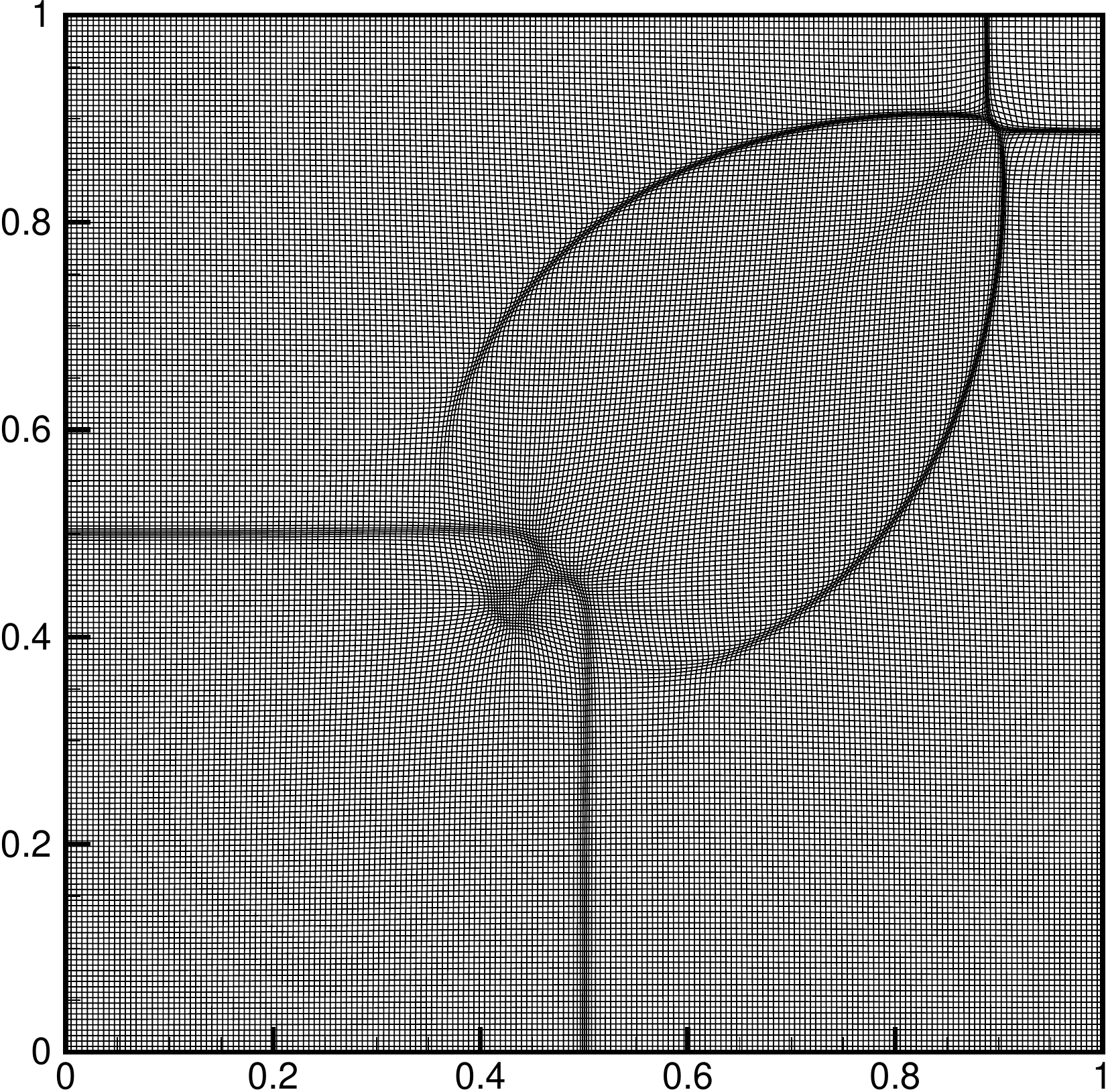}
	 			\caption{  {\tt MM-WENOMR} ($ 200\times 200$)}
	 			\end{subfigure}
	 			\begin{subfigure}[b]{0.32\textwidth}
	 			\centering
	 			\includegraphics[width=1.0\linewidth]{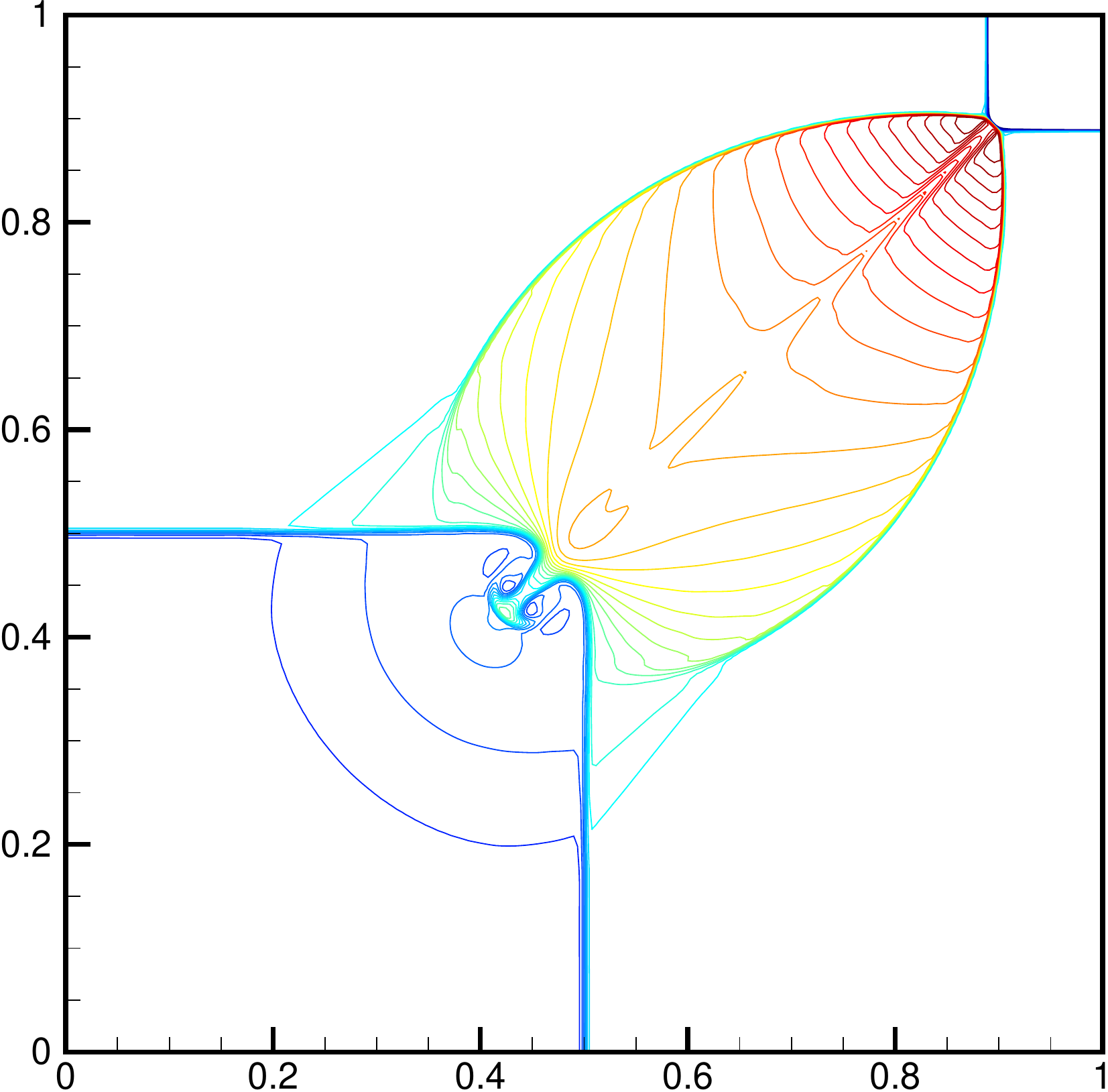}
	 			\caption{{\tt MM-WENOMR} ($ 200\times 200$)}
	 		\end{subfigure}	 	
 			\begin{subfigure}[b]{0.32\textwidth}
 			\centering
 			\includegraphics[width=1.0\linewidth]{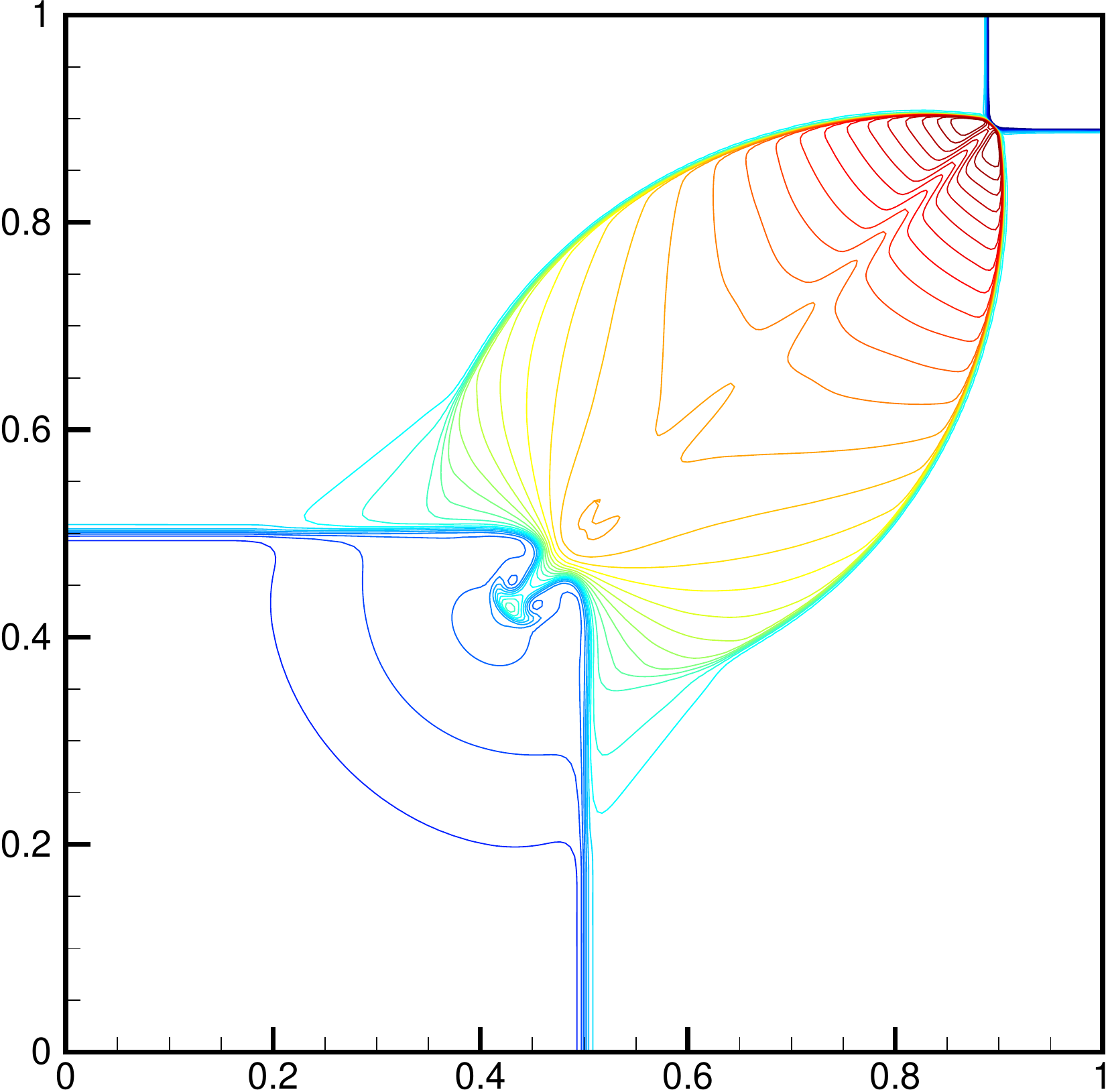}
 			\caption{{\tt MM-WENOJS} ($ 200\times 200$)}
 		\end{subfigure}

 		\begin{subfigure}[b]{0.32\textwidth}
 			\centering
 			\includegraphics[width=1.0\linewidth]{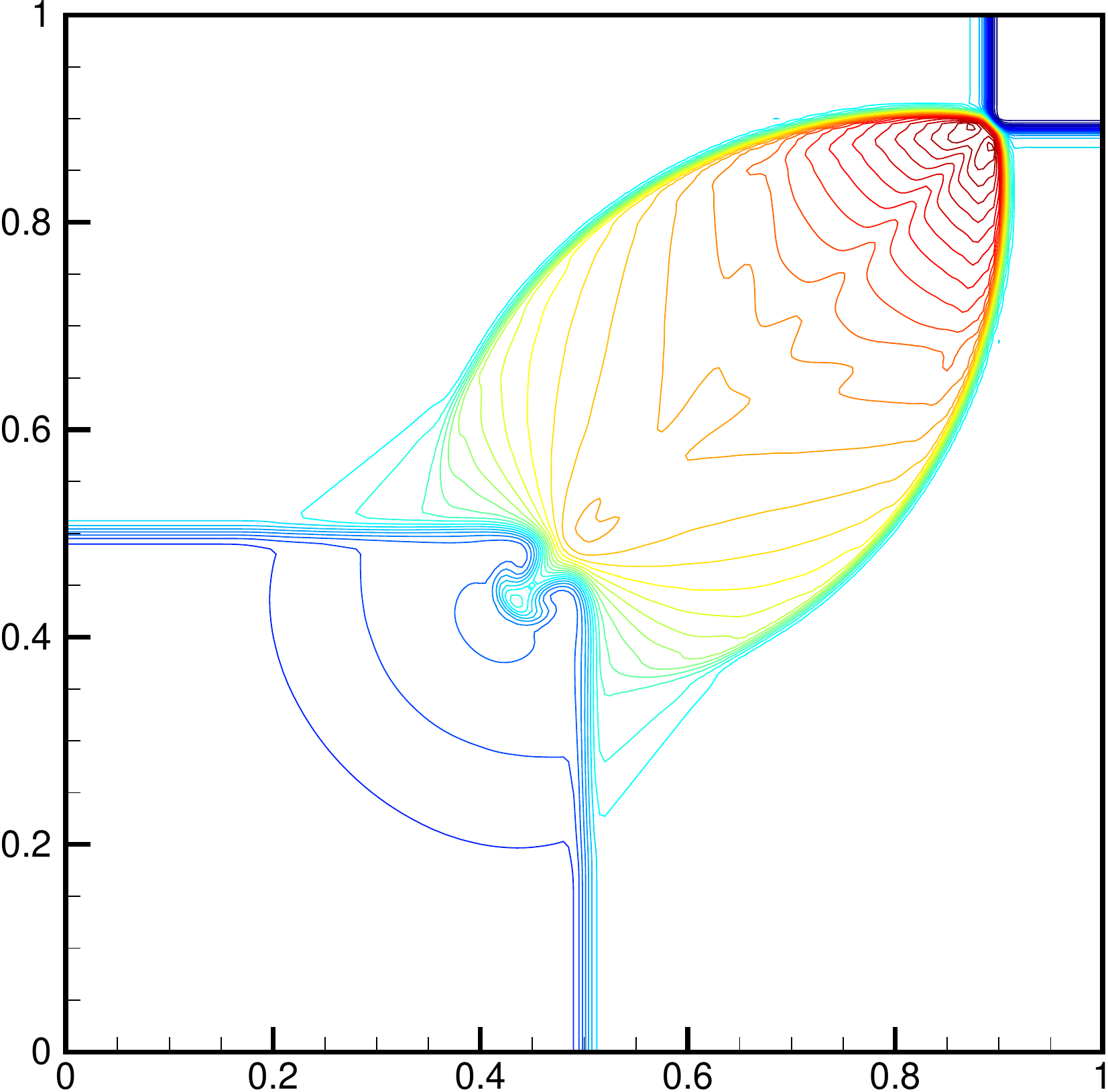}
 			\caption{ {\tt UM-WENOMR} ($ 200\times 200$)}
 		\end{subfigure}
	 			\begin{subfigure}[b]{0.32\textwidth}
	 			\centering
	 			\includegraphics[width=1.0\linewidth]{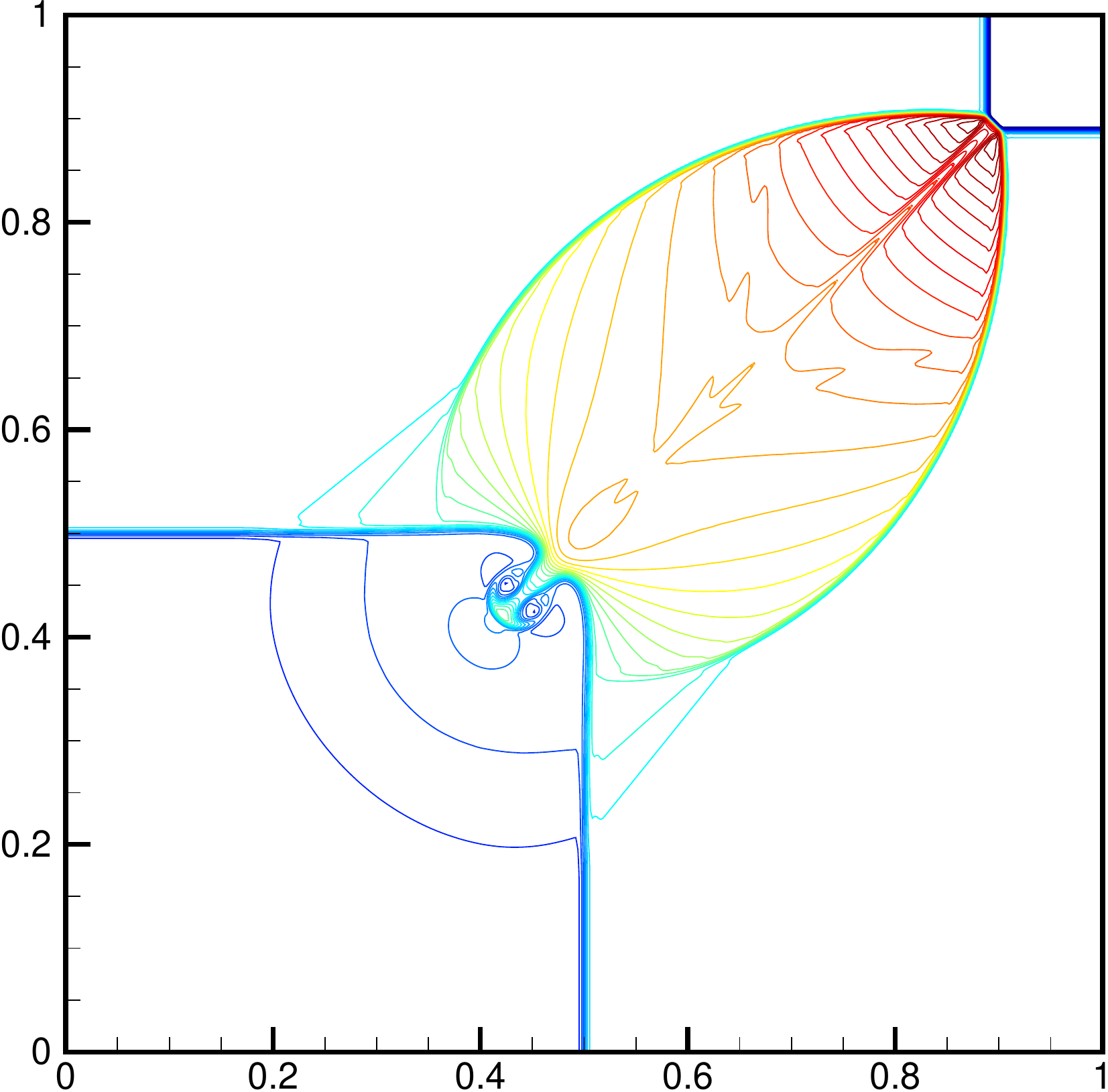}
	 			\caption{{\tt UM-WENOMR} ($ 500\times 500$)}
	 				\end{subfigure}
 					\begin{subfigure}[b]{0.316\textwidth}
 					\centering
 					\includegraphics[width=1.0\linewidth]{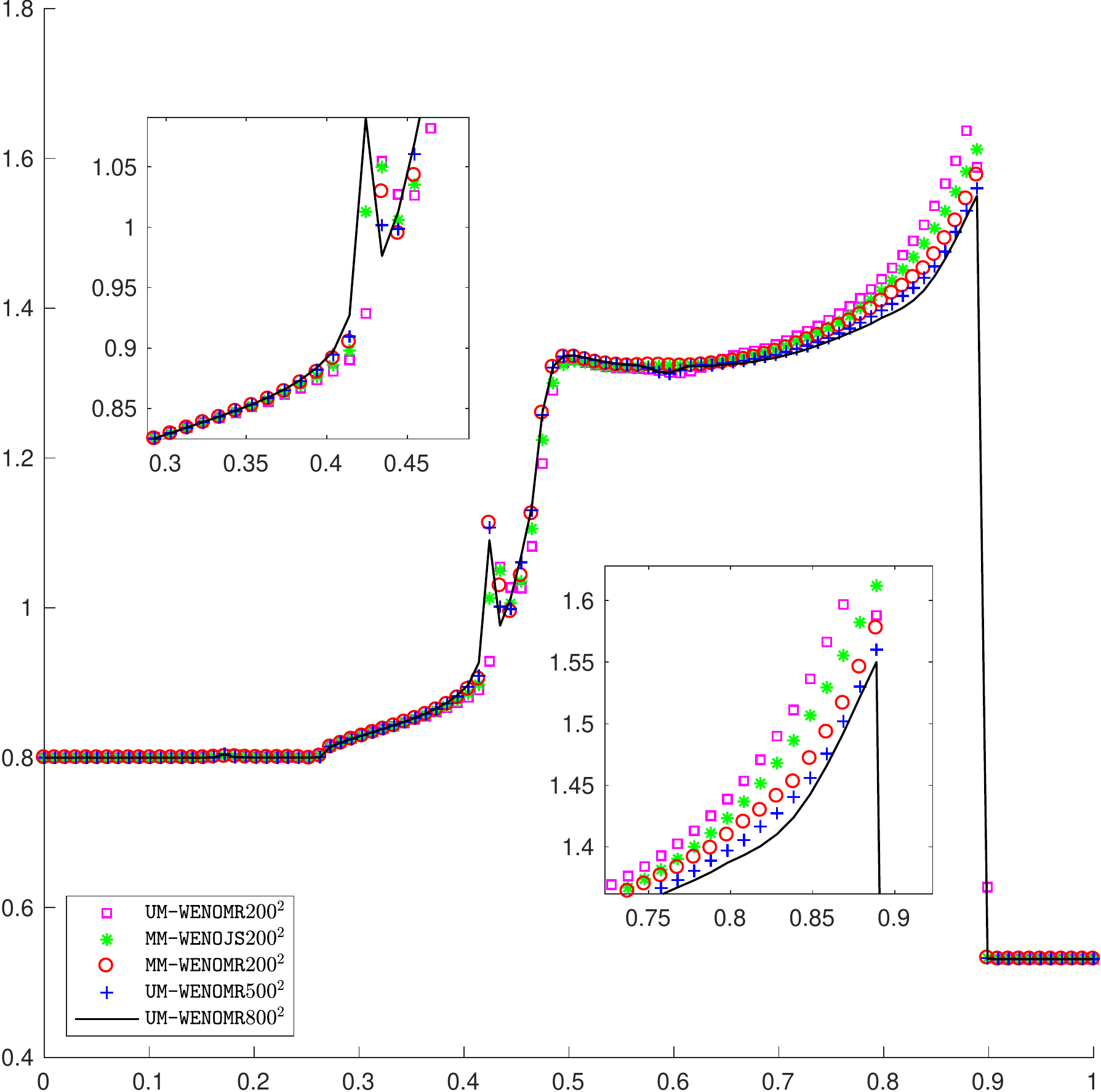}
 					\caption{$\rho_1$ along  $x_1 = x_2$}
 					\label{fig:RP1_Comp}
 				\end{subfigure}
	 		\caption{Example  \ref{ex:RP1}. Adaptive mesh of {\tt MM-WENOMR} with $200 \times 200$ cells,
 density contours (40 equally spaced
contour lines) of {\tt MM-WENOMR}, {\tt MM-WENOJS}, and
 {\tt UM-WENOMR}, and   densities along $x_1 = x_2$  at $t = 0.25$.}
	 		\label{fig:RP1}
	 	
	 \end{figure}
 \end{example}
	
	\begin{example}[2D Riemann problem \uppercase\expandafter{\romannumeral2}]\label{ex:RP2}\rm
	  The initial data are \cite{LaxLiu1998}
	 $$
	 \left(\rho_1, v_{1}, v_{2}, p\right)=\left\{\begin{array}{ll}
	 (1, 0.75, -0.5, 1), & x_1>0.5, x_2>0.5 ,\\
	 (2, 0.75, 0.5, 1), & x_1<0.5, x_2>0.5, \\
	 (1, -0.75, 0.5, 1),& x_1<0.5, x_2<0.5, \\
	 (3, -0.75, - 0.5, 1), & \text{otherwise,}
	 \end{array}\right.
	 $$
	 with $p_{\infty,1} = 0.$
	It describes the interaction of four contact discontinuities.
	
	 \begin{figure}[!ht]
		\centering
		\begin{subfigure}[b]{0.32\textwidth}
			\centering
			\includegraphics[width=1.0\linewidth]{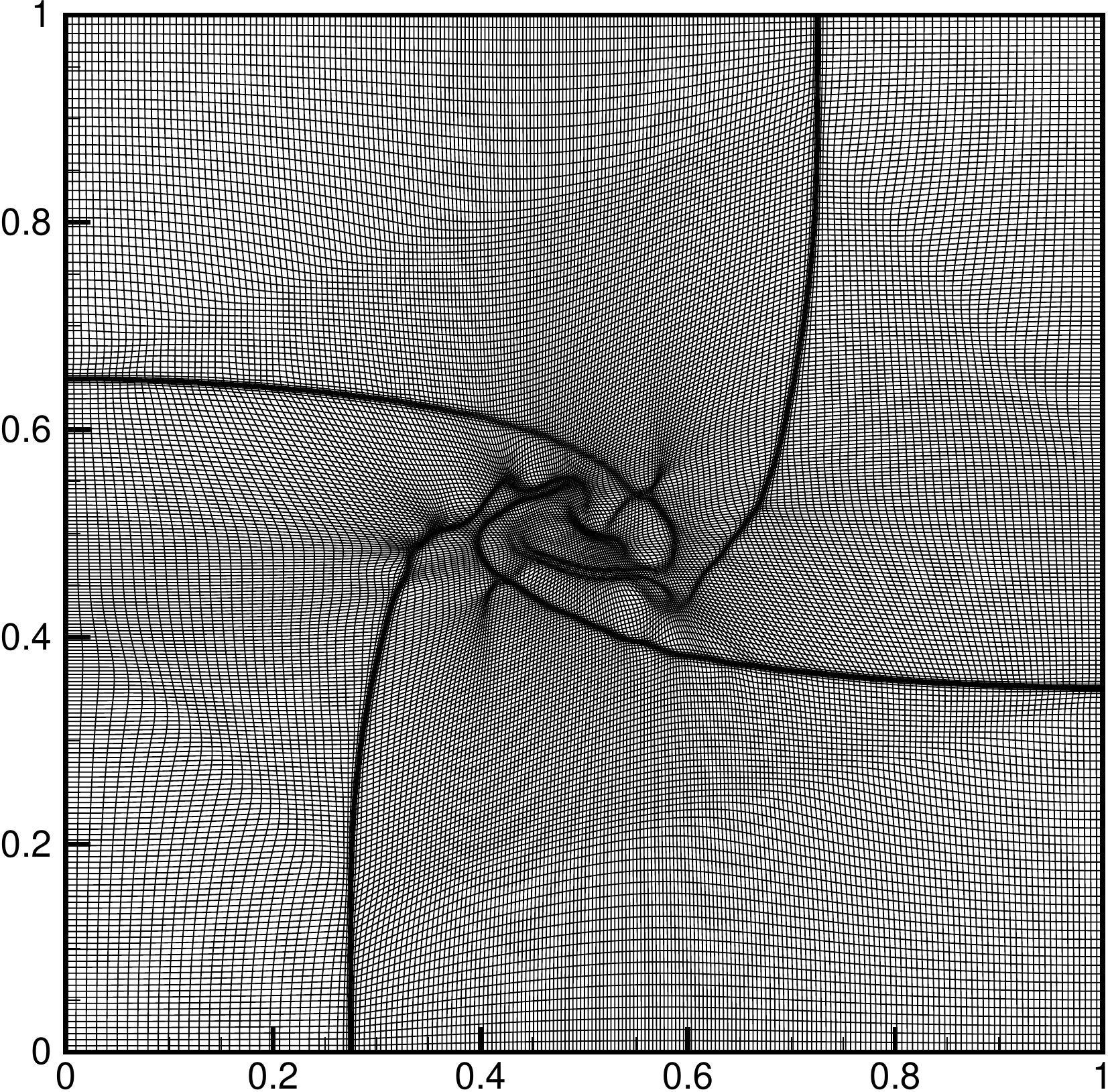}
			\caption{{\tt MM-WENOMR} ($ 200\times 200$)}
		\end{subfigure}
		\begin{subfigure}[b]{0.32\textwidth}
			\centering
			\includegraphics[width=1.0\linewidth]{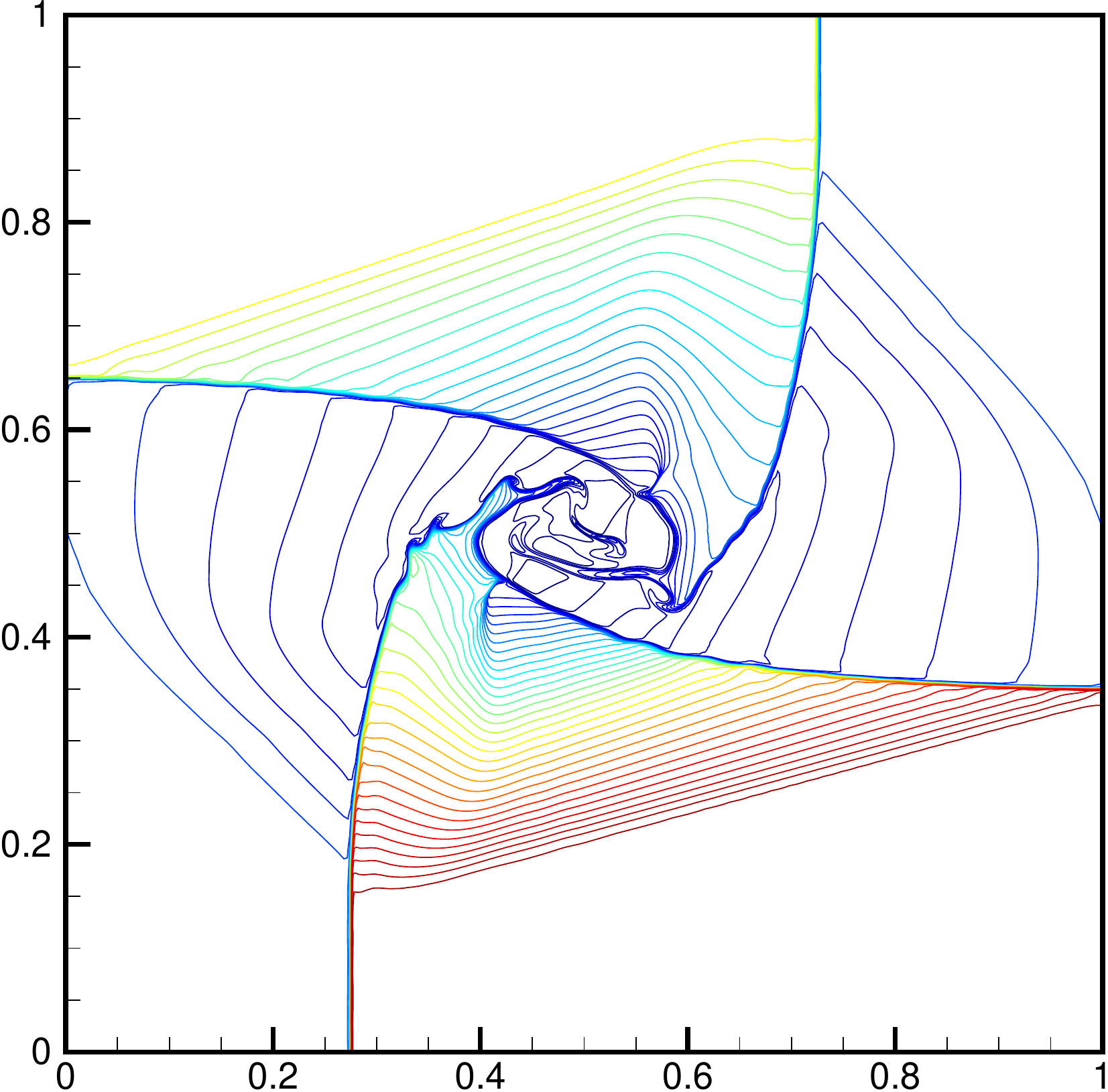}
			\caption{{\tt MM-WENOMR} ($ 200\times 200$)}
		\end{subfigure}
	\begin{subfigure}[b]{0.32\textwidth}
		\centering
		\includegraphics[width=1.0\linewidth]{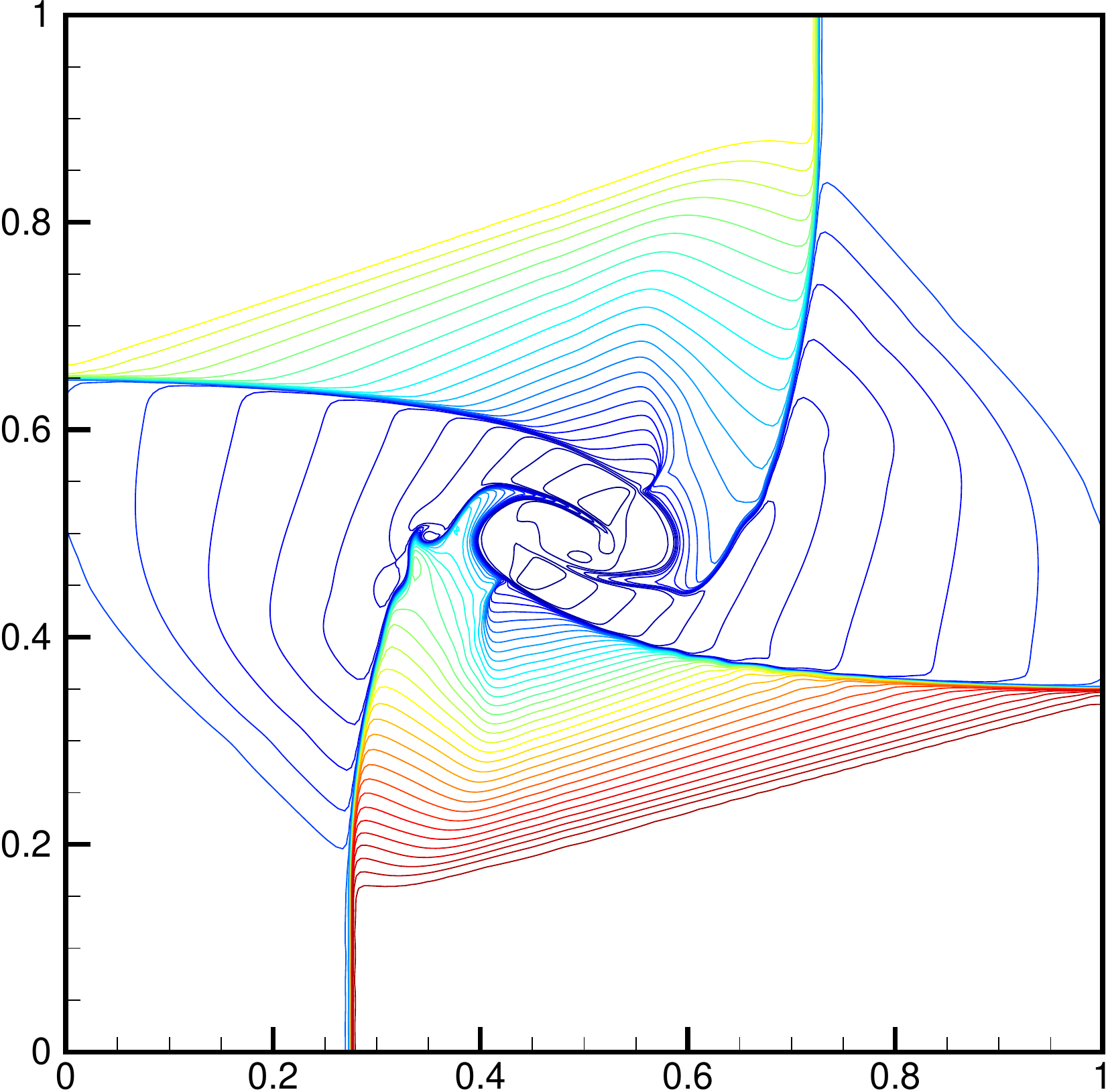}
		\caption{{\tt MM-WENOJS} with ($ 200\times 200$)}
	\end{subfigure}
		
		\begin{subfigure}[b]{0.32\textwidth}
			\centering
			\includegraphics[width=1.0\linewidth]{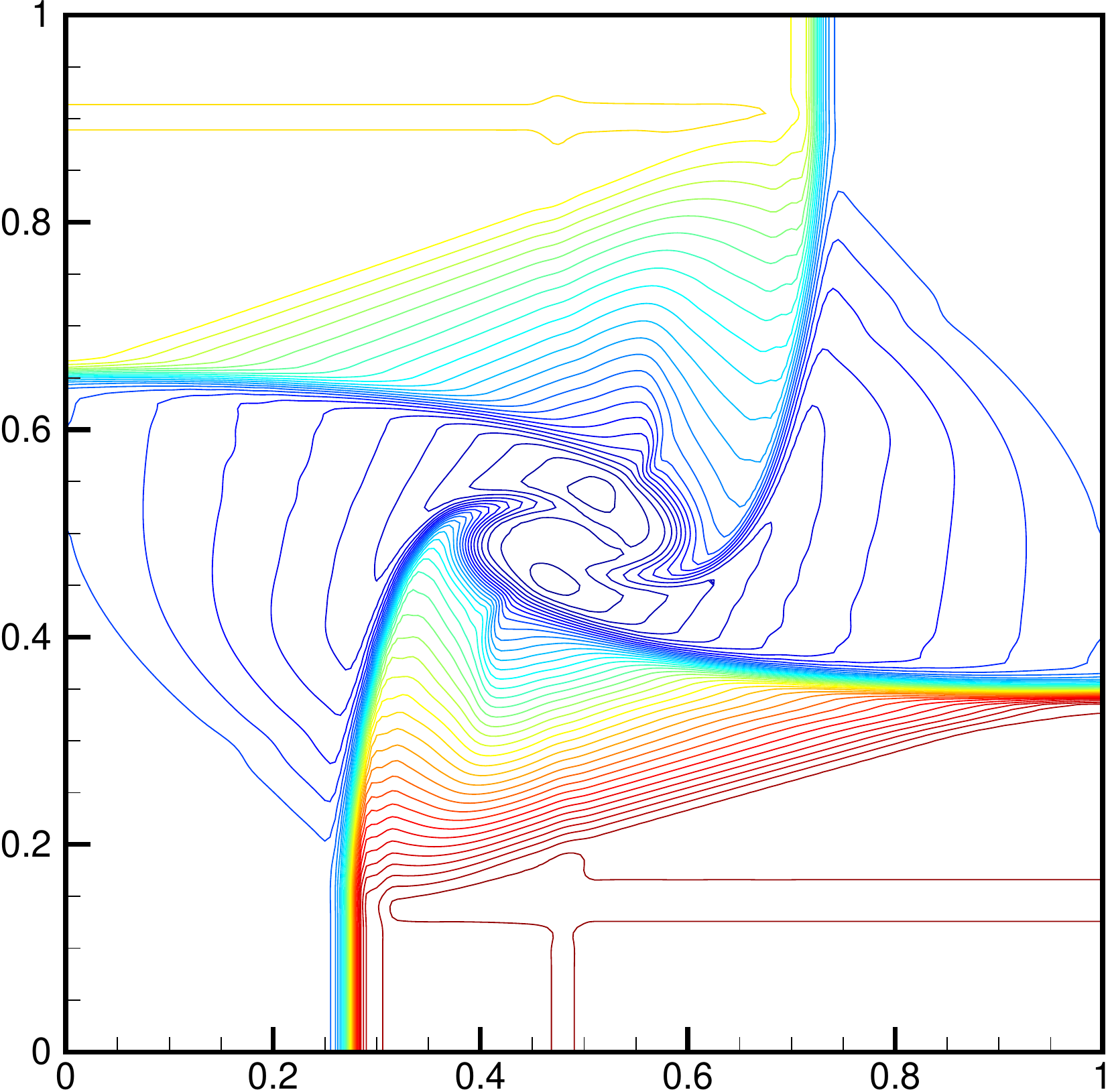}
			\caption{ {\tt UM-WENOMR} ($ 200\times 200$)}
		\end{subfigure}
		\begin{subfigure}[b]{0.32\textwidth}
			\centering
			\includegraphics[width=1.0\linewidth]{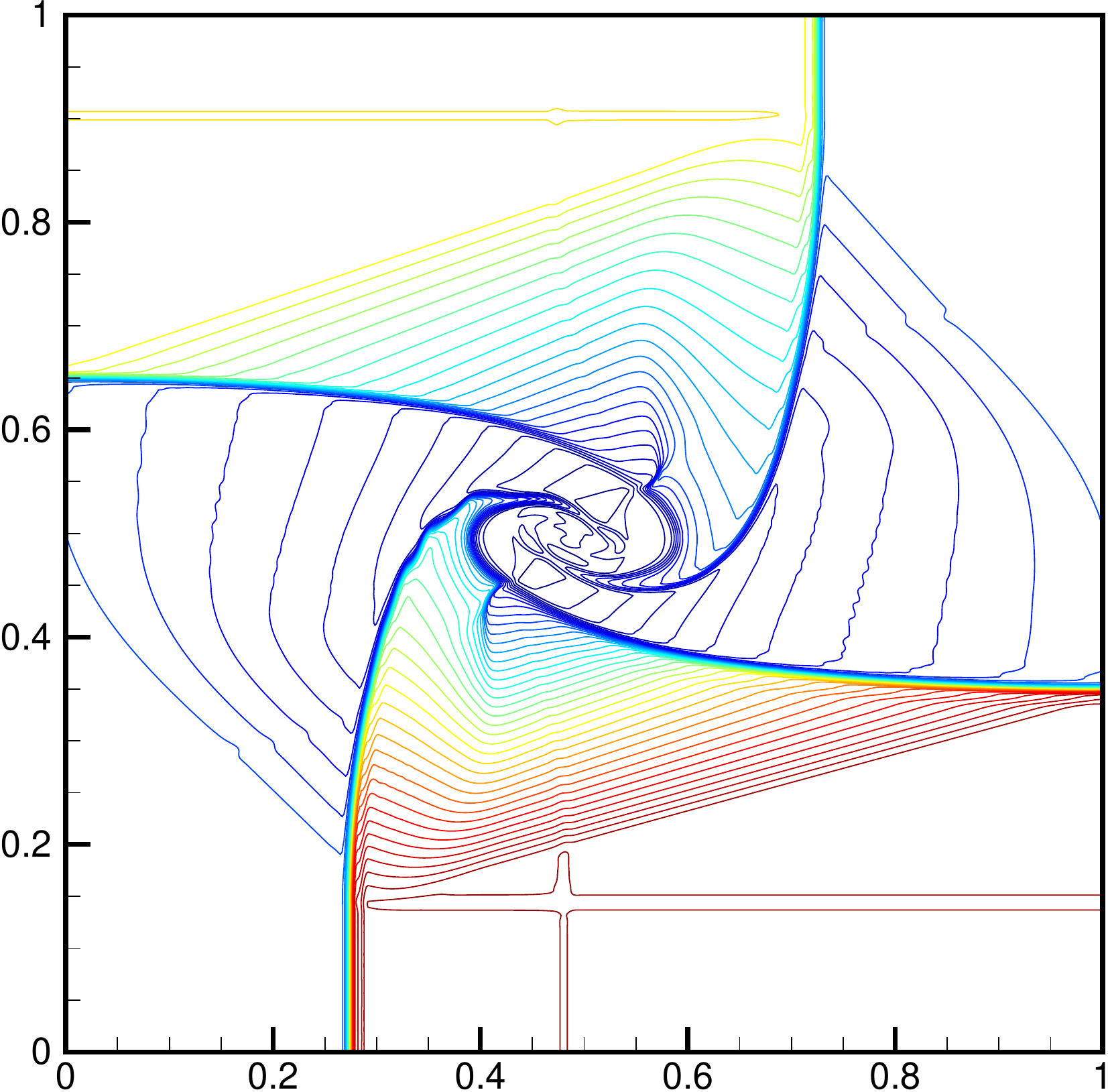}
			\caption{ {\tt MM-WENOMR} ($ 600\times 600$)}
		\end{subfigure}
			\begin{subfigure}[b]{0.316\textwidth}
			\centering
			\includegraphics[width=1.0\linewidth]{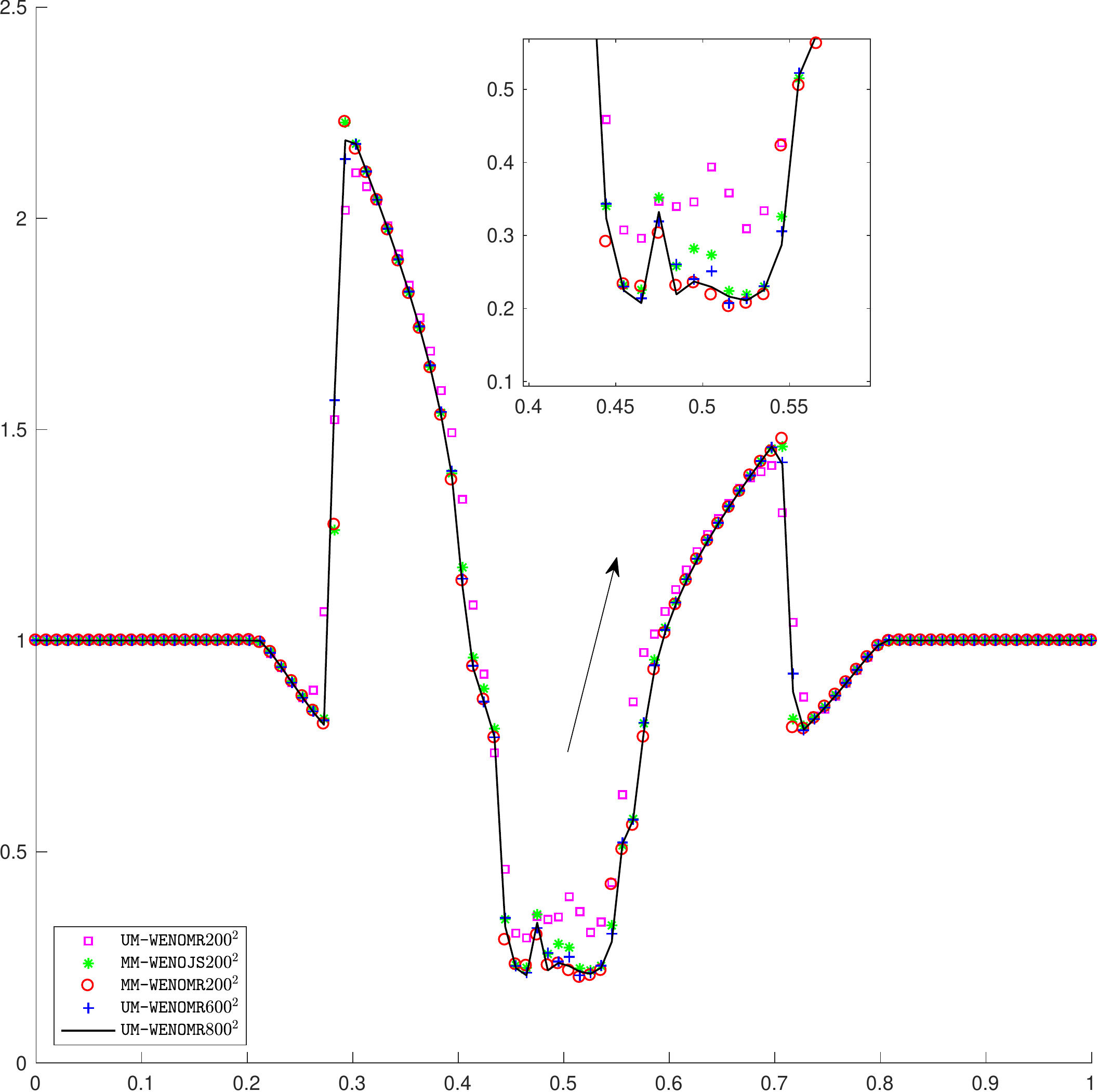}
			\caption{$\rho_1$ along  $x_1 = x_2$}
		\end{subfigure}
	\caption{Example  \ref{ex:RP2}.
 Adaptive mesh of {\tt MM-WENOMR} with $200 \times 200$ cells,
 density contours (40 equally spaced
contour lines) of {\tt MM-WENOMR}, {\tt MM-WENOJS}, and
 {\tt UM-WENOMR}, and   densities along $x_1 = x_2$  at $t = 0.3$.
}

		\label{fig:RP2_t03}
	\end{figure}
%
%
The monitor function and the linear weights of the multi-resolution WENO reconstruction  are the same as those used in Example \ref{ex:2DQP}.
	  Figure \ref{fig:RP2_t03} 
	  shows  the adaptive mesh of {\tt MM-WENOMR} with $200 \times 200$ cells and the densities at $t = 0.3$.
	  It is seen that the four initial contact discontinuities interact with each other to form a spiral with the low density around the center of the domain as time increases. The moving mesh methods capture the contact discontinuities and the rarefaction wave well,
and {\tt MM-WENOMR} displays more small scale structures and roll up of the slip lines than {\tt UM-WENOMR} with $600 \times 600$ cells and {\tt MM-WENOJS} with $200 \times 200$ cells.
	  The CPU times in Table \ref{CPU} clearly highlight the efficiency of {\tt MM-WENOMR}, which takes only  $21.5\%$  CPU time of {\tt UM-WENOMR} with $600 \times 600$ cells and $7.6\%$  time more than {\tt MM-WENOJS} with $200 \times 200$ cells.
	 \end{example}

\begin{example}[3D isentropic vortex]\label{ex:3DVortex}\rm
	This example is to check the accuracy of the 3D ES adaptive moving mesh method via the 3D isentropic vortex problem  describing a smooth isentropic vortex {moving} in a certain direction. The initial data are similar to that in \cite{BOSCHERI2014484} except for that the cylindrical vortex is rotated to the diagonal of  the domain $\Omega_p=[-10,10]^3$, 
	given by
	\begin{align*}
& T = 1 -\frac{ \epsilon^{2}}{8 \Gamma_1 \pi^{2}} e^{1-r^{2}},\quad
 \rho_1 = T^{1/(\Gamma_1-1)},\quad p = (\Gamma_1-1)\rho_1T,\\
& \bm{v} = \frac{1}{\sqrt{6}}\left(\tilde{v}_1 - \sqrt{3}\tilde{v}_2 + \sqrt{2}, \tilde{v}_1 + \sqrt{3}\tilde{v}_2 + \sqrt{2}, -2\tilde{v}_1+ \sqrt{2}\right)^{\rm{T}},
	\end{align*}
	where
	\begin{align*}
& p_{\infty,1} = 0, \quad c_{v,1} = 1, \quad\epsilon = 5,
	\quad r^2 = \tilde{x}_1^2 + \tilde{x}_2^2,\\
&	(\tilde{x}_1, \tilde{x}_2) = \dfrac{1}{\sqrt{6}}(\hat{x}_1 +2\hat{x}_2, \sqrt{3}\hat{x}_1),\quad (	\tilde{v}_1, \tilde{v}_2) =
(1, 1) + \frac{\epsilon}{2 \pi} e^{\frac{1-r^{2}}{2}}(\tilde{x}_2, -\tilde{x}_1),\\
&(\hat{x}_1,	\hat{x}_2) = \left(-x_1 +x_2 +20k_1, x_1-x_3+20k_2\right), \quad (\hat{x}_1,	\hat{x}_2)\in[-10, 10]^2, \quad k_1,k_2\in \mathbb{Z}.
	\end{align*}

Figure \ref{fig:3DVortex} shows the
 $\ell^1$- and $\ell^\infty$-errors in $\rho_1$ at $t = 0.1$, the orders of convergence obtained by using {\tt MM-WENOMR},
 and {the time-evolutions of the discrete} total entropy $\sum_{\bm{i}} J_{\bm{i}} \eta\left(\boldsymbol{U}_{\bm{i}}\right) \Delta \xi_{1} \Delta \xi_{2} \Delta\xi_{3}$ obtained by the EC adaptive moving mesh scheme and {\tt MM-WENOMR}.
The monitor function is  the same as that used  in  Example \ref{ex:2DVortex}, and the boundary points  move adaptively according to the periodic boundary conditions.
 The results show that {\tt MM-WENOMR} can  achieve the expected convergence orders, and the EC scheme almost keeps the total
entropy conservative whereas the total entropy of the ES scheme decays in time.

	\begin{figure}[!ht]
	\centering
	\begin{subfigure}[b]{0.4\textwidth}
		\centering
		\includegraphics[width=1.0\linewidth]{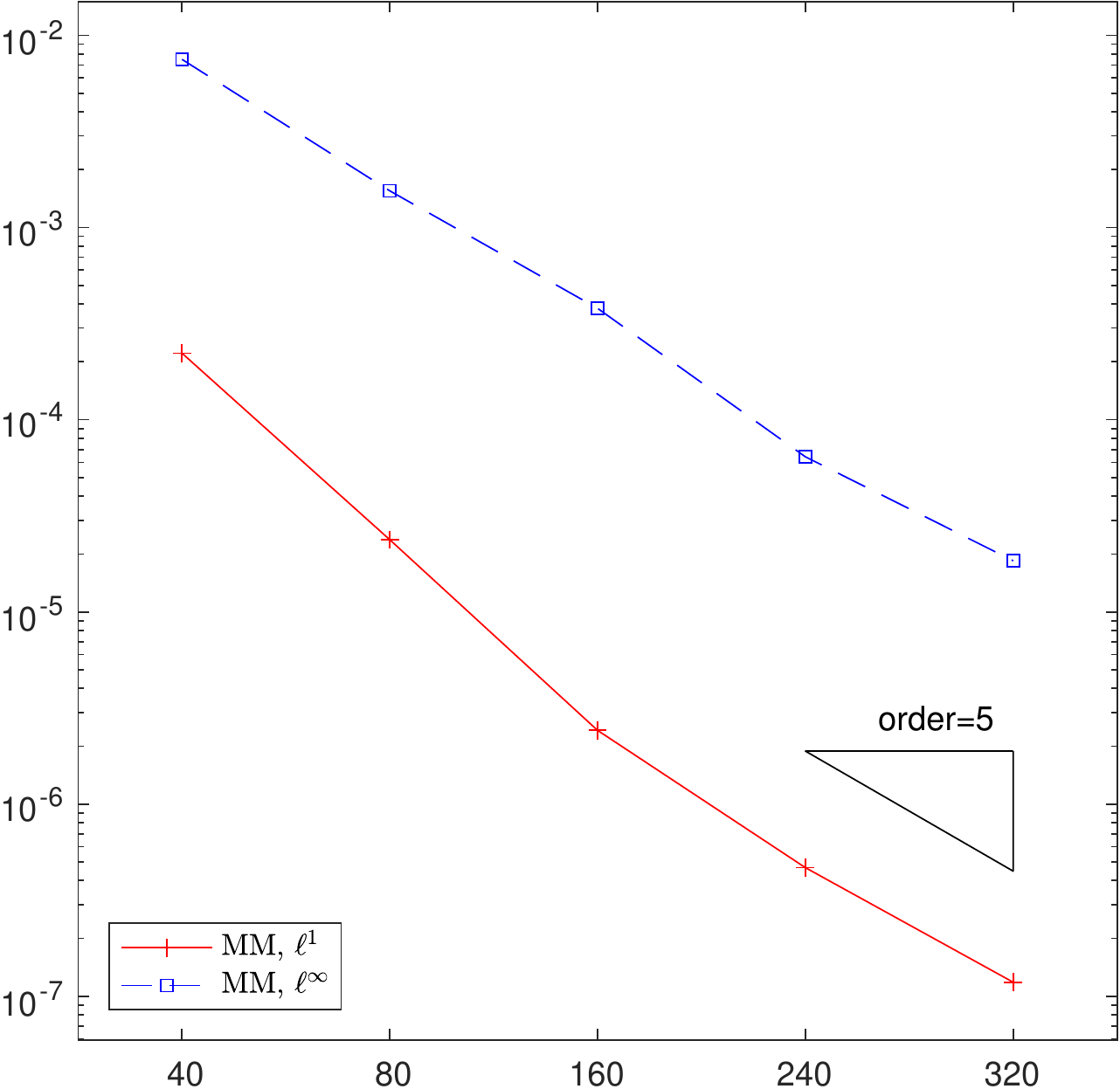}
	\end{subfigure}
	\begin{subfigure}[b]{0.4\textwidth}
		\centering
		\includegraphics[width=1.0\linewidth]{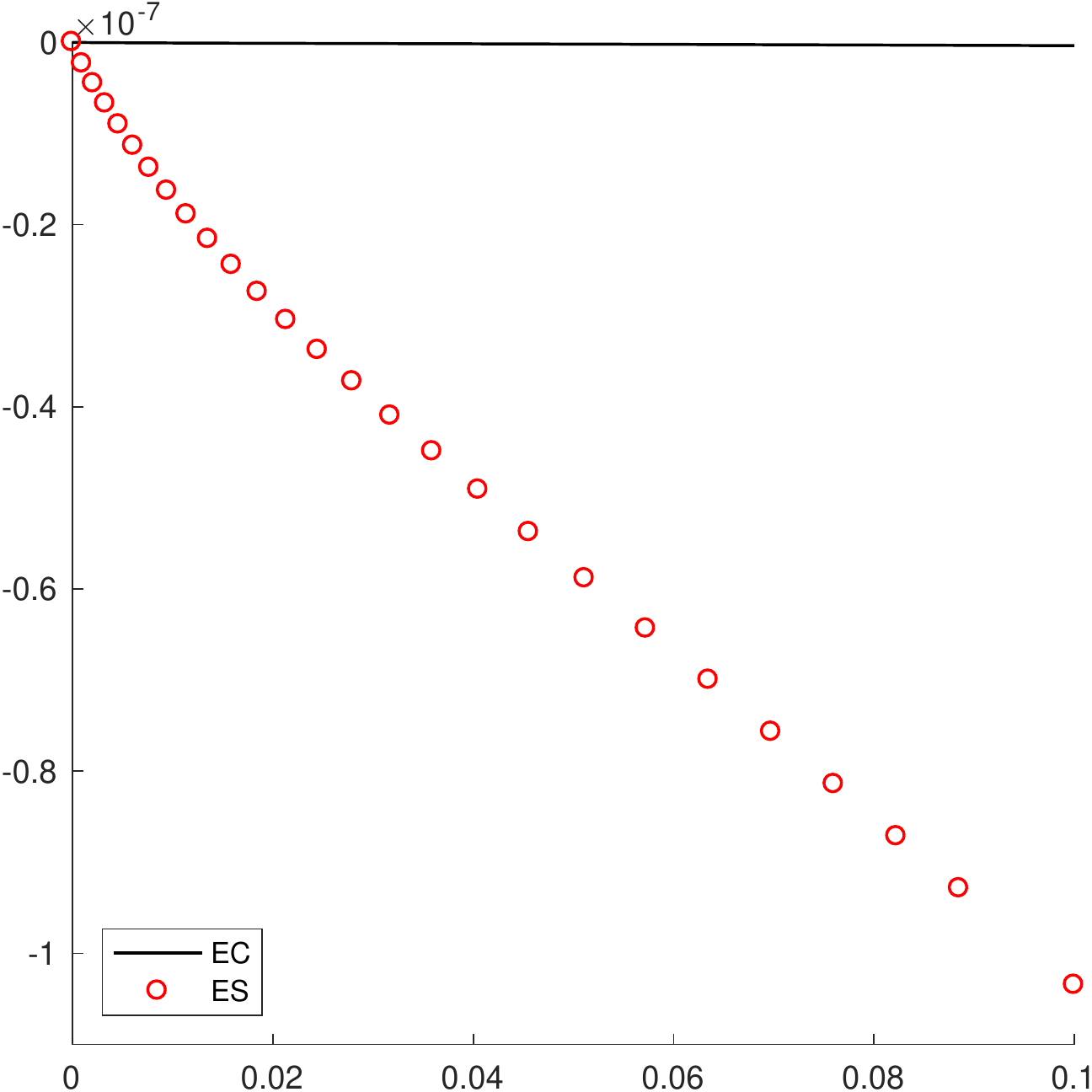}
	\end{subfigure}
	\caption{Example {\ref{ex:3DVortex}}.
 Left: $\ell^1$- and $\ell^\infty$-errors in $\rho_1$ at $t=0.1$ and the orders of convergence; right:  the total entropy  with respect to $t$ with $N_1 = 160$. 
 }
	\label{fig:3DVortex}
\end{figure}

\end{example}

 \begin{example}[3D spherical symmetric shock tube]\label{ex:3DSod}\rm
	The initial data are
	$$
	\left(\rho_1, v_{1}, v_{2},  v_{3}, p\right)=\left\{\begin{array}{ll}
	(1.0,0,0,0,1.0), & \sqrt{x_1^2+x_2^2+x_3^2} < 0.5,\\
	(0.125,0,0, 0,0.1), &\text {otherwise,}
	\end{array}\right.
	$$
	with $ p_{\infty,1} = 0$, and
	the  domain $\Omega_p$ is taken as $[0,1]^{3}$.
 The monitor function is chosen as \eqref{eq:monitor} with
	$\kappa = 1, \sigma_1 = \rho_1$ and $ \alpha_1 = 500$.
Figure \ref{fig:3DSod_Mesh} shows that
 the	 mesh points adaptively concentrate near the large gradient area of the density.
  Figure \ref{fig:3DSod_Rho} gives the densities $\rho_1$ along the line connecting $(0,0,0)$ and $(1,1,1)$, where
    the solid line denotes the reference solution obtained by  a second-order TVD scheme using uniform mesh of $8000$ cells in the 1D spherical coordinates.
It is seen that	  {\tt MM-WENOMR}  with $100^3$ cells  is better than {\tt UM-WENOMR} with $100^3$ cells near the rarefaction wave, the contact discontinuity and the shock wave, and  the adaptive moving mesh method can precisely capture the flow features. Table \ref{CPU_78} shows that {\tt MM-WENOMR}  with $100^3$ cells only takes $37.2\%$ CPU time of {\tt UM-WENOMR} with $200^3$ cells, but it gives comparable results, verifying the efficiency of {\tt MM-WENOMR}, and   the solution of {\tt MM-WENOMR} with $100^3$ cells is as good
	as  that of  {\tt MM-WENOJS} with  comparable CPU time.
	\begin{figure}[!ht]
		\centering
		\begin{subfigure}[b]{0.47\textwidth}
			\centering
			\includegraphics[width=1.0\linewidth, trim=20 10 2 1, clip]{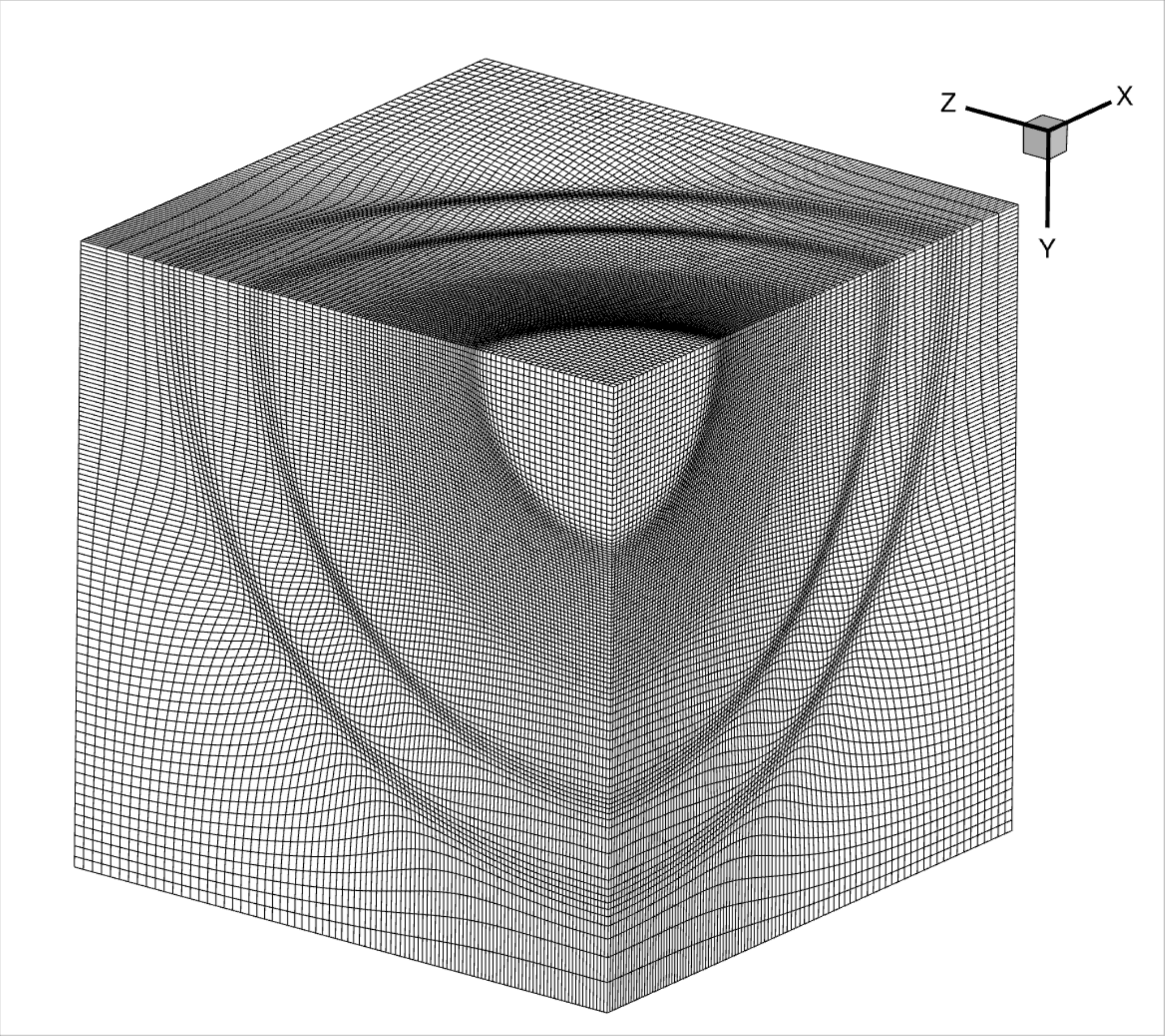}
			\caption{Adaptive mesh of $100^3$ cells}
			\label{fig:3DSod_Mesh}
		\end{subfigure}
		\begin{subfigure}[b]{0.39\textwidth}
			\centering
			\includegraphics[width=1.0\linewidth]{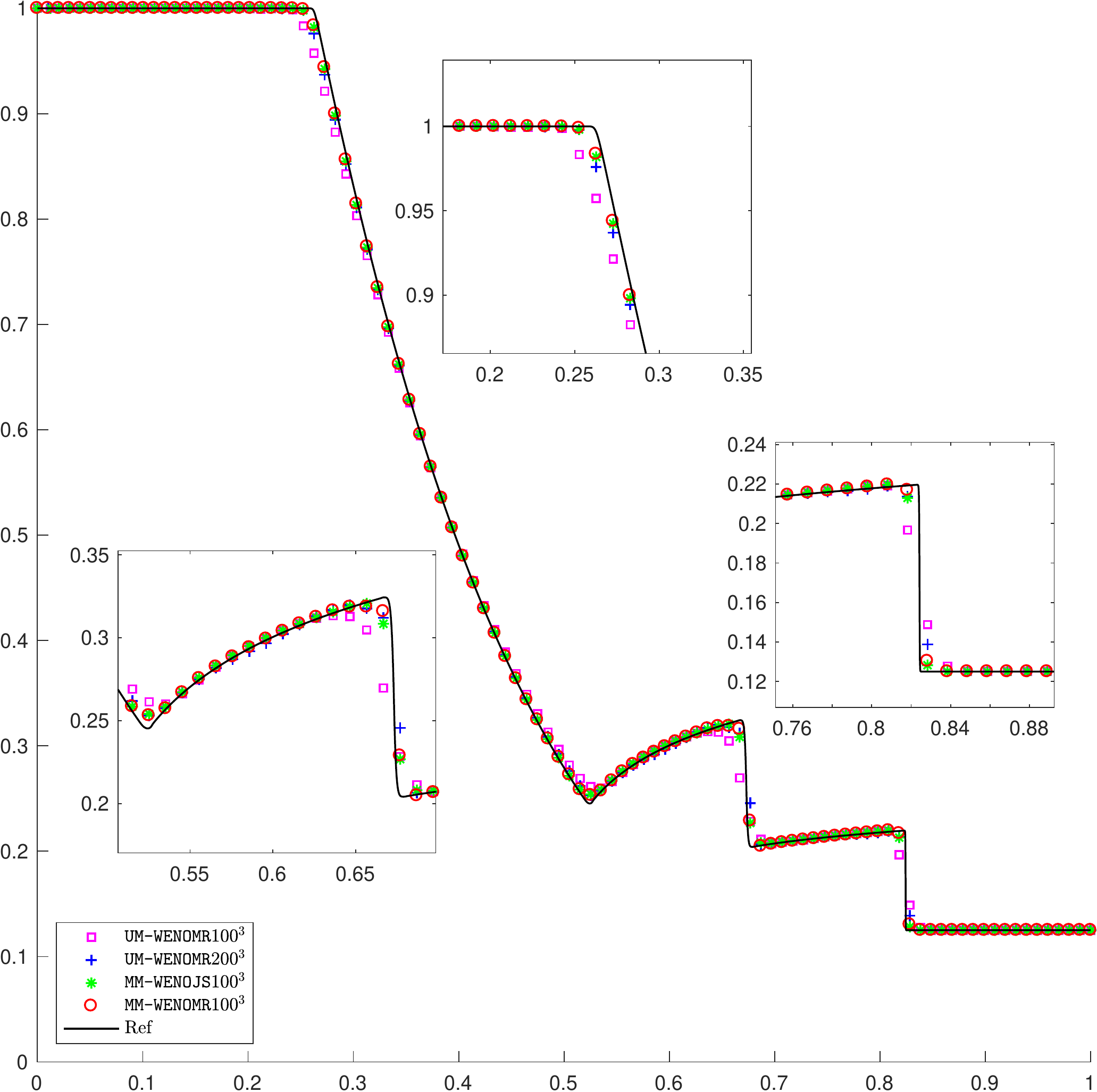}
			\caption{$\rho_1$}
			\label{fig:3DSod_Rho}
		\end{subfigure}
		\caption{Example {\ref{ex:3DSod}}. Adaptive mesh  of {\tt MM-WENOMR} and $\rho_1$ along the line connecting $(0,0,0)$ and $(1,1,1)$ at $t=0.2$.}
		\label{fig:3DSod}
	\end{figure}
	\begin{table}
		\centering
		\resizebox{.95\columnwidth}{!}{
			\begin{tabular}{l|cccc}
				\hline &  {\tt MM-WENOMR}    & {\tt MM-WENOJS} & {\tt UM-WENOMR}  & {\tt UM-WENOMR} \\
				\hline Example \ref{ex:3DSod} & $13 \mathrm{m} 53\mathrm{s}$  ($100^3$ cells) & $13 \mathrm{m} 52\mathrm{s}$  ($100^3$ cells)& $2\mathrm{m} 45\mathrm{s}$ ($100^3$ cells)& $37\mathrm{m} 19 \mathrm{s}$  ($200^3$ cells)\\
				\hline Example \ref{ex:tri} & $ 3 \mathrm{m} 21 \mathrm{s}$ ($350\times 150$ cells)  & $ 3 \mathrm{m} 3 \mathrm{s}$ ($350\times 150$ cells)& $ 49 \mathrm{s}$ ($350\times 150$ cells)  & $18 \mathrm{m} 7 \mathrm{s}$ ($1050\times450$ cells)\\
				\hline
			\end{tabular}
		}
		\caption{CPU times of Examples \ref{ex:3DSod}-\ref{ex:tri} ($32$ cores).}
		\label{CPU_78}
	\end{table}
\end{example}

\subsection{Two-component compressible Euler equations ($N=2$)}
This section solves the 2D and 3D two-component compressible Euler equations ($N=2$) with the ideal and stiffened EOS.

\begin{example}[2D tri-point problem]\label{ex:tri}\rm
It corresponds to a 2D three-state  Riemann problem in a rectangular domain $\Omega_p$, illustrated in Figure \ref{Tri_Dom},  and  has been widely used in  testing the high-resolution numerical schemes \cite{GALERA20105755}. Initially, $\Omega_p=[0,7]\times[0,3]$ is split into  three sub-domains $\Omega_1 = [0, 1]\times[0, 3]$, $\Omega_2 = [1, 7]\times[1.5, 3]$, and
$ \Omega_3 = [1, 7]\times[0, 1.5]$, and
the initial data are
$$
\label{Tri}
\left(\rho_{1}, \rho_{2}, v_{1}, v_{2}, p\right)=\left\{\begin{array}{ll}
(1-\epsilon,\epsilon,0,0,1), & (x_1, x_2) \in \Omega_1, \\
(0.125-\epsilon,\epsilon,0,0,0.1), & (x_1, x_2) \in \Omega_2,\\
(\epsilon,1-\epsilon,0,0,0.1), & (x_1, x_2) \in \Omega_3,\\
\end{array}\right.
$$
where $\epsilon = 10^{-5}$, $\Gamma_{1}=1.5$, $\Gamma_{2}=1.4$, $p_{\infty, 1} = p_{\infty,2} =0$, and $c_{v, 1}=c_{v, 2}=1$.

The adaptive meshes and densities at $t = 3.5$ and $5$ are plotted in Figures \ref{fig:2Dtri_35} and  \ref{fig:2Dtri_50},
where the monitor function is chosen as \eqref{eq:monitor} 	with
 $\kappa =1, \sigma_1 = \rho$ and $\alpha_1 = 1200$.
The densities  along the line connecting $(2,0)$ and $(7, 3)$ at $ t = 3.5$ and $5$ are shown in Figure \ref{fig:2Dtri_Cmp}.
One can see that   {\tt MM-WENOMR}  with $350 \times 150$ cells outperforms {\tt UM-WENOMR}  with $1050 \times 450$ cells and  {\tt MM-WENOJS} with $350 \times 150$ cells,
and the small scale structures and the Kelvin-Helmholtz instability can be clearly observed by using {\tt MM-WENOMR}.
 The  CPU times in  Table \ref{CPU_78} show that
{\tt MM-WENOMR} with   $350 \times 150$ cells   only takes $18.9\%$ CPU time of  {\tt MM-WENOMR} with $1050 \times 450$ cells and $9.8\%$  CPU time more than {\tt MM-WENOJS} with $350 \times 150$ cells.
	  \begin{figure}[!ht]
	\centering \includegraphics[width=0.6\linewidth]{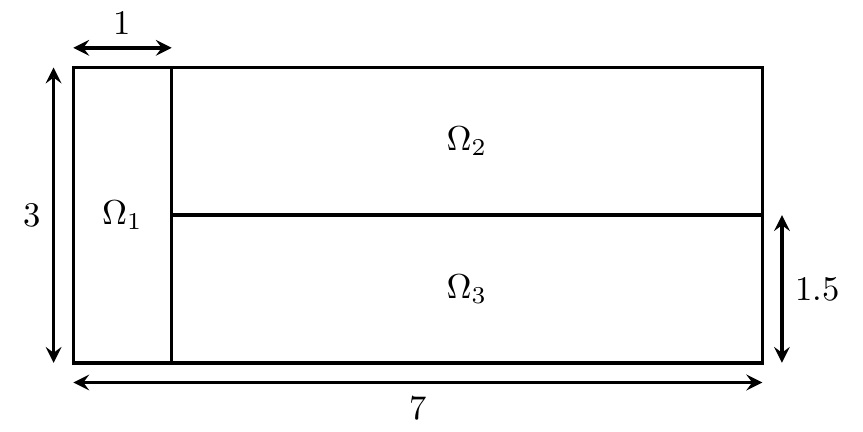}
	\caption{Example  \ref{Tri}.  Initial decomposition of $\Omega_p$.
}
	\label{Tri_Dom}
	
\end{figure}

\begin{figure}[!ht]
	\centering
	\begin{subfigure}[b]{0.32\textwidth}
		\centering
		\includegraphics[width=1.0\linewidth]{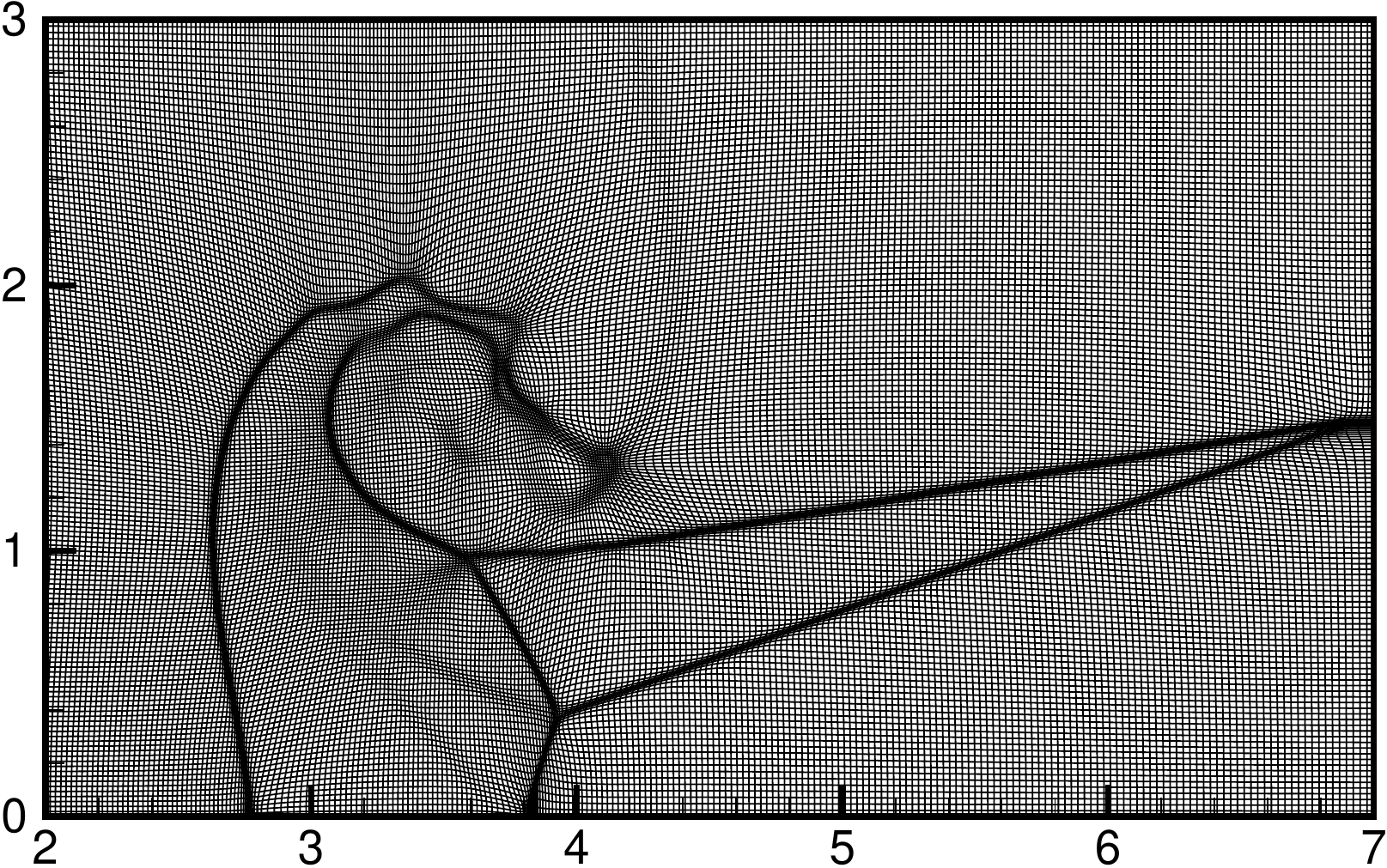}
		\caption{{\tt MM-WENOMR} ($350 \times 150$)}
	\end{subfigure}
	\begin{subfigure}[b]{0.32\textwidth}
		\centering
		\includegraphics[width=1.0\linewidth]{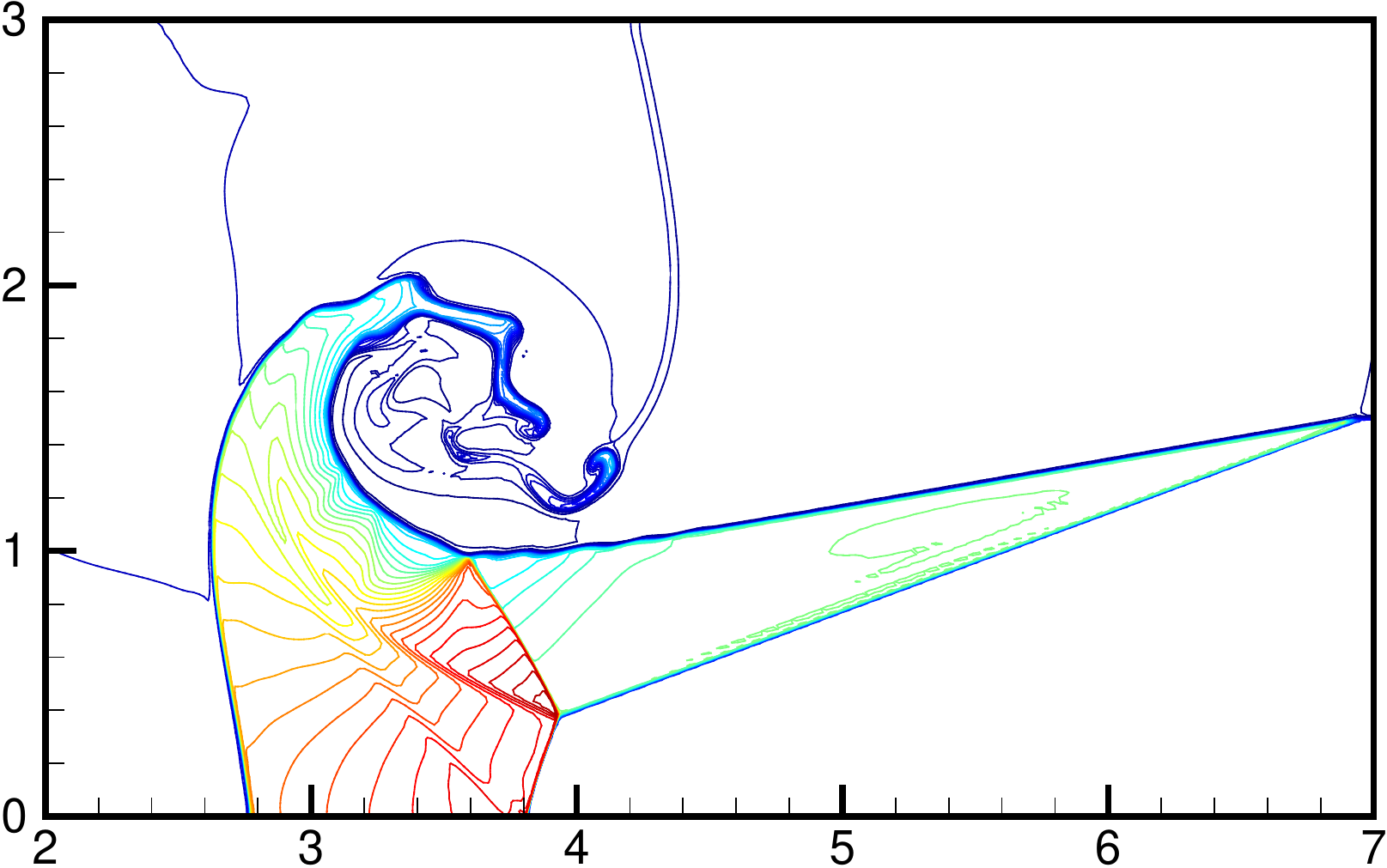}
		\caption{{\tt MM-WENOMR} ($ 350 \times 150$)}
	\end{subfigure}
	\begin{subfigure}[b]{0.32\textwidth}
	\centering
	\includegraphics[width=1.0\linewidth]{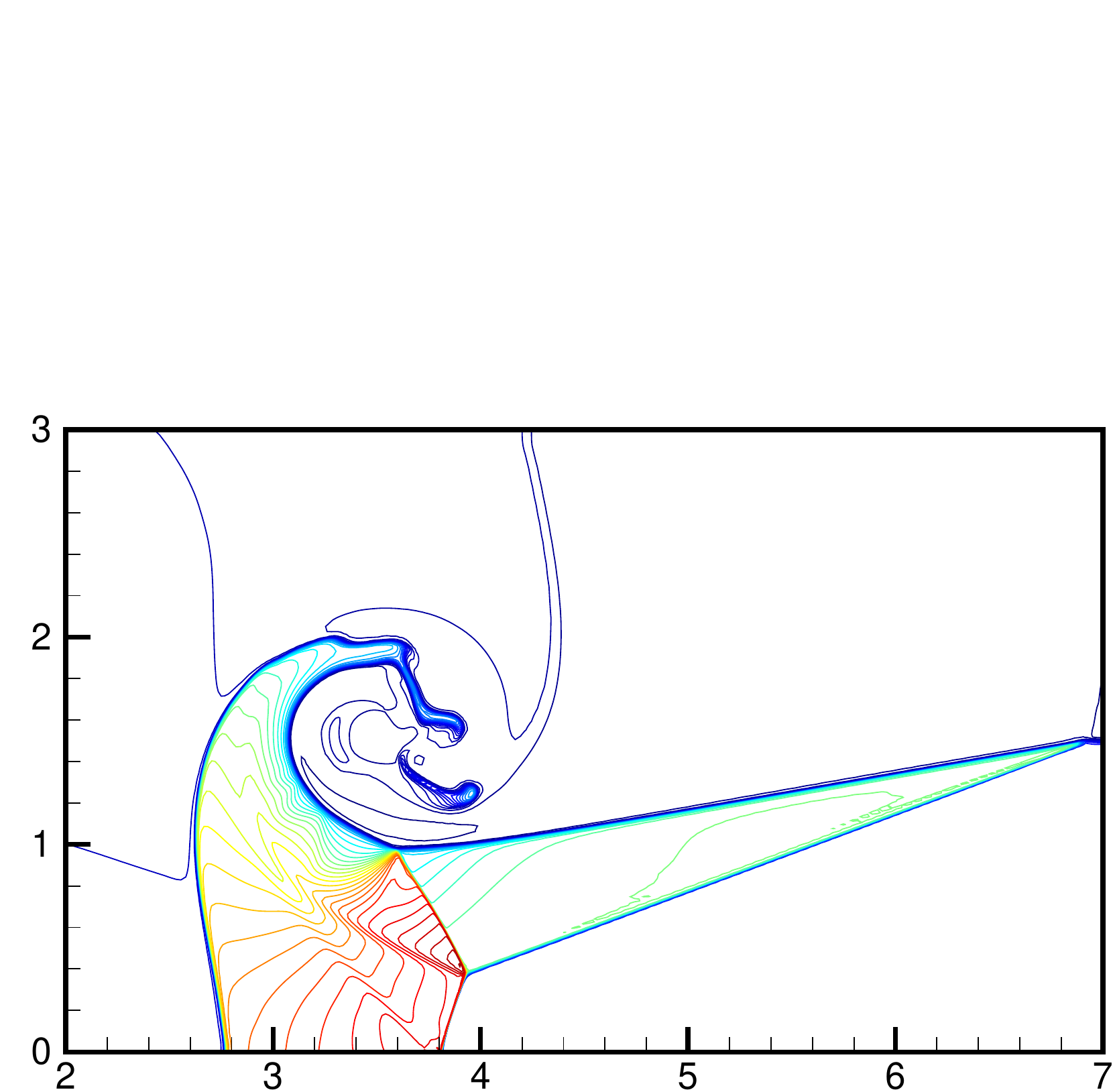}
	\caption{{\tt MM-WENOJS} ($ 350 \times 150$)}
\end{subfigure}
	
	\begin{subfigure}[b]{0.32\textwidth}
		\centering
		\includegraphics[width=1.0\linewidth]{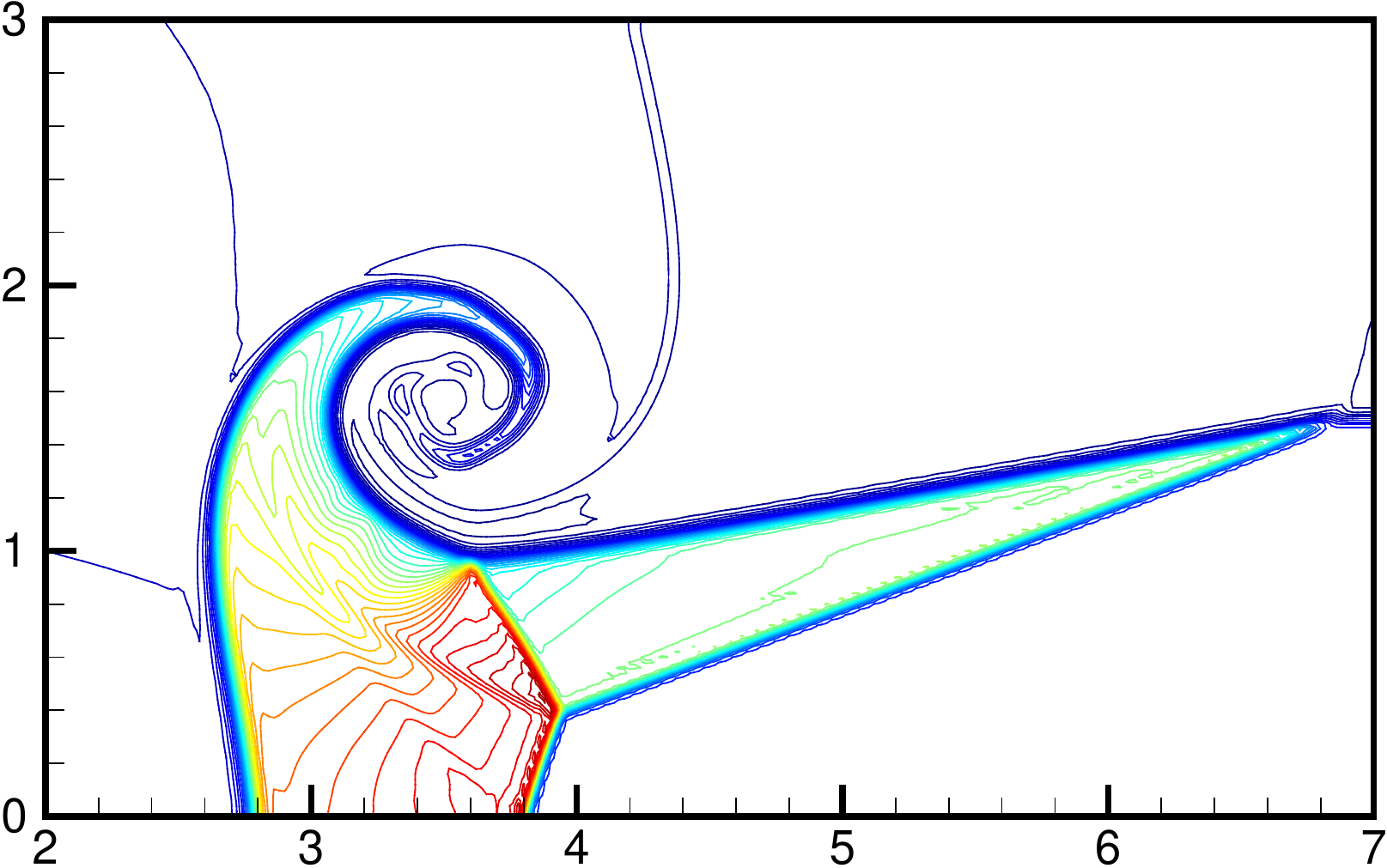}
		\caption{{\tt UM-WENOMR} ($ 350 \times 150$)}
	\end{subfigure}
	\begin{subfigure}[b]{0.32\textwidth}
		\centering
	\includegraphics[width=1.0\linewidth]{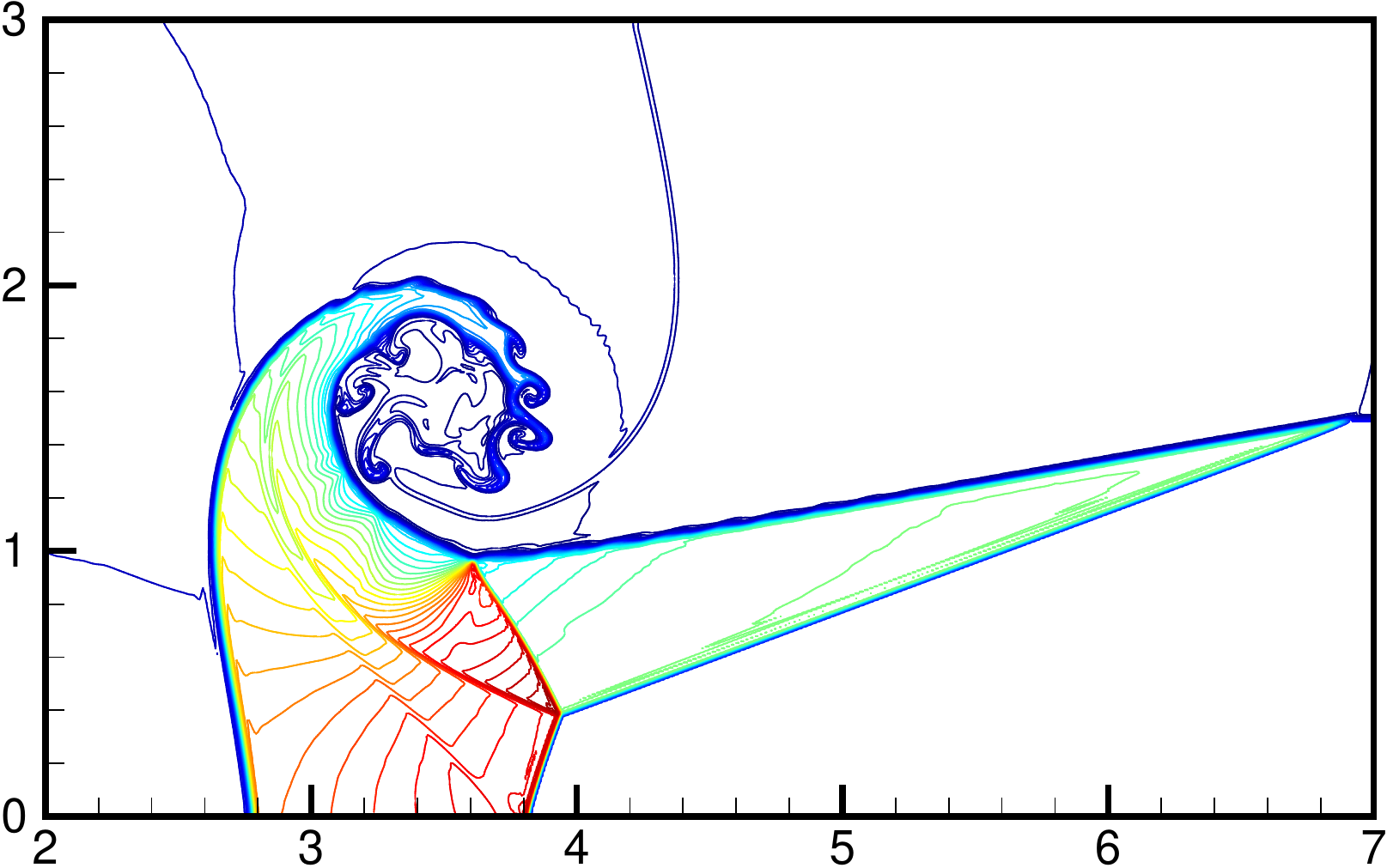}
		\caption{{\tt UM-WENOMR} ($ 1050 \times 450$) }
	\end{subfigure}
	\caption{Example  \ref{ex:tri}. Adaptive mesh of {\tt MM-WENOMR} with $350 \times 150$ cells and   density contours (40 equally spaced contour lines) at $t = 3.5$.}
	\label{fig:2Dtri_35}
	
\end{figure}
\begin{figure}[!ht]
	\centering
	\begin{subfigure}[b]{0.32\textwidth}
		\centering
		\includegraphics[width=1.0\linewidth]{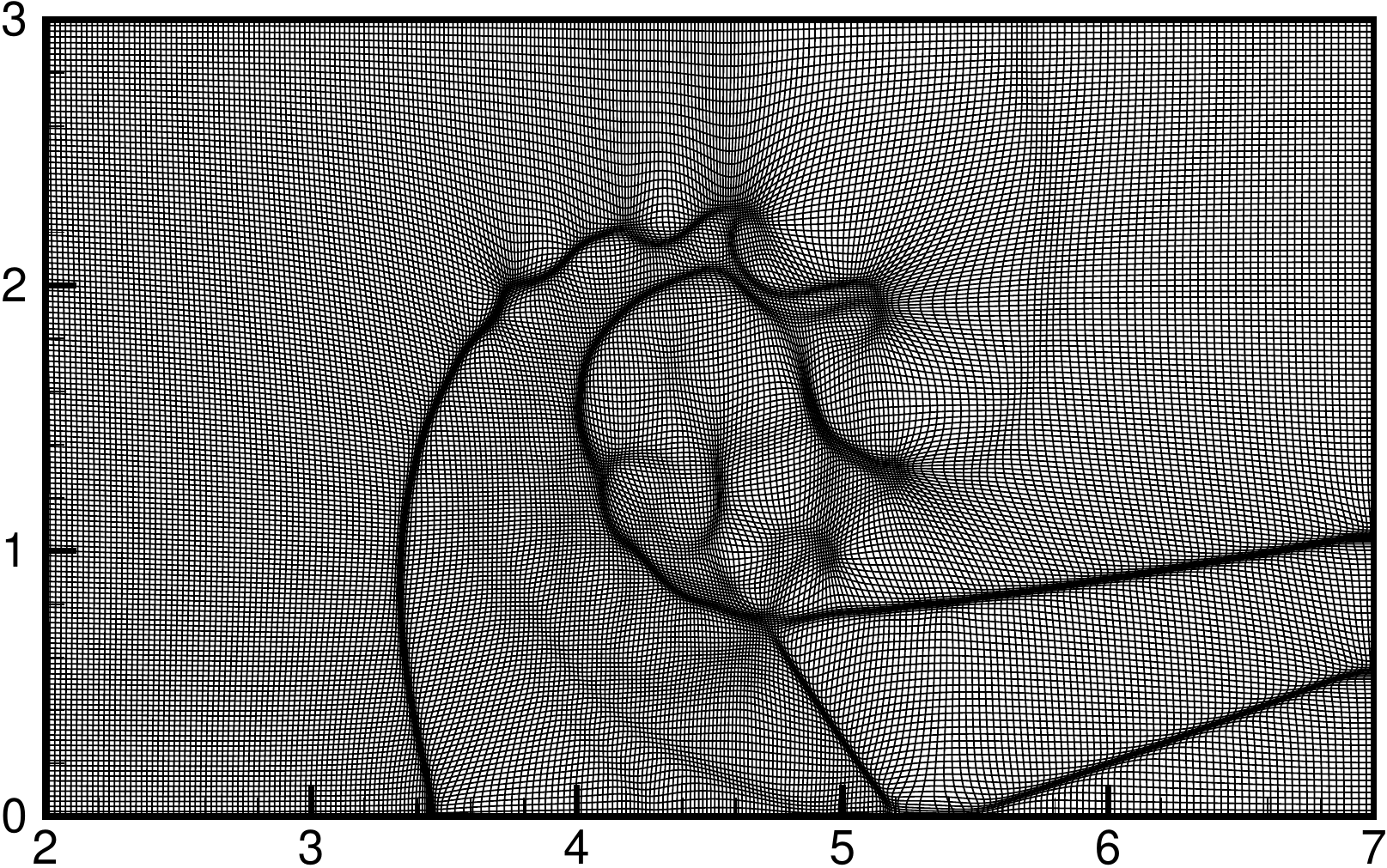}
		\caption{{\tt MM-WENOMR} ($350 \times 150$)}
	\end{subfigure}
	\begin{subfigure}[b]{0.32\textwidth}
		\centering
		\includegraphics[width=1.0\linewidth]{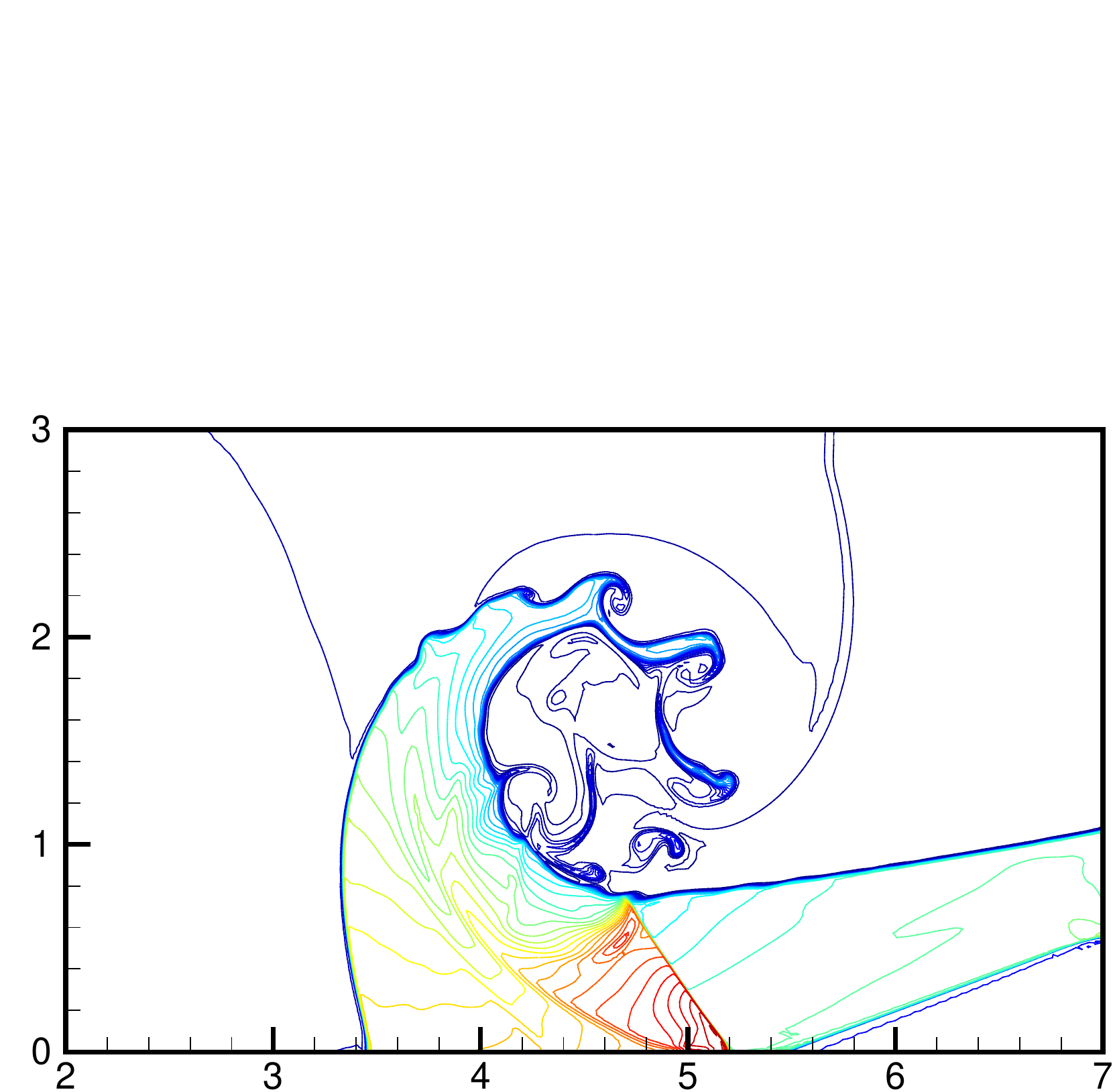}
		\caption{{\tt MM-WENOMR} ($ 350 \times 150$)}
	\end{subfigure}
	\begin{subfigure}[b]{0.32\textwidth}
		\centering
		\includegraphics[width=1.0\linewidth]{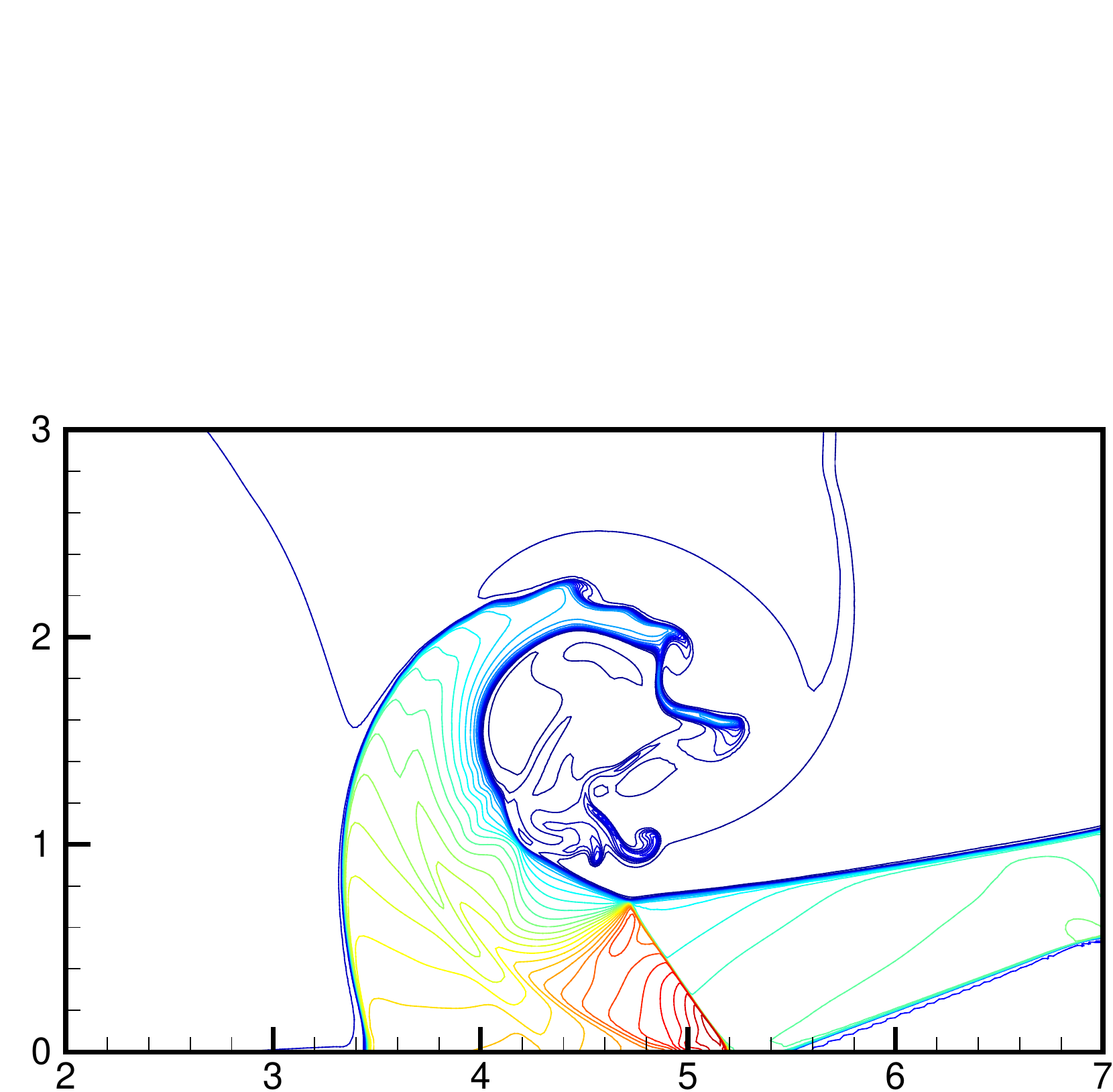}
		\caption{ {\tt MM-WENOJS} ($ 350 \times 150$)}
	\end{subfigure}
	
	\begin{subfigure}[b]{0.32\textwidth}
		\centering
		\includegraphics[width=1.0\linewidth]{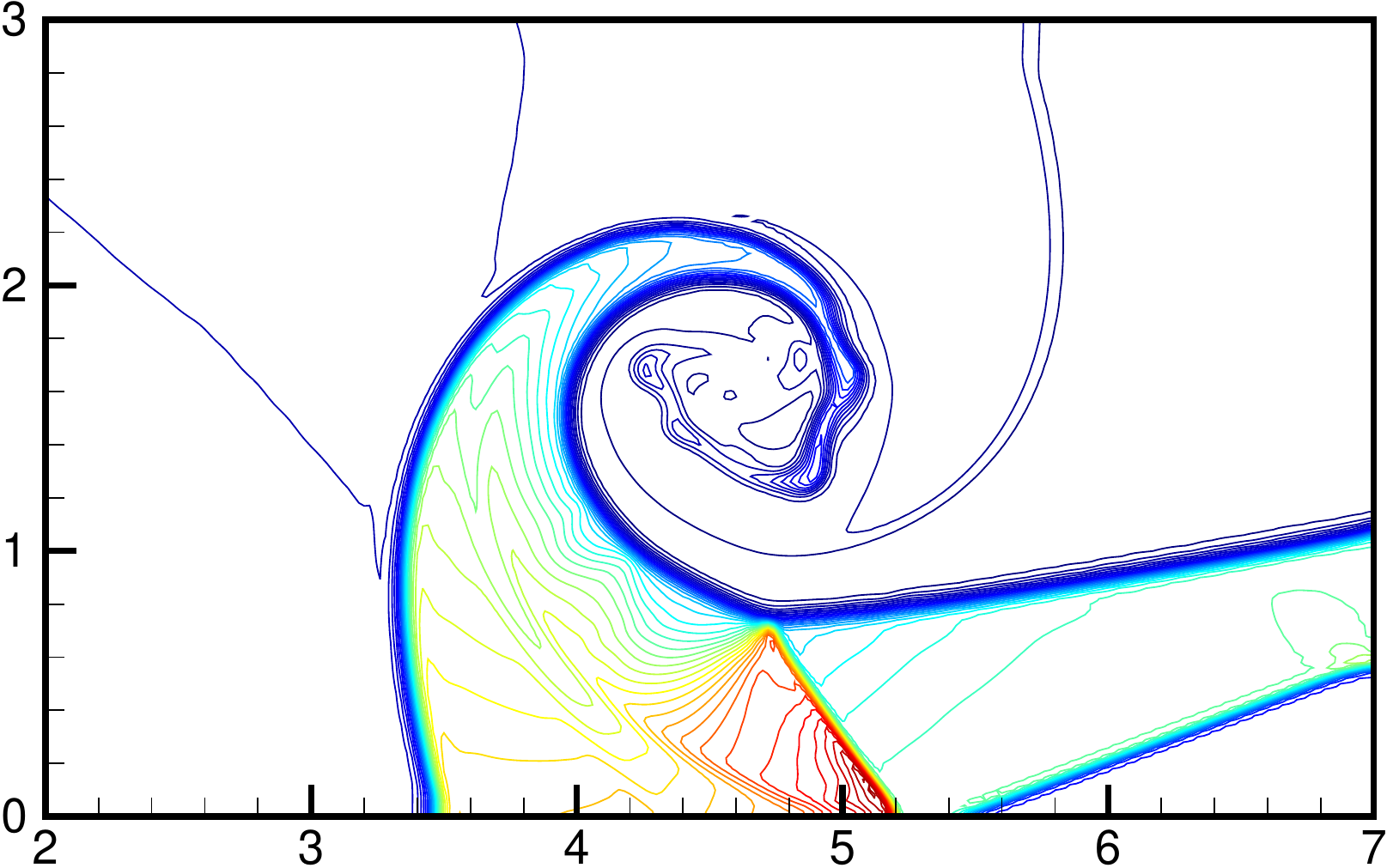}
		\caption{  {\tt UM-WENOMR} ($ 350 \times 150$)}
	\end{subfigure}
	\begin{subfigure}[b]{0.32\textwidth}
		\centering
		\includegraphics[width=1.0\linewidth]{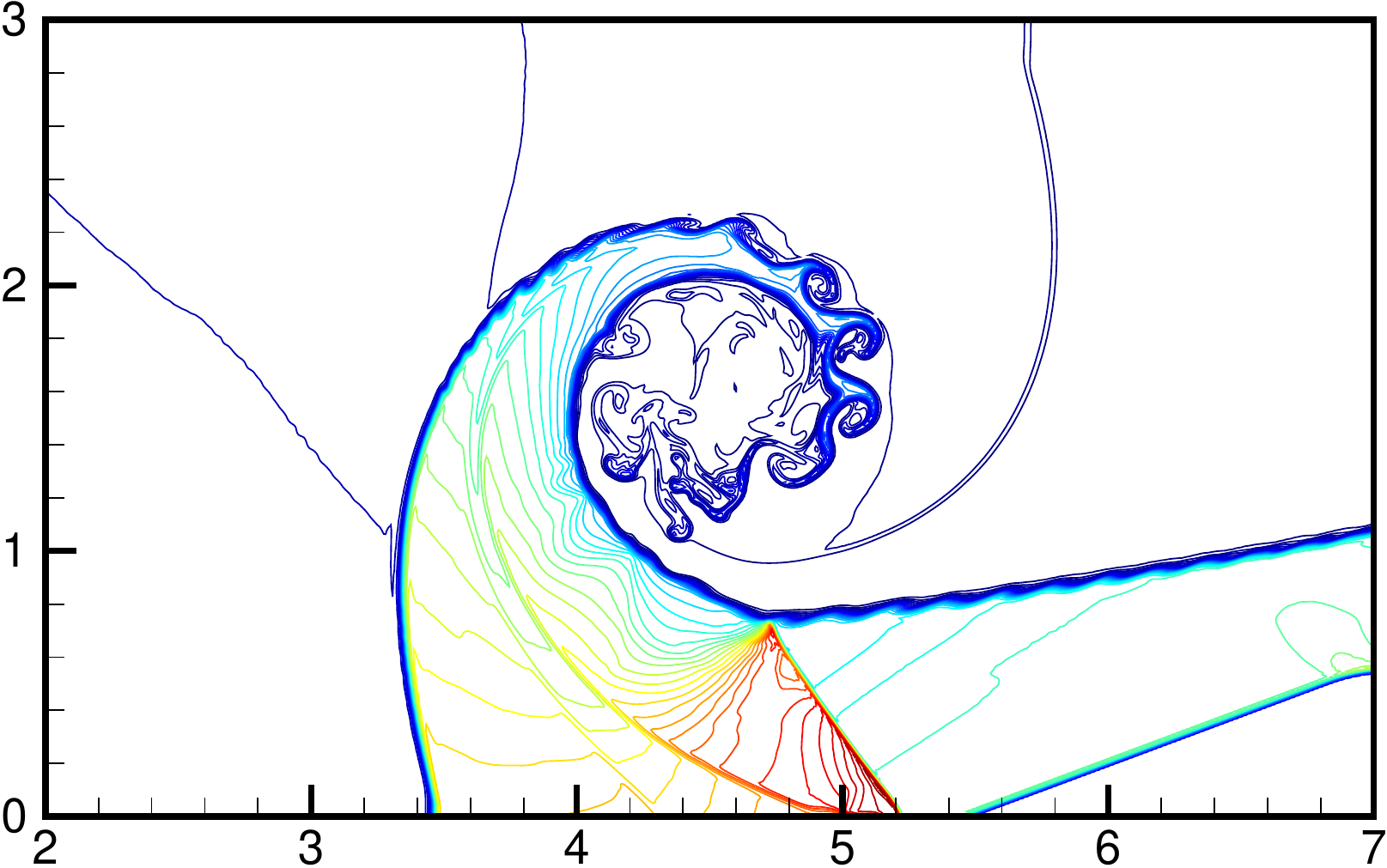}
		\caption{ {\tt UM-WENOMR} ($ 1050 \times 450$) }
	\end{subfigure}
	\caption{Same as Figure \ref{fig:2Dtri_35}, except for $t = 5$.
}
	\label{fig:2Dtri_50}
	
\end{figure}

\begin{figure}[!ht]
	\centering
		\begin{subfigure}[b]{0.33\textwidth}
		\centering
		\includegraphics[width=1.0\linewidth]{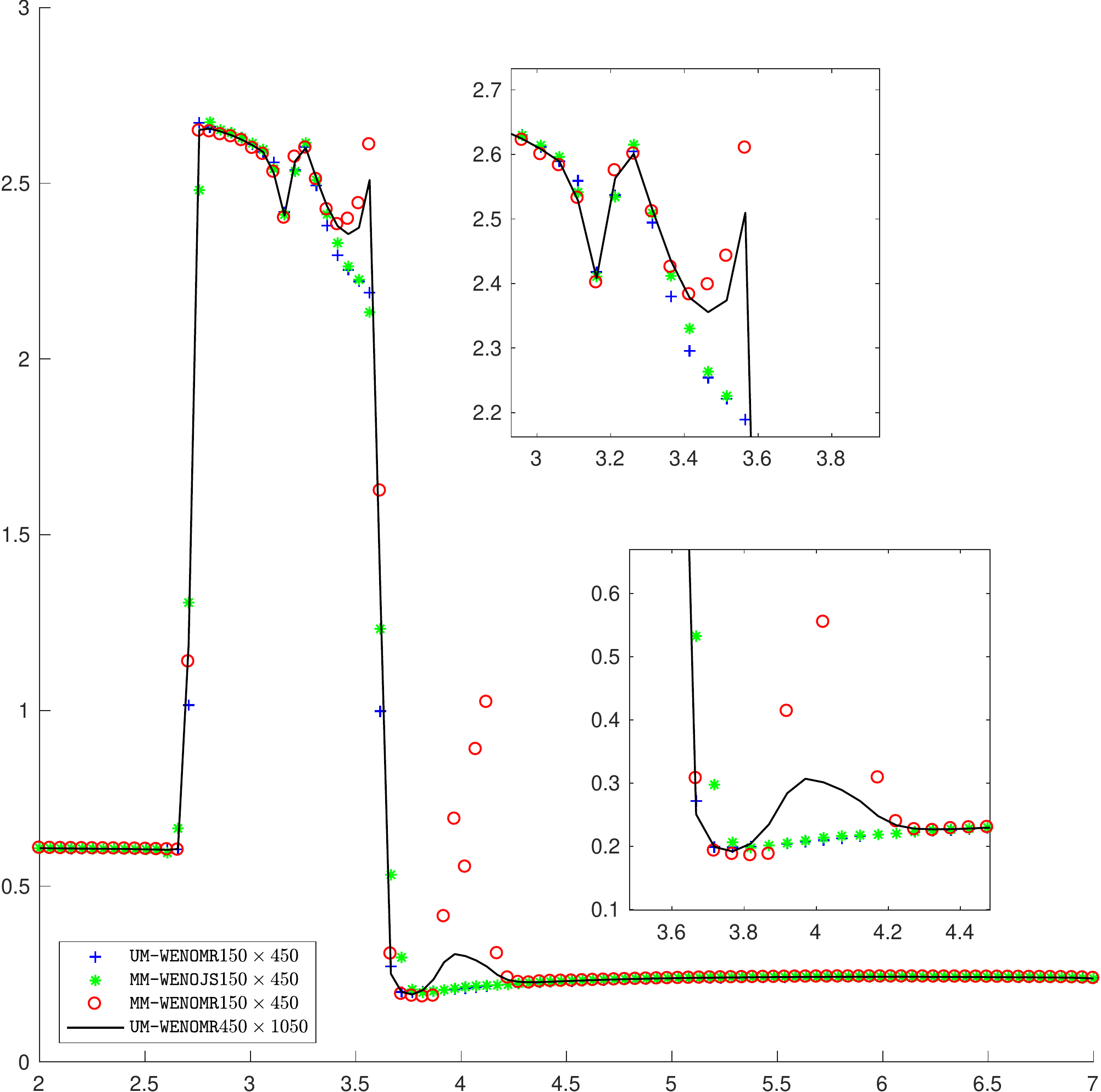}
		\caption{$t = 3.5$}
	\end{subfigure}
	\begin{subfigure}[b]{0.33\textwidth}
		\centering
		\includegraphics[width=1.0\linewidth]{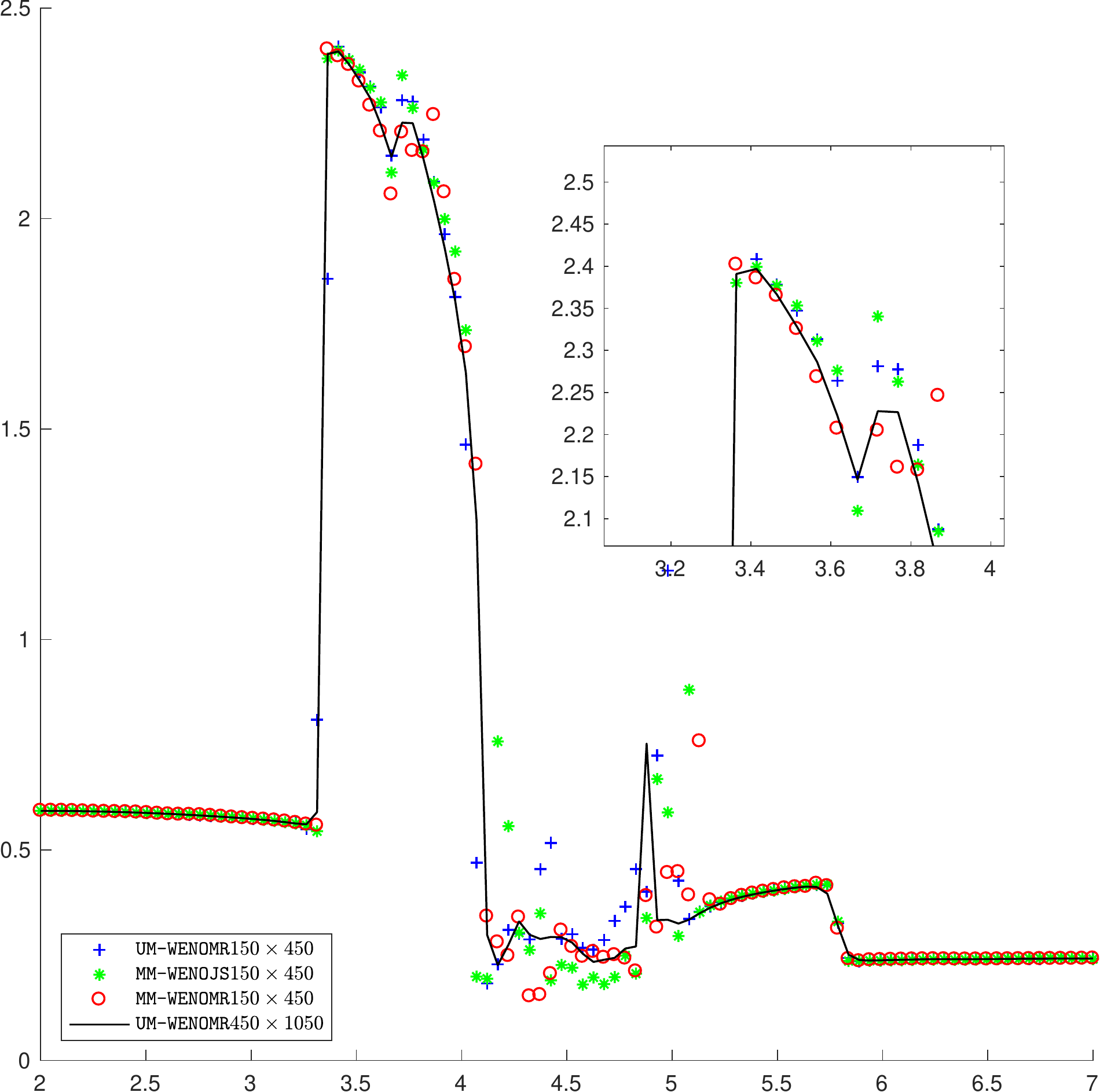}
		\caption{$t = 5$}
	\end{subfigure}
	\caption{Example  \ref{ex:tri}. Densities $\rho$ along the line connecting $(2,0)$ and $(7,3)$ at $t = 3.5$ and $ 5$, respectively.}
	\label{fig:2Dtri_Cmp}
	
\end{figure}

\end{example}

	\begin{example}[2D shock-bubble interaction \uppercase\expandafter{\romannumeral1}]\label{ex:ShockBubble}\rm
	This test is about the interaction of a   shock wave with a helium cylindrical bubble \cite{haas_sturtevant_1987}, and has been frequently applied to  numerical computations \cite{quirk_karni_1996}.
 Initially, the domain $\Omega_p=[0, 445]\times[-44.5, 44.5] $  is decomposed into three sub-domain as illustrated in Figure \ref{ShockBubble_Dom}, and
	a Mach $M_s = 1.22$ shock wave, positioned at $x_1 = 275$, moves 
	through the quiescent air  and will eventually meet a cylindrical helium bubble, centered at $(x_1,x_2) = [225, 0]$ of radius 25, filled with the helium contaminated with $28\%$ of air. 
	Reflecting boundary conditions are specified on the top  and bottom boundaries, while
	 outflow and inflow boundary conditions are applied on the left and right boundaries. The density of the bubble is determined based on the assumption that the regions $\Omega_1$ and $\Omega_2$ are in pressure and temperature equilibrium.
Specially,	the initial data are
	$$
	\left(\rho_{1}, \rho_{2}, v_{1}, v_{2}, p\right)=\left\{\begin{array}{ll}
	\left(\epsilon,1.225\left(R_{1} / R_{2}\right) - \epsilon, 0,0,101325\right), & (x_1, x_2) \in \Omega_1, \\
	(1.225 - \epsilon,\epsilon,0,0,101325), & (x_1, x_2) \in \Omega_2, \\
	(1.6861 - \epsilon,\epsilon,-113.5243,0,159060), & (x_1, x_2) \in \Omega_3,
	\end{array}\right.
	$$
	with
	$\epsilon = 0.03,
     \Gamma_{1}=1.4, \Gamma_{2}=1.647, p_{\infty,1} = p_{\infty,2} = 0, R_{1}=0.287$, and $R_{2}=1.578.
	$ 	
		  \begin{figure}[!ht]
		\centering
		\includegraphics[width=0.7\linewidth]{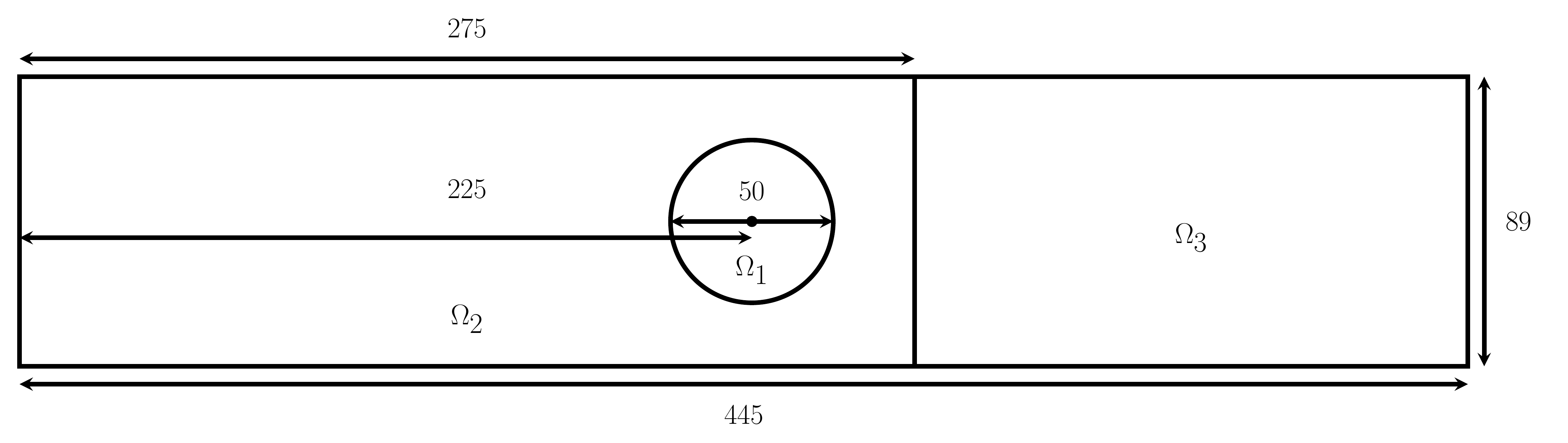}
		\caption{Example  \ref{ex:ShockBubble}. Initial decomposition of   $\Omega_p$.} 
		\label{ShockBubble_Dom}
		\end{figure}

%
 Figures \ref{fig:0}-\ref{fig:2} show the adaptive meshes and   the schlieren images   by the function
	\begin{align}
	\label{eq:schlieren}
	\Phi=\exp \left(-\Psi|\nabla \rho| /|\nabla \rho|_{\max }\right),
	\end{align}
	 with $\Psi=\left(30 \rho_{1}+150 \rho_{2}\right) / \rho$
	obtained by {\tt MM-WENOMR} and {\tt UM-WENOMR}   %
	  at $t=0.02$,  $0.052$, $0.076$,  $0.26$,
	 $0.452$, and $t=0.676$ (after the bubble is first hit by the incident shock wave),
where
 	the monitor function is chosen as
	 \eqref{eq:monitor}
	 	with
	 $\kappa = 1, \sigma_1 = \rho$ and $ \alpha_1 = 1000$.
   Notice that the top and bottom half parts of the  schlieren image are the results obtained
    respectively by {\tt MM-WENOMR} and {\tt UM-WENOMR},
and
 the	velocity of the shock wave is $-415.16$, thus it takes about   $t=0.06$ for the shock wave to meet the
	bubble. 
 We see that  the  mesh {points adaptively} concentrate near the large gradient area of the density
 and {\tt MM-WENOMR} captures the sharp bubble interfaces   and some small wave structures well. Table \ref{CPU_MC_Shock_bubble} tells us that
      {\tt MM-WENOMR} costs  $26.2
	\%$ CPU time of {\tt UM-WENOMR} with a finer mesh, when it obtains  even better results.
	\begin{table}
	\centering
	 \resizebox{.95\columnwidth}{!}{
	\begin{tabular}{l|ccc}
			\hline & {\tt MM-WENOMR}  & {\tt UM-WENOMR}& {\tt UM-WENOMR}  \\
		\hline Example \ref{ex:ShockBubble}&  $1 \mathrm{h} 4 \mathrm{m} $ ($800\times 160$ cells) 
		 & $ 9 \mathrm{m} 12\mathrm{s}$ ($800\times 160$ cells)  & $4 \mathrm{h} 4 \mathrm{m} $ ($2400\times480$ cells)\\
		 \hline Example \ref{ex:ShockBubble2}&  $ 9 \mathrm{m} 28\mathrm{s}$ ($800\times 160$ cells)    & $  1\mathrm{m} 47\mathrm{s}$ ($800\times 160$ cells)  & $40 \mathrm{m} 11\mathrm{s}$ ($2400\times480$ cells)\\
		\hline
		Example \ref{ex:3DShockBubble} &$ 11\mathrm{h} 6\mathrm{m}$   ($400 \times 80 \times 80$ cells)
		& $2\mathrm{h} 34\mathrm{m}$ ($400 \times 80 \times 80$ cells) & $38\mathrm{h} 2\mathrm{m}$ ($800 \times 160 \times 160$ cells) \\ \hline
		 Example \ref{ex:3DShockBubble2} &$ 2\mathrm{h} 9\mathrm{m}$  ($400 \times 80 \times 80$ cells)&$  46\mathrm{m} 32\mathrm{s}$   ($400 \times 80 \times 80$ cells)   & $11\mathrm{h} 58\mathrm{m}$ ($800 \times 160 \times 160$ cells) \\ \hline
	\end{tabular}
}
	\caption{CPU times of Examples
			\ref{ex:ShockBubble}-\ref{ex:3DShockBubble2} ($32$ cores).
	}
	\label{CPU_MC_Shock_bubble}
\end{table}
\begin{figure}[!ht]
	\centering
	
	\begin{subfigure}[b]{0.32\textwidth}
		\centering
		\includegraphics[width=1.8in,height=1.78in]{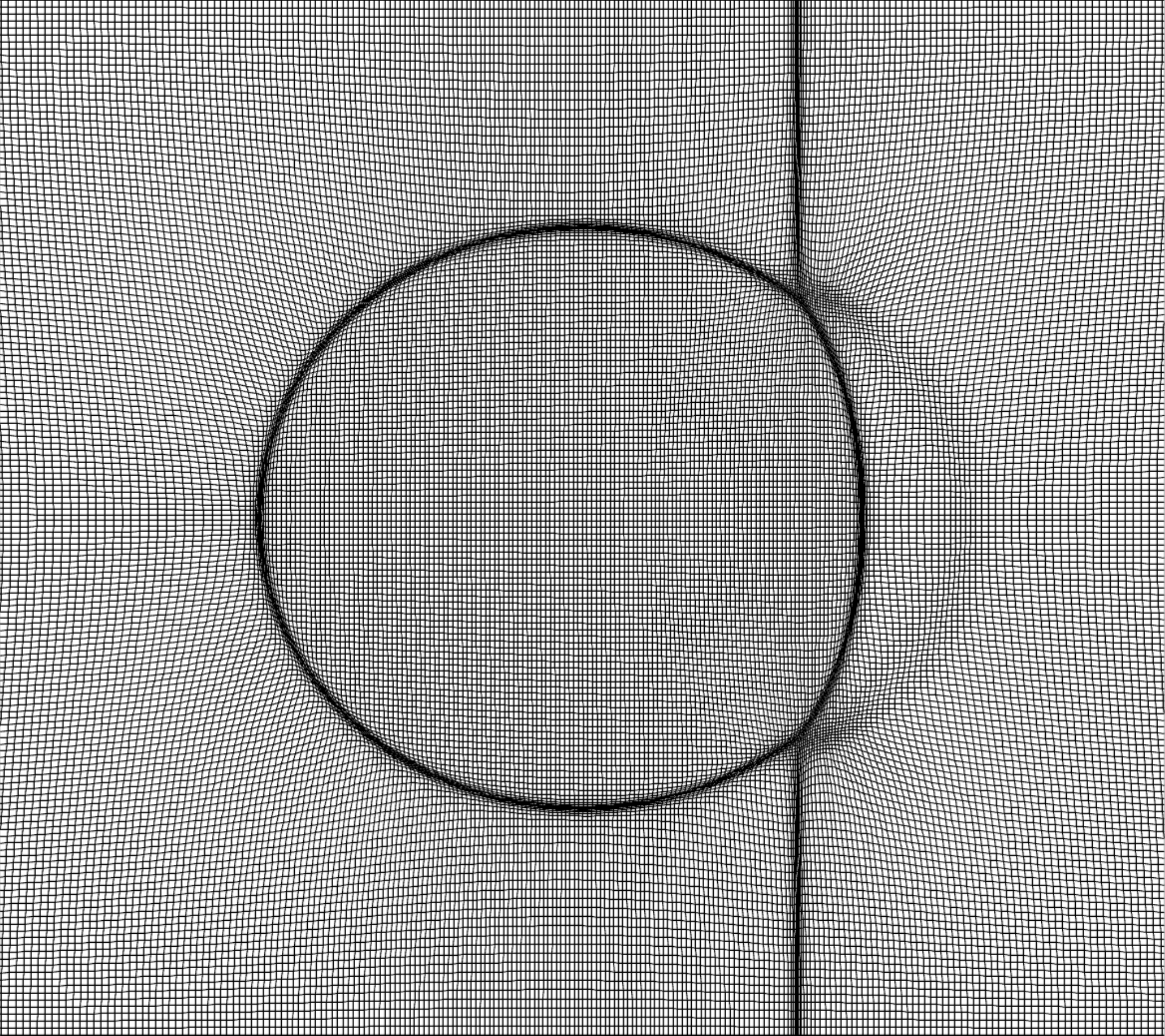}
	\end{subfigure}
	\begin{subfigure}[b]{0.32\textwidth}
		\centering
		\includegraphics[width=1.8in,height=1.78in]{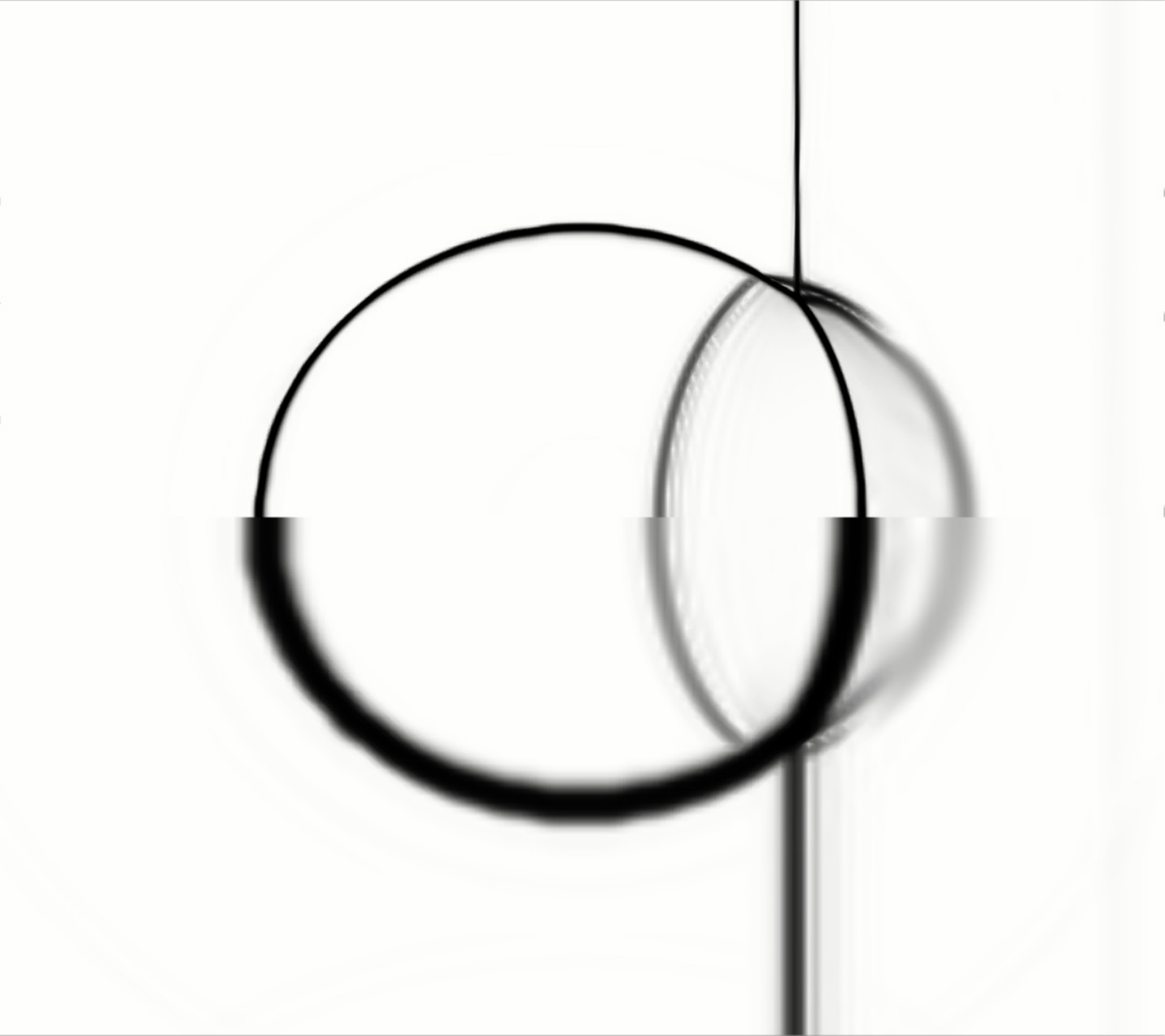}
	\end{subfigure}
	\begin{subfigure}[b]{0.32\textwidth}
		\centering
		\includegraphics[width=1.8in,height=1.78in]{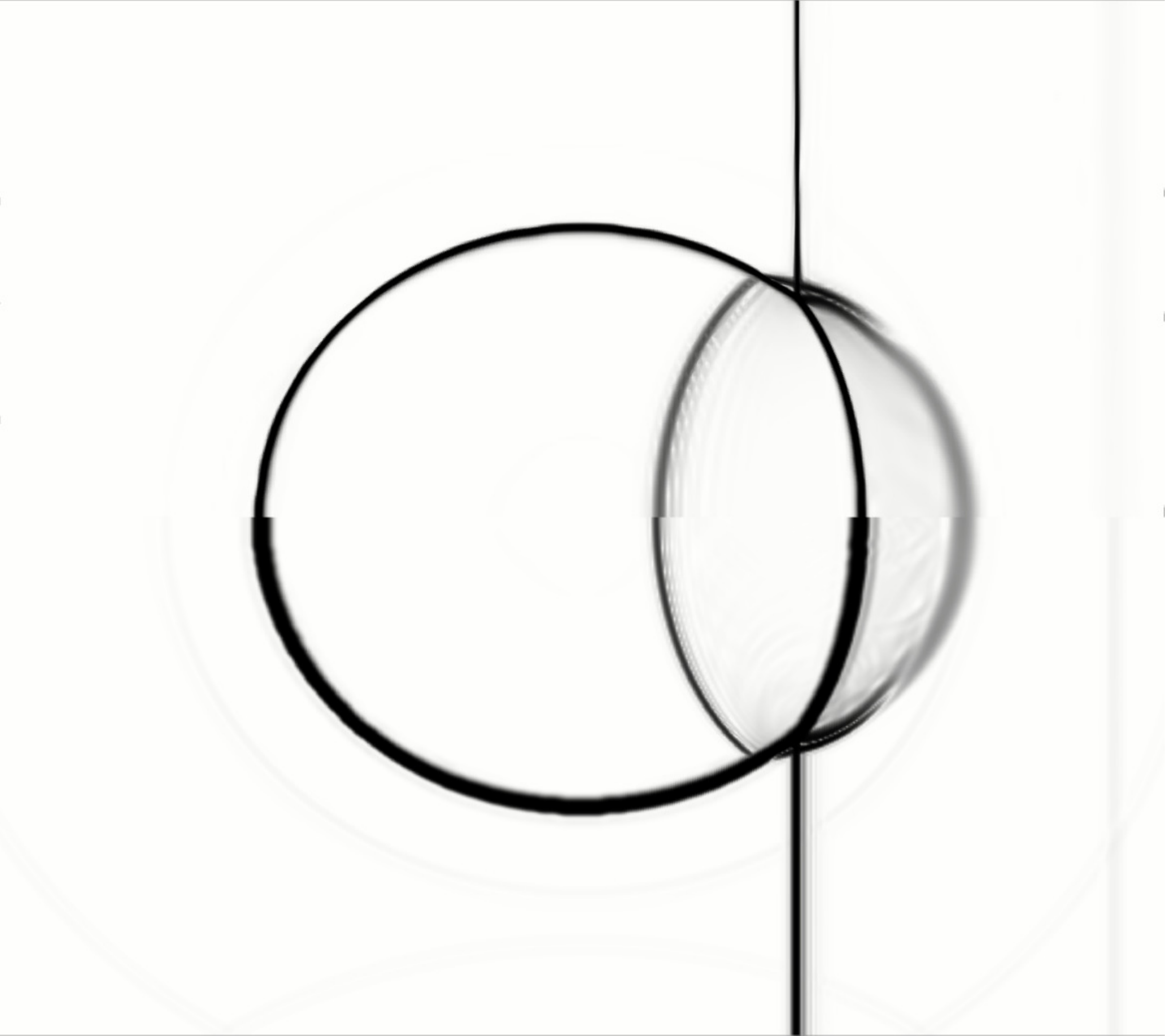}
	\end{subfigure}

	
	\begin{subfigure}[b]{0.32\textwidth}
		\centering
		\includegraphics[width=1.8in,height=1.78in]{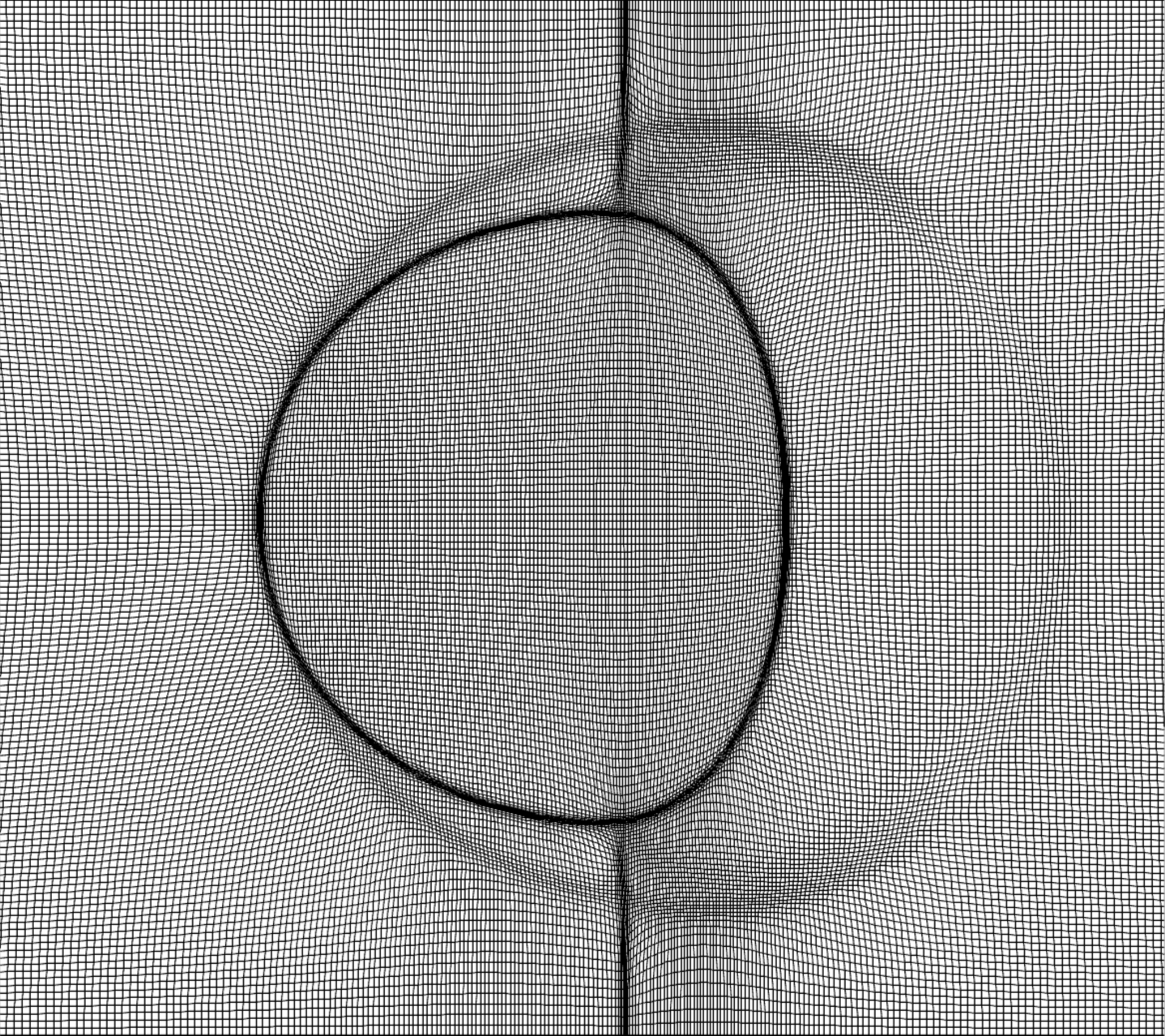}
	\end{subfigure}
	\begin{subfigure}[b]{0.32\textwidth}
		\centering
		\includegraphics[width=1.8in,height=1.78in]{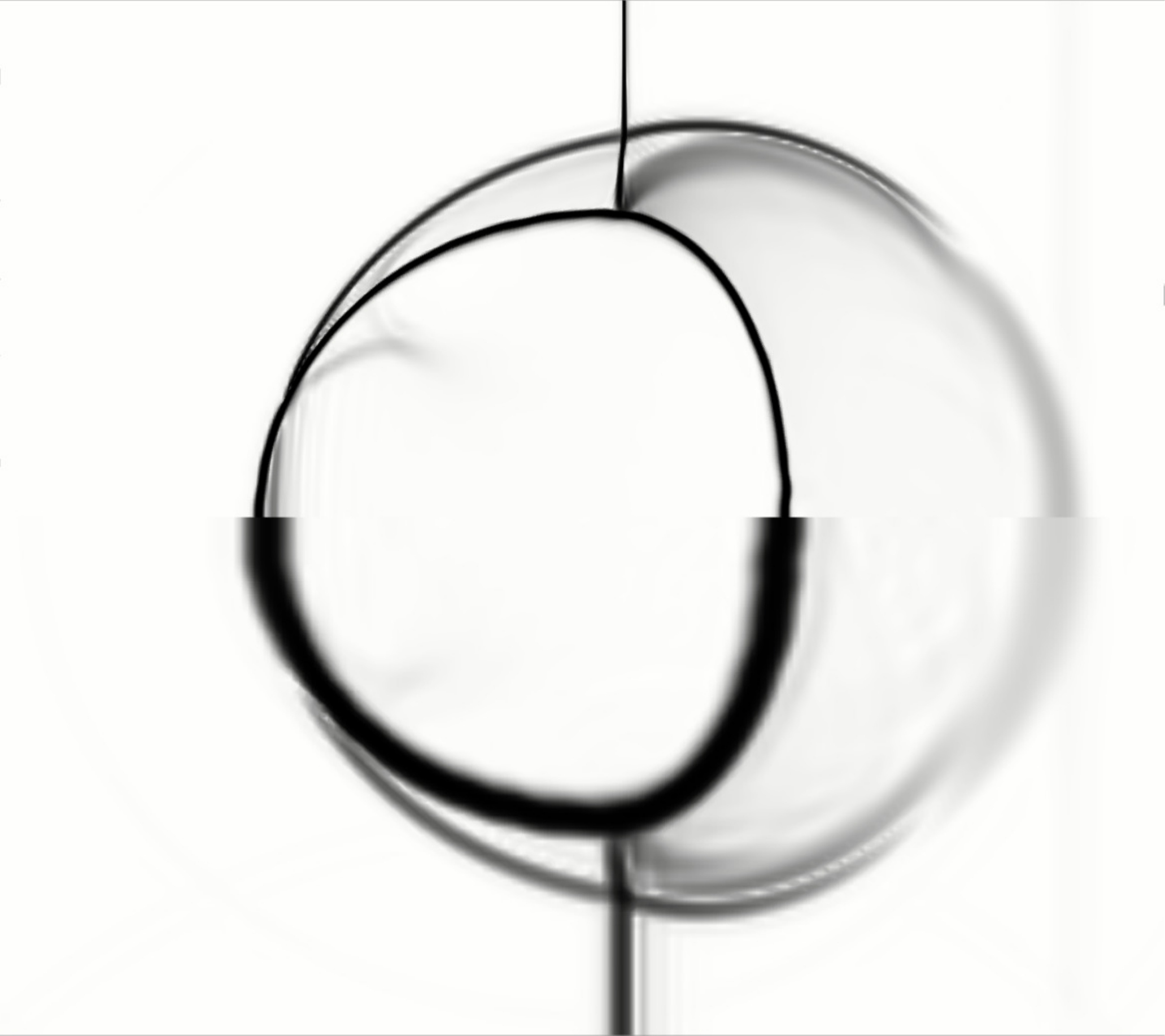}
	\end{subfigure}
	\begin{subfigure}[b]{0.32\textwidth}
		\centering
		\includegraphics[width=1.8in,height=1.78in]{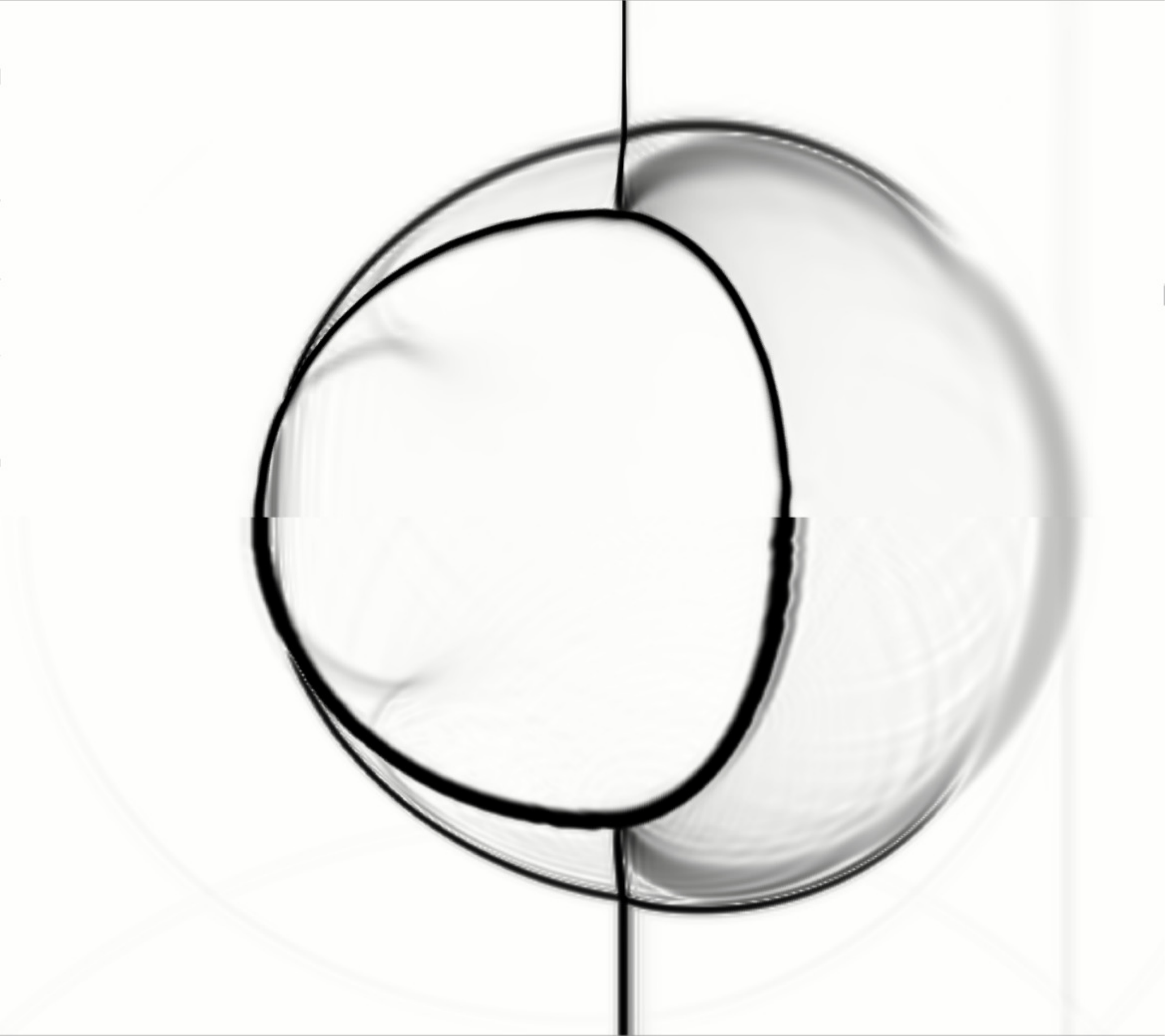}
	\end{subfigure}

\begin{subfigure}[b]{0.32\textwidth}
	\centering
	\includegraphics[width=1.8in,height=1.78in]{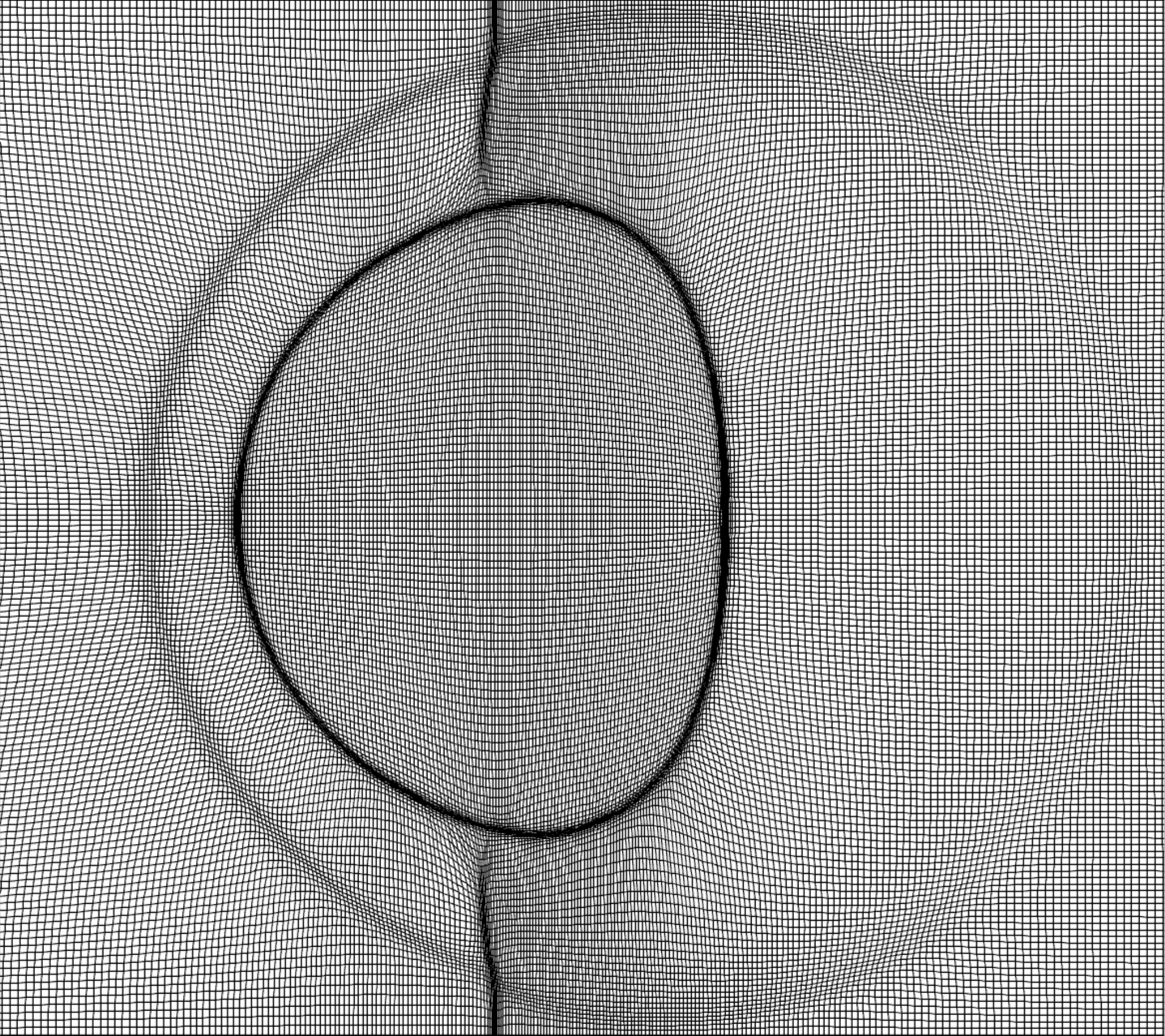}
\end{subfigure}
\begin{subfigure}[b]{0.32\textwidth}
	\centering
	\includegraphics[width=1.8in,height=1.78in]{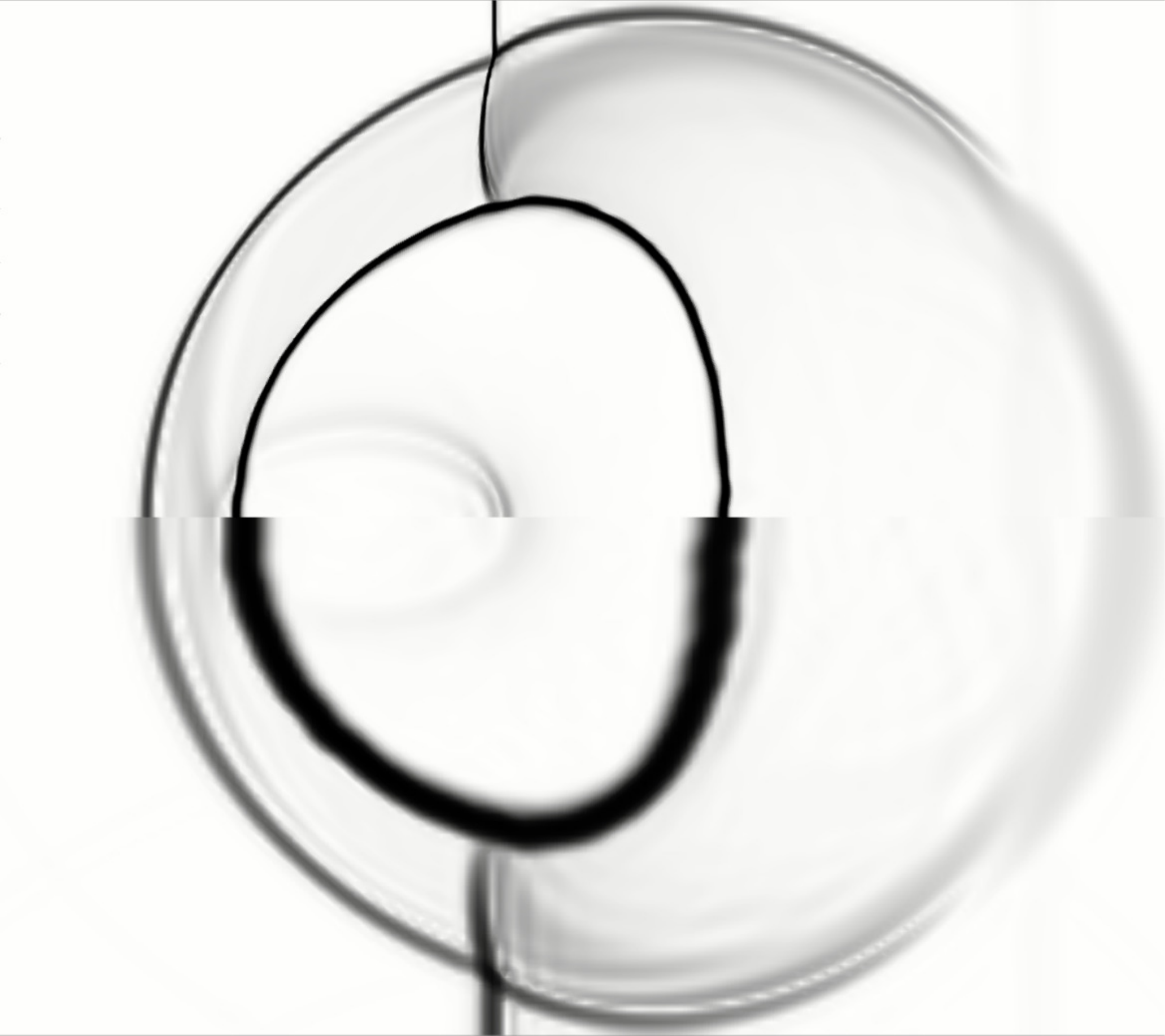}
\end{subfigure}
\begin{subfigure}[b]{0.32\textwidth}
	\centering
	\includegraphics[width=1.8in,height=1.78in]{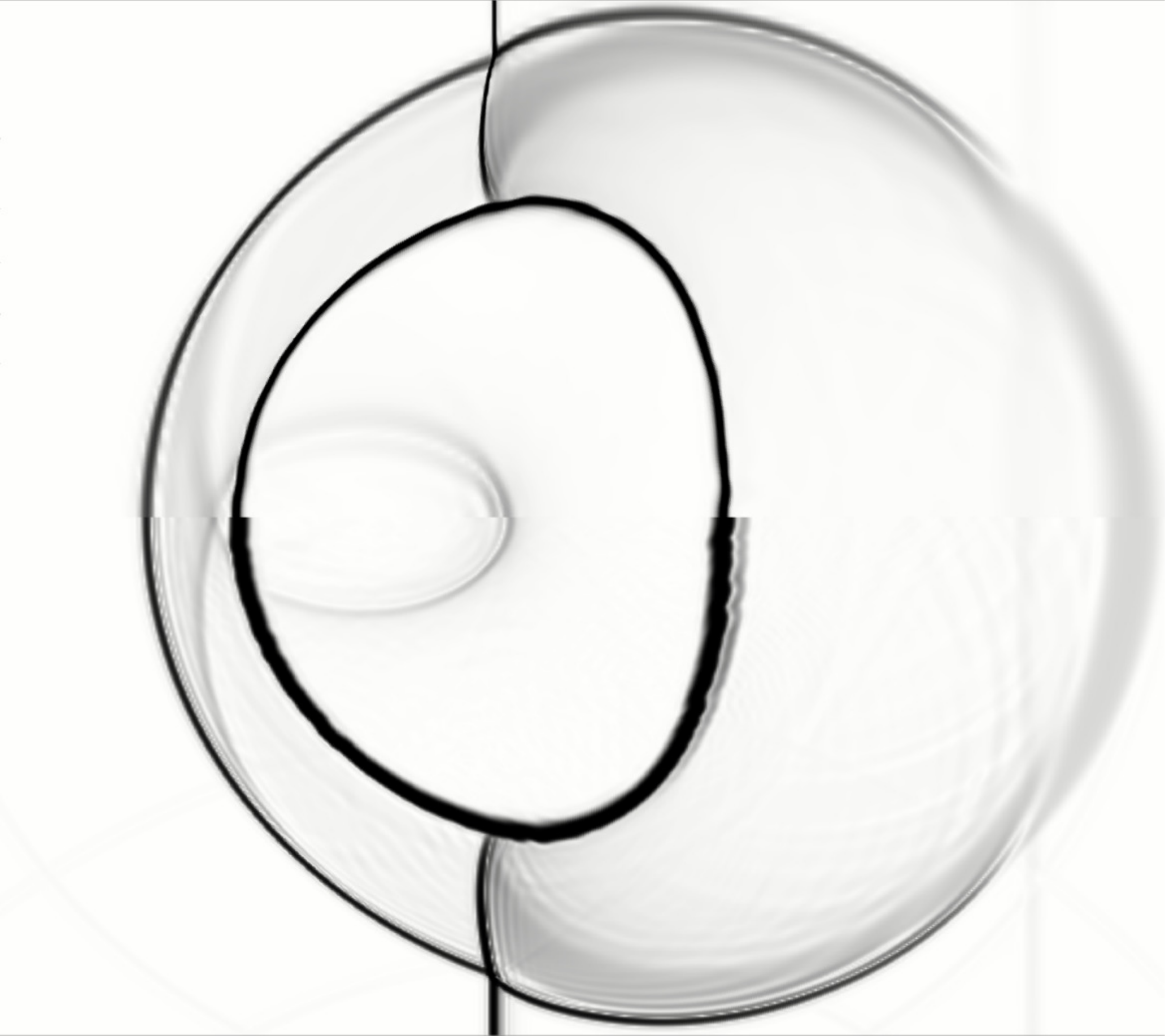}
\end{subfigure}
   \caption{Example \ref{ex:ShockBubble}. From left to right: adaptive meshes of $800 \times 160$ cells,  schlieren images  obtained by   {\tt MM-WENOMR} with $800 \times 160$ cells (top half) and   {\tt UM-WENOMR} with $800 \times 160$ cells (bottom half), and
   schlieren images by {\tt MM-WENOMR}  with $800 \times 160$ cells (top half) and   {\tt UM-WENOMR} with $2400 \times 480$ cells (bottom half). From top to bottom:  $t = 0.02, 0.052, 0.076$.}
   \label{fig:0}
	
\end{figure}

\begin{figure}[!ht]
	\centering

	\begin{subfigure}[b]{0.32\textwidth}
		\centering
		\includegraphics[width=1.8in,height=1.78in]{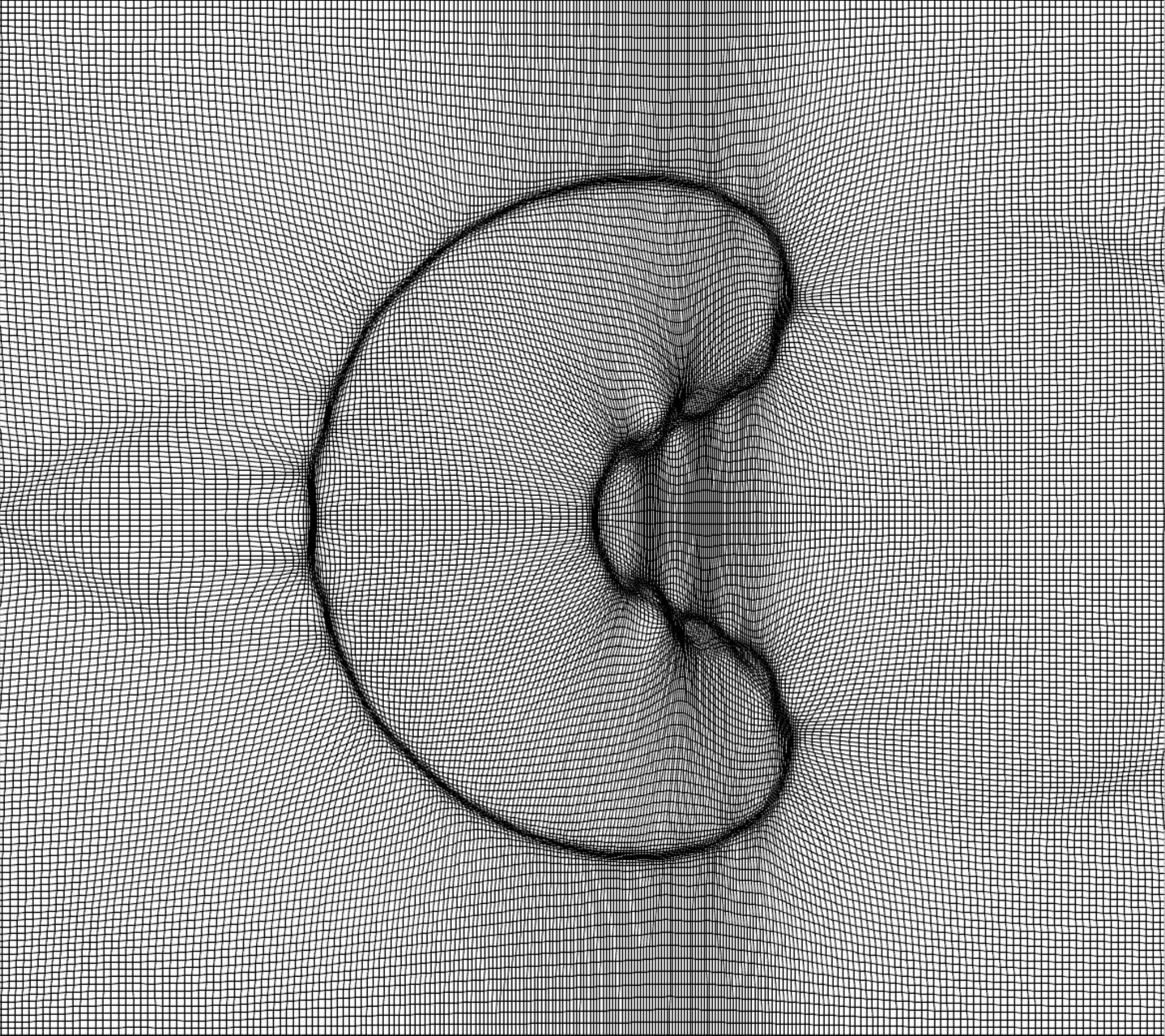}
	\end{subfigure}
	\begin{subfigure}[b]{0.32\textwidth}
		\centering
		\includegraphics[width=1.8in,height=1.78in]{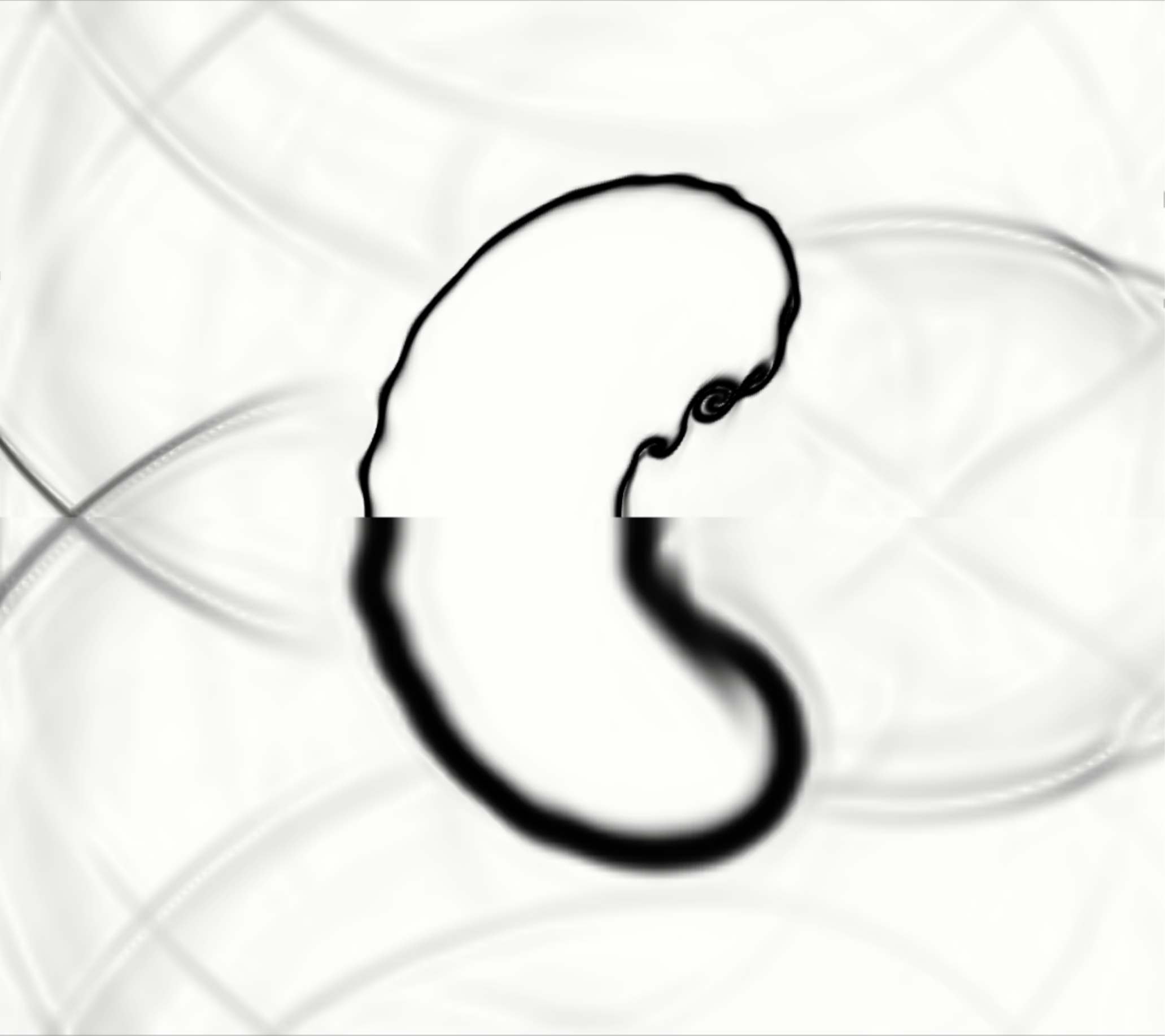}
	\end{subfigure}
	\begin{subfigure}[b]{0.32\textwidth}
		\centering
		\includegraphics[width=1.8in,height=1.78in]{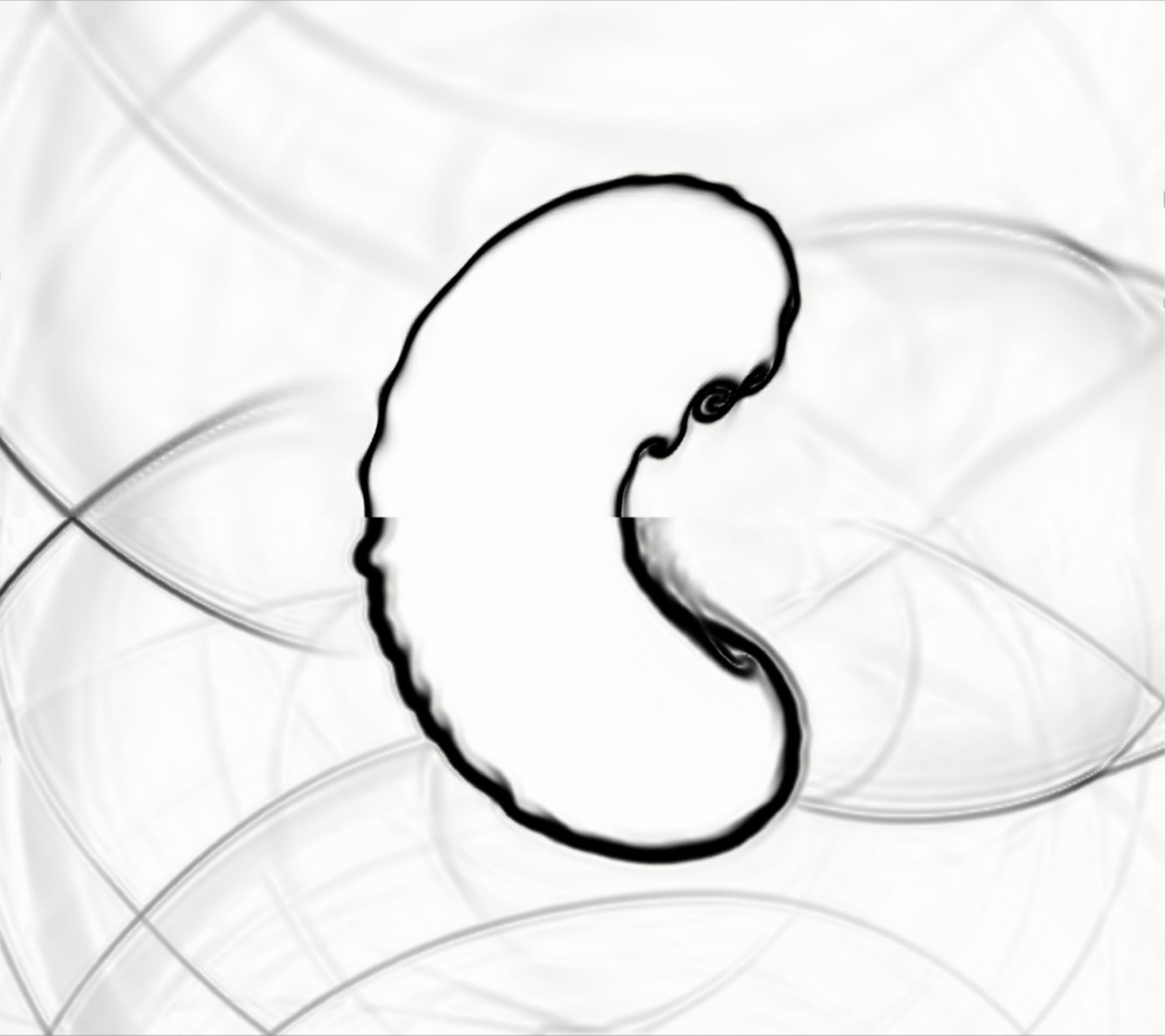}
	\end{subfigure}

	\begin{subfigure}[b]{0.32\textwidth}
		\centering
		\includegraphics[width=1.8in,height=1.78in]{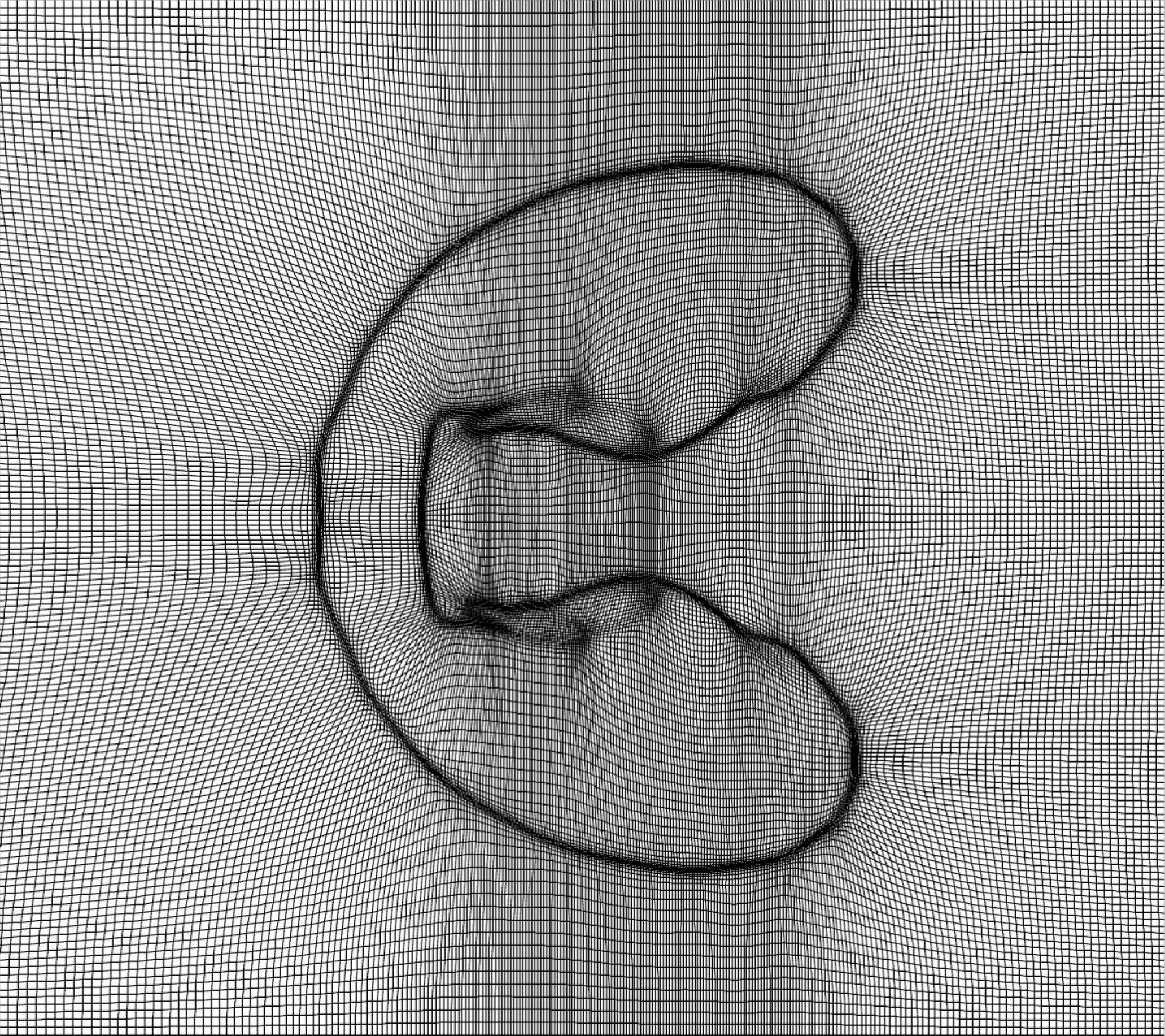}
	\end{subfigure}
	\begin{subfigure}[b]{0.32\textwidth}
		\centering
		\includegraphics[width=1.8in,height=1.78in]{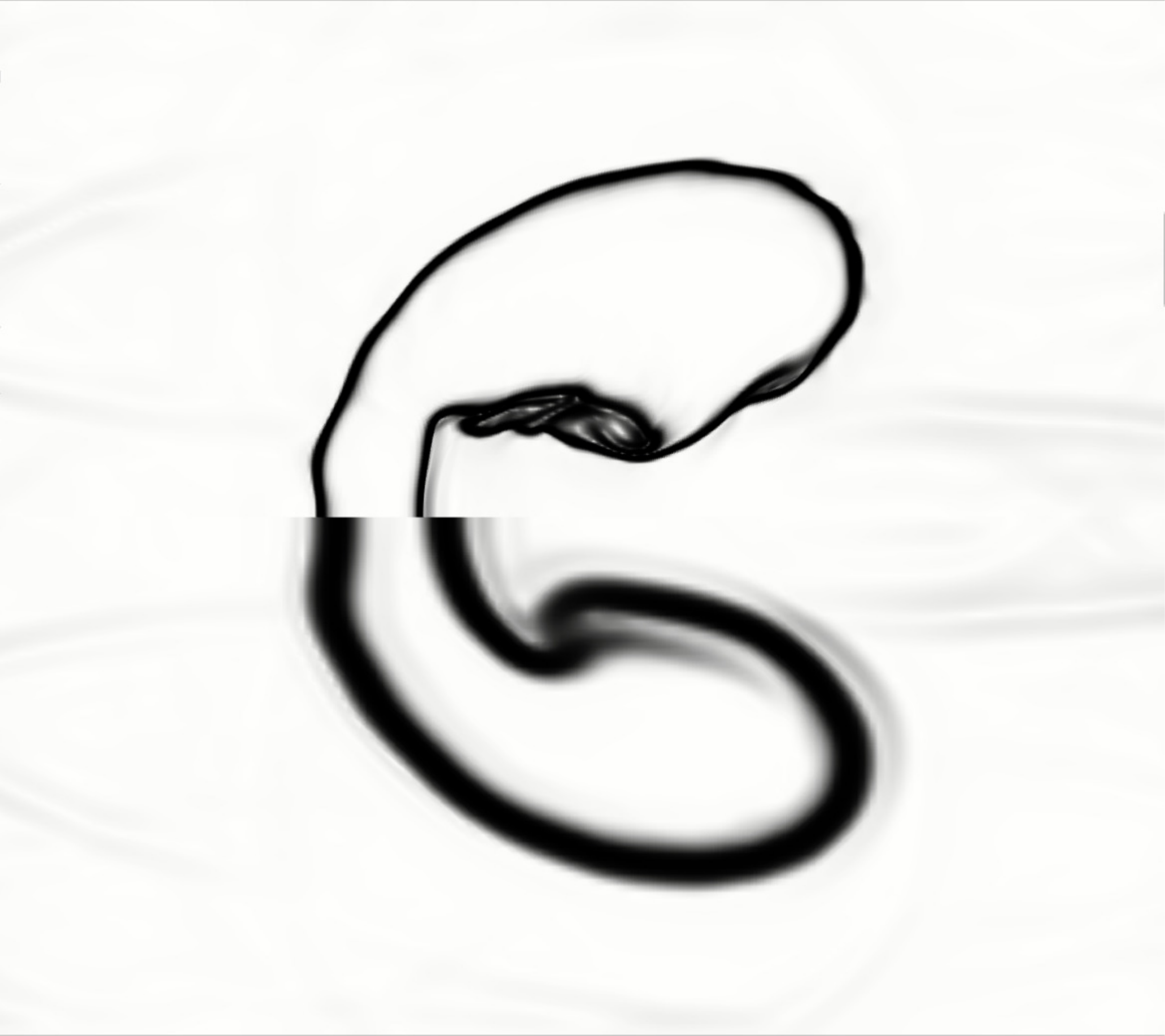}
	\end{subfigure}
	\begin{subfigure}[b]{0.32\textwidth}
		\centering
		\includegraphics[width=1.8in,height=1.78in]{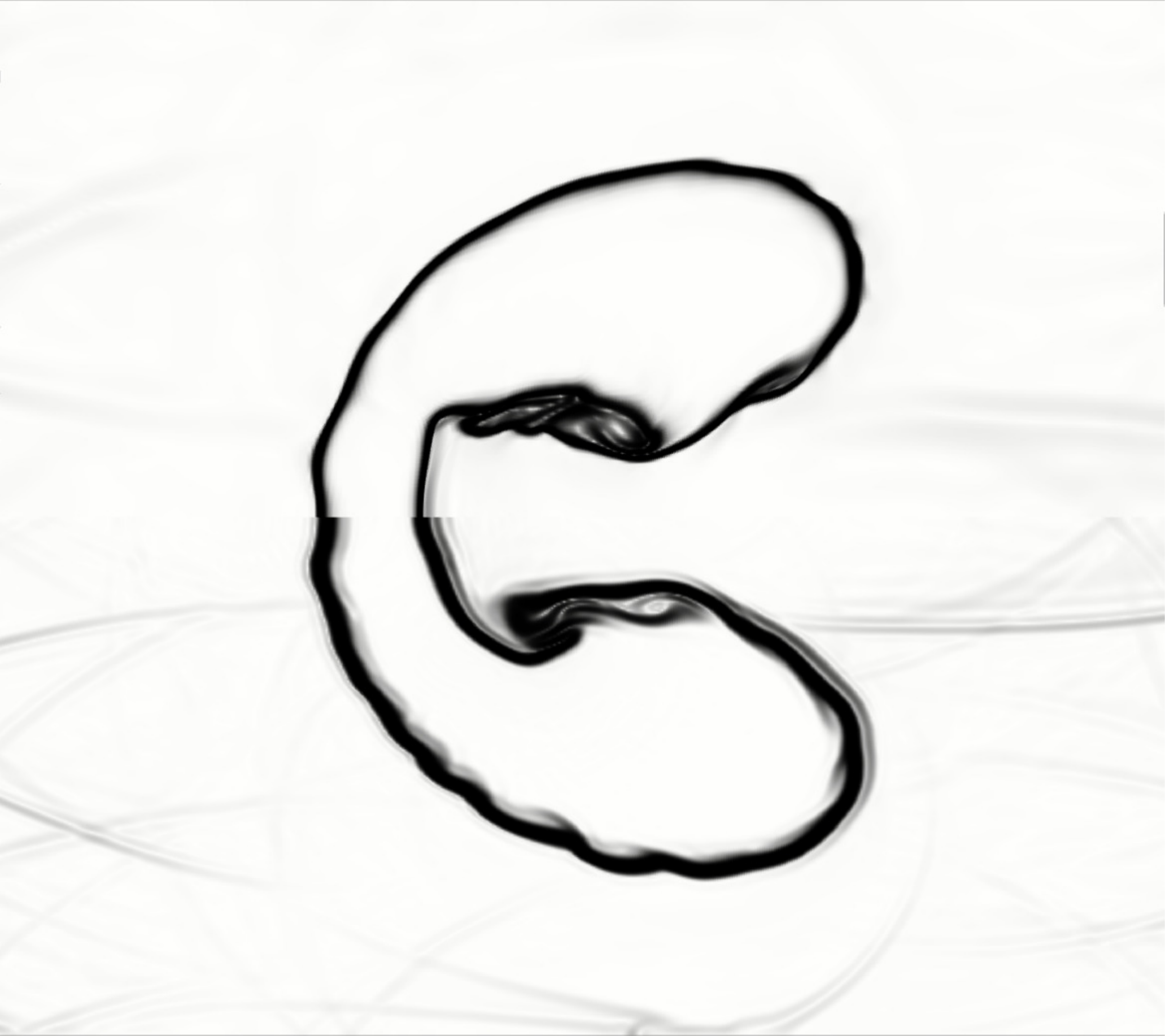}
	\end{subfigure}

\begin{subfigure}[b]{0.32\textwidth}
	\centering
	\includegraphics[width=1.8in,height=1.78in]{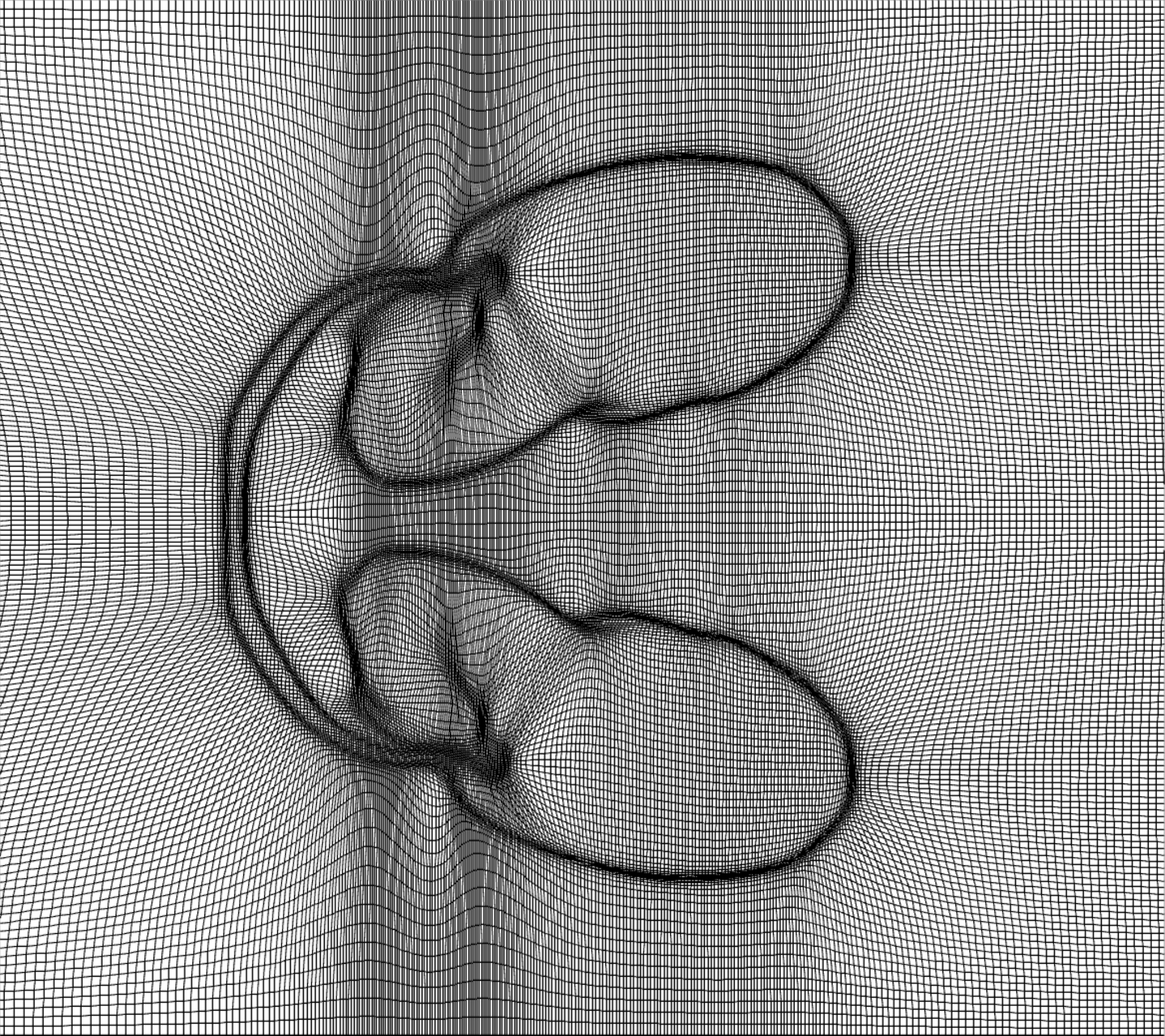}
\end{subfigure}
\begin{subfigure}[b]{0.32\textwidth}
	\centering
	\includegraphics[width=1.8in,height=1.78in]{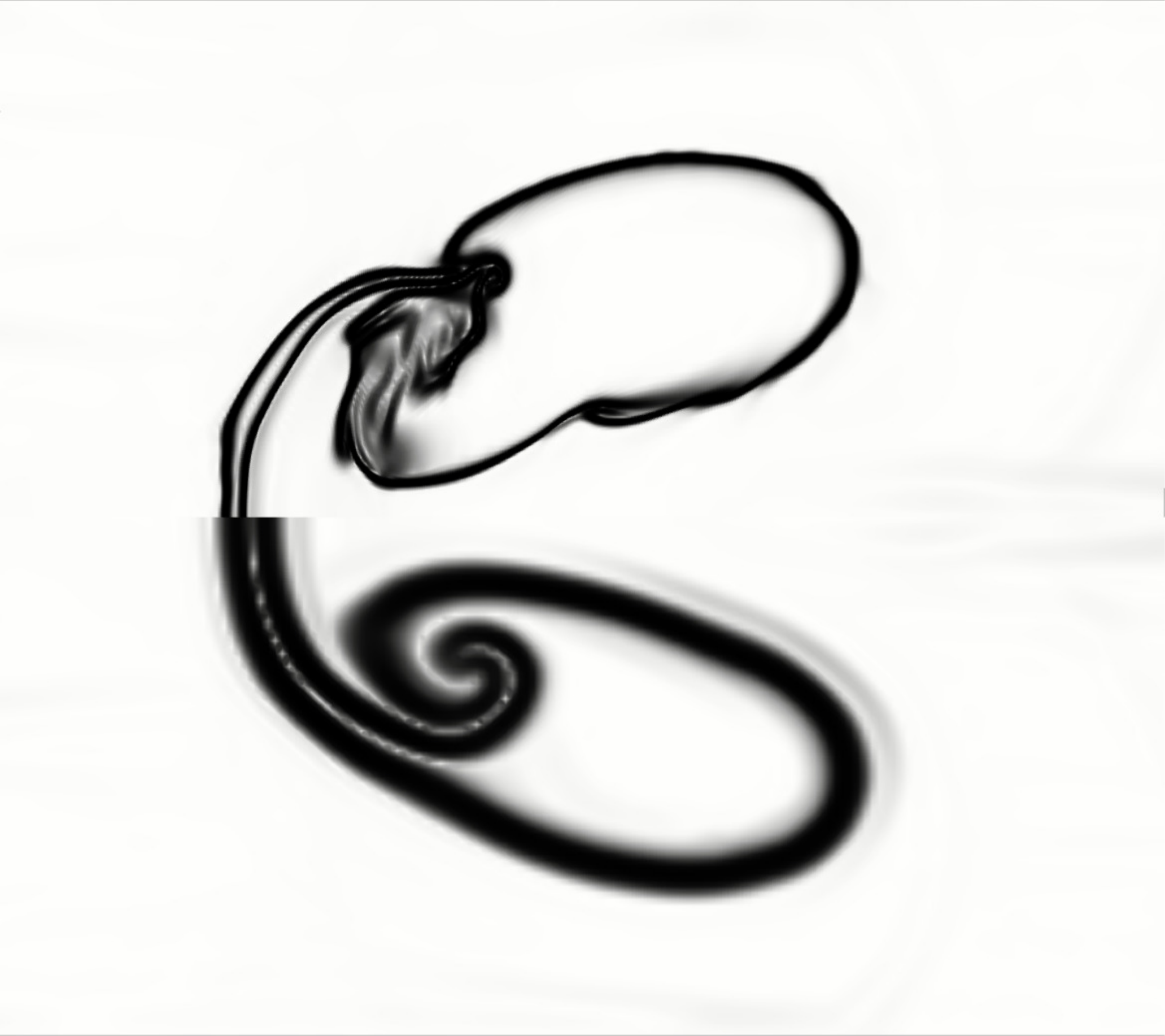}
\end{subfigure}
\begin{subfigure}[b]{0.32\textwidth}
	\centering
	\includegraphics[width=1.0\linewidth]{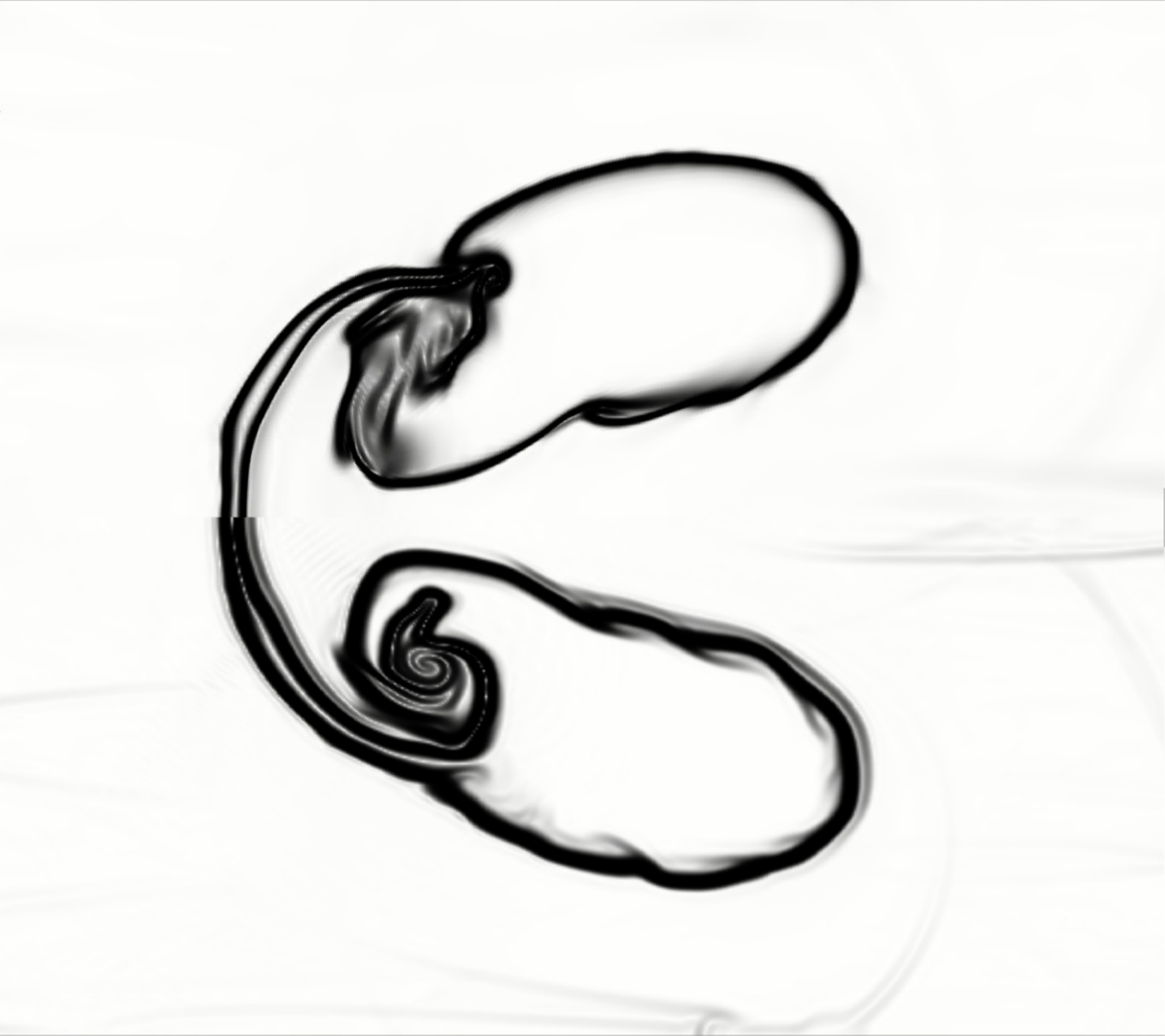}
\end{subfigure}
\caption{Same as Figure \ref{fig:0} except for $t = 0.26, 0.452, 0.676$.}
\label{fig:2}
\end{figure}
\end{example}

	\begin{example}[2D shock-bubble interaction \uppercase\expandafter{\romannumeral2}]\label{ex:ShockBubble2}\rm The shock-bubble interaction problem is extended to {the} stiffened gas.
	The  domain and the boundary conditions are the same as those in Example  \ref{ex:ShockBubble}.
Initially,
the regions $\Omega_2$ and $\Omega_3$ are filled with the stiffened gas, while  $\Omega_1$  is filled with the ideal gas,
and a Mach $M_s = 10$ shock wave at $x_1 = 275 $ moves to
	a cylindrical bubble centered at $(x_1,x_2) = [225, 0]$.
Specially,
	the initial data are
	$$
	\left(\rho_{1}, \rho_{2}, v_{1}, v_{2}, p\right)=\left\{\begin{array}{ll}
	\left(\epsilon,5 - \epsilon, 0,0,100\right), & (x_1, x_2) \in \Omega_1, \\
	(1 - \epsilon,\epsilon,0,0,100), & (x_1, x_2) \in \Omega_2, \\
	(1.980198 - \epsilon,\epsilon,-121.2497,0,29800), & (x_1, x_2) \in \Omega_3,
	\end{array}\right.
	$$
	with
	$\epsilon = 0.05,
  \Gamma_{1}=3.0, p_{\infty, 1} = 100,  \Gamma_{2}=1.4, p_{\infty,2} = 0, c_{v,1} = c_{v,2} = 1.
	$ 	

Figure \ref{fig:Shockbubble2Mesh} presents the adaptive mesh obtained by {\tt MM-WENOMR} with $800 \times 160$ cells at $t = 0.8$, and Figure \ref{fig:Shockbubble2} gives  the schlieren images  at $t = 0.8$,
where
	the monitor function is chosen as
	\eqref{eq:monitor}
	with
	$\kappa = 5$, $(\sigma_1,\sigma_2, \sigma_3,\sigma_4, \sigma_5) = (\Phi, u,\rho_1, \rho_2, p)$,  $(\alpha_1, \alpha_2, \alpha_3,\alpha_4, \alpha_5)  = (600,  500, 1200, 1200, 1200)$,
and the schlieren function is given by  \eqref{eq:schlieren} with $\Psi=\left(120 \rho_{1}+20 \rho_{2}\right) / \rho$.
We see that the mesh points adaptively concentrate near the
	shock wave and  the bubble interfaces,
  {\tt MM-WENOMR} with $800 \times 160$ cells is better  than {\tt UM-WENOMR} with the same number of cells,
 {\tt MM-WENOMR} can detect  the sharp bubble interfaces well,
 and {\tt MM-WENOMR} is efficient since
	the CPU time of {\tt MM-WENOMR} with $800 \times 160$ cells is $23.6\%$   of  {\tt UM-WENOMR} with $2400 \times 480$ cells shown in Table \ref{CPU_MC_Shock_bubble}.
	\begin{figure}[!ht]
		\centering
		\includegraphics[width=3.0in,height=1.8in]{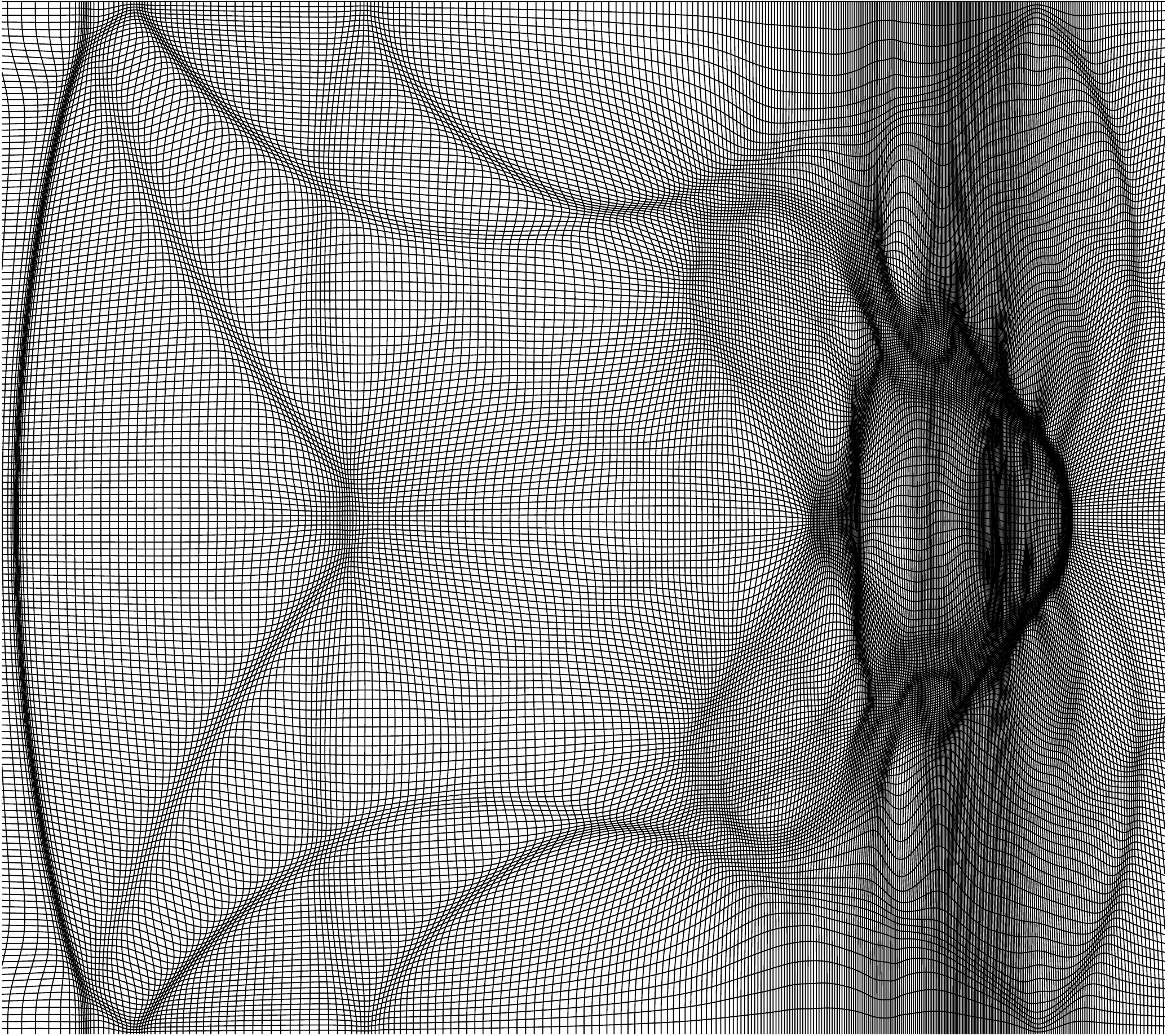}
	\caption{Example \ref{ex:ShockBubble2}. Adaptive mesh obtained by {\tt MM-WENOMR} with $800 \times 160$ cells  at $t = 0.8$.
}
\label{fig:Shockbubble2Mesh}
\end{figure}

	\begin{figure}[!ht]
		\centering
		
%
%

	\begin{subfigure}[b]{0.45\textwidth}
	\centering
	\includegraphics[width=2.5in,height=1.8in]{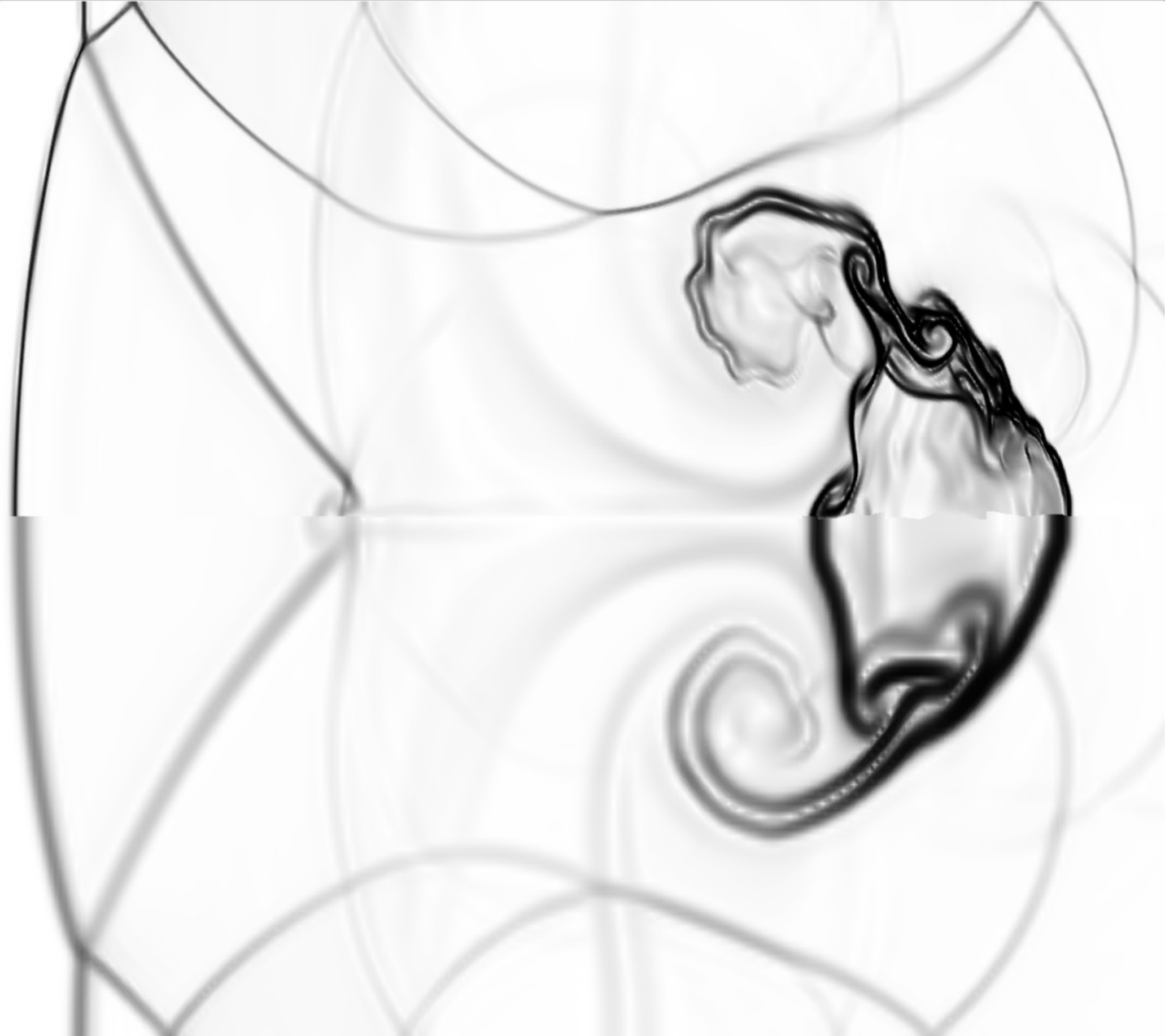}
\end{subfigure}
\begin{subfigure}[b]{0.45\textwidth}
	\centering
	\includegraphics[width=2.5in,height=1.8in]{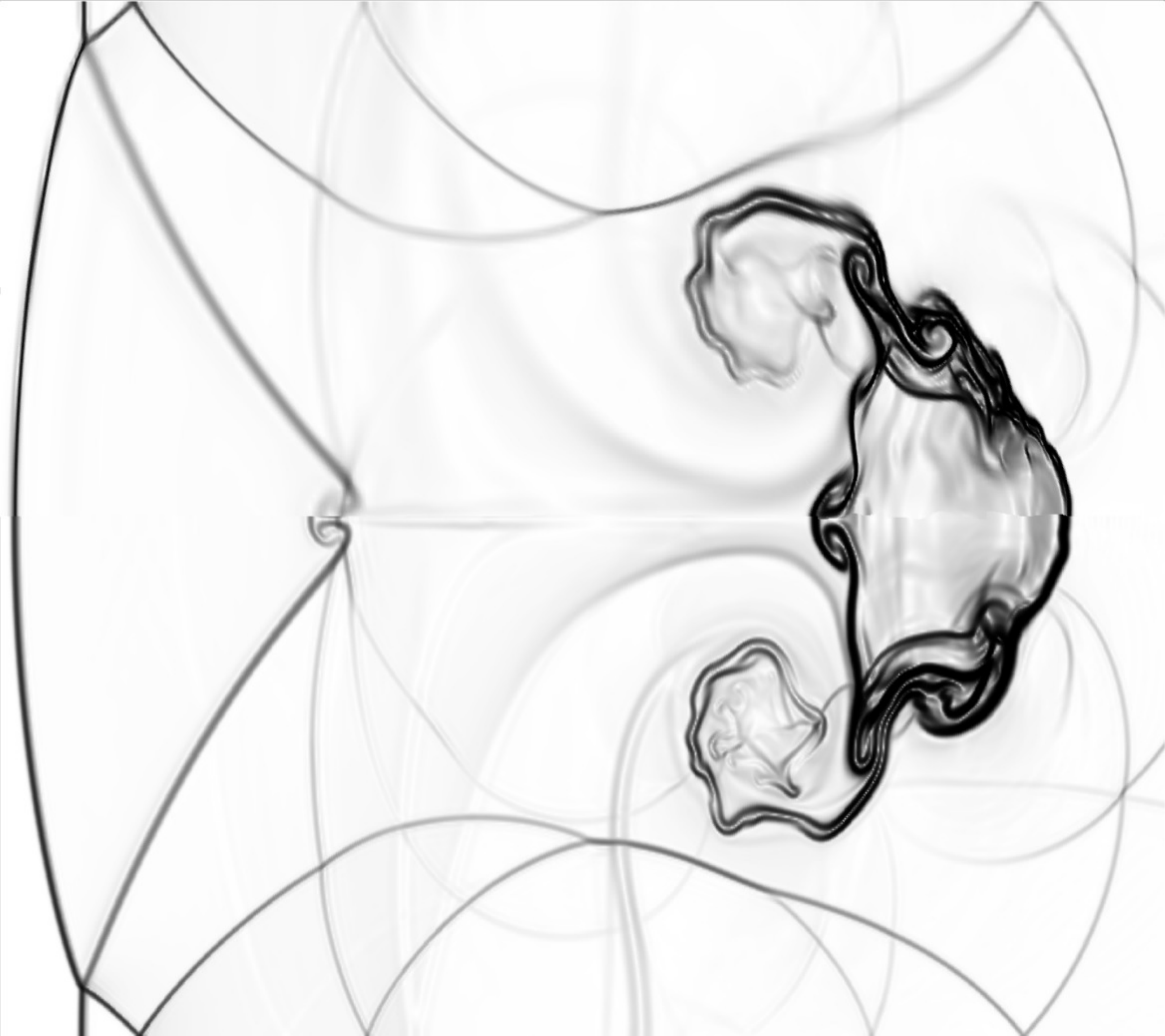}
\end{subfigure}
	
		\caption{Example \ref{ex:ShockBubble2}. Schlieren images of $\Phi$ at $t = 0.8$.  Left: {\tt MM-WENOMR} with $800 \times 160$ cells (top half)) and {\tt UM-WENOMR} with $800 \times 160$ cells (bottom half); right:  {\tt MM-WENOMR} with $800 \times 160$ cells (top half) and  {\tt UM-WENOMR} with $2400 \times 480$ cells (bottom half). 
		}
		\label{fig:Shockbubble2}
	\end{figure}

\end{example}

\begin{example}[3D shock-bubble interaction
	\uppercase\expandafter{\romannumeral1}]\label{ex:3DShockBubble}\rm
	This is an extension of  Example \ref{ex:ShockBubble}, and considers  a planer shock wave interacting  with a helium  bubble in the   domain  $[0, 445]\times[-44.5, 44.5] \times[-44.5, 44.5]$.
	The initial pre- and post-shock states are
	$$
	\left(\rho_{1}, \rho_{2}, v_{1}, v_{2}, v_{3}, p\right)=\left\{\begin{array}{ll}
	(1.225 - \epsilon, \epsilon,0,0,0,101325), & x_1 < 275,\\
	(1.6861 - \epsilon, \epsilon, -113.5243,0,0,159060), & x_1 > 275,
	\end{array}\right.
	$$
	and the bubble state is
	$$
		\left(\rho_{1}, \rho_{2}, v_{1}, v_{2}, v_{3}, p\right)=
		\left(\epsilon,1.225\left(R_{1} / R_{2}\right) - \epsilon, 0,0,0,101325\right),  ~\sqrt{(x_1-225)^2 + x_2^2+ x_3^2} < 25,
		$$
	with
		$ \epsilon = 0.03, p_{\infty, 1} = p_{\infty,2} = 0,
	 \Gamma_{1}=1.4, \Gamma_{2}=1.647, R_{1}=0.287, R_{2}=1.578$.
	
	 Figure \ref{fig:3DMesh} gives the close-up of the adaptive mesh, the iso-surface of $\rho = 0.66$, three offset 2D slices and two surface meshes near the bubble at $t = 0.72$, where
 the monitor function is chosen as
	 \eqref{eq:monitor}
	with $\kappa=1, \sigma_1=\rho$ and $\alpha_1=1200$,
and the linear weights of the multi-resolution WENO reconstruction are taken as $\chi_1 = 0.95, \chi_2 = 0.045$ and $\chi_3 = 0.005$.  Figure \ref{fig:3DRHO} shows    the schlieren images on the slice $x_3 = 0$ given by \eqref{eq:schlieren}  with $\Psi=\left(10 \rho_{1}+30 \rho_{2}\right) / \rho $, where the top half parts   are the results obtained by {\tt MM-WENOMR} with $400 \times 80 \times 80$ cells, while the left and right bottom half parts are those obtained by {\tt UM-WENOMR} with $400 \times 80 \times 80$ cells and  $800 \times 160 \times 160$ cells, respectively.
	One can  see that the solution of {\tt MM-WENOMR} with $400 \times 80 \times 80$ cells is comparable to that of {\tt UM-WENOMR} with $800 \times 160 \times 160$ cells,
but the CPU time of {\tt MM-WENOMR} with $400 \times 80 \times 80$ cells is $29.2\%$   of {\tt UM-WENOMR} with $800 \times 160 \times 160$ cells, see Table \ref{CPU_MC_Shock_bubble}.

\begin{figure}[!ht]
	\centering
	
	\begin{subfigure}[b]{0.45\textwidth}
		\centering
		\includegraphics[width=1.0\linewidth,trim=40 10 25 90, clip ]{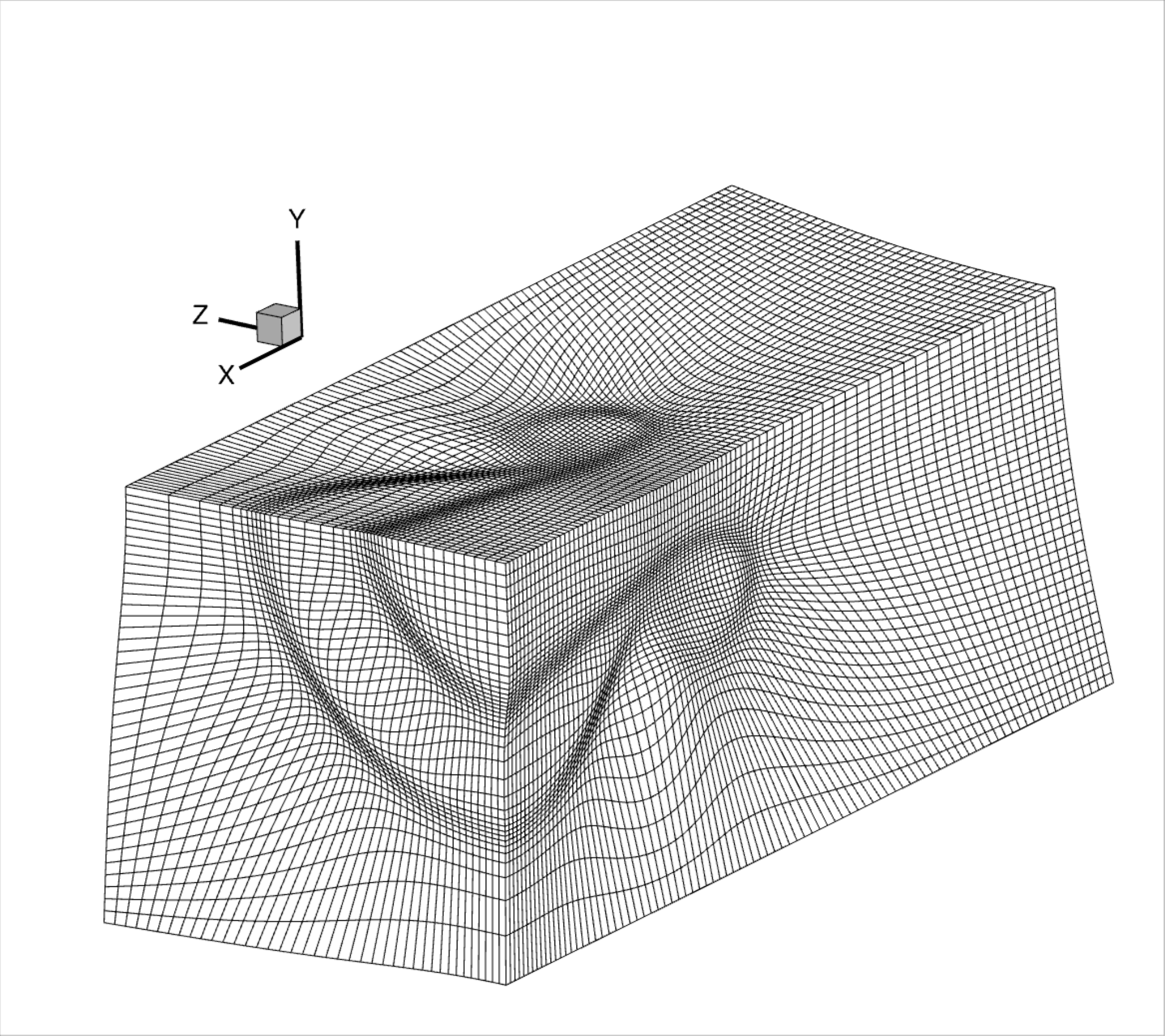}
	\end{subfigure}
\begin{subfigure}[b]{0.48\textwidth}
	\centering
	\includegraphics[width=1.0\linewidth,trim=45 30 20 40, clip ]{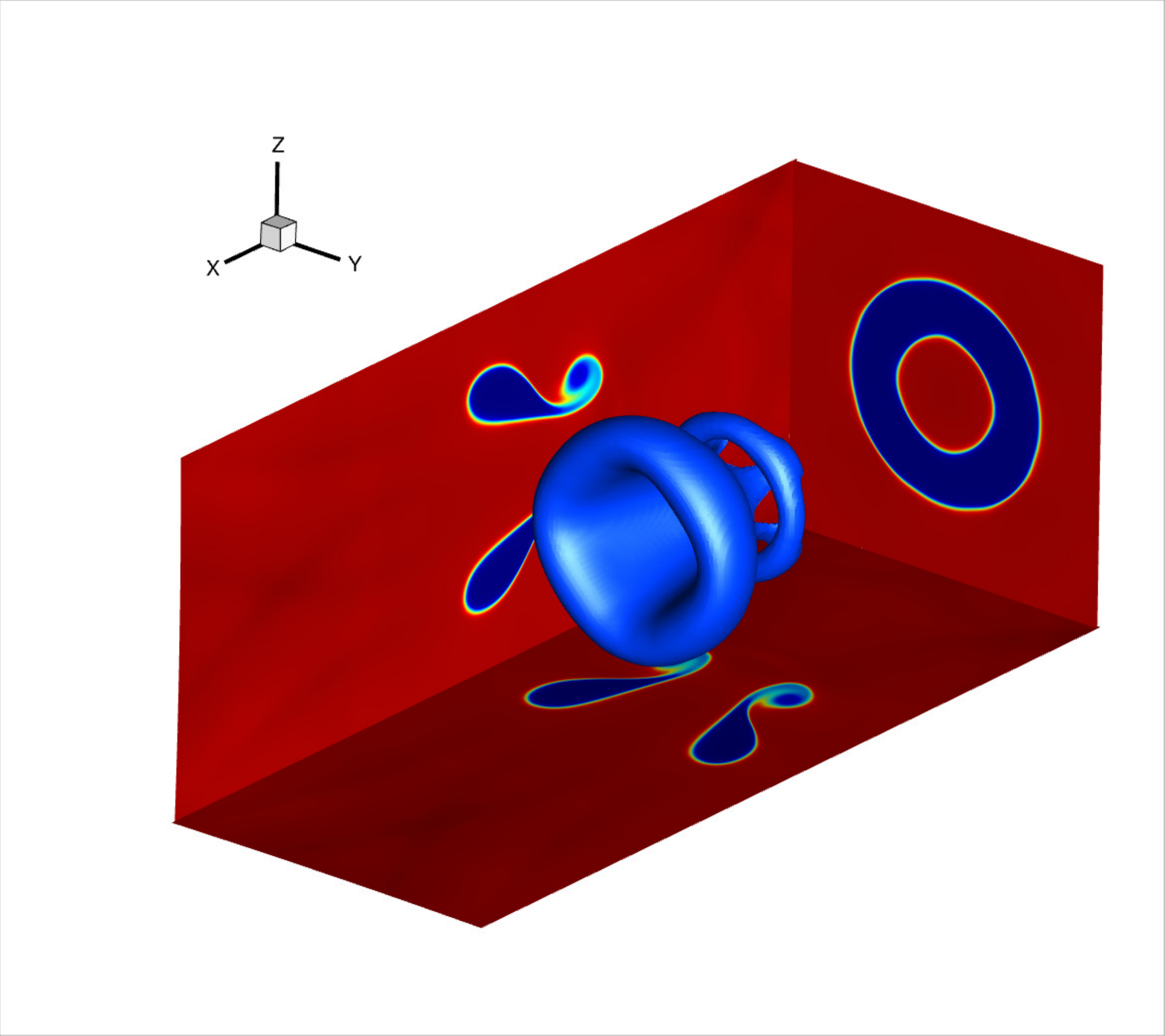}
\end{subfigure}

\begin{subfigure}[b]{0.44\textwidth}
	\centering
	\includegraphics[width=1.0\linewidth,trim=40 0 35 0, clip]{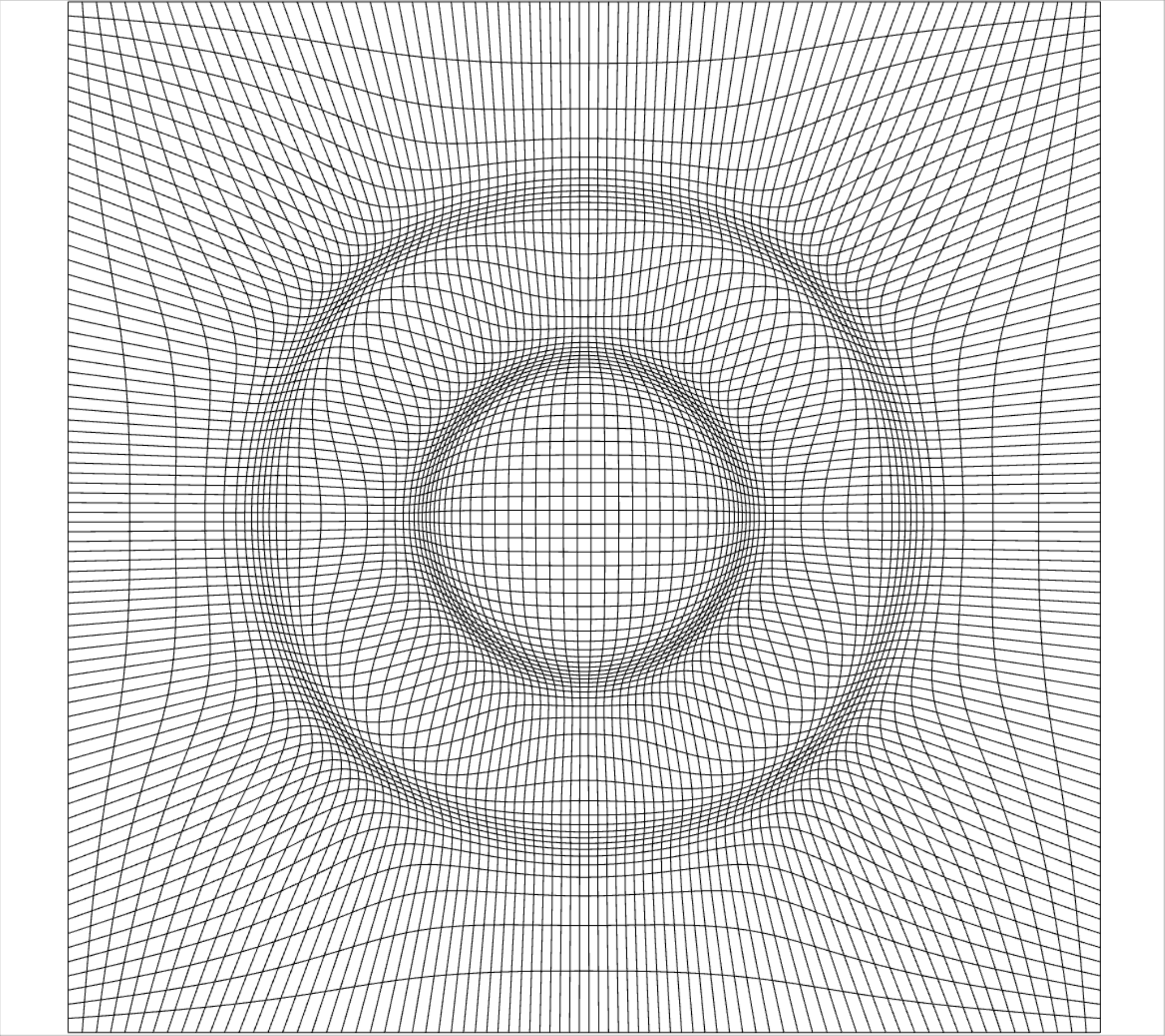}
\end{subfigure}
	\begin{subfigure}[b]{0.4435\textwidth}
		\centering
		\includegraphics[width=1.0\linewidth,trim=40 0 35 1, clip]{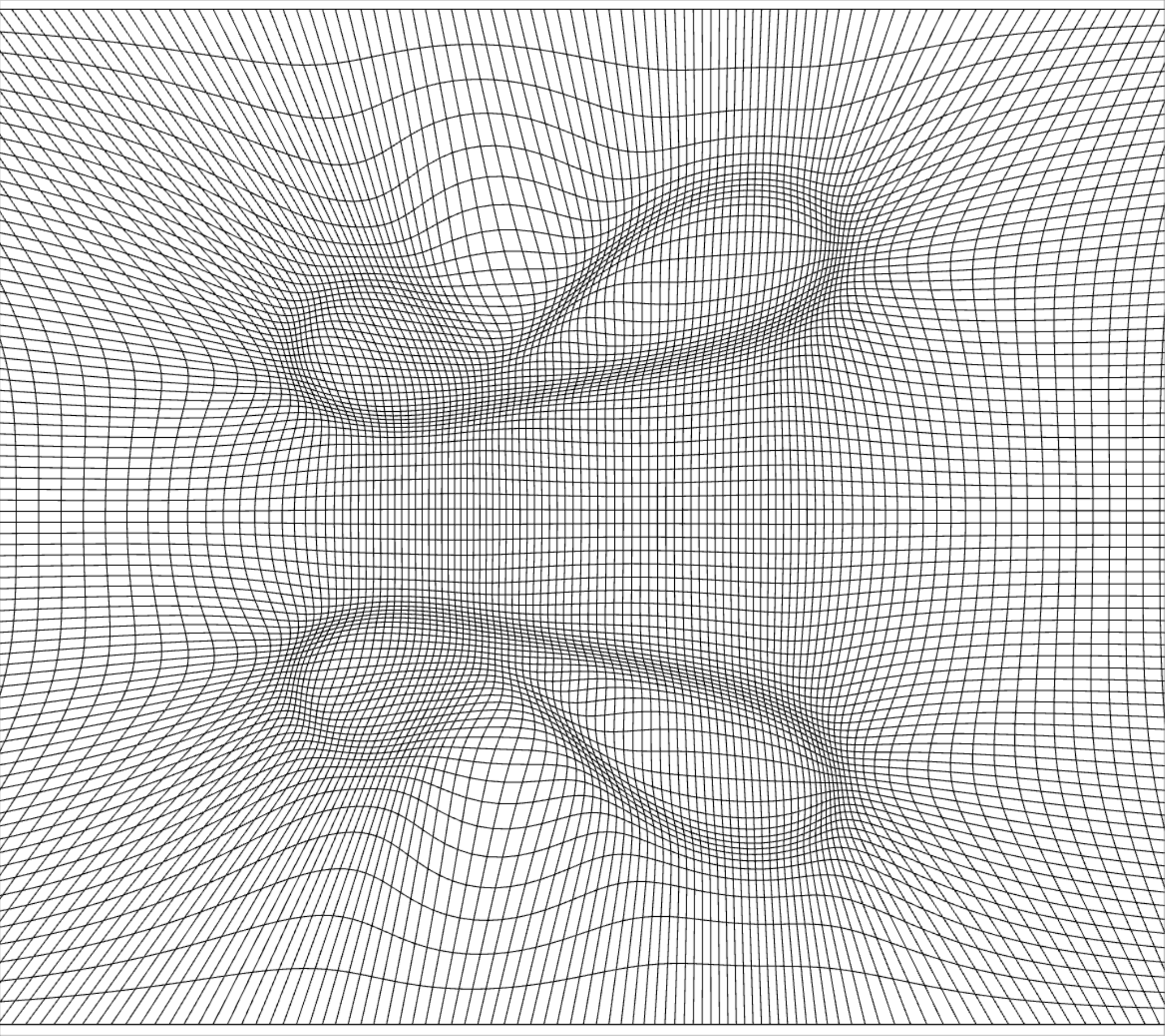}
	\end{subfigure}
\caption{Example \ref{ex:3DShockBubble}. Adaptive meshes  and $\rho$ at $t = 0.72$. Top left: close-up of the adaptive mesh, $i_1 \in [40, 130], i_2 \in [1, 40], i_3 \in [40, 80]$;  top right: the iso-surface of  $\rho = 0.66$ and  three offset 2D slices taken at $x_1 =
	136$, $x_2 = 0$, $x_3 = 0$;
	bottom left: the surface mesh with $i _1= 120$; bottom right:  the surface mesh with $i_3 = 40$.
}
\label{fig:3DMesh}
\end{figure}

\begin{figure}[!ht]
	\centering
	\vspace{-1.5cm}
	\begin{subfigure}[b]{0.45\textwidth}
		\centering
		\includegraphics[width=1.8in,height=1.78in,  trim=1 2 1 1, clip ]{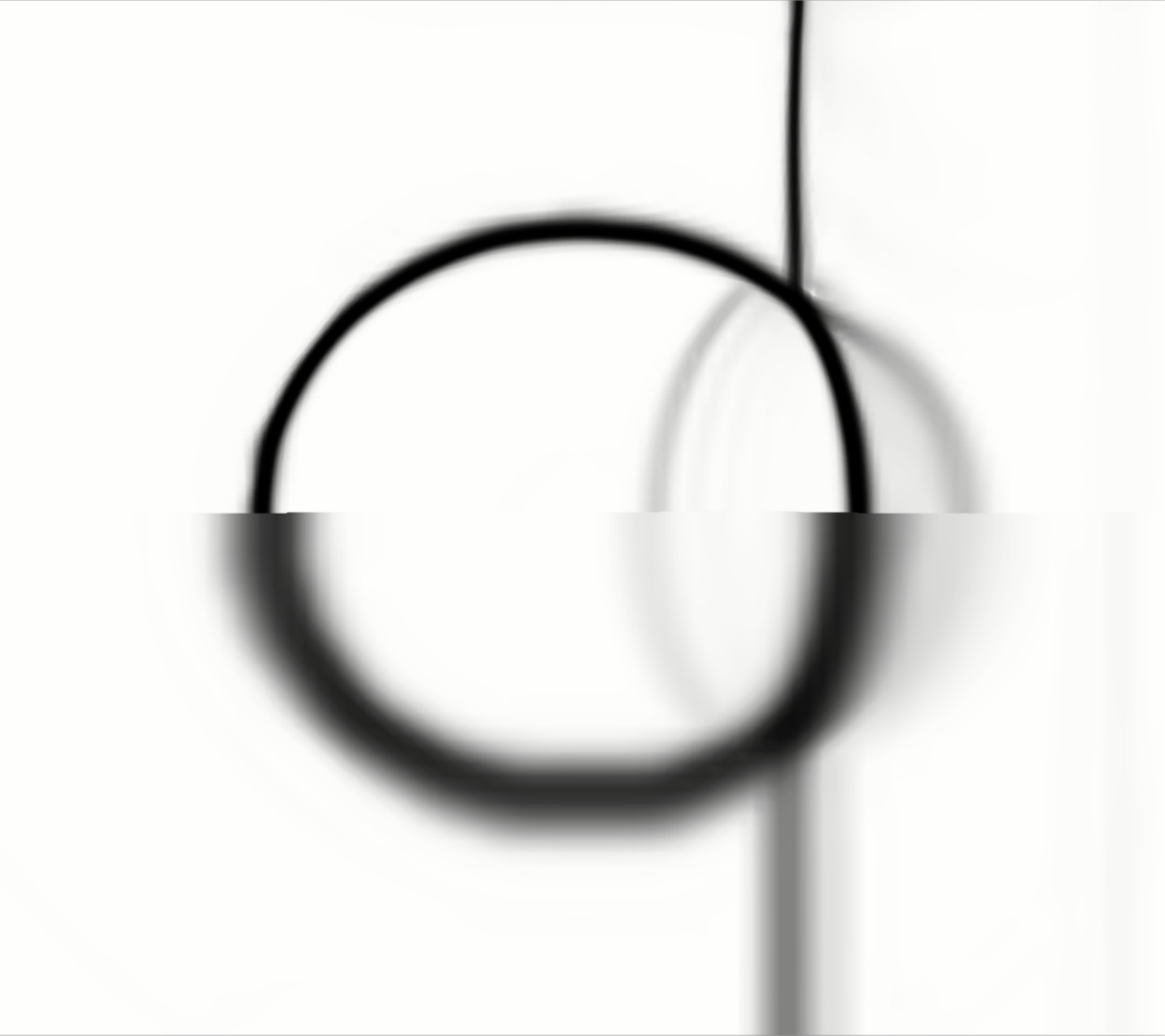}
	\end{subfigure}
	\begin{subfigure}[b]{0.45\textwidth}
		\centering
		\includegraphics[width=1.8in,height=1.78in, trim=1 2 1 1, clip]{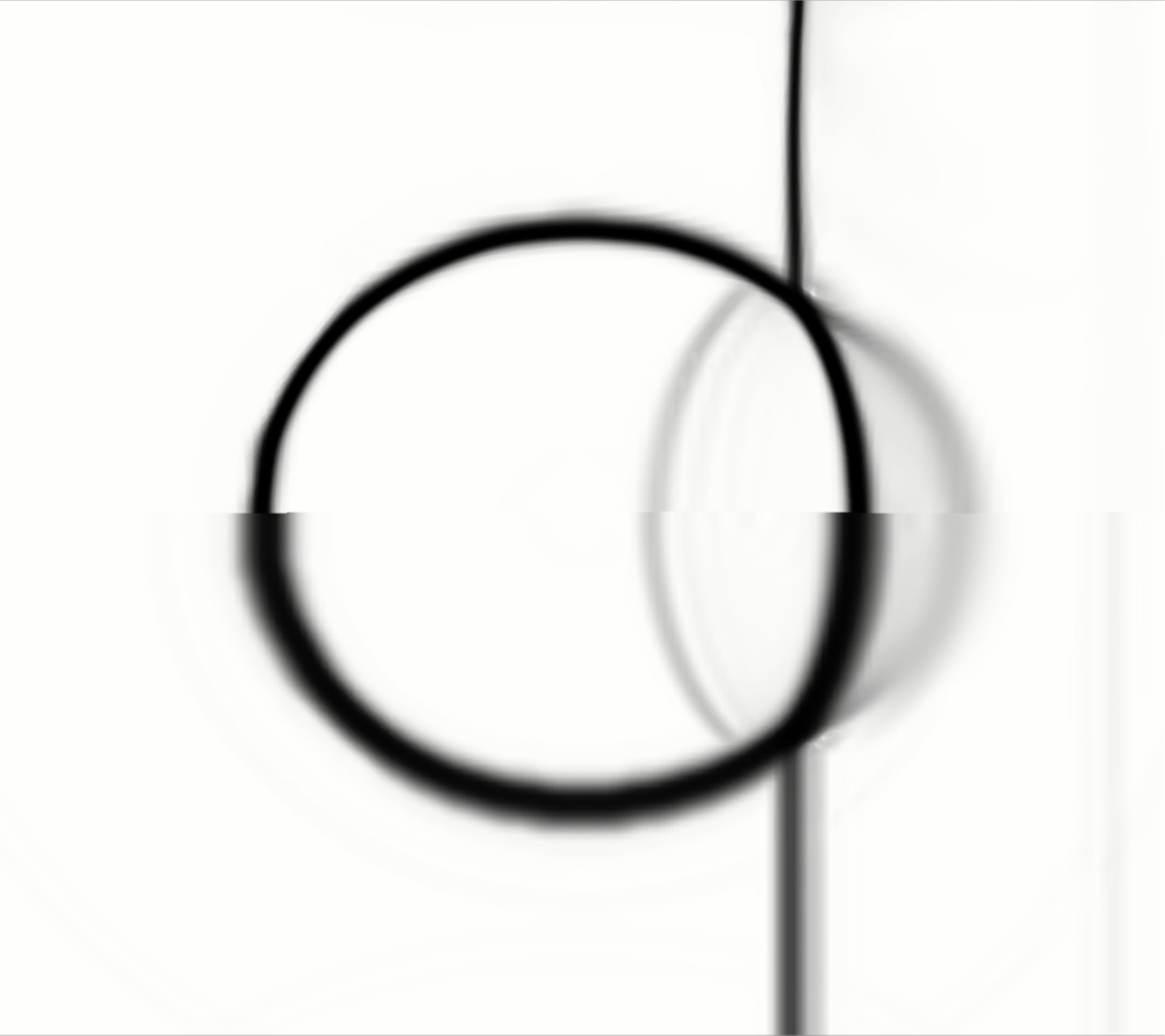}
	\end{subfigure}
	
		\begin{subfigure}[b]{0.45\textwidth}
		\centering
		\includegraphics[width=1.8in,height=1.78in,  trim=1 2 1 1, clip ]{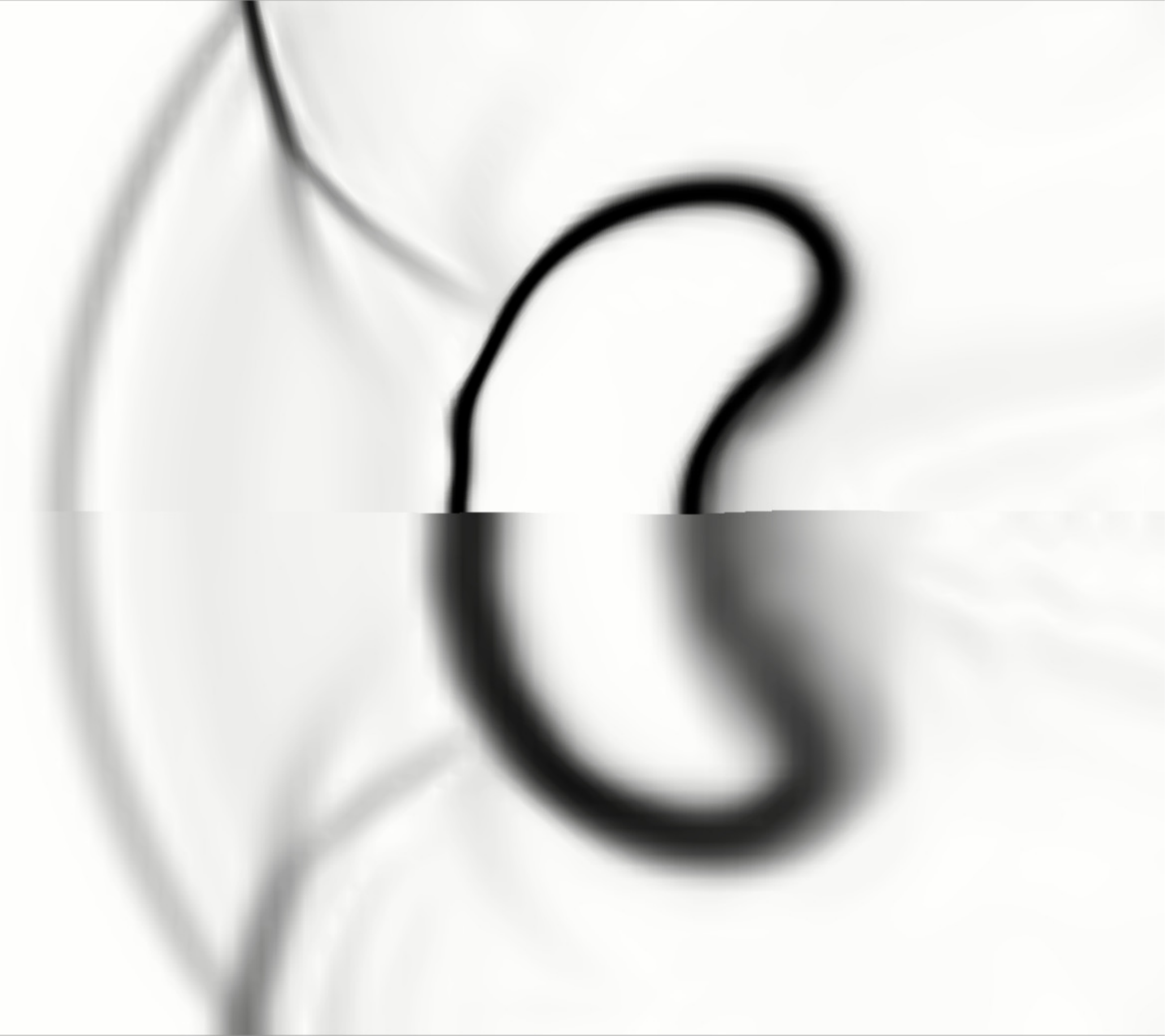}
	\end{subfigure}
	\begin{subfigure}[b]{0.45\textwidth}
		\centering
		\includegraphics[width=1.8in,height=1.78in, trim=1 2 1 1, clip]{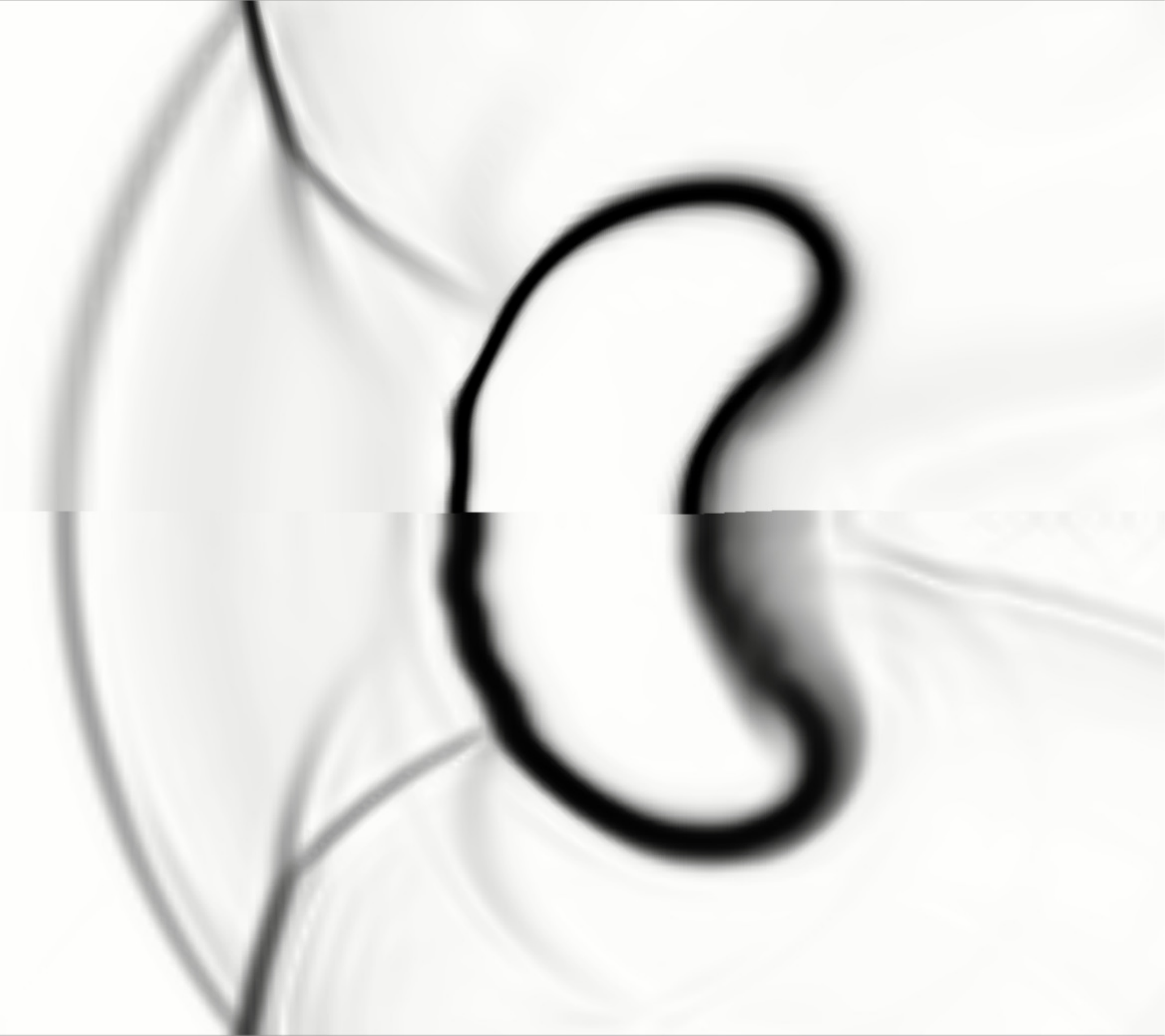}
	\end{subfigure}

	\begin{subfigure}[b]{0.45\textwidth}
	\centering
	\includegraphics[width=1.8in,height=1.78in,  trim=1 2 1 1, clip ]{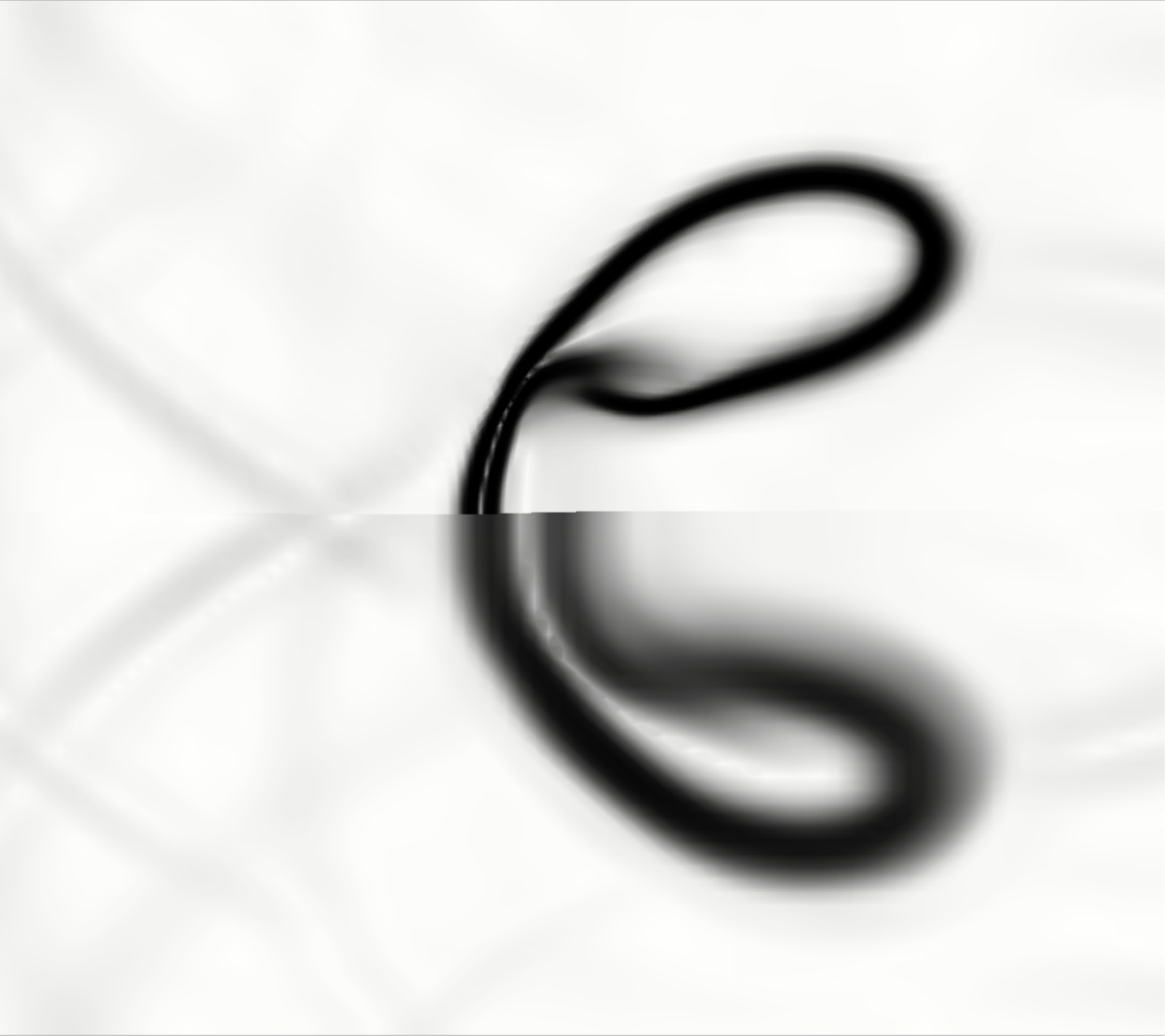}
\end{subfigure}
\begin{subfigure}[b]{0.45\textwidth}
	\centering
	\includegraphics[width=1.8in,height=1.78in, trim=1 2 1 1, clip]{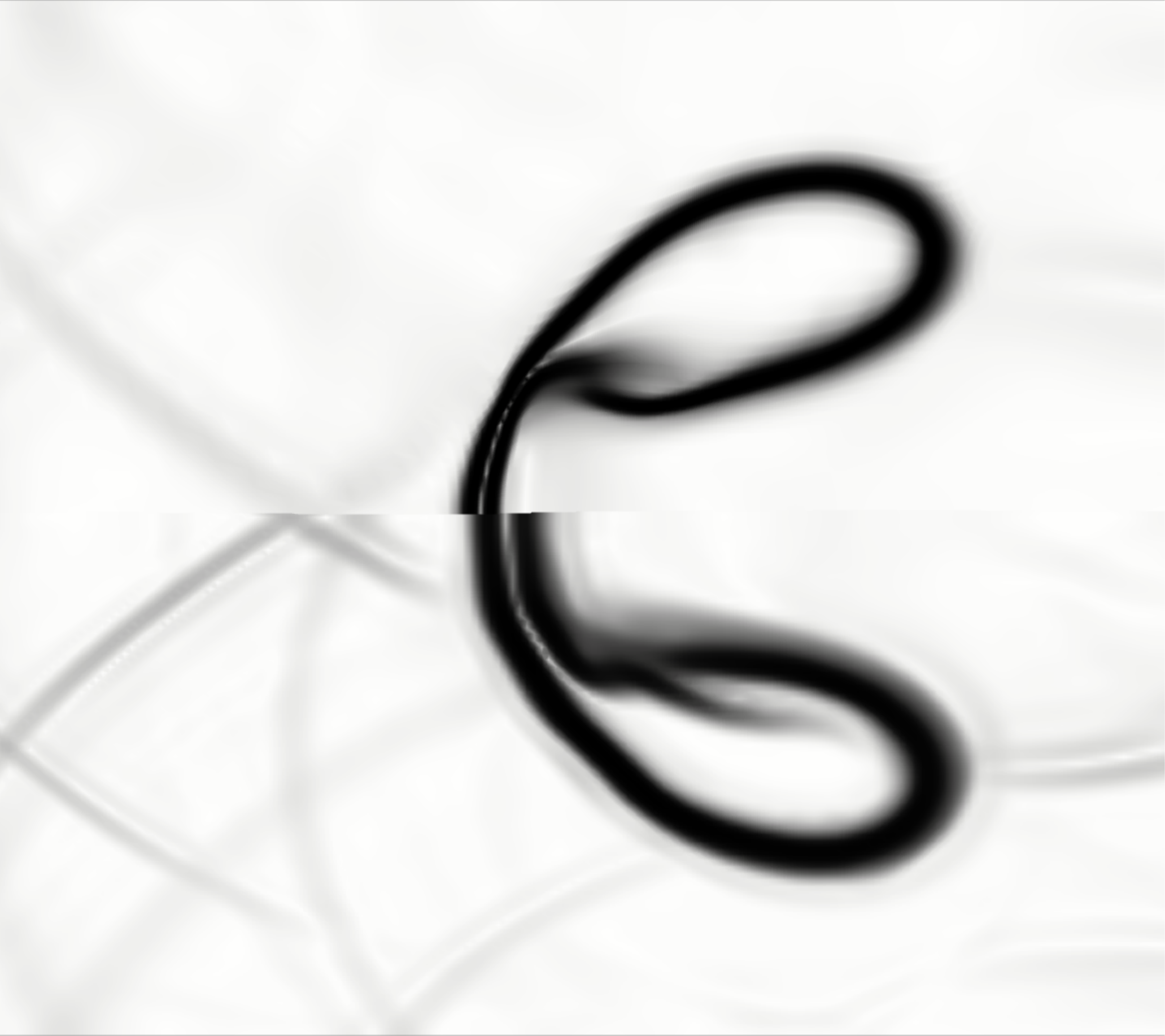}
\end{subfigure}

	\begin{subfigure}[b]{0.45\textwidth}
	\centering
	\includegraphics[width=1.8in,height=1.78in,  trim=1 2 1 1, clip ]{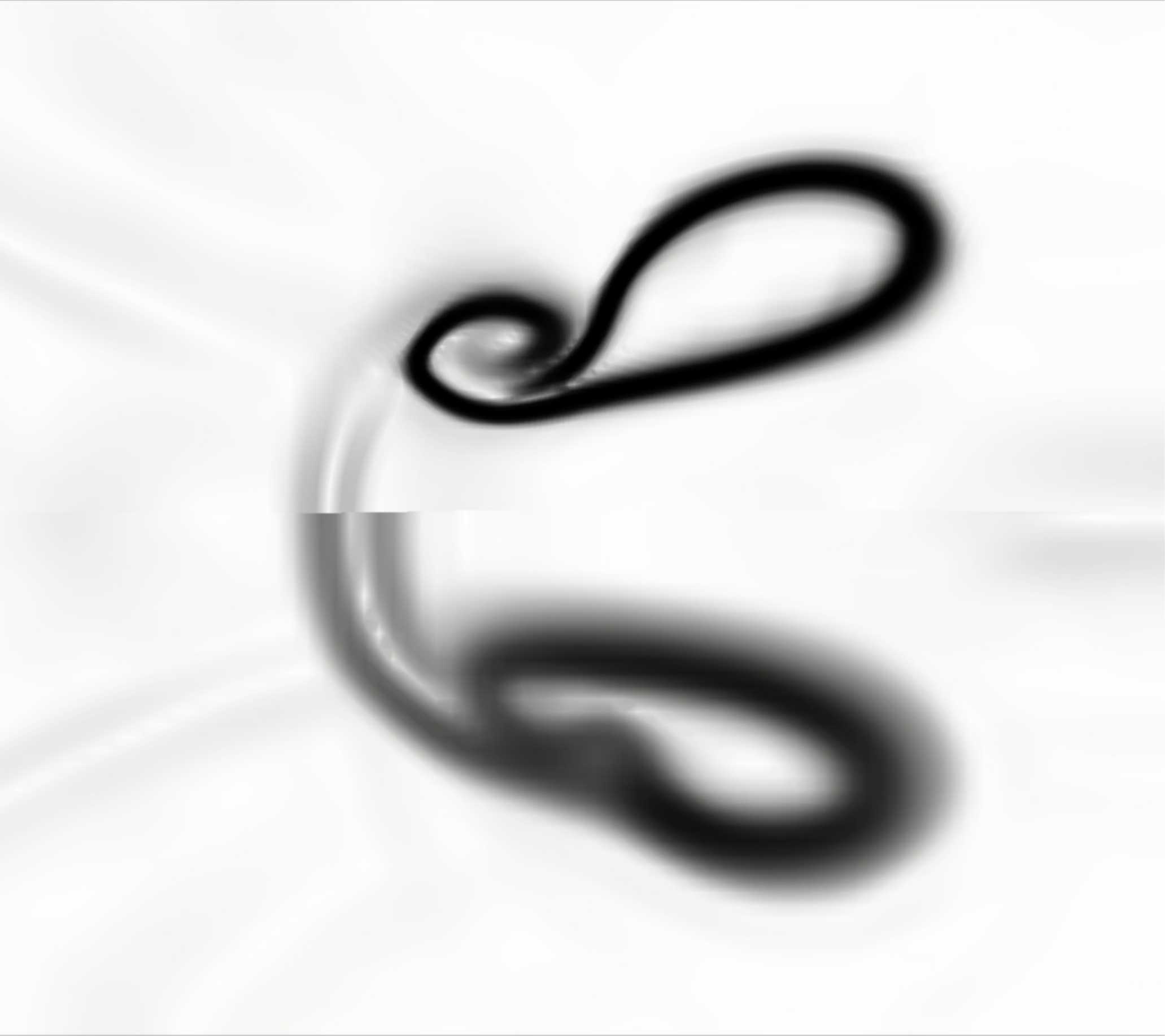}
\end{subfigure}
\begin{subfigure}[b]{0.45\textwidth}
	\centering
	\includegraphics[width=1.8in,height=1.78in, trim=1 2 1 1, clip]{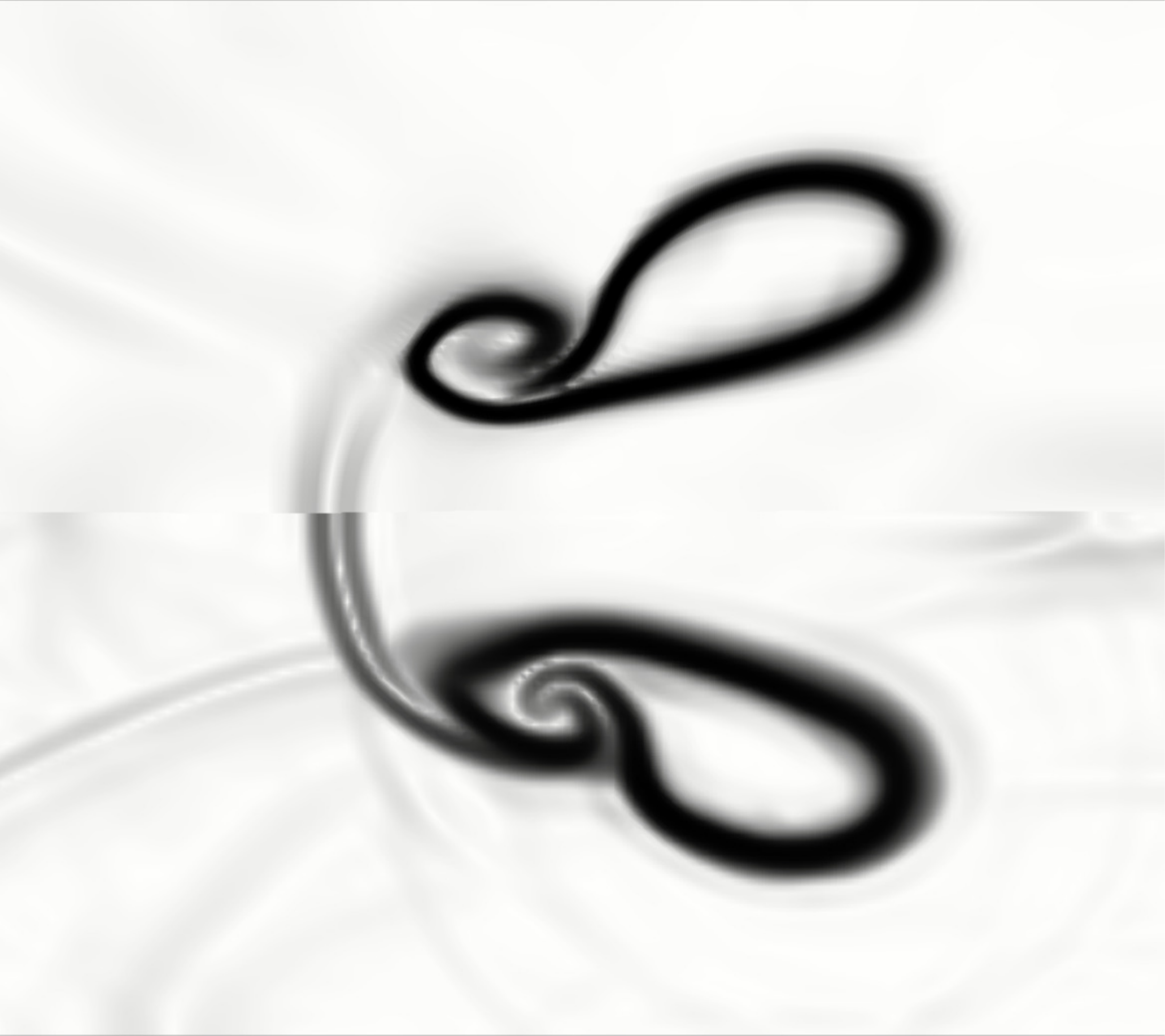}
\end{subfigure}

	\begin{subfigure}[b]{0.43\textwidth}
	\centering
	\includegraphics[width=1.8in,height=1.78in,  trim=1 2 1 1, clip ]{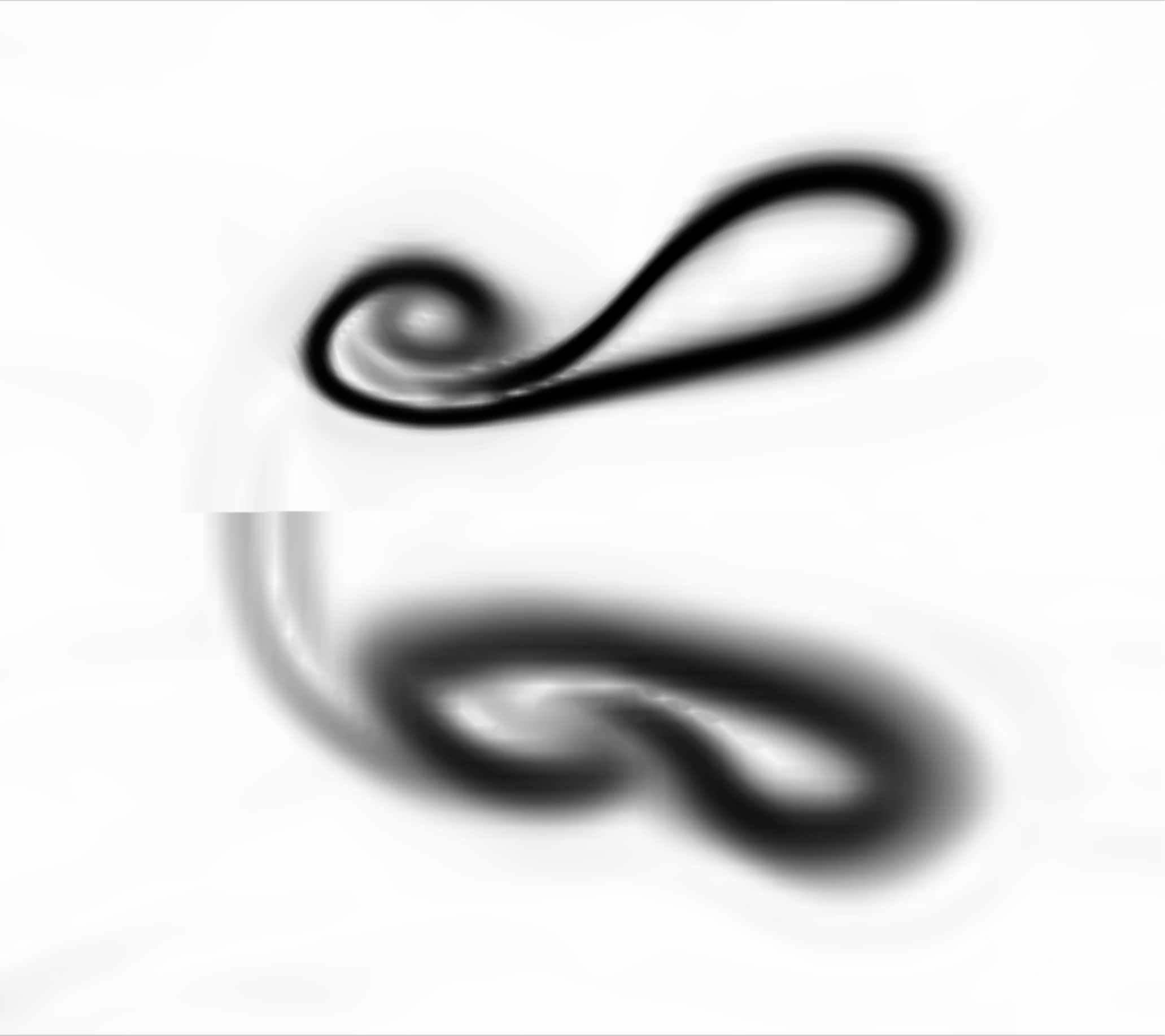}
\end{subfigure}
\begin{subfigure}[b]{0.43\textwidth}
	\centering
	\includegraphics[width=1.8in,height=1.78in, trim=1 2 1 1, clip]{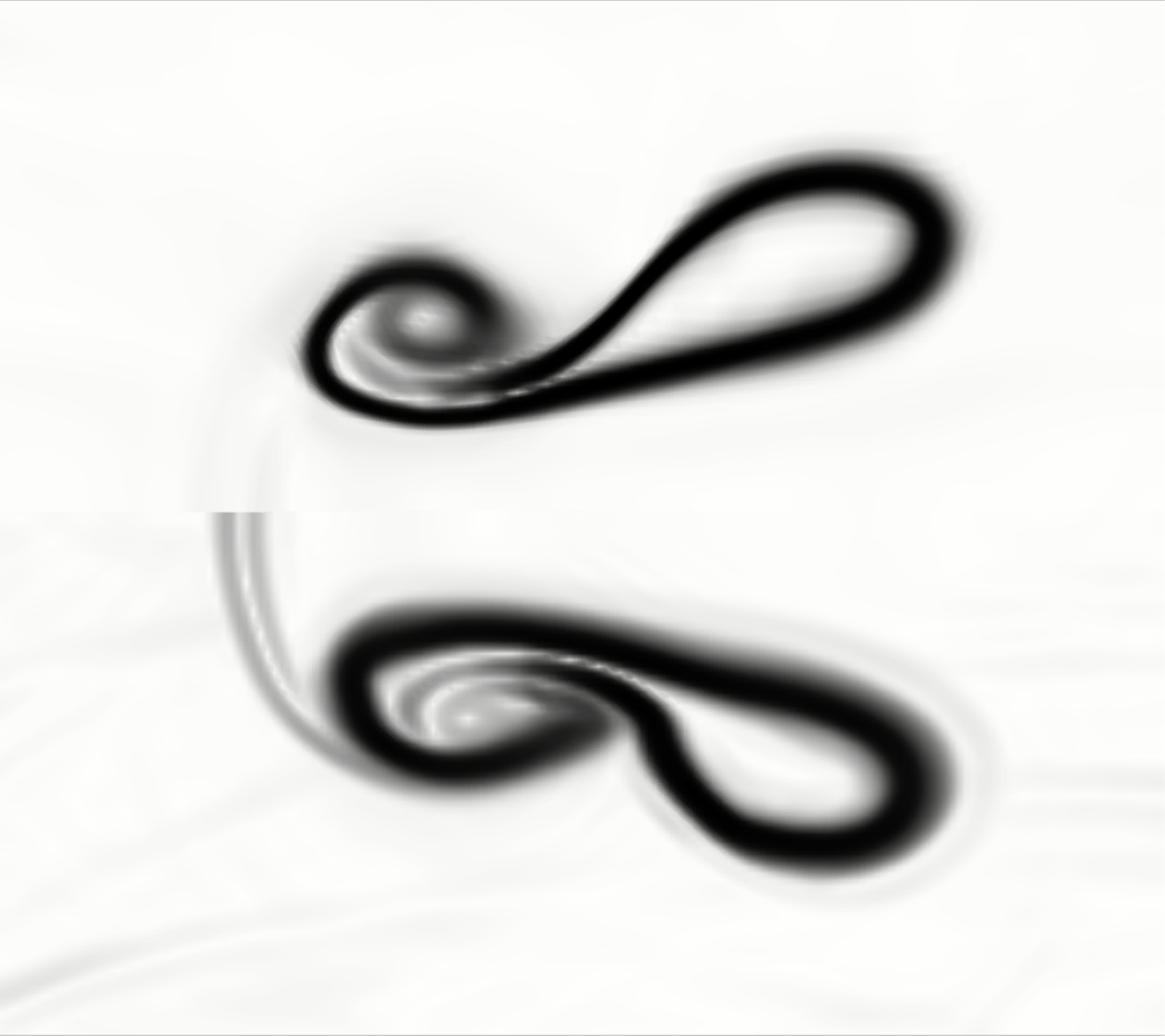}
\end{subfigure}
	\caption{Example \ref{ex:3DShockBubble}. Schlieren images $\Phi$ on the slice $x_3 = 0$ at $t = 0.08, 0.24, 0.40, 0.56, 0.72$ (from top to bottom). Left:  {\tt MM-WENOMR} (top half) and  {\tt UM-WENOMR} (bottom half) with $400 \times 80 \times 80$ cells;
right:
 {\tt MM-WENOMR} with $400 \times 80 \times 80$ cells  (top half)
 and {\tt UM-WENOMR} with $800 \times 160 \times 160$ cells (bottom half). }
	\label{fig:3DRHO}
\end{figure}

\end{example}

\begin{example}[3D shock-bubble interaction
	\uppercase\expandafter{\romannumeral2}]\label{ex:3DShockBubble2}\rm
	This test considers  an extended 3D problem of the interaction of the shock wave and  the bubble in the   domain  $[0, 445]\times[-44.5, 44.5] \times[-44.5, 44.5]$. A bubble of ideal gas initially located at $(x_1,x_2,x_3) = [225,0,0]$ with radius $r = 25$ is surrounded by a stiffened gas.
	The initial pre- and post-shock states are
	$$
	\left(\rho_{1}, \rho_{2}, v_{1}, v_{2}, v_{3}, p\right)=\left\{\begin{array}{ll}
	(1- \epsilon, \epsilon,0,0,0,100), & x_1 < 275,\\
	(1.980198 - \epsilon, \epsilon, -121.2497,0,0,29800), & x_1 > 275,
	\end{array}\right.
	$$
	and the bubble state is
	$$
	\left(\rho_{1}, \rho_{2}, v_{1}, v_{2}, v_{3}, p\right)=
	\left(\epsilon,5 - \epsilon, 0,0,0,100\right),  ~\sqrt{(x_1-225)^2 + x_2^2+ x_3^2} < 25,
	$$
	where
	$\epsilon = 0.05, c_{v,1} = c_{v,2}=1,
	\Gamma_{1}=3.0, p_{\infty, 1} = 100,  \Gamma_{2}=1.4, p_{\infty,2} = 0.
	$ 	
	
 Figure \ref{fig:3DMesh2} shows the close-up of the adaptive mesh, the iso-surfaces of $\rho = 11.11$, three offset 2D slices and two surface meshes near the bubble at $t = 0.64$, and  Figure \ref{fig:3DShockbubble2} gives   the schlieren images given by the function in \eqref{eq:schlieren} with $\Psi=\left(20 \rho_{1}+5 \rho_{2}\right) / \rho $ at different {times}, where
 	the monitor function is chosen as
	\eqref{eq:monitor}
	with
	$(\sigma_1,\sigma_2) = (\rho_1, p)$ and $ (\alpha_1, \alpha_2)  = (500,  500)$, and the linear weights of {the} multi-resolution WENO reconstruction are the same as those in Example \ref{ex:3DShockBubble}.
 It is shown that
 the solution given by  {\tt MM-WENOMR}  is superior to  the solution given by {\tt UM-WENOMR} with finer mesh, while {\tt MM-WENOMR} is efficient since the CPU time of {\tt MM-WENOMR} with $400 \times 80 \times 80$ cells is $18.0\%$ of {\tt UM-WENOMR} with $800 \times 160 \times 160$ cells, see Table \ref{CPU_MC_Shock_bubble}.

	\begin{figure}[!ht]
		\centering
		
		\begin{subfigure}[b]{0.45\textwidth}
			\centering
			\includegraphics[width=1.0\linewidth,trim=30 60 35 50, clip ]{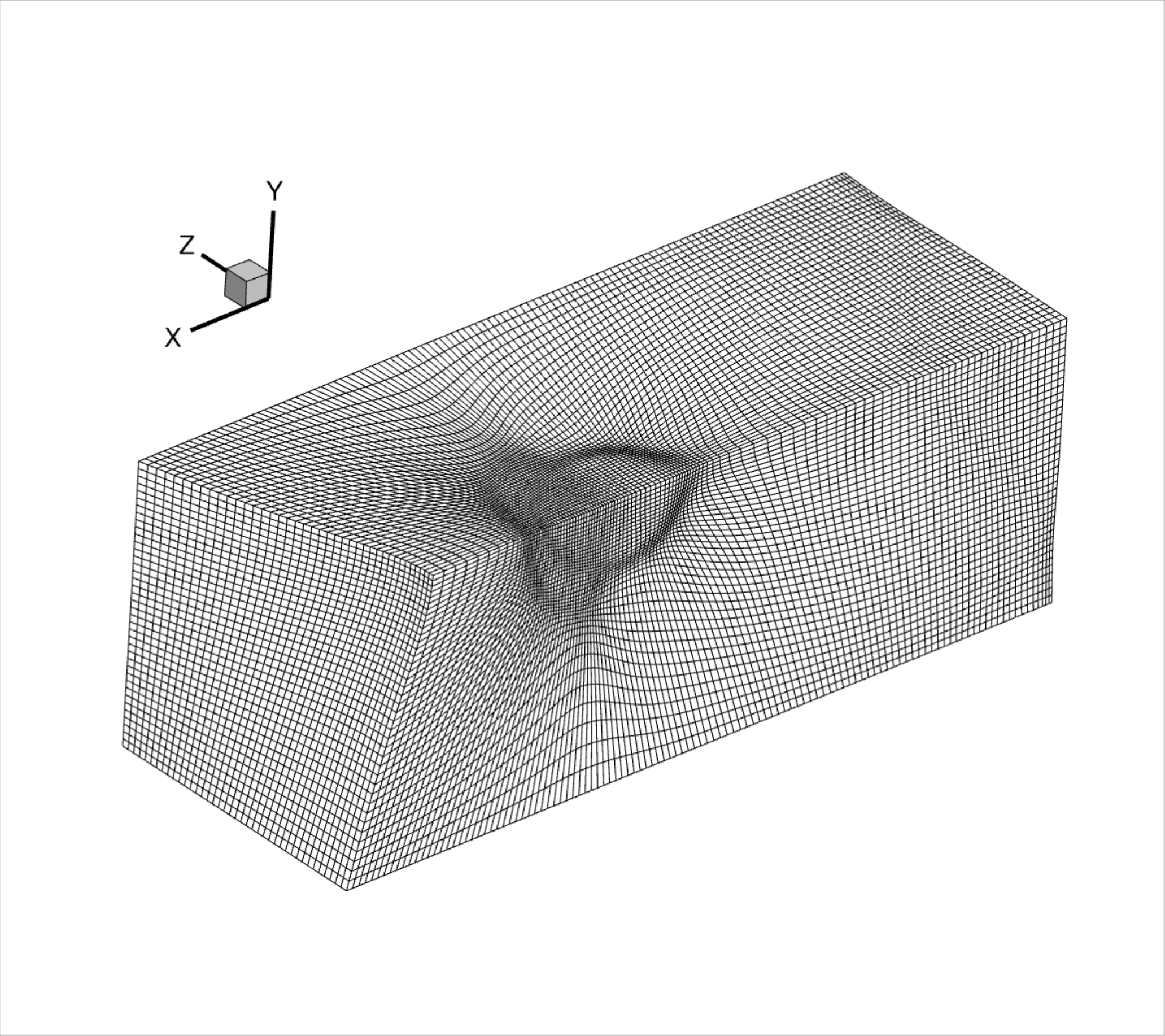}
		\end{subfigure}
		\begin{subfigure}[b]{0.45\textwidth}
			\centering
			\includegraphics[width=1.0\linewidth,trim=40 0 20 20, clip ]{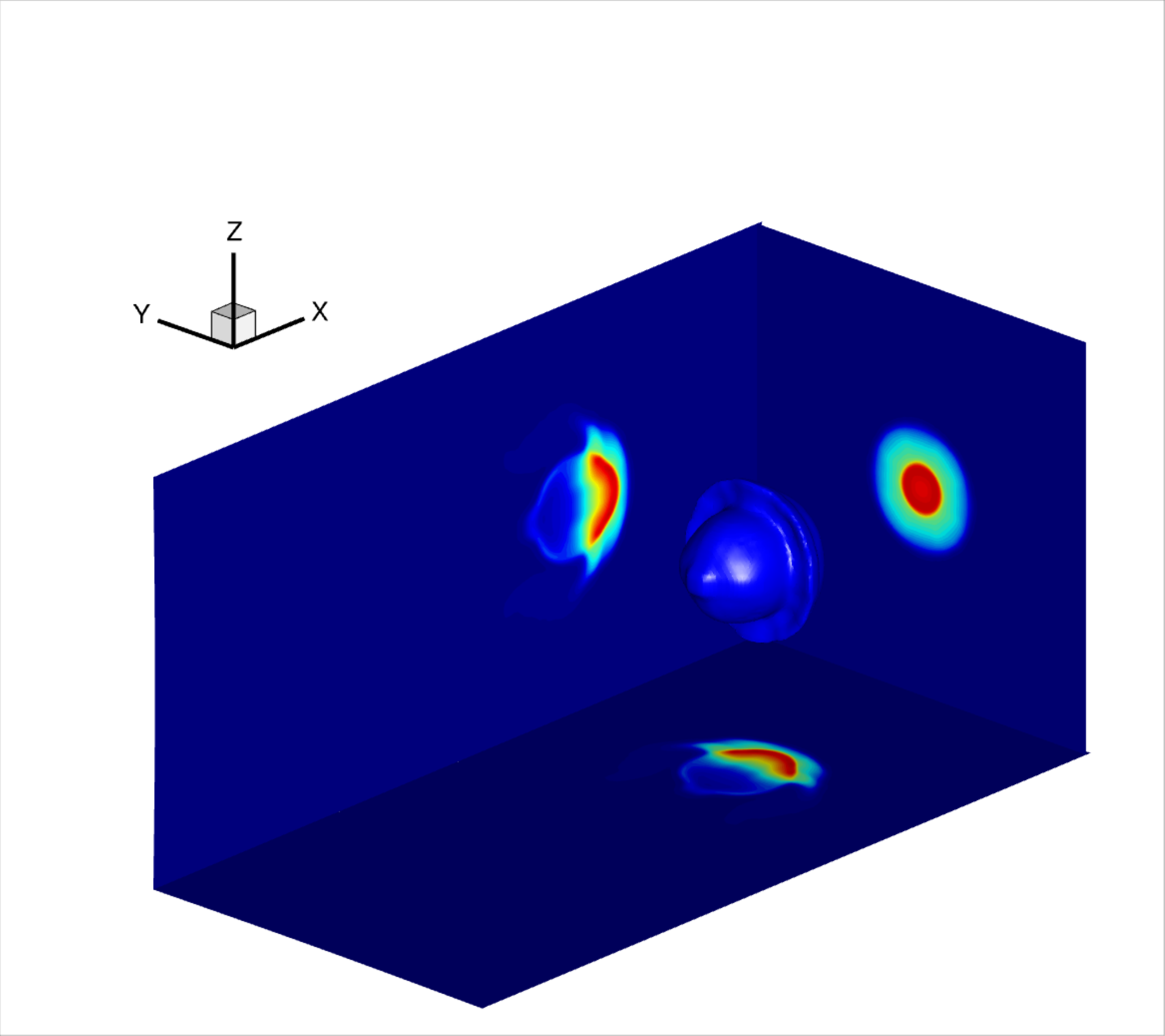}
		\end{subfigure}
		
		\begin{subfigure}[b]{0.42\textwidth}
			\centering
			\includegraphics[width=1.0\linewidth,trim=40 0 35 0, clip]{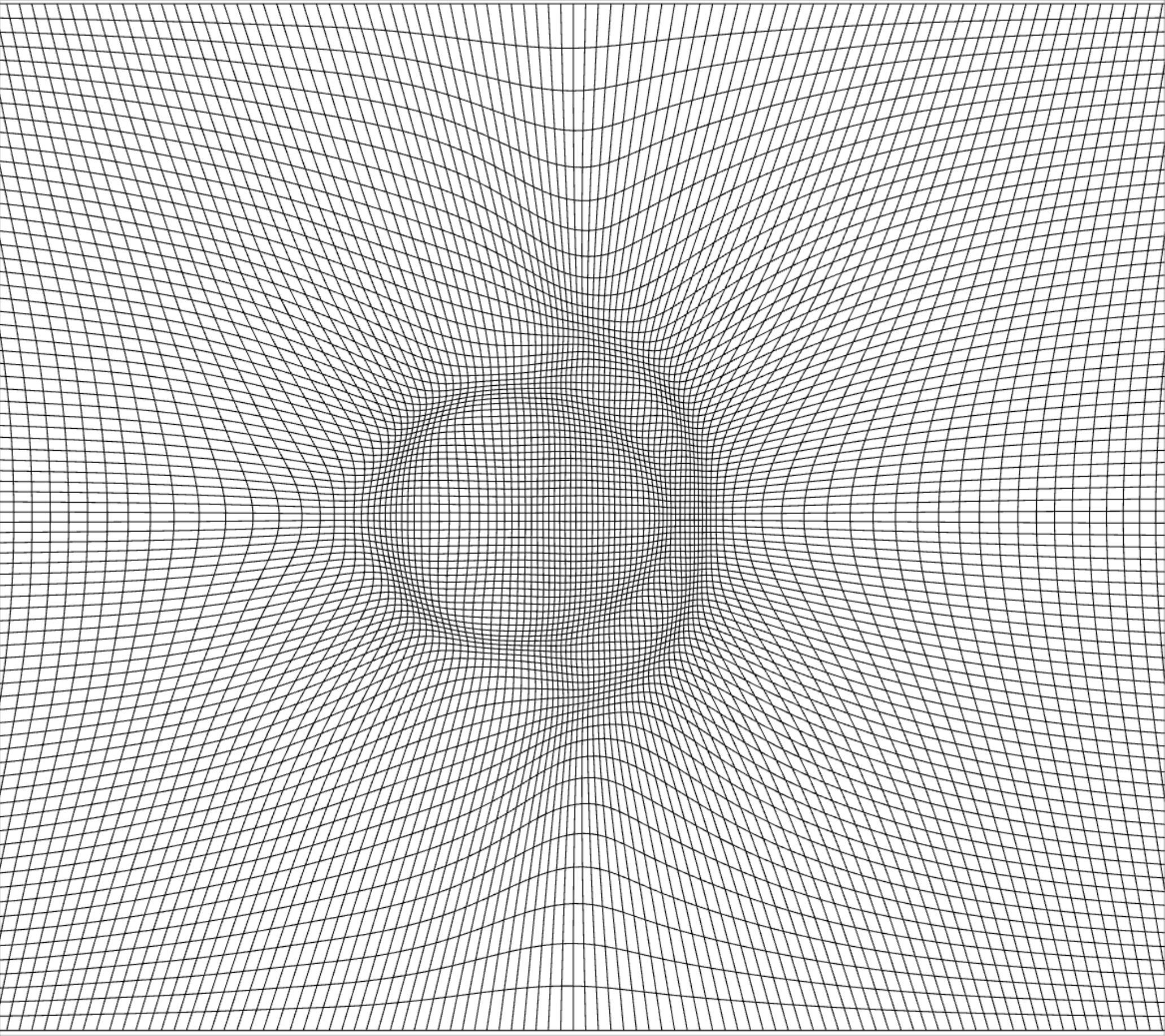}
		\end{subfigure}
		\begin{subfigure}[b]{0.42\textwidth}
			\centering
			\includegraphics[width=1.0\linewidth,trim=40 0 35 0, clip]{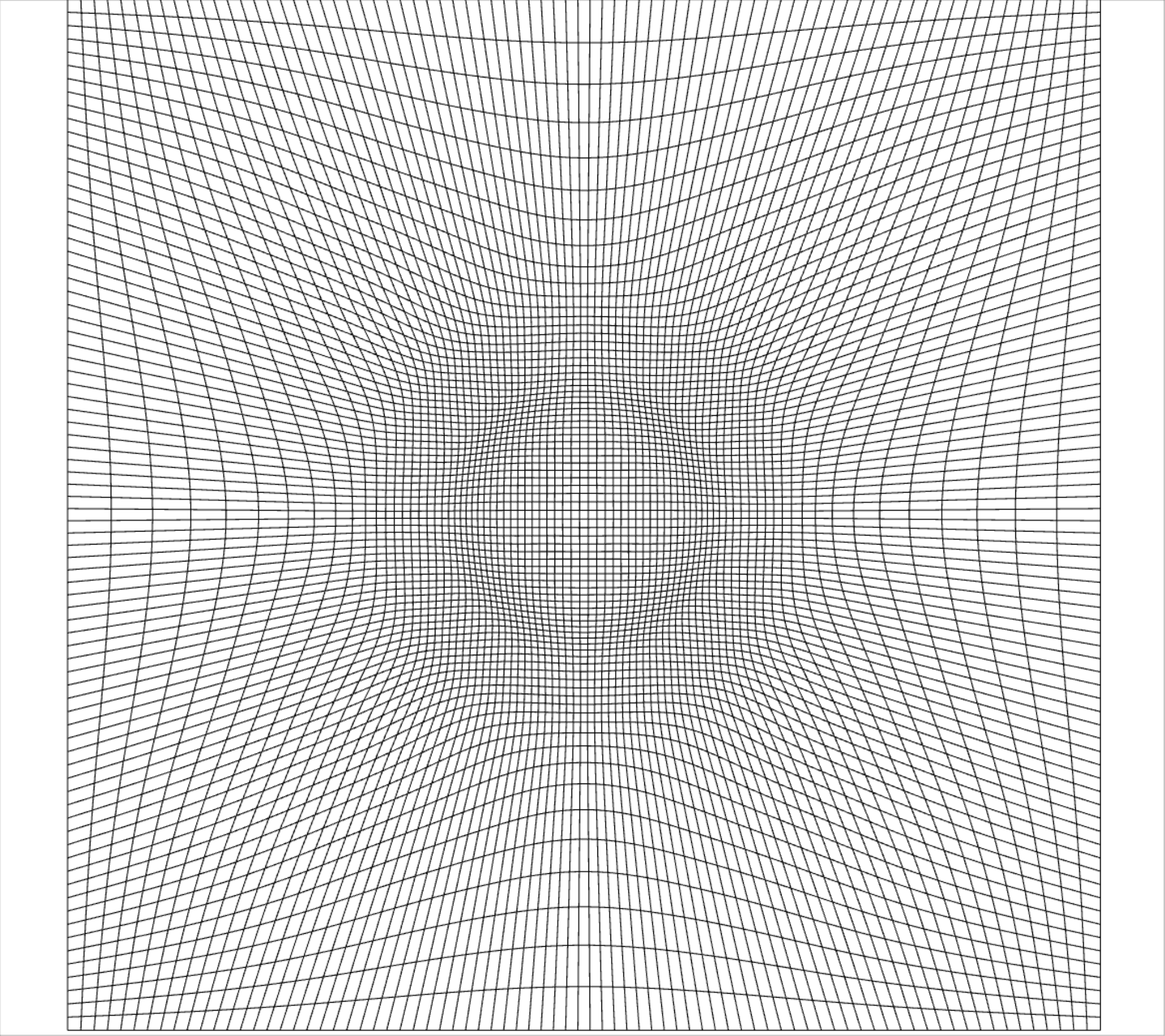}
		\end{subfigure}
		\caption{Example \ref{ex:3DShockBubble2}. Adaptive meshes  and $\rho$ at $t = 0.64$. Top left: close-up of the adaptive mesh, $i_1 \in [100, 230], i_2 \in [1, 40], i_3 \in [40, 80]$;  top right:
			the iso-surface of  $\rho = 11.11$ and  three offset 2D slices taken at $x_1 =
			214$, $x_2 = 0$, $x_3 = 0$;
			bottom left: the surface mesh with $i _1= 180$; bottom right:  the surface mesh with $i_2 = 40$.
		}
		\label{fig:3DMesh2}
	\end{figure}
	\begin{figure}[!ht]
		\centering
		
		\begin{subfigure}[b]{0.45\textwidth}
			\centering
			\includegraphics[width=2.6in,height=1.78in]{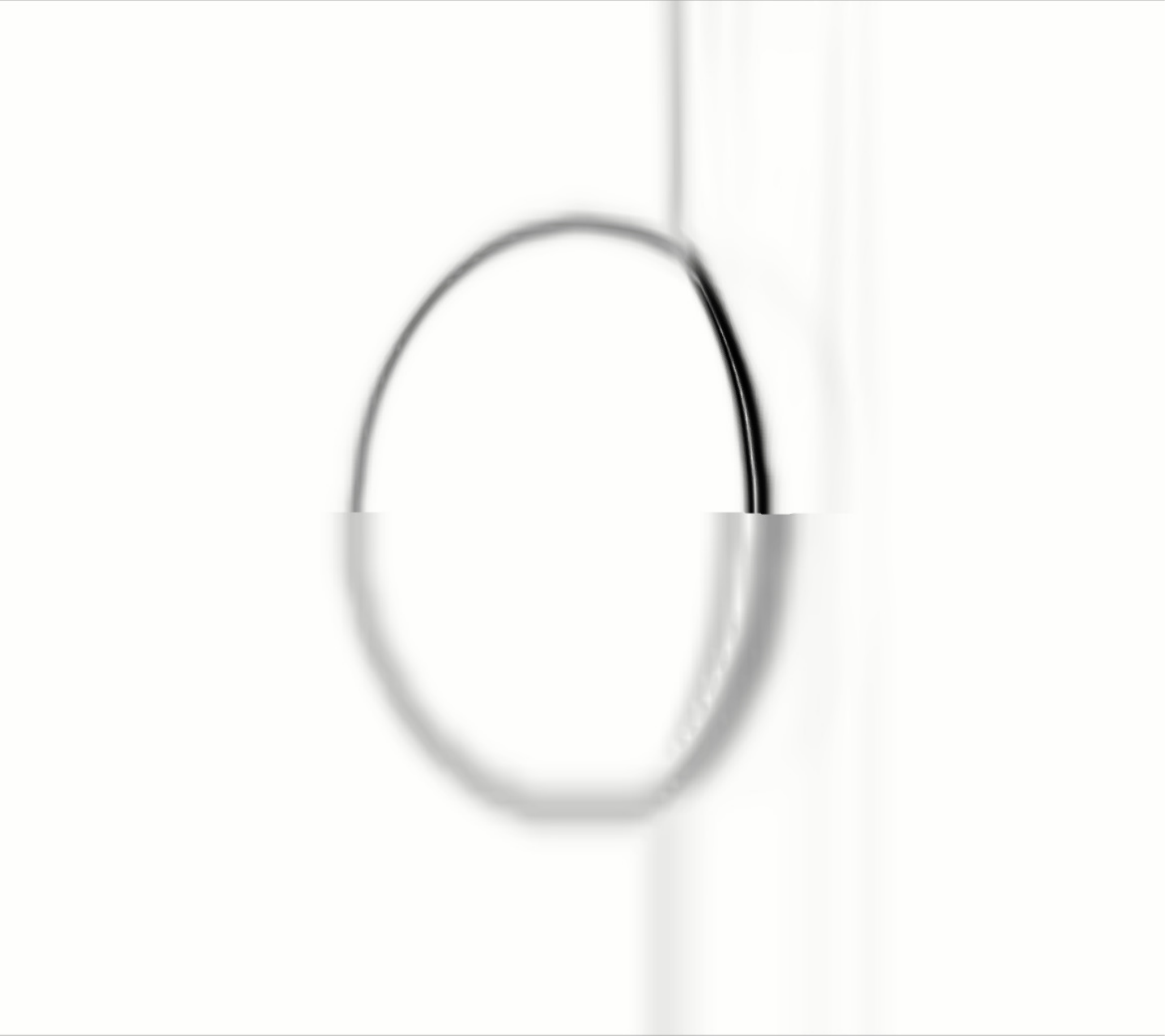}
		\end{subfigure}
		\begin{subfigure}[b]{0.45\textwidth}
			\centering
			\includegraphics[width=2.6in,height=1.78in]{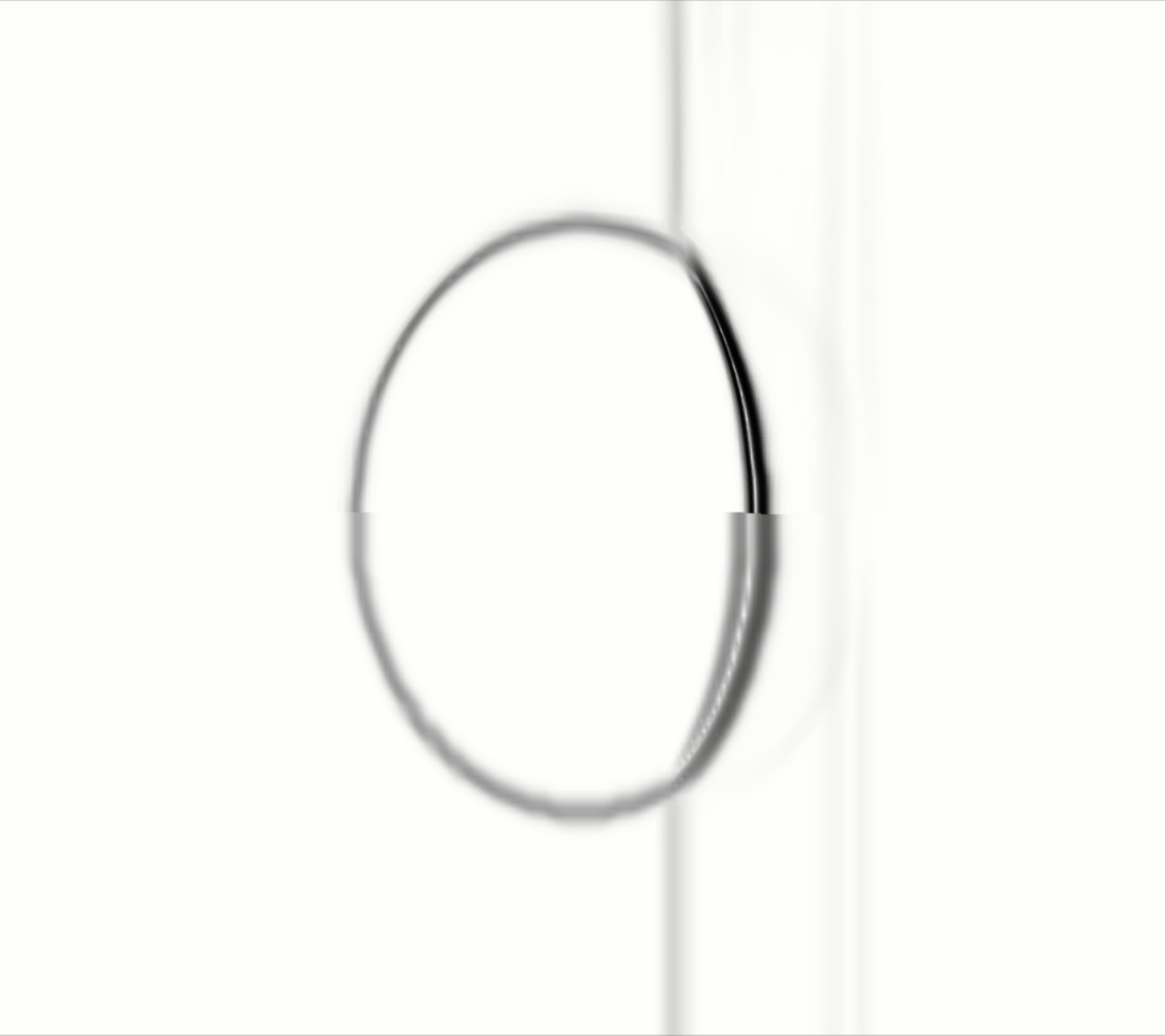}
		\end{subfigure}
		
		\begin{subfigure}[b]{0.45\textwidth}
			\centering
			\includegraphics[width=2.6in,height=1.78in]{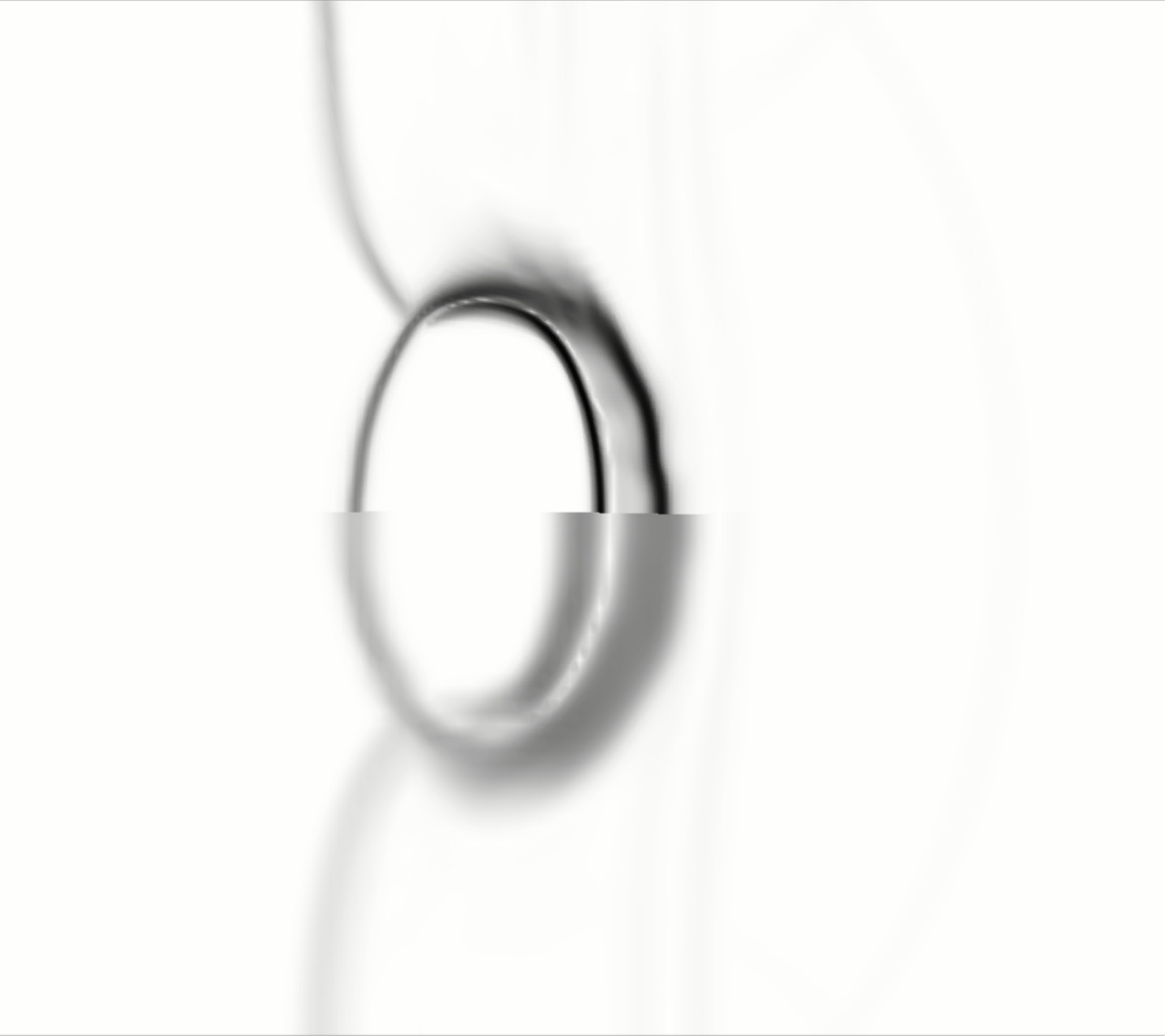}
		\end{subfigure}
		\begin{subfigure}[b]{0.45\textwidth}
			\centering
			\includegraphics[width=2.6in,height=1.78in]{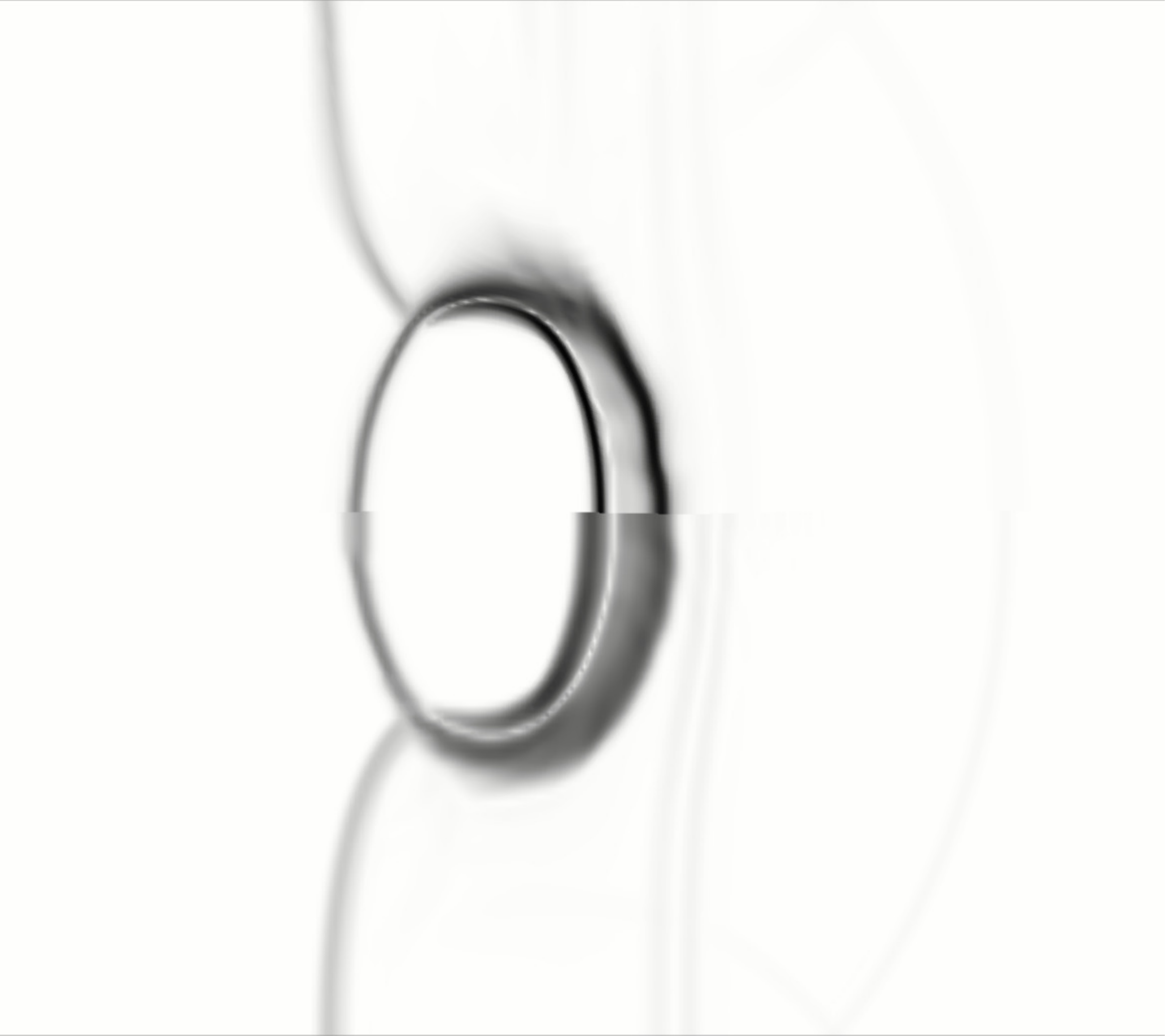}
		\end{subfigure}
		
		\begin{subfigure}[b]{0.45\textwidth}
			\centering
			\includegraphics[width=2.6in,height=1.78in]{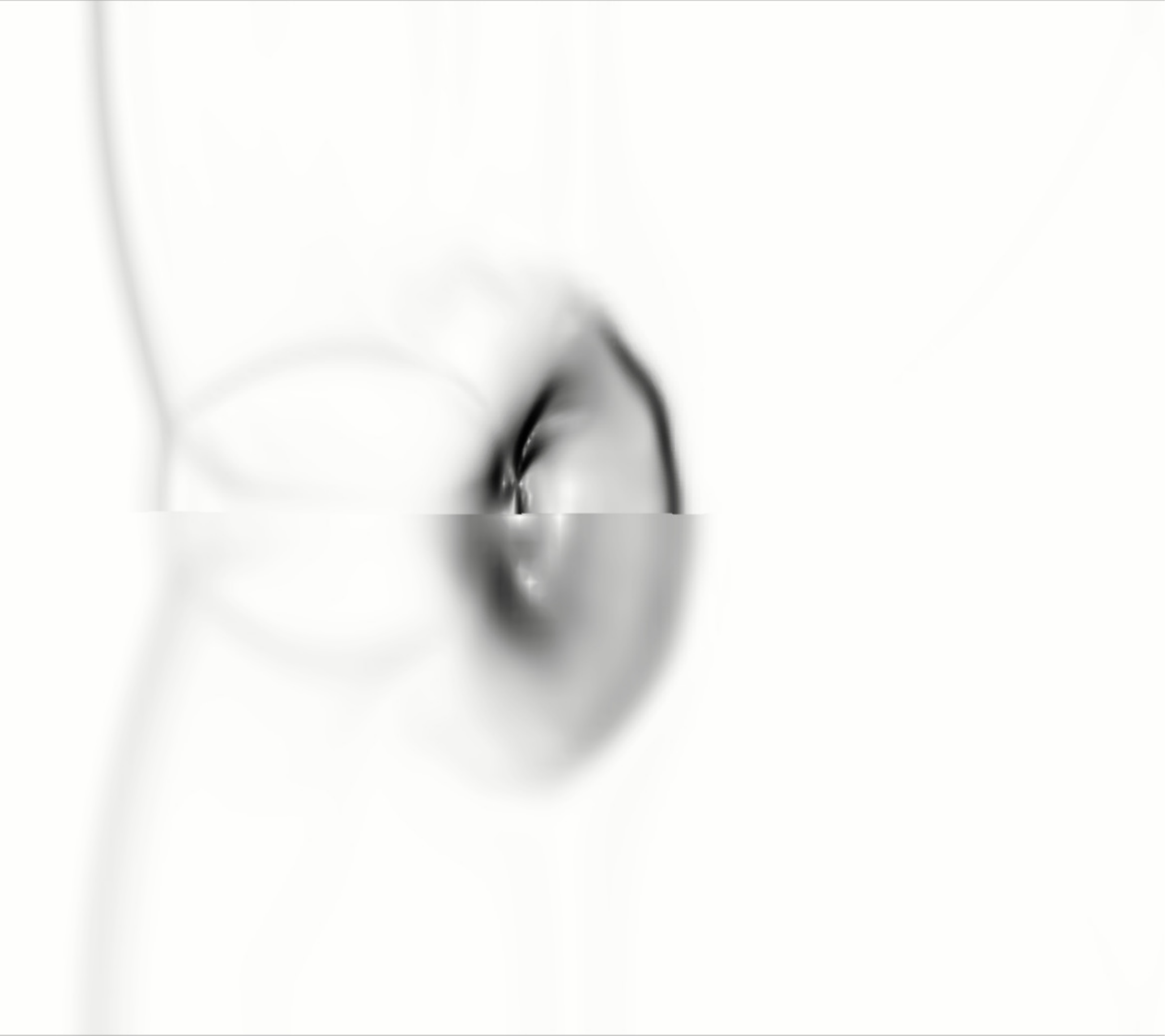}
		\end{subfigure}
		\begin{subfigure}[b]{0.45\textwidth}
			\centering
			\includegraphics[width=2.6in,height=1.78in]{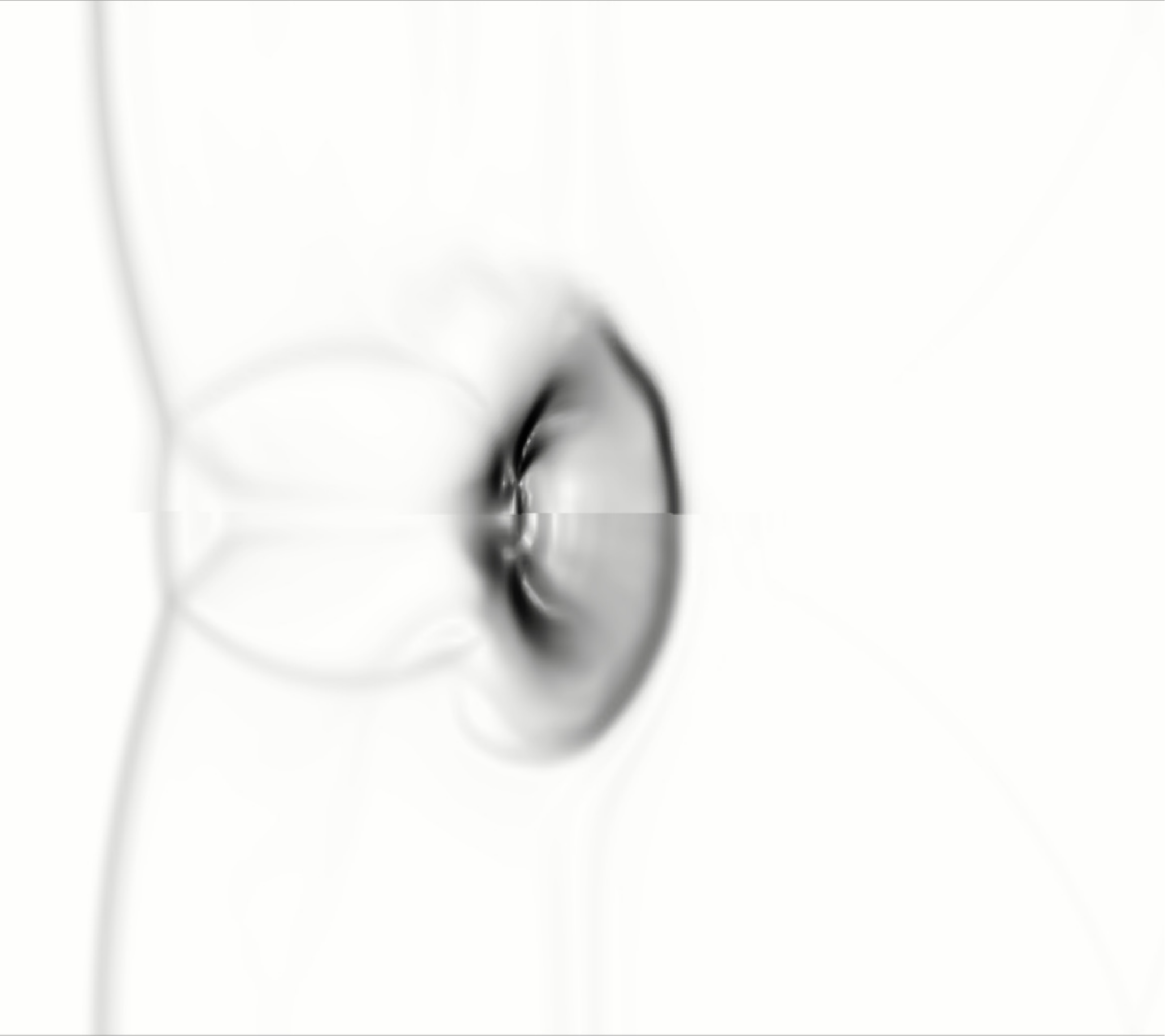}
		\end{subfigure}
		
		\begin{subfigure}[b]{0.45\textwidth}
			\centering
			\includegraphics[width=2.6in,height=1.78in]{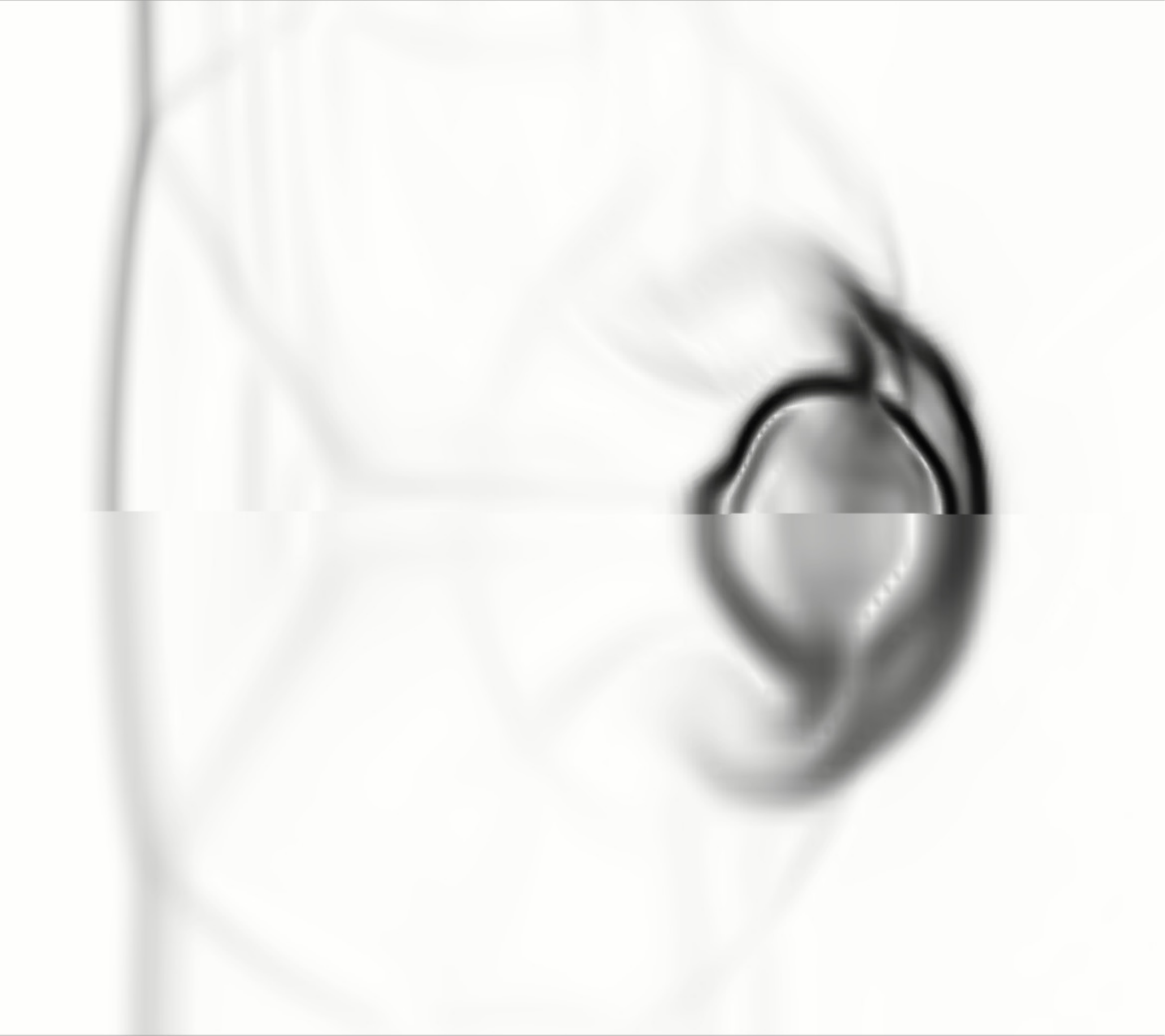}
		\end{subfigure}
		\begin{subfigure}[b]{0.45\textwidth}
			\centering
			\includegraphics[width=2.6in,height=1.78in]{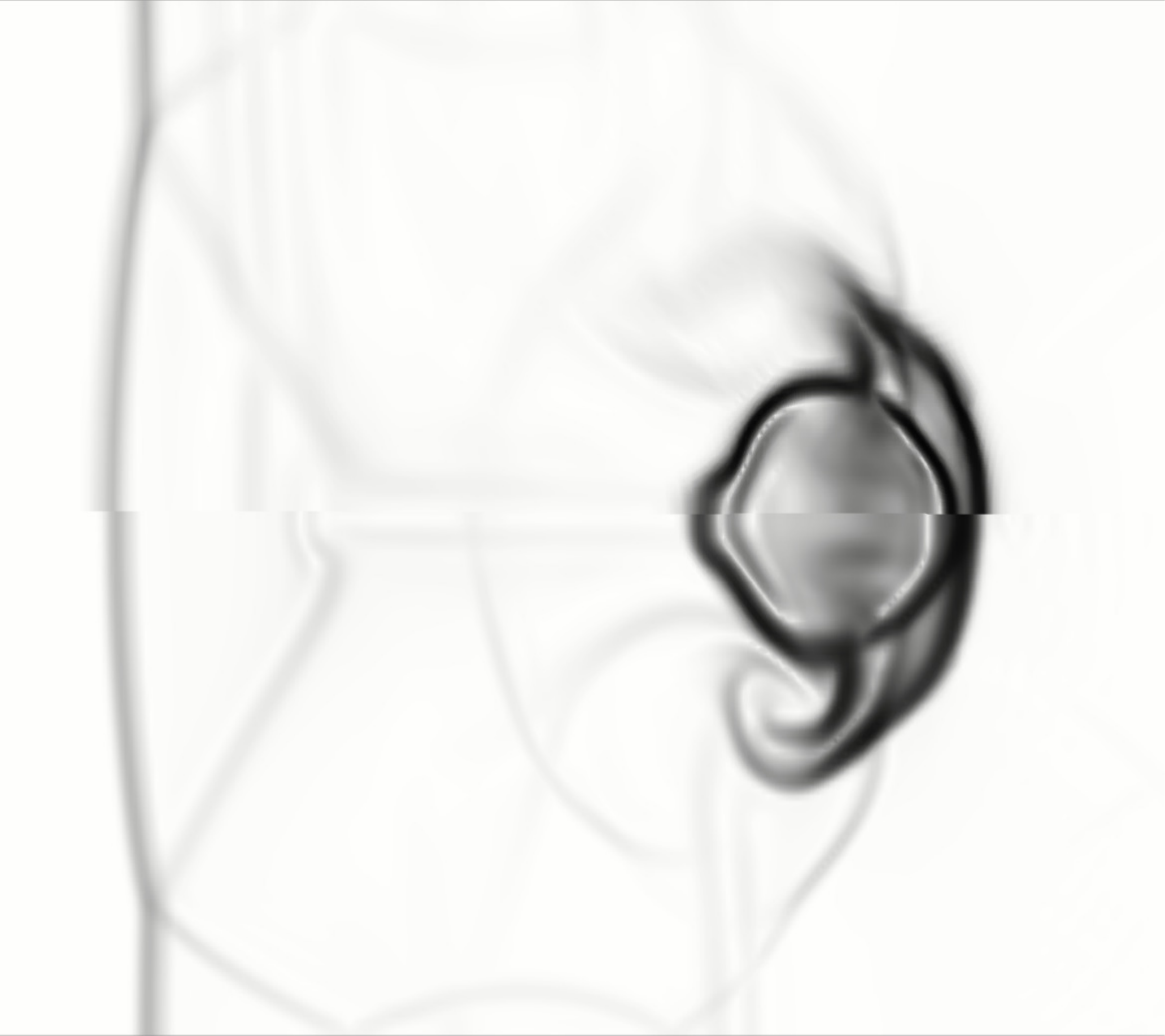}
		\end{subfigure}
		
		\caption{Example \ref{ex:3DShockBubble2}. Schlieren images $\Phi$ on $x_3 = 0$ at $t = 0.16, 0.32, 0.48, 0.64$ (from top to bottom). Left:    {\tt MM-WENOMR} (top half) and  {\tt UM-WENOMR} (bottom half) with $400 \times 80 \times 80$  cells;
right:    {\tt MM-WENOMR}  with $400 \times 80 \times 80$  cells (top half) and
 {\tt UM-WENOMR} with $800 \times 160 \times 160$ cells (bottom half). }
		\label{fig:3DShockbubble2}
	\end{figure}

\end{example}

%% file: Conc.tex
\section{Conclusion}\label{section:Conc}
This paper extended the  high-order entropy stable (ES) adaptive moving mesh finite difference schemes developed in 
		\cite{duan2021highorder}
		to the two- and three-dimensional (multi-component) compressible Euler equations
    	with   the stiffened equation of state.
	The two-point entropy conservative (EC) flux in the curvilinear coordinates was  first  constructed  according to the sufficient condition, and then combined with  the high-order  discretization of  the GCLs to derive
the high-order semi-discrete EC schemes.
 To ensure the entropy stability and to suppress the numerical oscillation, the high-order ES schemes were constructed by adding some suitable high-order  dissipation terms based on the multi-resolution WENO reconstruction for the scaled entropy variables to the high-order EC schemes.
 The high-order explicit strong-stability-preserving Runge-Kutta time discretization was implemented to obtain the fully-discrete high-order schemes,
  which could preserve the free-stream states and where the   mesh points were adaptively redistributed  by  solving the mesh iterative redistribution equations    with an
  appropriately chosen monitor function.  Several  2D and 3D numerical tests  were conducted on the parallel computer system with the MPI programming to validate the accuracy and the ability  to  resolve   effectively the localized structures of the present schemes.

\section*{Acknowledgments}
The authors were partially supported by the National Numerical Windtunnel project,
the National Key R\&D Program of China (Project Number 2020YFA0712000),
the National Natural Science Foundation of China (No.  12126302 \& 12171227),
and High-performance Computing Platform of Peking University.

%% file: WENO.tex
\appendix
\section{Multi-resolution WENO reconstruction}\label{section:WENO}
This appendix outlines  the 1D fifth-order multi-resolution WENO reconstruction proposed in \cite{WANG2021105138}. Different from the classical WENO reconstruction \cite{Jiang1996Efficient}, it uses three unequal-sized central stencils to get the WENO reconstructed polynomial with any positive linear weights whose sum is one,
 and can {achieve} the expected accuracy in the smooth regions and  the  good performance in  capturing sharp shock transitions.
For the given cell average values $\{W_{j}\}$, the fifth-order reconstruction process for $W^{\text{\tt\tiny WENO}}_{j\pm\frac12}$ is described as follows.
\begin{enumerate}	
	\item For each $j$, select three central spatial stencils $[\xi_{j-m+\frac12}, \xi_{j+m-\frac12}]$, $m = 1,2,3$, and then reconstruct 
$(2m-2)$th {degree} polynomials
 $h_m(\xi)$, $m = 1,2,3$, respectively, by
	\begin{align*}
	&	\frac{1}{\Delta \xi} \int_{\xi_{l-\frac12}}^{\xi_{l+\frac12}} h_{1}(\xi) d \xi=W_{l},\quad l = j,\\
	&	\frac{1}{\Delta \xi} \int_{\xi_{l-\frac12} }^{\xi_{l+\frac12}} h_{2}(\xi) d \xi=W_{l}, \quad l=j-1, j, j+1,\\
	&	\frac{1}{\Delta \xi} \int_{\xi_{l-\frac12}}^{\xi_{l+\frac12}} h_{3}(\xi) d \xi=W_{l}, \quad l=j-2, j-1, j, j+1, j+2.
	\end{align*}
	Those  polynomials $\{h_m(\xi)\}$
can be explicitly given by
	\begin{align}\label{basis}\begin{aligned}
&	h_{1}(\xi)=W_{j},   \quad	h_{2}(\xi)=a_{1} L_{1}(\xi)+a_{2} L_{2}(\xi)+a_{3} L_{3}(\xi), \\
&	h_{3}(\xi)=b_{1} L_{1}(\xi)+b_{2} L_{2}(\xi)+b_{3} L_{3}(\xi)+b_{4} L_{4}(\xi)+b_{5} L_{5}(\xi),
	\end{aligned}\end{align}
where $L_\ell$, $\ell=1,\cdots,5$, are  the (local) Legendre basis functions
			\begin{align*}
	L_{1}(\xi)&=1,\quad
	L_{2}(\xi)=\frac{\xi-\xi_{j}}{\Delta \xi},\quad
	L_{3}(\xi)=\left(\frac{\xi-\xi_{j}}{\Delta \xi}\right)^{2}-\frac{1}{12},\\
	L_{4}(\xi)&=\left(\frac{\xi-\xi_{j}}{\Delta \xi}\right)^{3}-\frac{3}{20}\left(\frac{\xi-\xi_{j}}{\Delta \xi}\right),\quad
	L_{5}(\xi)=\left(\frac{\xi-\xi_{j}}{\Delta \xi}\right)^{4}-\frac{3}{14}\left(\frac{\xi-\xi_{j}}{\Delta \xi}\right)^{2}+\frac{3}{560},
	\end{align*}
	and the coefficients $a_m, b_\ell$, $m=1,\cdots,3, \ell=1,\cdots,5$, are calculated by
	\begin{align*}
	&a_{1}=W_{j},\quad a_{2}=\frac{W_{j+1}-W_{j-1}}{2}, \quad a_{3}=\frac{W_{j-1}-2 W_{j}+W_{j+1}}{2},\\
   &b_{1}=W_{j},\quad b_{2}=\frac{11 W_{j-2}-82 W_{j-1}+82 W_{j+1}-11 W_{j+2}}{120},\\
&	b_{3}=\frac{-3 W_{j-2}+40 W_{j-1}-74 W_{j}+40 W_{j+1}-3 W_{j+2}}{56},\\
&b_{4}=\frac{-W_{j-2}+2W_{j-1}-2 W_{j+1}+W_{j+2}}{12},\\
&b_{5}=\frac{W_{j-2}-4 W_{j-1}+6 W_{j}-4 W_{j+1}+W_{j+2}}{24}.
	\end{align*}
	\item Compute the smoothness indicators $\beta_m$, $m = 1,2,3$,
by
\begin{align*}
\beta_{1}=\min \left(\left(W_{j}
-W_{j-1}\right)^{2},\left(W_{j+1}-W_{j}\right)^{2}\right),
\end{align*}
and
	\begin{align}\label{smoothness indicators}
	\beta_{m}=\sum_{\alpha=1}^{2m-2} \int_{\xi_{j-m+\frac{1}{2}}}^{\xi_{j+m-\frac{1}{2}}}(\Delta \xi)^{2 \alpha-1}\left(\frac{d^{\alpha} h_{m}(\xi)}{d\xi^{\alpha}}\right)^{2} d \xi,\quad m=2,3.
    \end{align}
Substituting \eqref{basis}  into \eqref{smoothness indicators} gives  
	\begin{align*}
\beta_{2}=a_{2}^{2}+\frac{13}{3} a_{3}^{2},\ \
\beta_{3}=\left(b_{2}+\frac{1}{10} b_{4}\right)^{2}+\frac{13}{3}\left(b_{3}+\frac{123}{455} b_{5}\right)^{2}+\frac{781}{20} b_{4}^{2}+\frac{1421461}{2275} b_{5}^{2}.
 \end{align*}
The above  $\beta_1$ avoids the zero value of $\beta_1$  calculated by \eqref{smoothness indicators} and the zeroth degree polynomial $h_1(\xi)$.

\item Choose the linear weights
 $\chi_{1} = 1/111, \chi_{2} = 10/111, \chi_{3} = 100/111$, unless otherwise stated,
 and compute the nonlinear weights $\{\hat{\omega}_{m}\}$ by
$$
\hat{\omega}_{m}=\frac{\bar{\omega}_{m}}{\sum_{\ell=1}^{3} \bar{\omega}_{\ell}}, \quad m=1,2,3,
$$
with
$$ \bar{\omega}_{m}=\chi_{m}\left(1+\frac{\tau}{\varepsilon+\beta_{m}}\right),
\quad
\tau=\left(\frac{\left|\beta_{3}-\beta_{1}\right|+\left|\beta_{3}
-\beta_{2}\right|}{2}\right)^{2},
$$
where $\varepsilon$ is taken as $10^{-10}$ to avoid the denominator becoming zero.
In practice, the linear weights can be any positive numbers whose sum is one. If the linear weight corresponding to the five-point sub-stencil takes a relatively large value,
then  sharper transitions in discontinuous regions can be obtained, but the numerical solution may tend to oscillate.

\item The reconstructed value $W^{\text{\tt\tiny WENO}}_{j+\frac12}$ by the fifth-order multi-resolution WENO reconstruction is computed by
$$
W^{\text{\tt\tiny WENO}}_{j+\frac12}=\hat{\omega}_{1} h_{1}\left(x_{j+\frac12}\right)+\hat{\omega}_{2} h_{2}\left(x_{j+\frac12}\right)+\hat{\omega}_{3}\left[\frac{1}{\chi_{3}} h_{3}\left(x_{j+\frac12}\right)-\sum_{m=1}^{2} \frac{\chi_{m}}{\chi_{3}} h_{m}\left(x_{j+\frac12}\right)\right],
$$
while   $W^{\text{\tt\tiny WENO}}_{j-\frac12}$ may  be obtained  by using the mirror-symmetry with respect to $x_{j}$.

\end{enumerate}

%% file: document.bbl
\begin{thebibliography}{100}

\bibitem{Bhoriya2020Entropy}
D.~Bhoriya and H.~Kumar, {Entropy-stable schemes for relativistic hydrodynamics
  equations}, \emph{Z. Angew. Math. Phys.}, 71 (2020),  1--29.

\bibitem{Biswas2018Low}
B.~Biswas and R.K. Dubey, {Low dissipative entropy stable schemes using third
  order WENO and TVD reconstructions}, \emph{Adv. Comput. Math.}, 44 (2018),
  1153--1181.

\bibitem{BOSCHERI2014484}
W.~Boscheri and M.~Dumbser, A direct {A}rbitrary-{L}agrangian–{E}ulerian
  {ADER-WENO} finite volume scheme on unstructured tetrahedral meshes for
  conservative and non-conservative hyperbolic systems in {3D}, \emph{J.
  Comput. Phys.}, 275 (2014),  484--523.

\bibitem{Brackbill1993An}
J.U. Brackbill, An adaptive grid with directional control, \emph{J. Comput.
  Phys.}, 108 (1993),  38--50.

\bibitem{Brackbill1982Adaptive}
J.U. Brackbill and J.S. Saltzman, Adaptive zoning for singular problems in two
  dimensions, \emph{J. Comput. Phys.}, 46 (1982),  342--368.

\bibitem{BRIO2001177}
M.~Brio, A.~Zakharian, and G.~Webb, Two-dimensional {R}iemann solver for
{E}uler equations of gas dynamics, \emph{J. Comput. Phys.}, 167 (2001),
177--195.


\bibitem{Budd2009Adaptivity}
C.J. Budd, W.Z. Huang, and R.D. Russell, Adaptivity with moving grids,
  \emph{Acta Numer.}, 18 (2009),  111--241.

\bibitem{CAO1999221}
W.M. Cao, W.Z. Huang, and R.D. Russell, An r-adaptive finite element method
  based upon moving mesh {PDE}s, \emph{J. Comput. Phys.}, 149 (1999),
  221--244.

\bibitem{Carpenter2014Entropy}
M.H. Carpenter, T.C. Fisher, E.J. Nielsen, and S.H. Frankel, {Entropy stable
  spectral collocation schemes for the Navier-Stokes equations: Discontinuous
  interfaces}, \emph{SIAM J. Sci. Comput.}, 36 (2014),  B835--B867.

\bibitem{CENICEROS2001609}
H.D. Ceniceros and T.Y. Hou, An efficient dynamically adaptive mesh for
  potentially singular solutions, \emph{J. Comput. Phys.}, 172 (2001),
  609--639.

\bibitem{Chen2020Review}
T.H. Chen and C.-W.~Shu, {Review of entropy stable discontinuous Galerkin methods
  for systems of conservation laws on unstructured simplex meshes}, \emph{CSIAM
  Trans. Appl. Math.}, 1 (2020),  1--52.

\bibitem{davis1992interface}
S.F. Davis, An interface tracking method for hyperbolic systems of conservation
  laws, \emph{Appl. Numer. Math.}, 10 (1992),  447--472.

\bibitem{Davis1982}
S.F. Davis and J.E. Flaherty, An adaptive finite element method for
  initial-boundary value problems for partial differential equations,
  \emph{SIAM J. Sci. Stat. Comput}, 3 (1982),  6--27.

\bibitem{duan2021highorder}
J.M. Duan and H.Z. Tang, High-order accurate entropy stable adaptive moving  mesh finite difference schemes for special relativistic
  (magneto)hydrodynamics, \emph{J. Comput. Phys.}, online (2022), 111038.

\bibitem{Duan2020RHD}
J.M. Duan and H.Z. Tang, {High-order accurate entropy stable finite difference
  schemes for one- and two-dimensional special relativistic hydrodynamics},
  \emph{Adv. Appl. Math. Mech.}, 12 (2020),  1--29.

\bibitem{Duan2020RMHD}
J.M. Duan and H.Z. Tang, {High-order accurate entropy stable nodal
  discontinuous Galerkin schemes for the ideal special relativistic
  magnetohydrodynamics}, \emph{J. Comput. Phys.}, 421 (2020),  109731.

\bibitem{DUAN2021109949}
J.M. Duan and H.Z. Tang, Entropy stable adaptive moving mesh schemes for 2{D}
  and 3{D} special relativistic hydrodynamics, \emph{J. Comput. Phys.}, 426
  (2021),  109949.

\bibitem{MR991358}
G.~Fernandez and B.~Larrouturou, Hyperbolic schemes for multi-component {E}uler
  equations, in \emph{Nonlinear {H}yperbolic {E}quations-{T}heory, {C}omputation
  {M}ethods, and {A}pplications ({A}achen, 1988)}, vol.~24 of \emph{Notes Numer.
  Fluid Mech.} (1989), 128--138.

\bibitem{Fjordholm2012Arbitrarily}
U.S. Fjordholm, S.~Mishra, and E.~Tadmor, {Arbitrarily high-order accurate
  entropy stable essentially non-oscillatory schemes for systems of
  conservation laws}, \emph{SIAM J. Numer. Anal.}, 50 (2012),  544--573.

\bibitem{GALERA20105755}
S.~Galera, P.H. Maire, and J.~Breil, A two-dimensional unstructured
  cell-centered multi-material {ALE} scheme using {VOF} interface
  reconstruction, \emph{J. Comput. Phys.}, 229 (2010),  5755--5787.

\bibitem{Gassner2013A}
G.J. Gassner, {A skew-symmetric discontinuous Galerkin spectral element
  discretization and its relation to SBP-SAT finite difference methods},
  \emph{SIAM J. Sci. Comput.}, 35 (2013),  1233--1253.

\bibitem{Gottlieb2001Strong}
S.~Gottlieb, C.-W. Shu, and E.~Tadmor, Strong stability-preserving high-order
  time discretization methods, \emph{SIAM Rev.}, 43 (2001),  89--112.

\bibitem{GOUASMI2020112912}
A.~Gouasmi, K.~Duraisamy, and S.M. Murman, Formulation of entropy-stable
  schemes for the multicomponent compressible {E}uler equations, \emph{Comput.
  Methods Appl. Mech. Engrg.}, 363 (2020),  112912.

\bibitem{haas_sturtevant_1987}
J.~Haas and B.~Sturtevant, {Interaction of weak shock waves with cylindrical
  and spherical gas inhomogeneities}, \emph{J. Fluid Mech.}, 181 (1987),
  41--76.

\bibitem{1987Explicit}
A.~Habbal, A.~Dervieux, H.~Guillard, and B.~Larrouturou, Explicit calculation
  of reactive flows with an upwind finite element hydrodynamical code,
  \emph{INRIA},  (1987).

  \bibitem{Hiltebrand2014Entropy}
  A.~Hiltebrand and S.~Mishra, {Entropy stable shock capturing space-time
  	discontinuous Galerkin schemes for systems of conservation laws},
  \emph{Numer. Math.}, 126 (2014),  103--151.

\bibitem{Ismail2009Affordable}
F.~Ismail and P.L. Roe, {Affordable, entropy-consistent Euler flux functions II
  : Entropy production at shocks}, \emph{J. Comput. Phys.}, 228 (2009),
  5410--5436.

\bibitem{Jiang1996Efficient}
G.S. Jiang and C.-W. Shu, {Efficient implementation of weighted WENO schemes},
  \emph{J. Comput. Phys.}, 126 (1996),  202--228.

\bibitem{LARROUTUROU199159}
B.~Larrouturou, How to preserve the mass fractions positivity when computing
  compressible multi-component flows, \emph{J. Comput. Phys.}, 95 (1991),
  59--84.

\bibitem{Larrouturou1989}
B.~Larrouturou and L.~Fezoui,
  On the equations of multi-component perfect of real gas inviscid flow,
   In \emph{Nonlinear Hyperbolic Problems} (edited by C. Carasso, P. Charrier, B. Hanouzet, and J.L. Joly), Lecture Notes in Mathematics, vol 1402, Springer, Berlin, Heidelberg, 1989, 69-98.

\bibitem{LaxLiu1998}
P.D. Lax and X.D. Liu, Solution of two-dimensional {R}iemann problems of gas
dynamics by positive schemes, \emph{SIAM J. Sci. Comput.}, 19 (1998),
319--340.

\bibitem{etde_20457592}
O.~Le~M\'{e}tayer, J.~Massoni, and R.~Saurel, Elaborating equations of state of
  a liquid and its vapor for two-phase flow models, \emph{Int. J. Thermal
  Sciences}, 43 (2004),  265--276.

\bibitem{Lefloch2002Fully}
P.G. LeFloch, J.M. Mercier, and C.~Rohde, {Fully discrete entropy conservative
  schemes of arbitraty order}, \emph{SIAM J. Numer. Anal.}, 40 (2002),
  1968--1992.


\bibitem{2007PLUTO}
A.~Mignone, G.~Bodo, S.~Massaglia, T.~Matsakos, O.~Tesileanu, C.~Zanni, and
  A.~Ferrari, {PLUTO}: {A} numerical code for computational astrophysics,
  \emph{Astrophys. J. Suppl. Ser.}, 170 (2007),  228--242.

\bibitem{miller1996high}
G.H. Miller and E.G. Puckett, A high-order godunov method for multiple
  condensed phases, \emph{J. Comput.Phys.}, 128 (1996),  134--164.

\bibitem{Miller1981}
K.~Miller, Moving finite elements. {II}, \emph{SIAM J. Numer. Anal.}, 18
  (1981),  1033--1057.

\bibitem{MULDER1992209}
W.~Mulder, S.~Osher, and J.A. Sethian, Computing interface motion in
  compressible gas dynamics, \emph{J. Comput. Phys.}, 100 (1992),  209--228.

\bibitem{puckett19923d}
E.~Puckett and J.~Saltzman, A 3{D} adaptive mesh refinement algorithm for
  multimaterial gas dynamics, \emph{Phys. D}, 60 (1992),  84--93.

\bibitem{quirk_karni_1996}
J.J. Quirk and S.~Karni, On the dynamics of a shock–bubble interaction,
  \emph{J. Fluid Mech.}, 318 (1996),  129--163.

\bibitem{Ren2000An}
W.Q. Ren and X.P. Wang, An iterative grid redistribution method for singular
  problems in multiple dimensions, \emph{J. Comput. Phys.}, 159 (2000),
  246--273.


\bibitem{Renac2021}
F. Renac, Entropy stable, robust and high-order DGSEM for the compressible multicomponent Euler equations,
\emph{J. Comput. Phys.}, 445 (2021), 110584.
\bibitem{SHYUE1998208}
K.M. Shyue, An efficient shock-capturing algorithm for compressible
  multicomponent problems, \emph{J. Comput.Phys.}, 142 (1998),  208--242.

\bibitem{Stockie2001}
J.M. Stockie, J.A. Mackenzie, and R.D. Russell, A moving mesh method for
  one-dimensional hyperbolic conservation laws, \emph{SIAM J. Sci. Comput}, 22
  (2001),  1791--1813.

\bibitem{Tadmor1987The}
E.~Tadmor, {The numerical viscosity of entropy stable schemes for systems of
  conservation laws, I}, \emph{Math. Comp.}, 49 (1987),  91--103.

\bibitem{Tadmor2003Entropy}
E.~Tadmor, {Entropy stability theory for difference approximations of nonlinear
  conservation laws and related time-dependent problems}, \emph{Acta Numer.},
  12 (2003),  451--512.

\bibitem{Tang2005Moving}
T.~Tang, Moving mesh methods for computational fluid dynamics, \emph{Contemp.
  Math.}, 383 (2005),  141--173.

\bibitem{VISBAL2002155}
M.R. Visbal and D.V. Gaitonde, On the use of higher-order finite-difference
  schemes on curvilinear and deforming meshes, \emph{J. Comput. Phys.}, 181
  (2002),  155--185.

\bibitem{Wang2004A}
D.S. Wang and X.P. Wang, A three-dimensional adaptive method based on the
  iterative grid redistribution, \emph{J. Comput. Phys.}, 199 (2004),
  423--436.

\bibitem{WANG2021105138}
Z.~Wang, J.~Zhu, L.~Tian, Y.~Yang, and N.~Zhao, An efficient fifth-order finite
  difference multi-resolution {WENO} scheme for inviscid and viscous flow
  problems, \emph{Comput. $\&$ Fluids}, 230 (2021),  105138.

\bibitem{Winslow1967Numerical}
A.M. Winslow, Numerical solution of the quasilinear {P}oisson equation in a
  nonuniform triangle mesh, \emph{J. Comput. Phys.}, 1 (1967),  149--172.

  \bibitem{wu2008general}
  H.~Wu and Z.J. Shen, A general {R}iemann solver for {E}uler equations,
  \emph{Int. J. Numer. Meth Fluids}, 56 (2008),  1993--2020.

\bibitem{Wu2020Entropy}
K.L. Wu and C.-W. Shu, Entropy symmetrization and high-order accurate entropy
  stable numerical schemes for relativistic {MHD} equations, \emph{SIAM J. Sci.
  Comput.}, 42 (2020),  A2230--A2261.

\bibitem{xu1997bgk}
K.~Xu, {BGK}-based scheme for multicomponent flow calculations, \emph{J.
  Comput.Phys.}, 134 (1997),  122--133.

\bibitem{Zhang1993Discrete}
H.~Zhang, M.~Reggio, J.Y. Tr$\acute{\text{e}}$panier, and R.~Camarero,
Discrete form of the {GCL} for moving meshes and its implementation in {CFD} schemes,
  \emph{Comput. {$\&$} Fluids}, 22 (1993),  9--23.




%
%

\end{thebibliography}
